\documentclass[11pt]{amsart}
\usepackage{graphicx}
\usepackage{epsfig}
\usepackage{amsmath}
\usepackage{stackrel, amssymb}
\usepackage{amscd}
\usepackage{latexsym}
\usepackage{tabularx}
\usepackage{a4wide}
\usepackage[usenames]{color}
\usepackage{subfigure}
\usepackage{url}
\usepackage{verbatim}
\usepackage{mathtools}
\usepackage{circuitikz}
\usepackage{mathabx}
\usepackage{tikz-cd}
\usepackage[T1]{fontenc}
\DeclareMathAlphabet{\mathdutchcal}{U}{dutchcal}{m}{n}

\newcommand{\dz}{\ssub{\mathdutchcal{z}}!}
\newcommand{\dq}{\ssub{\mathdutchcal{q}}!}

\newcommand{\Gm}{\ssub{\ssub[-1pt]{\G}!_{\mathbf m}}!}

\newcommand{\car}[1]{\mathrm{char}(#1)}
\newcommand{\dpat}{\mathdutchcal{p}}
\definecolor{cadmiumgreen}{rgb}{0.0, 0.42, 0.24}
\usepackage{enumitem} 
\setlist[enumerate]{label=\textnormal{(\arabic*)}}

\ctikzset{bipoles/resistor/height=0.15}
\ctikzset{bipoles/resistor/width=0.4}
\usepackage[
colorlinks, citecolor=cadmiumgreen, linkcolor=blue,
pdfauthor={Omid Amini, Eduardo Esteves, Eduardo Garcez}, 
pdftitle={limit-canonical-series},
pdfstartview ={FitV},
]{hyperref}
\usepackage{tikz, float} 
\usetikzlibrary {positioning}
\usetikzlibrary{calc,decorations.markings}
\usetikzlibrary{shapes,snakes}
\numberwithin{equation}{section}
\tikzstyle{Cwhite}=[scale = .8,circle, fill = white, minimum size=3mm] 
\tikzstyle{Cgray}=[scale = .4,circle, fill = gray, minimum size=3mm] 
\tikzstyle{Cblack2}=[scale = .4,circle, fill = black, minimum size=5mm] 
\tikzstyle{Cblack}=[scale = .7,circle, fill = black, minimum size=3mm]
\tikzstyle{C0}=[scale = .9,circle, fill = black!0, inner sep = 0pt, minimum size=3mm]
\tikzstyle{C1}=[scale = .7,circle, fill = black!0, inner sep = 0pt, minimum size=3mm]
\tikzstyle{Cred}=[scale = .4,circle, fill = red, minimum size=3mm] 
\usepackage[matrix,arrow,tips,curve]{xy}
\usepackage{pb-diagram,pb-xy}
\usepackage{verbatim}
\usepackage{mathrsfs}
\usepackage{color}
\usepackage[leqno]{amsmath}
\usepackage{scalerel}
\usepackage{marginnote}

\RequirePackage[capitalize, noabbrev]{cleveref}
\newtheorem{thm}{Theorem}[section]

\newtheorem{lemma}[thm]{Lemma}
\newtheorem{problem}[thm]{Problem}

\newtheorem{prop}[thm]{Proposition}

\newtheorem{cor}[thm]{Corollary}

\theoremstyle{definition}
\newtheorem{defii}[thm]{Definition}
\newenvironment{defi}
  {\pushQED{\qed}\defii}
  {\popQED\enddefii}
\newtheorem{remm}[thm]{Remark}
\newenvironment{remark}
  {\pushQED{\qed}\remm}
  {\popQED\endremm}

\newtheorem{exx}[thm]{Example}
\newenvironment{example}
  {\pushQED{\qed}\exx}
  {\popQED\endexx}
\numberwithin{equation}{section}

\newcommand{\G}{\mathbf G}
\newcommand{\K}{\ssub{\mathbf K}!}
\renewcommand{\k}{\ssub{\mathbf k}!}
\newcommand{\R}{\ssub{\mathbb R}!}
\renewcommand{\P}{\ssub[-1pt]{\mathbf P}!}

\newcommand{\Q}{\ssub{\mathbf Q}!}

\newcommand{\Z}{\ssub{\mathbb Z}!}
\renewcommand{\:}{\colon}
\newcommand{\C}{{\mathscr C}}

\newcommand{\sshateta}{\ssub{\hat\eta}!}

\newcommand{\g}{\mathfrak{g}}
\newcommand{\red}{\ssub{\mathrm{red}}!}
\newcommand{\res}{\ssub{\mathrm{res}}!}
\newcommand{\Res}{\ssub{\mathrm{Res}}!}
\newcommand{\bRes}{\ssub{\mathbf{Res}}!}

\newcommand{\valuation}{\mathfrak v}
\newcommand{\val}{\mathrm{val}}
\newcommand{\trop}{\mathrm{trop}}
\newcommand{\an}{{^\mathrm{an}}}
\usepackage{todonotes}

\newcommand{\rmint}{{\mathrm{int}}}

\NewDocumentCommand{\ssub}{O{0pt} O{.9} m t! e{_^}}{
  #3%
  \IfValueT{#5}{
    \IfBooleanTF{#4}{\sb{\hspace{#1}\scaleobj{#2}{#5}}}{\sb{#5}}
  }
  \IfValueT{#6}{
  \IfBooleanTF{#4}{\sp{\hspace{#1}\scaleobj{#2}{#6}}}{\sp{#6}}
}
}

\NewDocumentCommand{\ssubb}{O{0pt} O{0pt} O{.9} m t! e{_^}}{
  #4%
  \IfValueT{#6}{
    \IfBooleanTF{#5}{\sb{\hspace{#1}\scaleobj{#3}{#6}}}{\sb{#6}}
  }
  \IfValueT{#7}{
  \IfBooleanTF{#5}{\sp{\hspace{#2}\scaleobj{#3}{#7}}}{\sp{#7}}
}
}
\ExplSyntaxOn
\NewDocumentCommand{\tossub}{o o m}{
  \expandafter\let\csname old\cs_to_str:N #3\endcsname#3
  \renewcommand#3%
  {\ssub[#1][#2]{\csname old\cs_to_str:N #3\endcsname}}
}
\ExplSyntaxOff

\newcommand{\stable}{\mathscr A}

\newcommand{\rquot}[2]{#1\bigl/#2}
\newcommand{\adfun}{\ssub{\mathscr F}!}

\newcommand{\cd}{\mathit{cd}}
\newcommand{\dd}{\mathit{d}}

\newcommand{\VS}{\ssub{\mathbf U}!}
\newcommand{\scrU}{\ssub{\mathscr U}!}
\newcommand{\scrL}{\ssub{\mathscr L}!}

\tossub\omega
\tossub\Omega
\tossub\alpha
\tossub\beta
\tossub\sigma
\tossub\mu
\tossub\nu
\tossub\varphi
\tossub\delta
\tossub\Delta
\tossub\Gamma
\tossub\gamma
\tossub\zeta
\tossub\eta
\tossub\stable
\tossub\C
\tossub\iota
\tossub\rho
\tossub\varrho

\tossub\g
\tossub\red
\tossub\cd
\tossub\dd
\tossub\varpi
\tossub\lambda
\tossub\epsilon
\tossub\varepsilon
\tossub\theta
\tossub\xi

\newcommand{\sssigma}{\ssub[-1pt]{\sigma}!}
\newcommand{\ssSigma}{\ssub{\Sigma}!}

\newcommand{\sszeta}{\ssub[-1pt]{\zeta}!}
\newcommand{\sstau}{\ssub[-1pt]{\tau}!}

\newcommand{\ssL}{\ssub{L}!}

\newcommand{\cT}{\ssub{\mathcal T}!}
\newcommand{\cO}{\ssub{\mathcal O}!}
\newcommand{\hatcO}{\ssub{\hat{\mathcal O}}!}

\newcommand{\ssN}{\ssub{N}!}
\newcommand{\ssI}{\ssub{I}!}
\newcommand{\ssJ}{\ssub{J}!}
\newcommand{\ssP}{\ssub{P}!}
\newcommand{\ssT}{\ssub{T}!}
\newcommand{\calJ}{\mathcal J}

\newcommand{\rmH}{\ssub{\mathrm{H}}!}

\newcommand{\sspi}{\ssub{\pi}!}
\newcommand{\ssG}{\ssub{G}!}
\newcommand{\ssD}{\ssub{D}!}

\newcommand{\ssW}{\ssub[-1pt]{\mathrm W}!}
\newcommand{\ssR}{\ssub[-1pt]{\mathrm R}!}
\newcommand{\ssn}{\ssub{n}!}
\newcommand{\ssq}{\ssub{q}!}
\newcommand{\ssh}{\ssub{h}!}

\newcommand{\ssE}{\ssub{E}!}

\newcommand{\sse}{\ssub{e}!}

\newcommand{\ssa}{\ssub{a}!}
\newcommand{\ssb}{\ssub{b}!}
\newcommand{\barb}{\ssub{\bar b}!}
\newcommand{\bara}{\ssub{\bar a}!}

\newcommand{\ssd}{\ssub{d}!}
\newcommand{\ssx}{\ssub{x}!}
\newcommand{\ssy}{\ssub{y}!}
\newcommand{\ssz}{\ssub{z}!}
\newcommand{\hatz}{\ssub{\hat z}!}
\newcommand{\barz}{\ssub{\bar z}!}
\newcommand{\wtz}{\ssub{\wt z}!}
\newcommand{\frakp}{\ssub{\mathfrak p}!}
\newcommand{\frakq}{\ssub{\mathfrak q}!}

\newcommand{\ssw}{\ssub{w}!}
\newcommand{\ssr}{\ssub{r}!}
\newcommand{\ba}{\ssub{\mathbf a}!}

\newcommand{\sst}{\ssub{t}!}

\newcommand{\ssX}{\ssub{X}!}
\newcommand{\scrX}{\ssub{\mathscr X}!}
\newcommand{\wtscrX}{\ssub{\wt{\mathscr X}}!}
\newcommand{\scrZ}{\ssub{\mathscr Z}!}

\newcommand{\ssboldPhi}{\ssub{\mathbf{\Phi}}!}
\newcommand{\ssPhi}{\ssub{\Phi}!}

\newcommand{\bphi}{\ssub{\mathbb{\phi}}!}
\newcommand{\ssm}{\ssub{m}!}
\newcommand{\ssp}{\ssub{p}!}

\newcommand{\ssS}{\ssub{S}!}
\newcommand{\N}{\mathbb N}
\newcommand{\ev}{\ssub{\mathrm{ev}}!}

\newcommand{\filt}{\ssub{\mathrm{F}}!}

\renewcommand{\setminus}{\smallsetminus}
\renewcommand{\emptyset}{\varnothing}

\newcommand{\Div}{\operatorname{Div}}
\renewcommand{\div}{\mathrm{div}}
\newcommand{\E}{\ssub{\mathbb E}!}
\newcommand{\A}{\ssub{A}!}   
\newcommand{\barA}{\ssub{\bar A}!}   

\newcommand{\ord}{\mathrm{ord}}
\newcommand{\he}{\ssub{\scaleto{\mathfrak{h}}{8pt}}!}  
\newcommand{\te}{\ssub{\scaleto{\mathfrak{t}}{6pt}}!}  

\makeatletter
\newsavebox\myboxA
\newsavebox\myboxB
\newlength\mylenA

\newcommand*\overbar[2][0.75]{%
    \sbox{\myboxA}{$\m@th#2$}%
    \setbox\myboxB\null
    \ht\myboxB=\ht\myboxA%
    \dp\myboxB=\dp\myboxA%
    \wd\myboxB=#1\wd\myboxA
    \sbox\myboxB{$\m@th\overline{\copy\myboxB}$}
    \setlength\mylenA{\the\wd\myboxA}
    \addtolength\mylenA{-\the\wd\myboxB}%
    \ifdim\wd\myboxB<\wd\myboxA%
       \rlap{\hskip 1\mylenA\usebox\myboxB}{\usebox\myboxA}%
    \else
        \hskip -0.5\mylenA\rlap{\usebox\myboxA}{\hskip 0.5\mylenA\usebox\myboxB}%
    \fi}

\newcommand{\comp}[1]{\overbar[.6]{#1}} 
\newcommand{\mg}{\mathscr M}
\newcommand{\mgbar}{\ssub[-1pt]{\comp\mg}!}

\newcommand{\wt}{\widetilde}

\newcommand{\Rat}{\mathrm{Rat}}

\newcommand{\slztwist}[1]{\partial_{_{ \hspace{-.04cm}#1}}}

\newcommand{\rest}[1]{\raisebox{-1pt}{$\vert$}_{#1}}

\newcommand{\proj}[1]{\theta_{_{ \hspace{-.07cm}#1}}}

\newcommand{\slp}{\mathrm{sl}}
\newcommand{\TT}{\mathrm{T}}

\newcommand{\cl}{\ssub{\mathscr D}!}  
\newcommand{\clcan}{\ssub{\mathscr K}!}  

\newcommand{\supp}{\ssub{\mathrm{supp}}!}

\newcommand{\st}{\bigm|} 

\newcommand{\Fun}{\mathscr E} 
\tossub\Fun 

\newcommand{\id}{\mathrm{Id}}
\tossub\id
\newcommand{\zero}{0}
\newcommand{\one}{\ssub{\mathbf 1}!}
\tossub\zero
\tossub\one

\newcommand{\ssrho}{\ssub{\rho}!}

\newcommand{\sspsi}{\ssub{\psi}!}

\tossub \chi

\newcommand{\varX}{\ssub{\mathbf X}!}
\newcommand{\varV}{\ssub[-1pt]{\mathbf V}!}
\newcommand{\compvarV}{\ssub[-1pt]{\overline{\mathbf V}}!}

\newcommand{\fX}{\mathfrak X}
\newcommand{\varC}{\mathbf C}

\newcommand{\varD}{\mathbf D}
\newcommand{\varH}{\mathbf H}
\newcommand{\varM}{\ssub{\mathbf M}!}

\newcommand{\varm}{\mathbf m}
\newcommand{\varR}{\mathbf R}

\newcommand{\varf}{\mathbf f}

\newcommand{\varalpha}{\ssub{\mathbf \alpha}!}
\newcommand{\cG}{\ssub{\mathcal G}!}
\newcommand{\scrG}{\ssub{\mathscr G}!}

\newcommand{\cW}{\ssub{\mathcal W}!}
\newcommand{\cC}{\ssub{\mathcal C}!}

\newcommand{\valgroup}{\Lambda}

\newcommand{\ssA}{\ssub{A}!}
\newcommand{\ssV}{\ssub{V}!}
\newcommand{\ssH}{\ssub{H}!}

\newcommand{\codim}{\mathrm{codim}}
\newcommand{\ssu}{\ssub{u}!}
\newcommand{\ssv}{\ssub{v}!}

\definecolor{cadmiumgreen}{rgb}{0.0, 0.42, 0.24}
\definecolor{darkred}{rgb}{.6,0,0}
\definecolor{byzant}{rgb}{0.74, 0.2, 0.64}
 \definecolor{pblue}{rgb}{0.11, 0.22, 0.73}
\definecolor{pgreen}{rgb}{0.0, 0.65, 0.58}
\definecolor{aqua}{rgb}{0.0, 1.0, 1.0}
\definecolor{lblue}{rgb}{0.0, 0.55, 1.0}
\definecolor{pblue}{rgb}{0.11, 0.22, 0.73}

\DeclareFontFamily{U}{MnSymbolA}{}
\DeclareSymbolFont{MnSyA}{U}{MnSymbolA}{m}{n}
\DeclareFontShape{U}{MnSymbolA}{m}{n}{
    <-6>  MnSymbolA5
   <6-7>  MnSymbolA6
   <7-8>  MnSymbolA7
   <8-9>  MnSymbolA8
   <9-10> MnSymbolA9
  <10-12> MnSymbolA10
  <12->   MnSymbolA12}{}
\DeclareMathSymbol{\dashedleftarrow}{\mathrel}{MnSyA}{98}
\DeclareMathSymbol{\dashedrightarrow}{\mathrel}{MnSyA}{96}

\makeatletter
\newcommand{\topreleft}[1]{%
  \vbox {\m@th\ialign{##\crcr
  \topreleftfill \crcr
  \noalign{\kern-\p@\nointerlineskip}
  $\hfil\displaystyle{#1}\hfil$\crcr}}}

\newcommand{\topreright}[1]{%
  \vbox {\m@th\ialign{##\crcr
  \toprerightfill \crcr
  \noalign{\kern-\p@\nointerlineskip}
  $\hfil\displaystyle{#1}\hfil$\crcr}}}

\def\topreleftfill{%
  $\m@th%
  \dashedleftarrowtip%
  \mkern-1mu%
  \xleaders\hbox{$\mkern2mu\shortbar\mkern-1mu$}\hfill%
  \mkern1mu%
  \shortbar%
  \mkern0.5mu%
$}

\def\toprerightfill{%
  $\m@th%
  \mkern.5mu%
  \shortbar%
  \mkern-1mu%
  \xleaders\hbox{$\mkern2mu\shortbar\mkern-1mu$}\hfill%
  \mkern1mu%
  ‌​\dashedrightarrowtip%
$}

\def\dashedleftarrowtip{%
  \raisebox{\z@}[4.0pt][0.0pt]{$\mathord{\dashedleftarrow}$}}

\def\dashedrightarrowtip{%
  \raisebox{\z@}[4.0pt][0.0pt]{$\mathord{\dashedrightarrow}$}}

\def\shortbar{%
  \smash{\scalebox{0.4}[1.0]{$-$}}}
\makeatother

\newcommand{\orient}{\ssub{\scaleto{\mathscr{O}}{5pt}}!}

\newcommand{\sss}{\ssub{s}!}

\newcommand{\ssvarphi}{\ssub{\varphi}!}

\newcommand{\Wexp}{\ssub{\mathrm{W}}!^{^{\scaleto{\mathrm{exp}}{3pt}}}}

\newcommand{\hatWexp}{\ssub{\widehat{\mathrm{W}}}!}

\newcommand{\bV}{\ssubb[-2pt]{\mathbf{\Upsilon}}!}

\DeclareMathAlphabet\mathbfcal{OMS}{cmsy}{b}{n}

\newcommand{\GlobSp}{\ssub{\mathbfcal{G}}!}

\renewcommand{\hat}{\widehat}

\renewcommand{\ssX}{\ssub{X}!}

\newcommand{\rmU}{\ssub[-1pt]{\mathrm{U}}!}

\newcommand{\Fl}{\mathcal F\hspace{-.05cm}\scaleto{\ell}{7pt}}

\newcommand{\scrC}{\mathscr C}

\newcommand{\bul}{\bullet}

\newcommand{\unde}{\underline e}

\newcommand{\ssK}{\ssub{K}!}
\newcommand{\ssZ}{\ssub{Z}!}

\makeatletter
\newcommand{\extp}{\@ifnextchar^\@extp{\@extp^{\,}}}
\def\@extp^#1{\mathop{\bigwedge\nolimits^{\!#1}}}
\makeatother

\newcommand{\FC}{\ssub[-.5pt]{\mathcal FC}!}

\renewcommand{\deg}{\mathrm{deg}}

\newcommand{\varGr}{\mathbf{Gr}}

\newcommand{\varT}{\mathbf T}
\newcommand{\grass}{\mathrm{Grass}}
\newcommand{\B}{\ssub{\mathbf B}!}
\newcommand{\Br}{\ssub{{\mathcal B}\hspace{-.02cm}\mathit{r}}!}
\newcommand{\PSL}{{\mathcal P}{\mathcal S}{}\mathcal L}
\newcommand{\PS}{\ssub{{\scaleto{{\mathfrak S}}{7pt}}}!}
\newcommand{\slsfunc}{\scaleto{{\mathfrak S}}{7pt}}

\newcommand{\altcone}{\ssub[-1pt]{\scaleto{\cC}{7pt}}!}
\newcommand{\intaltcone}{\ssub[-1pt]{{\ring{\scaleto{\cC}{7pt}}}}!}

\newcommand{\cone}{\ssub[-1pt]{{\sigma}}!}
\newcommand{\intcone}{\ssub[-1pt]{{\ring\sigma}}!}

\newcommand{\subface}{\ssub{\prec}!}
\newcommand{\subfaceq}{\ssub{\preceq}!}
\newcommand{\supface}{\ssub[-1pt]{\succ}!}
\newcommand{\supfaceq}{\ssub{\succeq}!}

\newcommand{\upmin}{\mathrm{UpMin}}
\newcommand{\downsum}{\mathrm{DownSum}}

\title{Limit canonical series}

\author{Omid Amini}
\author{Eduardo Esteves}
\author{Eduardo Garcez}

\date{\today}
\address{CNRS - CMLS, \'Ecole Polytechnique, Palaiseau, France}
\email{omid.amini@polytechnique.edu}

\address{Instituto Nacional de Matem\'atica Pura e Aplicada, Rio de Janeiro, Brazil}
\email{esteves@impa.br}

\address{Universidade Estadual do Cear\'a, Aracati, Brazil}
\email{jemgarcez@gmail.com}

\begin{document}

\begin{abstract} We describe the limits of canonical series along families of curves degenerating to a nodal curve which is general for its topology, in the weak sense that the branches over nodes on each of its components are in general position. We define a fan structure on the space of edge lengths on the dual graph of the limit curve, and construct a projective variety parameterizing the limits, organized in strata associated to the cones of this fan. This extends 
to all topologies the works by Eisenbud--Harris (Invent. Math. 87: 496$-$515, 1987) on curves of compact type and Esteves--Medeiros (Invent. Math. 149: 267$-$338, 2002) on two-component curves.
\end{abstract}

\maketitle

\setcounter{tocdepth}{1}
\tableofcontents

\section{Introduction}

\subsection{Overview} Let $X$ be a stable curve of genus $g$ on the boundary of the Deligne--Mumford compactification $\mgbar_{g}$ of the moduli of smooth curves $\mg_g$ over an algebraically closed field $\k$ of characteristic zero. In this paper, we address the following problem. 

\begin{problem}\label{prob:main} Describe the limits of the spaces of Abelian differentials 
along any family of smooth curves degenerating to $X$. 
\end{problem}  

Denote by $G=(V,E)$ the dual graph of $X$, with vertices $v\in V$ in bijection with the components $\varC_{v}$ of the normalization of $X$, and edges $e\in E$ in bijection with the nodes $\ssp^e$ of $X$. We say that $X$ is \emph{general in its topology}, in the weak sense, if the collection of branches on $\varC_v$ over nodes of $X$ is general for each vertex $v\in V$. 

Our main results, Theorems~\ref{thm:characterization-FC}~and~\ref{thm:variety-lcs}, give a complete solution to the problem stated above for a nodal curve $X$ which is general in its topology, the first theorem describing and the second parameterizing the limits.
 
In a nutshell, the solution is given as follows. For each vertex $v\in V$, denote by $\Omega!_v$ the space of meromorphic differentials of $\varC_v$, and put $\Omega\coloneqq\bigoplus\Omega!_v$. A degeneration of smooth curves to $X$ gives rise, via tropicalization, to an edge length function $\ell \colon E\to \R_{>0}$ on the dual graph of $X$. This function corresponds to the singularity degrees of the total space of the degeneration at the nodes $\ssp^e$. Additionally, we obtain a collection of $g$-dimensional subspaces $\ssW_h\subset\Omega$ indexed by real-valued functions $h\colon V\to \R$,  arising as limits of spaces of Abelian differentials along the degeneration; see Section~\ref{sec:tropicalization}. These limit spaces exhibit a finite generation property: only finitely many of them are needed to determine all the others. Geometrically, there is a tiling of the standard simplex $\Delta!_g\coloneqq\{q\in\R^V_{\geq 0}\,|\,\sum q(v)=g\}$, each $\ssW_h$ contributing with a tile, and the spaces associated to the full-dimensional tiles generate the others. We say these spaces form the \emph{fundamental collection of limit canonical series associated to the degeneration}. 

 We consider all degenerations to $X$ with a given edge length function $\ell\colon E\to\R_{>0}$, and give a complete description of the resulting fundamental collections of limit canonical series; see Theorems~\ref{step2sing}~and~\ref{thm:characterization-FC}. Moreover, as $\ell$ varies, we show that the answer does not vary for $\ell$ within the relative interior of certain rational polyhedral cones in $\R^E_{\geq 0}$, that we describe explicitly; see Section~\ref{sec:fan-structure}. 

We prove that these cones form a rational fan $\Sigma$ on the space of edge lenghts; see Section~\ref{sec:fan-structure-fc}. Furthermore, we show that the set of orbits $\Gm^V\ssW_h$, for $\ssW_h$ in the fundamental collection associated to a degeneration to $X$, under the componentwise action of the algebraic torus $\Gm^V$ on $\Omega$, corresponds to a point on a projective variety $\varV=\bigcup_{\sigma\in \Sigma}\varV_{\sigma}$, stratified according to $\Sigma$; see Theorem~\ref{thm:variety-lcs}. This projective variety $\varV$, which we call \emph{the variety of limit canonical series}, is our solution to Problem~\ref{prob:main}. The constructions are all effective and give an explicit description of the limits. 

\vskip0.4cm

Limits of linear series on curves of compact type were studied by Eisenbud and Harris in the 80's \cite{EH86} as a tool for understanding the moduli of curves $\mg_g$. Of special importance is the case of canonical series. Eisenbud and Harris \cite[Thm.~2.4]{EH87a} were able to describe limits of canonical series on curves of compact type which are general in their topology. Going beyond compact type, there is initial work by Coppens and Gatto \cite{CG} for two-component curves, and a more in-depth study by the second author and Medeiros \cite{EM} in the case the two-component curve is general in its topology. The techniques in \cite{EM} were applied in \cite{ES07} to describe limits of canonical series along families with constant edge length function degenerating to curves with several components, but with the unexpected condition that each two components intersect. As in \cite{EH87a}, the limit curve was assumed to be general in its topology. 

The work in \cite{EM} was specially thorough in the sense that all possible degenerations to a given two-component curve were considered, and a variety parameterizing all limits of canonical series was constructed and described in terms of a polyhedral decomposition of the space of singularity degrees. It is with this thoroughness that we are able in the present paper to extend the work in \cite{EM} to all nodal curves $X$, with an arbitrary number of components, meeting in whatever ways, as long as the branches over nodes on each component of $X$ are in general position. Moreover, it is worth noting that, even for two-component curves $X$, we go beyond what was done in \cite{EM}, as we parameterize all limits of canonical series, not just those which are \emph{aspects} with focus on the components of $X$, following a terminology introduced in \cite{EH86} and employed in \cite{EM, ES07}.

\vskip0.4cm

There are two recent developments that play a key role in our approach. In reverse chronological order:

\vskip0.2cm

$\bullet$ The first is the polyhedral approach developed in \cite{AE-modular-polytopes} to address the question of describing limits of linear series. Making use of this approach in the present paper,  we associate to each degeneration to a nodal curve $X$ of genus $g$, a certain tiling of $\Delta!_g$ by base polytopes of integer polymatroids. The full-dimensional tiles in the tiling are in correspondence with the limit canonical series $\ssW_h$ in the fundamental collection. They generate all the other limits by a finite generation theorem proved in \cite{AE-modular-polytopes}. 

As the degeneration varies, so does the corresponding tiling. However, all such tilings are coarsenings of a fundamental tiling of $\Delta!_g$ by what we call \emph{brick polytopes} or simply \emph{bricks}; see Section~\ref{sec:bricks-fans}. Thus, rather than parameterizing the limit canonical series by the multidegrees of their divisors, as was done in previous works~\cite{EH87a, EM}, we parameterize them by the bricks their associated polytopes contain.

This new approach allows us to associate to each family a multiset of fixed size, consisting of the elements of the fundamental collection counted with multiplicities, each element being counted precisely according to the number of bricks its associated polytope contains. This novel perspective proves to be highly effective in organizing the arising data when considering all degenerations to $X$ at the same time, across all possible edge length functions. 

Indeed, once the characterization of the fundamental collection is established in the paper (Theorems~\ref{step2sing}~and~\ref{thm:characterization-FC}), and the fan structure on the space of edge lengths is constructed (Section~\ref{sec:fan-structure}), the approach outlined above leads to the definition of the variety $\varV$ of limit canonical series as a subvariety of a second variety expressed as a product, indexed by the bricks, of certain torus quotients of Grassmannians. These are GIT quotients, the condition that the base polytope of the subspace contains the brick being a new interpretation of a GIT stability condition. This interpretation is formalized in Theorem~\ref{GIT-projective} in Appendix~\ref{sec:torusgrass}.

 In~\cite[\S 9.2]{AE-modular-polytopes}, we established a connection between our polyhedral approach to the study of linear series and the theory of Chow quotients of Grassmannians, as initiated by Kapranov~\cite{Kap-chow, GKZ} and Lafforgue~\cite{Laf03}, and further developed by Giansiracusa and Wu~\cite{GW22}. In contrast, the construction of the variety $\varV$ is not tied to a Chow quotient but rather to a Mumford quotient.  This construction shares some common features with the work by Thaddeus \cite{Tad99} on complete collineations, although he worked in a distinct and more restricted setting. For further details, see Section~\ref{sec-var-lim-can-series} and Appendix~\ref{sec:torusgrass}.
\vskip0.2cm

$\bullet$ The second development is recent work by Bainbridge, Chen, Gendron, Grushevsky and M\"oller \cite{BCGGM18}, who study degenerations of single Abelian differentials (one-dimensional subspaces of the $g$-dimensional space of Abelian differentials) on families of Riemann surfaces, and construct a compactification of their moduli. This is done by considering tuples of meromorphic differentials on the components of stable curves $X$ that satisfy certain pole, zero, and residue conditions. Most notably, they discovered the \emph{global residue conditions} by an application of Stokes Formula; see Section~\ref{sec:residue-spaces-and-conditions}. The pole and zero conditions, which were previously known, impose constraints on the limit canonical series $\ssW_h$ in the fundamental collection for each degeneration to a stable curve $X$. These constraints are maximal when the branches on the components of $X$ over the nodes are in general position, but turn out to be not enough to determine the $\ssW_h$. After a careful study of the residue conditions in our companion work \cite{AEG-residue-polytope}, we came to expect that adding the global residue conditions under the general position hypothesis could be the missing piece of puzzle to handle Problem~\ref{prob:main} in this level of generality. This led us on a long journey to establish what are now our Theorems~\ref{step2sing}~and~\ref{thm:characterization-FC}. 

 We observe that in \cite{BCGGM18}, the authors work one differential at a time, the goal being to find a meaningful compactification of the projective Hodge bundle over $\mg_g$. 
In particular, the approach does not distinguish whether two limit tuples of differentials on a stable curve $X$ arise from the same degeneration to $X$. It would be interesting to investigate linear series of canonical divisors of intermediate rank, thereby bridging the works done in \emph{loc.~cit.}~and here.

\smallskip

Beyond the works and novel ideas mentioned above, this paper introduces two key innovations that are essential for establishing our results, enabling us to go substantially beyond what was known before.  

\smallskip

$\bullet$ The first is the use of polymatroids in analyzing dimension counts, a critical step in proving our characterization theorems. Polymatroids provide a geometric framework for organizing the data, resulting in an elegant and comprehensible setup, which we believe will be of independent interest for further applications of combinatorial and tropical techniques in degeneration problems. 
Without them, handling the intricate combinatorial arguments would have been either unfeasible or extraordinarily challenging. 

\smallskip

$\bullet$ The second is the canonical fan defined on the space of edge lengths of the dual graph of a stable curve. We associate a fan to each brick, and deduce the canonical fan as a common refinement of these fans. From the perspective of tropical geometry, each space of edge lengths on a dual graph corresponds to a stratum in the moduli space of tropical curves of the given genus; see~\cite{ACP}. By combining the fan structures on the various strata associated with different dual graphs, we can construct a new fan structure on the tropical moduli space. Our results show that the structure of the variety of limit canonical series remains unchanged for all edge lengths lying in the interior of a given cone. We conjecture that this new fan structure on the tropical moduli space represents the tropicalization of a novel compactification of the moduli space of curves, in the spirit of \emph{loc.~cit.}, over which the varieties of limit canonical series reside. A quick analysis suggests that this compactification will be a modification of $\mgbar_{g}$ with many more divisors, described in terms of the fan structure. This connects well with the broader open question posed by Eisenbud and Harris in \cite{EHBull} of extending their theory in \cite{EH86} to all stable curves. To the best of our knowledge, this is the first instance of a fan structure on the space of edge lengths of a dual graph in relation to this open problem. The fan we obtain seems to be a big refinement of the one by Abreu and Pacini \cite{AP}, who studied a tropical analogue of Mainò's moduli space of enriched curves~\cite{Maino}. 

\vskip0.2cm
 
 Finally, we note that one of the main goals in \cite{EM, ES07, EH87a} was to describe limits of Weierstrass points. The description in the full generality considered in our setting is  obtained in Section~\ref{sec:applications}, where other applications are suggested.

\vskip0.2cm

In the rest of this introduction, we present our main constructions and theorems.  We use the framework of tropical geometry, which proves particularly convenient in our setting for handling all degenerations to a nodal curve with variable edge length functions. We refer to the survey papers~\cite{BJ} and~\cite{JP} and the references therein for an introduction to the topic and a sample of results.

\subsection{Submodular functions} They appear throughout this paper. For each finite set $V$, let $\ssub{2}!^V$ be the collection of subsets of $V$. A function $\varphi \colon \ssub{2}!^V \to \R$ with $\varphi(\emptyset)=0$ is called \emph{submodular} if for each two subsets $\ssI_1,\ssI_2\subseteq V$, we have 
\[
\varphi(\ssI_1)+\varphi(\ssI_2)\geq \varphi(\ssI_1\cup \ssI_2)+ \varphi(\ssI_1\cap \ssI_2).
\]
Denote by $\R^V$ the space of real-valued functions $q\colon V \to \R$. For each $I\subseteq V$ and $q\in\R^V$, set $q(I)\coloneqq\sum_{v\in I}q(v)$. To each submodular function $\varphi$, we associate 
 \[
\P_{\varphi}\coloneqq \left\{q\in \R^V\,\st\,q(I)\leq \varphi(I)\text{ for each } I\subseteq V,\text{ with equality if }I=V\right\},
\]
the \emph{base polytope} of the polymatroid on the ground set $V$ defined by $\varphi$. It is easy to see that $\P_{\varphi}$ has maximum dimension $|V|-1$ if and only if $\varphi$ is \emph{simple}, that is,  
$\varphi(I) + \varphi(\ssI^c) >\varphi(V)$ for every proper nonempty subset $I\subset V$, where $\ssI^c\coloneqq V\setminus I$ (see \cite[Prop.~2.7]{AE-modular-polytopes}). 

Given a field $\k$ and a collection $(\VS_v)_{v\in V}$ of vector spaces over $\k$, each finite-dimensional subspace $\ssW \subseteq \bigoplus_{u\in V} \VS_u$ gives rise to a submodular function $\nu!^*_\ssW$, satisfying $\nu!^*_\ssW(I)=\dim_{\k} \proj{I}(\ssW)$ for each $I\subseteq V$, where $\proj{I}\colon \bigoplus_{u\in V} \VS_u\to \bigoplus_{u\in I} \VS_u$ is the projection map. Submodular functions of this form arise naturally from tropicalization of linear series~\cite{AE-modular-polytopes}.
 
\subsection{Tropicalization and limit canonical series} \label{sec:tropicalization-introduction} Let $\K$ be an algebraically closed field which is complete for a nontrivial non-Archimedean valuation $\valuation$. Denote  by $\varR$ its valuation ring, by $\varm$ the maximal ideal of $\varR$, by $\k$ the residue field, and by $\valgroup$ the value group. Let $\varX$ be a smooth proper curve of genus $g$ over $\K$, and denote by $\varX^{\an}$ the Berkovich analytification of $\varX$. The points on $\varX^{\an}$ are in bijection with the union of the closed points on $\varX$ and the valuations on the function field $\K(\varX)$ extending $\valuation$.

Let $\Gamma$ be a skeleton of $\varX^{\an}$ associated to a semistable vertex set $V$. It is a metric graph. Let $G=(V,E)$ be the underlying graph and $\ell \colon E \to \R_{>0}$ be the edge length function. When the curve $\varX$ is defined over a discrete valued subfield $\K_0\subseteq\K$, there exists a stable model $\fX\to \mathrm{Spec}(\varR_1)$ for $\varX$ over a finite extension $\K_1$ of $\K_0$, where $\varR_1$ is the valuation ring of $\K_1$, and we can take for $G$ the dual graph of the stable reduction $X$ of $\varX$, and for $\ell$ the function assigning to each $e\in E$ the singularity degree of $\fX$ at the corresponding node, normalized by the degree of the field extension $\K_1/\K_0$. The metric graph associated to $(G, \ell)$ is well-defined and naturally embeds in $\varX^{\an}$ as a skeleton.

In any case, for each vertex $v$ of $G$, the residue field of the corresponding valuation on $\K(\varX)$, denoted $\k(v)$, has transcendental degree one over $\k$. Let $\varC_v$ be the smooth proper curve over $\k$ with field of functions $\k(v)$. Each edge $e$ of $G$ gives a point on $\varC_u$ and another on $\varC_v$, where $u$ and $v$ are the vertices connected by $e$. Identifying these points for each edge produces a nodal curve $X$ with dual graph $G$. We say $(X, \Gamma)$ \emph{is} (\emph{obtained by}) \emph{the tropicalization of $\varX$}. 

Let $\varH$ be the space of Abelian differentials of $\varX$, also refered to as the \emph{canonical series of} $\varX$. For each $v\in V$, denote by $\Omega!_v$ the space of meromorphic differentials of $\varC_v$, and put $\Omega \coloneqq\bigoplus_{v\in V}\Omega!_v$. To each function $h \colon V \to \valgroup$ we associate a $g$-dimensional subspace $\ssW_h \subset \Omega$ defined by reduction of $\varH$ relative to $h$, and called the \emph{limit canonical series} associated to $h$; see Section~\ref{sec:tropicalization}. If $\varX$ is defined over a discrete valued field $\K_0$, and admits a stable model $\fX \to \mathrm{Spec}(\varR_0)$ over the valuation ring $\varR_0$ of $\K_0$ with stable reduction $X$, and $h$ is integer-valued, the subspace $\ssW_h \subset \Omega$ can be identified with the restriction to $X$ of the space of global sections of 
\[
\omega!_{\fX/\mathrm{Spec}(\varR_0)}\big(-\sum_{v\in V}h(v)\ssX_v\big),
\]
the relative canonical sheaf of the stable model twisted by the vertical divisor $-\sum h(v)\ssX_v$ on $\fX$. Here, $\ssX_v$ is the component of $X$ corresponding to $v$ (whose normalization is $\varC_v$).   

For each function $h\colon V\to\valgroup$, we denote by  $\nu!_h^* \colon \ssub{2}!^V \to \Z$ the submodular function associated to $\ssW_h \subset \Omega = \bigoplus_{v\in V} \Omega!_v$ and by $\P_{\nu!^*_h}$ the corresponding base polytope. We say that the subspace $\ssW_h\subseteq\Omega$ is \emph{simple} if $\nu!_h^*$ is simple. We denote by $\FC(\varX)$ the collection of simple subspaces $\ssW_h\subset\Omega$, and call it the  \emph{fundamental collection of limit canonical series} on $X$ associated to $\varX$. The following theorem implies that $\FC(\varX)$ is finite, and~\cite[Thm.~9.6]{AE-modular-polytopes} implies that the spaces in $\FC(\varX)$ generate all the other limit canonical series $\ssW_h$.

\newtheorem*{thm:W-tiling}{\cref{thm:W-tiling}}
\begin{thm:W-tiling} \emph{The elements $\ssW_h$ of $\FC(\varX)$ are precisely those limit canonical series whose polytope $\P_{\nu!_h^*}$ is full-dimensional. They and their faces form a tiling of the standard simplex $\Delta!_g$ in $\R^V$.} 
\end{thm:W-tiling}

\subsection{Level structures and their residue spaces}\label{sec:ls-rs}  Each element of $\ssW_h$ is a tuple of meromorphic differentials on $\varC_v$, for $v\in V$. It is subject to certain residue conditions:

\newtheorem*{thm:ImWh}{\cref{thm:ImWh}}
\begin{thm:ImWh} \emph{Assume $\car\k=0$. Then, for each $h\colon V\to\valgroup$, the residue map $\Res\colon\Omega\to\k^{\E}$ takes $\ssW_h$ into the $\pi$-residue subspace $\GlobSp_\pi\subseteq\k^{\E}$, where $\sspi$ is the level structure on $G$ given by $h$.}
\end{thm:ImWh}

We explain the statement.  First, $\E$ is the set of arrows of $G$. Each edge of $G$ can be oriented in two ways, giving rise to two arrows. Each $a\in\E$ is thus supported on a unique $e\in E$, the reverse arrow denoted $\bar a$. We write $a=uv$ if $a=u\to v$ connects $u$ to $v$, in which case we call $u$ the \emph{tail} and $v$ the \emph{head} of $a$, and write $e=\{u,v\}$. For $v\in V$, let $\E_v\subseteq\E$ be the set of arrows with tail $v$.

Second, to each $a=uv\in\E$ we associate the branch $\ssp^a$ on $\varC_u$ above the node $p^e$ corresponding to the underlying edge $e=\{u,v\} \in E$. Denote by $\res_a \colon \Omega!_u \to \k$ the residue map at $\ssp^a$. Collecting the maps $\res_a$ for $a\in \E$ yields a \emph{residue map}:
\begin{equation}
\Res \colon \Omega \to \k^{\E}.
\end{equation}
Third, a partition of $V$ is a collection of  pairwise disjoint nonempty subsets whose union is $V$. A \emph{level structure} on $G$ is the data of an ordered partition of $V$, i.e., a partition of $V$ equipped with a total order on its elements. Equivalently, a level structure is the data $\pi$ of a set $\ssV_\pi$ endowed with a total order $\subface_{\pi}$ and a surjection $h=\ssh_\pi\colon V \to \ssV_\pi$. For each $u,v\in V$, abusing notation, we write $u\subface_\pi v$ if $h(u) \subface_\pi h(v)$. We refer to $(G, \pi)$ as a \emph{level graph}. Each function $h\colon V \to \R$ gives naturally rise to a level structure on $G$ which we denote by $\sspi_h$. 

\begin{figure}[ht]
    \centering
\begin{tikzpicture}[scale=.8pt]

\draw[line width=0.4mm] (5,-2) -- (6,0);
\draw[line width=0.4mm] (7,-2) -- (6,0);
\draw[line width=0.4mm] (7,-2) -- (7,0);
\draw[line width=0.4mm] (5,-2) -- (7,0);
\draw[line width=0.4mm] (8.5,-2) -- (7,0);
\draw[line width=0.4mm] (7,-2) -- (8.5,-2);
\draw[line width=0.4mm] (8.5,-2) -- (7,-2);

\draw[line width=0.1mm] (6,0) circle (0.8mm);
\filldraw[aqua] (6,0) circle (0.7mm);
\draw[line width=0.1mm] (7,0) circle (0.8mm);
\filldraw[aqua] (7,0) circle (0.7mm);
\draw[line width=0.1mm] (5,-2) circle (0.8mm);
\filldraw[aqua] (5,-2) circle (0.7mm);
\draw[line width=0.1mm] (7,-2) circle (0.8mm);
\filldraw[aqua] (7,-2) circle (0.7mm);
\draw[line width=0.1mm] (8.5,-2) circle (0.8mm);
\filldraw[aqua] (8.5,-2) circle (0.7mm);

\draw[black] (6,0.25) node{$\ssu_4$};
\draw[black] (7,0.25) node{$\ssu_5$};
\draw[black] (5,-2.3) node{$\ssu_1$};
\draw[black] (7,-2.3) node{$\ssu_2$};
\draw[black] (8.5,-2.3) node{$\ssu_3$};

\end{tikzpicture} 
\caption{A level graph $(G,\pi)$ with two levels:  $\ssV_\pi=\{1,2\}$ with order $1<2$. The arrows $\ssu_1\ssu_4, \ssu_1\ssu_5, \ssu_2\ssu_4, \ssu_2\ssu_5, \ssu_3\ssu_5$ are upward. The edge $\{\ssu_2,\ssu_3\}$ is horizontal.}
\label{fig:digraph-0}
\end{figure}
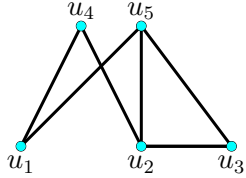

It remains to explain the subspace $\GlobSp_\pi\subseteq\k^{\E}$. Its introduction and the  terminology are motivated by the study of residues of meromorphic differentials in \cite{BCGGM18} and~\cite{TT22}. Given a level structure $\pi$ on $G$, we call $e=\{u,v\}\in E$ \emph{vertical} if $\ssh_\pi(u)\neq \ssh_\pi(v)$ and \emph{horizontal} otherwise. An arrow $a=u\to v\in\E$ is called \emph{compatible with $\pi$} or \emph{upward} if $u \supface_\pi v$. (Visually, both $\to$ and $\supface_\pi$ point to $v$.) Its reverse arrow $\bar a$ is called \emph{downward}. All the other arrows are called \emph{horizontal}. We denote by $\ssA_{\pi}$ the set of upward arrows and by $\barA_\pi$ the set of downward arrows. Figure~\ref{fig:digraph-0} shows a level graph with two levels. 

Put $\bV\coloneqq\k^{\E}$. For each $\psi \in \bV$ and $a\in \E$, denote by $\psi_a$ the $a$-coordinate of $\psi$. Define the \emph{$\pi$-residue space} of $G$, denoted $\GlobSp_\pi$, as the subspace of $\bV$ consisting of those functions $\psi$ which verify the residue conditions~\ref{cond:R1} through \ref{cond:R4}, below:

\begin{enumerate}[label=(R\arabic*)]
    \item \emph{Vanishing along downward arrows:}\label{cond:R1}
    \[\sspsi_a=0  \qquad \text{for every arrow $a \in\barA_{\pi}$}.\]

    \item \emph{Local residue conditions:}\label{cond:R2}
    \[\sum_{a\in\E_v}\sspsi_a=0 \qquad \text{for every $v\in V$}.\]

    \item \emph{Rosenlicht conditions:}\label{cond:R3}
    \[
    \sspsi_a+\sspsi_{\bar a}=0 \qquad \text{ for each horizontal arrow $a$ for $\pi$}.  
    \]

    \item \emph{Global residue conditions:}\label{cond:R4}
    \[
    \sum_{a\in \ssA^{\Xi}_{\pi, n}}\psi_a=0
    \]
\end{enumerate}
for each level $n\in \ssV_\pi$ and each connected component $\Xi$ of the subgraph $G[\ssV_{h \subface_\pi n}]$ of $G$ induced on the set of vertices $\ssV_{h \subface_\pi n} \subseteq V$ of level smaller than $n$, with $\A_{\pi,n}^{\Xi}\subseteq\A_{\pi}$ denoting the set of upward arrows with tail of level $n$ and head in $\Xi$. 

In Figure~\ref{fig:digraph-0}, taking $n=2$, the set $\ssV_{h \subface_\pi 2}$ consists of the vertices $\ssu_4$ and $\ssu_5$, with the induced graph $G[\ssV_{h \subface_\pi 2}]$ having no edges. Taking $\Xi$ to be the singleton with vertex set $\{\ssu_5\}$, the set $\ssA^\Xi_{\pi,2}$ has three arrows: $\ssu_1\ssu_5, \ssu_2\ssu_5, \ssu_3\ssu_5$. The corresponding global residue condition is the equation $\psi_{\ssu_1\ssu_5}+\psi_{\ssu_2\ssu_5}+\psi_{\ssu_3\ssu_5}=0$.

The main content of the theorem is \ref{cond:R4}.  Tropicalization of $\alpha \in \varH$ gives rise to a piecewise integral affine map $\trop(\varalpha) \colon \Gamma \to \R$ and reduction of $\alpha$ relative to the tropicalization gives an element $\wt\alpha=(\wt\alpha_v)_{v\in V} \in \Omega$. By~\cite{BCGGM18} for Riemann surfaces, and~\cite{TT22} in our more general setting, $\Res(\wt\alpha)$ verifies \ref{cond:R4}.

Given $h\colon V\to\valgroup$, each $\beta=(\beta_v)_{v\in V}\in \ssW_h$ comes from an $\alpha\in\varH$ with $\trop(\alpha)\rest{V}\geq h$, with $\beta_v=\wt\alpha_v$ if $\trop(\alpha)(v)=h(v)$, and $\beta_v=0$ otherwise.

If there is $\alpha\in\varH$ with $\trop(\alpha)\rest{V}= h$, then $\wt\alpha \in \ssW_h$, from which we deduce \ref{cond:R4} for each element of $\Res(\ssW_h)$. If there is no such $\alpha$, we use \cite[Thm.~9.6]{AE-modular-polytopes} to obtain the result from the previous special case (see Proposition~\ref{prop:compatibility-residue}).

\subsection{Residue polytopes}\label{sec:rp-intro} For the trivial ordered partition of $V$ consisting of a single level, Conditions \ref{cond:R1} and \ref{cond:R4} are vacuum, and the corresponding residue space is the first homology group $\ssH_1(G, \k)$, which is of dimension $g(G) = |E|-|V|+1$. We prove in~\cite[Thm.~1.1]{AEG-residue-polytope} that this holds more generally: $\dim_{\k}\GlobSp_\pi=g(G)$ for each ordered partition $\pi$ of $V$. Moreover, we show that  $\GlobSp_{\pi}$ is obtained from $\ssH_1(G, \k)$ by a process we call \emph{splitting}, and it arises as a certain limit of  $\ssH_1(G, \k)$ in the Grassmannian of $g(G)$-dimensional subspaces of $\bV$, see~\cite[Thm.~6.1]{AEG-residue-polytope}. This gives a new take on the global residue conditions.

For each level structure $\pi$ on $G$, denote by 
\[
\gamma!_\pi \colon \ssub{2}!^V \to \Z, \qquad \gamma!_\pi(I) \coloneqq \dim_{\k}\proj{I}\left(\GlobSp_\pi\right)  \qquad \textrm{for each $I\subseteq V$},
\]
the submodular function associated to the $\pi$-residue space $\GlobSp_\pi \subseteq \bV$, for the decomposition $\bV = \bigoplus_{v\in V} \bV_v$, where $\bV_v\coloneqq\k^{\E_v}$ for each $v\in V$. By~\cite[Thm.~1.2]{AEG-residue-polytope}, $\P_{\gamma_\pi}$ are the faces of the polytope associated to the trivial ordered partition of $V$ that we call the \emph{residue polytope of $G$}. The residue polytope of the complete graph on four vertices is depicted in azure in  Figure~\ref{fig:k4}. 

\subsection{Characterization Theorem} We explained tropicalization and the way it yields a nodal curve and certain spaces of meromorphic differentials on it. We change the perspective now. Fix a nodal (connected projective) curve $X$ of arithmetic genus $g>0$ over an algebraically closed field $\k$ of characteristic zero and let $(G,E)$ be its dual graph. Consider smooth proper curves $\varX$ over non-Archimedean valued fields $\K$ with residue field $\k$ that tropicalize to $(X, \Gamma)$, for $\Gamma$ a metric graph with underlying graph $G$. For each $v\in V$, let $\varC_v$ be the normalization of the component of $X$ associated to $v$, and $\Omega!_v$ be the space of meromorphic differentials on $\varC_v$. Put $\Omega\coloneqq\bigoplus_{v\in V}\Omega!_v$. Our aim is to characterize the collections of $g$-dimensional subspaces of $\Omega$ that are of the form $\FC(\varX)$ for some $\varX$, see Theorem~\ref{thm:characterization-FC} below.

A \emph{slope function on $G$} is a function $s\colon\E\to\Z$ which verifies
\[
  s(a)+s(\bar a) \in \{-1, 0\}\quad\text{for each arrow }a\in\E.
\]
Given a slope function $s$, we denote by $\ssA_s$ the set of arrows $a\in \E$ with $s(\bar a)<0$ (equivalently, either $s(a)>0$, or $s(a)=0$ and $s(a)+s(\bar a)=-1$), and call the arrows in $\ssA_s$ \emph{upward} for $s$. The reverse arrows are called \emph{downward} for $s$; their set is denoted $\barA_s$.  If $s(a)+s(\bar a)=0$, we say that $a$ and the underlying edge are \emph{integer} for $s$. Denote by $\ssE_s^{\rmint}$ (resp.~$\ssA_s^{\rmint}$) the set of integer vertical edges (resp.~integer upward arrows).

Each function $h\colon V \to \R$ yields a slope function $\sss_h\colon\E\to\Z$ defined by
\[
\sss_h(a)\coloneqq\slztwist{\ell}h(a)\coloneqq \Big\lfloor \frac{h(u)-h(v)}{\ell_a}\Big\rfloor \qquad \text{ for each }a=uv\in\E,
\]
where $\ell_a$ is the length of the edge underlying $a$. We have $\ssA_{\sss_h}=\ssA_{\sspi_h}$,  for $\sspi_h$ the level structure induced by $h$.

A pair $(s,\pi)$ consisting of a slope function $s$ and a level structure $\pi$ on $G$ is called a \emph{slope-level pair} provided that $\ssA_s=\ssA_{\pi}$. Given $h\colon E \to \R$, the pair $(\sss_h,\sspi_h)$ is thus a slope-level pair.

To each slope function $s$ we associate the function $\zeta!_s \colon \ssub{2}!^V \to \Z$ defined by
\[
\zeta!_s \colon \ssub{2}!^V \to \Z, \qquad \zeta!_s(I) = |\ssA_s(\ssI^c, \ssI)|+ \sum_{a\in\E(I,\ssI^c)}s(a) \qquad \text{for all }\, I\subseteq V,
\]
where $\E(I,\ssI^c)$ is the set of arrows with tail in $I$ and head in $\ssI^c$, and where $\ssA_s(\ssI,\ssI^c)= \E(I,\ssI^c) \cap \ssA_s$. Our Proposition~\ref{prop:zetah} states that $\zeta!_s$ is submodular.

To the slope-level pair $(s,\pi)$ on $G$, we associate the submodular function $\eta_{s,\pi}\colon 2^V\to\Z$ defined by 
\[
\eta!_{s,\pi}\coloneqq\gamma!_{\pi}+\g+\zeta!_s
\]
where, as before, $\gamma!_{\pi}$ is the submodular function of $\GlobSp_{\pi}\subseteq\bV$, and  $\g \colon V \to \Z$ is the \emph{genus function} of $X$ associating to each vertex $v$ the genus $\g(v)$ of $\varC_v$.

In addition, for each slope function $s\colon E\to\Z$ and $v\in V$, we put
\begin{equation}\label{eq:shf-Lh-intro}
\ssL_{s,v}\coloneqq \omega!_{v}\Big(\sum_{\substack{a\in\E_v}} \big(1+s(a)\big) \ssp^a\Big),
\end{equation}
where $\omega!_v$ is the sheaf of Abelian differentials on $\varC_v$, and $p^a$ is the branch on $\varC_v$ over the node of $X$ corresponding to the underlying edge for each arrow $a\in\E_v$. Viewing $\omega!_v$ as a subsheaf of the constant sheaf on $\Omega!_v$, we view $\ssL_{s,v}$ as a sheaf of meromorphic differentials whose divisors of poles and zeros are constrained by the sum of $(1+s(a))p^a$ over $a\in\E_v$.

By a \emph{partial gluing data for the $\ssL_{s,v}$} we mean a collection $\varrho=(\varrho!_a)_{a\in\ssA_s^{\rmint}}$ of isomorphisms 
\[
\varrho!_a\colon\ssL_{s,u}\rest{p^a}\xrightarrow{\cong}\ssL_{s,v}\rest{p^{\bar a}}\quad\text{for }a=uv\in\ssA_s^{\rmint}.
\]
Then, we let $\ssL_s^{\varrho}$ be the subsheaf of the sum on $X$ of the $\ssL_{s,v}$ defined by imposing equality away from the nodes $\ssp^e$ corresponding to integer edges $e$ for $s$~\ref{Exp1}, Rosenlicht condition \ref{cond:R3} at the node $\ssp^e$ for horizontal $e$~\ref{Exp2}, and gluing
of the two sheaves $\ssL_{s,u}$ and $\ssL_{s,v}$ at the points $\ssp^a$ and $\ssp^{\bar a}$ for each $a=uv\in \ssA_s^{\rmint}$ given by $\varrho!_a$~\ref{Exp3}; see Section~\ref{sec:submodular-Wh}.

Given a partial gluing $\varrho$ as above, we set
\[
\Wexp_{s,\pi}(\varrho) \coloneqq\Res^{-1}(\GlobSp_{\pi})\cap H^0(X,\ssL_{s}^{\varrho})\subset\Omega.
\]

Our Theorem~\ref{thm:characterization-FC} gives a complete characterization of $\FC(\varX)$ for the $\varX$ whose tropicalization is $(X,\Gamma)$ for a fixed $\Gamma$, or equivalently, a fixed edge length function $\ell\colon E\to\R$. Moreover, we show lated that $\FC(\varX)$ depends only on the associated set $\PS_{\ell}$ of slope-level pairs, rather than on $\ell$ itself, see Section~\ref{sec:canonica-fan-intro}.

Assume that $\car\k$ is zero and $X$ is general in its topology. Then, we prove that for each $h\colon V\to\valgroup$, the subspace $\ssW_h\subset\Omega$ is simple if and only if the submodular function $\eta!_{\sss_h,\sspi_h}$ is positive and simple. In this case, there is a unique partial gluing data $\varrho!_h=\varrho!_h(\varX)$ for the sheaves $\ssL_{\sss_h,v}$ on $\varC_v$ defining the sheaf $\ssL_{\sss_h}^{\varrho!_h}$ on $X$, such that $\ssW_h\subseteq H^0(X,\ssL_{\sss_h}^{\varrho!_h})$ and $\ssW_h=\Wexp_{\sss_h,\sspi_h}(\varrho!_h)$, see Theorem~\ref{thm:characterization-FC}.

We will explain below how the data of the partial gluing $\varrho!_h=\varrho!_h(\varX)$ can be extracted from $\varX$. Having done so, 
the theorem allows us to find the elements $\ssW_h$ of the fundamental collection by just computing the submodular function $\eta_{\sss_h,\sspi_h}$ and checking whether it is positive and simple. Furthermore, if this happens, then $\ssW_h$ is completely determined by the residue conditions, Theorem~\ref{thm:ImWh}, by pole and zero conditions derived from the inclusion $\ssW_h\subseteq\bigoplus H^0(\varC_v,\ssL_{\sss_h,v})$, and the compatibility condition given by the gluing data $\varrho!_h$. That the elements in $\ssW_h$ satisfy all these conditions is one thing, the main ingredient being Theorem~\ref{thm:ImWh}. That the conditions are enough requires the machinery developed in Section~\ref{genpos} to handle general position conditions, with the main theorems being:
\smallskip

\noindent $\bullet$ the first Realizability Theorem~\ref{thm:global}, which allows us to compute the submodular function of the subspace
\[
\hatWexp_{s,\pi} \coloneqq \Res^{-1}(\GlobSp_{\pi}) \bigcap \left(\bigoplus_{v\in V} H^0(\varC_v,\ssL_{s,v})\right) \subseteq \Omega
\]
for slope-level pairs $(s,\pi)$ such that $\eta!_{s,\pi}$ is positive (see Theorem~\ref{step1sing});

\noindent $\bullet$ the second Realizability Theorem~\ref{thm:eta-gen-2-hyper}, which allows us to handle the gluing conditions as we pass from $\hatWexp_{s,\pi}$ to $\Wexp_{s,\pi}(y)$, and establish that the submodular function of the subspace $\Wexp_{s,\pi}(y)\subseteq\Omega$ is the \emph{UpMin transform} $\upmin(\eta!_{s,\pi})$ of $\eta!_{s,\pi}$, for any partial gluing data $y$ for the $\ssL_{s,v}$ (see Theorem~\ref{step2sing-expected}).

 Once these theorems are established, the equality $\ssW_h=\Wexp_{\sss_h,\sspi_h}(\varrho!_h)$ boils down to the fact that both spaces have dimension $g$.

Submodular functions play a crucial role in our dimension counts, and the UpMin transform introduced in Section~\ref{sec:upmin} is key. For each function $\varphi \colon \ssub{2}!^V \to \R$, the UpMin transform of $\varphi$ is the function 
\[
\upmin(\varphi)(I) =\min_{K\supseteq I} \varphi(K) \qquad \text{for all } I\subseteq V.
\]
It transforms nonnegative submodular functions to nonnegative submodular functions (Proposition~\ref{prop:min-operation}), can be iterated (Proposition~\ref{prop:xi-upmin}), and preserves and detects simpleness (Proposition~\ref{prop:simpleness-min-operation}), fundamental properties we use over and over. To boot, it admits a beautiful geometric meaning, see Proposition~\ref{prop:beautiful-geometric-meaning}.

\vskip0.2cm

It remains to explain what the $(\varrho!_{h})_{h}$ appearing in the theorem are precisely. In order to do this, we fix isomorphisms $\mathcal O_{\varC_v}(\ssp^a)\rest{\ssp^a}\cong\k$ for each $v\in V$ and $a\in\E_v$. Let $\omega!_X$ be the canonical sheaf on $X$, and fix isomorphisms $\omega!_X\rest{\ssp^e}\cong\k$ for each $e\in E$. Using them, we define isomorphisms $\ssL_{s,v}\rest{\ssp^a}\cong\k$ for each slope function $s\colon\E\to\Z$, each $v\in V$ and $a\in\E_v$. By using these isomorphisms, and abusing of notation, we may view any partial gluing data 
as and element of $\prod_{a\in\ssA^{\rmint}_s}\k^{\times}$. The partial gluing $\varrho!_h=\varrho!_h(\varX)$ appearing in Theorem~\ref{thm:characterization-FC} is then given by 
\[
\varrho!_{h,a} = \ssrho_e^{\sss_h(a)} \qquad \qquad  \text{for all } a\in \ssA_s^\rmint
\]
where $e$ is the underlying edge of $a$ and $\ssrho_e$ is the \emph{leading coefficient} of the smoothing of the node $\ssp^e$ in $\varX$, see Appendix~\ref{sec:smoothing-and-gluing}. It follows moreover from Theorem~\ref{thm:appendix-C} and Theorem~\ref{thm:characterization-FC} that the data 
\[
(\varrho!_{h,a})_{h\in\R^V,a\in\ssA_{\sss_h}^{\rmint}}
\]
appearing in the description of the fundamental collection
for $\varX$ varying among all smooth curves over valued fields $\K$ tropicalizing to $(X,\Gamma)$ is the same as 
\[
(\ssrho_a^{\sss_h(a)})_{h\in\R^V,a\in\ssA_{\sss_h}^{\rmint}}
\]
for $\rho$ varying in $\Gm^E(\k) = \ssub{\bigl(\k^\times\bigr)}!^E$. 

\subsection{The canonical fan of an augmented graph} \label{sec:canonica-fan-intro} Let $G=(V,E)$ be a finite connected graph and $\g\colon V\to\Z_{\geq 0}$ a function. We call $(G,\g)$ an \emph{augmented graph}. For each edge length function $\ell\colon E\to\R_{>0}$, we consider the collection $\PS_\ell$ of slope-level pairs $(\sss,\sspi)$ on $G$ such that:
\begin{enumerate}[label=($\roman*$)]
    \item \label{sp:1} the submodular function $\eta!_{s,\pi}=\gamma!_\pi+\g+\zeta!_s$ is positive and simple, and
    \item \label{sp:2} there exists $h\colon V\to\R$ such that $(\sss,\sspi)=(\sss_h,\sspi_h)$.
\end{enumerate}
Our next main result, Theorem~\ref{thm:Sigma}, asserts that those $\ell'\in\R^E_{>0}$ for which $\PS_{\ell'}=\PS_\ell$ form the relative interior of a rational polyhedral cone $\sssigma_{\PS_{\ell}}$ in $\R_{\geq 0}^E$, and that the collection $\Sigma=\Sigma(G, \g)$ of these cones, and their faces appearing in the boundary $\R_{\geq 0}^E\setminus\R_{>0}^E$, is a rational fan on $\R_{\geq0}^E$. We call $\Sigma$ the \emph{canonical fan of $(G, \g)$}.  

We prove this theorem as follows. First, we associate to each slope-level pair $(s,\pi)$ on $G$ the \emph{cone of edge length functions} $\intcone_{s,\pi}$ and its closure  $\cone_{s,\pi}$ in $\R^E$, defined by imposing Condition~\ref{sp:2}:
\begin{align*}
\intcone_{s,\pi} &\coloneqq \Bigl\{\ell \in \R_{>0}^E \,\st\, \exists\,\, h\colon V \to \R \textrm{ with }\sss_h=s \textrm{ and } \sspi_h=\pi \Bigr\}, \qquad \cone_{s,\pi} = \overline{\intcone_{s,\pi}}.
\end{align*}
Our Theorem~\ref{thm:characterization-slope-level-cone} gives an explicit description of $\cone_{s,\pi}$ in terms of the inequalities in~\eqref{eq:slope-level-compatible}. In particular, this shows that $\cone_{s,\pi}$ is rational polyhedral. Furthermore,  Theorem~\ref{thm:face-cone} describes the faces of $\cone_{s,\pi}$ as cones $\cone_{s',\pi'}$ associated to slope-level pairs $(\sss', \sspi')$ obtained from $(s,\pi)$ by a \emph{squashing process}, see Section~\ref{sec:squashing}.

Then, we treat Condition~\ref{sp:1} as follows. For each slope-level pair $(s,\pi)$ on $G$, let $\P_{s,\pi}$ be the polytope associated to the submodular function $\upmin(\eta!_{s,\pi})$. Proposition~\ref{prop:simpleness-min-operation} yields that $\eta!_{s,\pi}$ is positive and simple if and only if $\P_{s,\pi}$ is full-dimensional, thus, if and only if $\P_{s,\pi}$ contains a full-dimensional \emph{brick}.

A brick is a polytope $\B\subseteq\Delta!_{g}$ defined as the set of all points $q \in \Delta!_g$ that verify the inequalities 
\[
\lfloor \beta(S)\rfloor \leq q(S) \leq \lceil \beta(S)\rceil\qquad \qquad \forall \,\, S \subseteq V
\]
for some $\beta\in\Delta!_{g}$. By \cite[Thm.~10.4]{AE-modular-polytopes}, the collection of bricks is a tiling of $\Delta!_{g}$. Furthermore, the base polytope $\P$ of any polymatroid with integral vertices lying in $\Delta!_{g}$ is tiled by the set of bricks included in $\P$. We let $\Br_g$ be the collection of full-dimensional bricks.

A slope-level pair $(s,\pi)$ on $G$ is called \emph{permissible} if there exists an $\ell\in \R_{>0}^E$ such that  $(s,\pi) \in\PS_\ell$, that is, 
\begin{itemize}
\item $\eta!_{s,\pi}$ is positive and simple, and
\item there is an $\ell \in \R_{>0}^E$ and an $h\colon V \to \R$ such that $(s,\pi) = (\sss_h, \sspi_h)$. 
\end{itemize}
For each permissible $(s,\pi)$, the polytope $\P_{s,\pi}$ is full-dimensional, it thus contains a brick in $\Br_g$. Let $\PSL=\PSL(G, \g)$ be the set of permissible slope-level pairs for  $(G, \g)$. We prove that $\PSL$ is finite. We can write 
\[
\PSL=\bigcup_{\B\in\Br_g}\PSL(\B) \quad \text{where} \quad \PSL(\B)\coloneqq\bigl\{(s,\pi)\in\PSL\,\st\, \P_{s,\pi}\supseteq\B\bigr\}.
\]

\newtheorem*{thm:fan-sigmaB}{\cref{thm:fan-sigmaB}}
\begin{thm:fan-sigmaB} \emph{Let $\B \in \Br_g$ be a full-dimensional brick in $\Delta!_g$. The collection of cones $\cone_{s,\pi}$ for $(s,\pi) \in \PSL(\B)$ and all their faces which lie on the boundary $\R_{\geq 0 }^E \setminus \R_{>0}^E$ form a rational fan with support $\R_{\geq 0}^E$.}
\end{thm:fan-sigmaB}
Denote by $\ssSigma_\B$ the rational fan with support $\R_{\geq 0}^E$ produced by this theorem. 

\newtheorem*{thm:Sigma}{\cref{thm:Sigma}}
\begin{thm:Sigma} Notation as above,  $\Sigma$ is a rational fan with support $\R_{\geq 0}^E$. It satisfies
\begin{equation}\label{eq:sigma-sigmaB}
\Sigma =  \Bigl\{\bigcap_{\B \in \Br_g} \sssigma_\B \, \st\, \sssigma_\B \text{ in } \ssSigma_\B \text{ for all } \B\in \Br_g\Bigr\}.
\end{equation}
In other words, $\Sigma=\bigwedge_{\B \in \Br_g} \ssSigma_{\B}$ is the meet of the fans $\ssSigma_{\B}$ associated to the full-dimensional bricks $\B$ in $\Delta!_g$.
\end{thm:Sigma}

\subsection{Variety of limit canonical series} Having understood by Theorem~\ref{thm:characterization-FC} what the set of $\FC(\varX)$ is, for smooth proper curves $\varX$ over valued fields $\K$ tropicalizing to $(X,\Gamma)$ for some $\Gamma$, and having understood by Theorem~\ref{thm:Sigma} how these sets vary with $\ell$, we can put all the $\FC(\varX)$ together in a parameter space.

Again, let $X$ be a nodal curve of arithmetic genus $g>0$ over an algebraically closed field $\k$ of characteristic zero, whose branches over nodes are in general position on the components of its normalization. Let $G=(V,E)$ be its dual graph and $\g$ its genus function. Recall the notation used above. We fix isomorphisms $\mathcal O_{\varC_{v}}(\ssp^a)\rest{\ssp^a}\cong\k$ for each $v\in V$ and $a\in \E_v$, and $\omega!_X\rest{\ssp^e}\cong\k$ for each $e\in E$. They induce isomorphisms $\ssL_{s,v}\rest{p^a}\cong\k$ for each slope function $s\colon\E\to\Z$, and each $v\in V$ and $a\in\E_v$.

Let $\Sigma$ be the canonical fan of $(G, \g)$. As $\PSL$ is finite, there is a finite-dimensional subspace $\scrU_v \subset \Omega!_v$ for each vertex $v\in V$ such that $H^0(\varC_v, \ssL_{s,v}) \subseteq \scrU_v$ for each $(s,\pi)\in \PSL$ and each $v\in V$. Put $\scrU\coloneqq\bigoplus\scrU_v$. Then,
\[
\Wexp_{s,\pi}(y) \subseteq H^0(X,\ssL_s^{y})\subseteq\bigoplus_{v\in V} H^0(\varC_v, \ssL_{s,v}) \subseteq \bigoplus_{v\in V}\scrU_v=\scrU
\]
for each $(s,\pi)\in \PSL$ and each partial gluing data $y$ for the $\ssL_{s,v}$.

Consider the Grassmannian $\varGr\coloneqq\grass\big(g,\scrU\big)$ of $g$-dimensional subspaces of $\scrU$. For each $\B \in \Br_g$, we define the \emph{brick locus} 
\[
\varGr^{\B}\coloneqq\grass^{\B}(g,\scrU)\coloneqq\Bigl\{\ssW \subseteq \scrU\,\st \,\P_{\nu!^*_\ssW}\supseteq \B\Bigr\},
\]
where $\P_{\nu!^*_\ssW}\subseteq \Delta!_g$ is the polytope associated to the submodular function $\nu!^*_{\ssW}$. 

The component-wise action of the torus $\Gm^V$ on $\scrU =\bigoplus_{v\in V}\scrU_v$ induces a natural action of $\Gm^V$ on $\varGr$, with the diagonal $\Gm \hookrightarrow \Gm^V$  acting trivially. Let $\varT \coloneqq\rquot{\Gm^V}{\Gm}$ be the quotient torus. The brick locus $\varGr^{\B}$ is invariant under the action of $\varT$. We interpret the bricks as stability chambers corresponding to the various linearizations of the action by $\varT$, and prove in Theorem~\ref{GIT-projective} that the categorical quotient $\rquot{\varGr^{\B}}{\varT}$, a GIT quotient, exists and is projective. 

\smallskip

For each cone $\sigma$ in $\Sigma$ that intersects $\R^E_{>0}$, and each full-dimensional brick $\B\in\Br_g$, denote by $\sssigma_{\B}$ the smallest cone in $\ssSigma_{\B}$ that contains $\sigma$. Define the map 
\[
\ssboldPhi_{\sigma}\colon\Gm^E \, \xrightarrow{\ssub{\bigl( \ssPhi^{\B}_{\sssigma_{\B}}\bigr)}!_{\B\in\Br_g}} \prod_{\B\in\Br_g}\frac{\varGr^{\B}}{\varT}
\]
whose $\B$-component $\ssPhi_{\sssigma_{\B}}^{\B} \colon\Gm^E \longrightarrow \frac{\varGr^{\B}}{\varT}
$ is given as follows: Let $(\sss,\pi)\in \PSL(\B)$ be the unique permissible slope-level pair such that $\sssigma_{s,\pi} = \sssigma_{\B}$. We view each element $\rho \in \Gm^E(\k)$ as a map $\rho \colon E \to \k^\times$, and consider its restriction to $\ssE^\rmint_s$, viewed as a map from $\ssA^{\rmint}_s$ to $\k^\times$ given by  $a\mapsto \ssrho_a$ with $\ssrho_a =\ssrho_e$  for $e$ the underlying edge of $a$. We associate to $\rho$ the partial gluing data $y=\ssy(s,\rho)$ for the $\ssL_{s,v}$, which corresponds to the data $(\ssy_a=\ssrho_a^{s(a)})_{a\in\ssA_{s}^{\rmint}}$ via the isomorphisms $\ssL_{s,v}\rest{p^a}\cong\k$ for $v\in V$ and $a\in\E_v$. Let $\ssL_s^{\ssy(s,\rho)}$ be the resulting subsheaf of the sum $\bigoplus\ssL_{s,v}$, that is,
\[
\ssL_{s}^{y(s,\rho)}= \bigcap_{a=uv\in\ssA^\rmint_s}\ker\Bigl(\bigoplus_{w\in V} \ssL_{s,w}\longrightarrow\ssL_{s,u}\rest{\ssp^a}\oplus\ssL_{s,v}\rest{p^{\bar a}}\cong\k\oplus\k \xrightarrow{(-\ssrho_a^{-s(a)},1)}\k\Bigr).
\]
Then,  $\ssboldPhi^{\B}_{\sssigma_\B}$ is defined by sending $\rho$ to the orbit under $\varT$ of $\Wexp_{s,\pi}(\ssy(s,\rho))$, where, 
\begin{equation*}
\Wexp_{s,\pi}(y(s,\rho)) = \Res^{-1}(\GlobSp_{\pi}) \bigcap H^0(X,\ssL_{s}^{\ssy(s,\rho)})\subseteq \bigoplus_{v\in V} H^0(\varC_v,\ssL_{s,v}) \subseteq \scrU. 
\end{equation*}

For each cone $\sigma \in \Sigma$ that intersects $\R_{>0}^E$, we set
\[
\varV_\sigma\coloneqq\mathrm{Im}(\ssboldPhi_{\sigma})\subseteq
\prod_{\B\in\Br_g}\frac{\varGr^{\B}}{\varT}, \quad \text{and}\quad  \varV\coloneqq\bigcup_{\sigma\in\Sigma} \varV_\sigma.
\]

Using the squashing process of Section~\ref{sec:fan-structure} and our study of the $\pi$-residue spaces in \cite{AEG-residue-polytope} (see Section~\ref{sec:rp-intro}), we prove the following theorem.

\newtheorem*{thm:variety-lcs}{\cref{thm:variety-lcs}}
\begin{thm:variety-lcs} \emph{Let $X$ be a nodal curve over an algebraically closed field $\k$ of characteristic zero with branches over nodes in general position on the components of its normalization, and let $G$ be its dual graph. Then:
\begin{enumerate}
    \item The variety of limit canonical series $\varV$ parametrizes all fundamental collections associated to any smooth proper curve  $\varX$ over a valued field $\K$ that tropicalizes to $(X, \Gamma)$ for any metric graph $\Gamma$ over $G$. 
    \smallskip
   \item For each cone $\tau$ of $\Sigma$ meeting $\R_{>0}^E$, we have that $\compvarV_{\tau}=\bigcup_{\sigma\supseteq\tau}\varV_{\sigma}$.  In particular, $\varV$ is a projective variety.
\end{enumerate}
}
\end{thm:variety-lcs}

\subsection{Applications}\label{sec:applications} Let $X$ be as before, and $G=(V,E)$ its dual graph. One of the main goals in \cite{EH87a, EM, ES07} was to describe limits of Weierstrass points along \emph{smoothings} of $X$, that is, one-parameter families of smooth curves degenerating to $X$, for curves $X$, respectively, of compact type, with two components, and whose each pair of components intersect (with the restriction that the smoothing be \emph{regular}, that is, with regular total space). We improve on all these results by describing these limits in the full generality considered in this paper. 

In fact, in order to do this, we need only part of the data of the  fundamental collection of simple limit canonical series. 

We use a refinement of \cite[Thm.~2.8]{EM}, which gives a formula for the limit Weierstrass divisor in a smoothing of $X$. The formula requires the determination of the limit canonical series which have foci on the components of $X$. It is actually enough to produce for each $v\in V$ a limit canonical series whose restriction to the component associated to $v$ has the same rank, and adjust the so-called correction numbers accordingly; cf.~\cite[Rmk.~2.9]{EM}. This is the case of the simple limit canonical series $\ssW_h$ whose polytopes contain a vertex of the standard simplex $\Delta_g$; see~\cite[Thm.~9.4]{AE-modular-polytopes}. The variety that parameterizes the limit Weierstrass divisors is thus the image of the projection of the subvariety
\[
\varV=\bigcup_{\sigma\in\Sigma} \varV_\sigma \subset \prod_{\B\in\Br_g} \frac{\varGr^{\B}}{\varT}
\]
onto the product indexed by the $|V|$ bricks $\B\in\Br_g$ which contain a vertex of $\Delta_g$. This variety provides a broad generalization for all curves of the variety constructed in \cite{EM} for two-component curves.

Our next theorem describes how each point on this image determines the limit Weierstrass divisor of any smoothing of $X$ which gives rise to it. In the algebraic setting, its proof uses the ideas above, and in the more general analytic setting considered in this paper, the proof uses a refined extension of \cite[Thm.~2.8]{EM} given in~\cite[\S{A}]{AGR23} using Berkovich geometry.

For each $v\in V$, we call  \emph{$v$-extremal} the vertex $q$ of $\Delta_g$ with $q(v)=g$. 
For each $a\in\E_v$, denote by $\ssm_a$ the maximum order of poles of the meromorphic differentials in $\scrU_v$ at $\ssp^a$. Put $\scrL_v\coloneqq\omega!_v(\sum_{a\in\E_v}\ssm_a \ssp^a)$. Then, $\scrU_v\subseteq H^0(\varC_v,\scrL_v)$. Finally, for each $e\in E$, let $\ssm_e\coloneqq\ssm_a+\ssm_{\bar a}$, with $a$ and $\bar a$ the arrows lying over $e$. 

\begin{thm} Let $\varX$ be a smooth proper curve over a valued field $\K$ which tropicalizes to $X$. For each $v\in V$, let $\ssW(\varX,v)$ be the limit canonical series in the fundamental collection $\FC(\varX)$ whose associated polytope contains the $v$-extremal vertex of $\Delta_g$, and let $\ssR(\varX,v)$ be the ramification divisor of the linear series of sections of $\scrL_v$ induced by the projection of $\ssW(\varX,v)$ in $\scrU_v$. Then, the limit Weierstrass divisor $\ssR(\varX)$ on $X$ is given by
\[
\ssR(\varX)=\sum_{v\in V}\ssR(\varX,v)+\sum_{e\in E}g(g+1-\ssm_e)p^e.
\]
\end{thm}

(We remark that the $\ssm_e$ depend on the choice of the $\scrU_v$ but the right-hand side in the formula above does not.)

In the above theorem, we used only the ``extremal'' limit canonical series in the fundamental collection. However, as observed already in \cite{EO}, additional limits are needed, in general, to describe the limits of canonical divisors.

In addition, by considering all limit canonical series in the fundamental collection, we can define a rational map from $X$ to a quiver Grassmannian, in the spirit of \cite{ESV2}. This construction can be used to describe a limit of the canonical map of $\varX$, extending the theory of admissible covers introduced by Harris and Mumford \cite{HM82} to the setting of canonical maps. Namely, each element $\ssW_h$ in $\FC(\varX)$ defines a rational map $X \dashrightarrow \mathbb{P}(\ssW_h^*)$, well-defined away from a finite set of points, precisely because $\ssW_h$ belongs to $\FC(\varX)$.  Combining them yields a rational map from $X$ to the product of the $\mathbb{P}(\ssW_h^*)$.  Remarkably, this map factors through a specific quiver Grassmannian, the one which parametrizes one-dimensional subrepresentations of the dual of the representation formed by the $\ssW_h$ of the quiver coming from the polyhedral decomposition of $\Delta_g$ by the polytopes $\P_{\nu!_h^*}$, given in Theorem~\ref{thm:W-tiling}. This factorized map can be viewed as a limit of the canonical map of $\varX$.

\subsection{Related work} \label{sec:related-work}

The literature on degenerations of line bundles over families of curves is extensive. We refer to \cite{DS, OS79, AK1, AK2, Cap, Pan96, Esteves01, MRV, KP}, among others. Part of the motivation for these works was to understand degenerations of linear series, which, however, proved to be significantly subtler, with only a few partial results known. Two theories are available today in restricted settings. The first, by Harris and Mumford \cite{HM82}, the theory of admissible covers, applies to degenerations of rank one linear series. The second, systematized in the pioneering work by Eisenbud and Harris \cite{EH86}, the theory of limit linear series, applies for any rank but only to limit curves of compact type. In \cite{EHBull}, Eisenbud and Harris posed the problem of extending their theory to all stable curves, a problem that remained largely unsolved to this day. Efforts were made, for instance in \cite{Ran85, Est98, Os-moduli, Os-pseudocompact}, but the ideas introduced were not further developed, likely because the combinatorial framework was rather unwieldy. 
Our framework, introduced in ~\cite{AE-modular-polytopes}, provides a new combinatorial and polyhedral foundation for a program aimed at tackling this problem, with the work presented in this paper as its first application.

On the combinatorial side, divisor theory for graphs was introduced in \cite{BN07} by Baker and Norine,  and later extended to tropical and hybrid settings. For a survey of these developments, we refer to Baker and Jensen~\cite{BJ} and Jensen and Payne~\cite{JP21}. The degeneration of linear series from a tropical perspective has been explored in~\cite{JP22, AG22}. Notably, tropical methods have been instrumental in studying the geometry of curves and their moduli spaces, with a significant application in the recent work by Farkas, Jensen, and Payne~\cite{FJP}. Recent work on the degeneration of tropical line bundles include \cite{AP20, AAPT, MMUV, CPS23}. Tropicalizations of the moduli space of admissible covers and the moduli space of weighted stable curves have been studied in~\cite{CMR16} and~\cite{Ulr15, CHMR16}.

\subsection*{Acknowledgment} The work that led to the present paper started in 2015 during a visit of the second author to ENS Paris. Over the past 10 years, we have benefited from the support and hospitality of various institutions. We warmly thank ENS, IMPA, \'Ecole Polytechnique, CNRS, the French-Brazilian network in mathematics, Math+ research center in Berlin, COFECUB-CAPES, CNPq and FAPERJ, among others. This work was partially supported by the "CAPES-COFECUB" programme (project number: Ma 1017/24), funded by the French Ministry for Europe and Foreign Affairs, the French Ministry for Higher Education and CAPES.


\section{Submodular functions}\label{sec:submodular} 
Let $V$ be a finite nonempty set. Denote by $\ssub{2}!^V$ the family of subsets of $V$. We study properties of functions $\varphi\colon \ssub{2}!^V\to\R$. 

Given two functions $\varphi,\mu \colon \ssub{2}!^V \to \R$, we write $\varphi \geq \mu$ if $\varphi(I) \geq \mu(I)$ for every $I\subseteq V$. Furthermore, we write $\varphi>\mu$ if the inequalities are strict for proper subsets $I$, i.e., equality might hold only if $I=\emptyset$ or $I=V$. (It will be often the case for us that $\varphi$ and $\mu$ agree on $\emptyset$ or $V$.) 

We call $\varphi\colon \ssub{2}!^V\to\R$ \emph{nonnegative} if $\varphi\geq 0$ and \emph{positive} if $\varphi>0$. 
We call $\varphi$ \emph{nonincreasing} if $\varphi(J)\geq\varphi(I)$ for $J\subseteq I\subseteq V$, and \emph{nondecreasing} if the reverse inequalities hold.

For each $I\subseteq V$, we denote by $\ssI^c$ the complement of $I$ in $V$.

\subsection{Submodular functions}
We say a function $\varphi \colon \ssub{2}!^V\to\R$
is \emph{submodular} if $\varphi(\emptyset)=0$ and 
\[
  \varphi(\ssI_1)+\varphi(\ssI_2) \geq \varphi(\ssI_1\cup \ssI_2)+ \varphi(\ssI_1\cap \ssI_2)
\quad\text{for each }\ssI_1,\ssI_2\subseteq V.
\]
 We call the quantity $\varphi(V)$ the \emph{range} of $\varphi$. Similarly, a function  $\mu \colon \ssub{2}!^V \to \R$ is called \emph{supermodular} if $\mu (\emptyset)=0$ and the inequalities above are all reversed; the quantity $\mu(V)$ is called the range of $\mu$. A function which is both supermodular and submodular is called \emph{modular}. Modular functions on $\ssub{2}!^V$ are in bijection with $\R^V$: each $q\in \R^V$ can be viewed as a modular function $q\colon\ssub{2}!^V\to\R$ by putting $q(I)\coloneqq \sum_{v\in I}q(v)$ for each $I\subseteq V$, with the convention $q(\emptyset)\coloneqq 0$.

\subsection{Adjoints} For each $\varphi \colon \ssub{2}!^V\to\R$ with $\varphi(\emptyset)=0$, let $\varphi!^* \colon \ssub{2}!^V\to\R$ be the \emph{adjoint to $\varphi$} defined by
\[
\varphi!^*(I)\coloneqq\varphi(V)-\varphi(V\setminus I)\quad\text{for each }I\subseteq V.
\]
It is easy to see that $\varphi$ is submodular, resp.~supermodular, if and only if $\varphi!^*$ is supermodular, resp.~submodular. Furthermore, $\varphi$ and $\varphi!^*$ have the same range, and $\ssub{(\varphi!^*)}!^* = \varphi$. For a modular function $q$, we have $\ssub{q}!^*=q$.

If $\varphi$ is submodular, then $\varphi^*\leq\varphi$. Indeed, applying the submodularity inequality to $I$ and $\ssI^c$, and using that $\varphi(\emptyset) =0$, we get that $\varphi^*(I)\leq\varphi(I)$ for each $I\subseteq V$; equality holds for $I=\emptyset$ and $I=V$. 

\subsection{Base polytopes} \label{sec:submodular-polytope} 
To each submodular function $\varphi\colon 2^V\to\R$ we associate the base polytope of the polymatroid associated to $\varphi$
\[
\P_{\varphi}\coloneqq\bigl\{q\in\R^V\,\st \, q(I)\leq\varphi(I)\,\quad \forall I\subseteq V, \text{ with equality if }I=V\bigr\},
\]
or simply, the (\emph{base}) \emph{polytope} of $\varphi$.

\subsection{The UpMin transform} \label{sec:upmin}
Given $\varphi \colon \ssub{2}!^V\to \R$, define $\chi\colon \ssub{2}!^V\to \R$ by 
\[
\chi(I) \coloneqq \min_{K\supseteq I} \varphi(K) \qquad \text{for each }I \subseteq V.
\]
We say that $\chi$ is the \emph{UpMin} transform of $\varphi$ and write $\chi = \upmin(\varphi)$.

\begin{prop}\label{prop:min-operation} Let $\varphi \colon \ssub{2}!^V\to \R$. Then, $\upmin(\varphi)$ is nondecreasing. Also, if $\varphi$ is submodular and nonnegative, then so is $\upmin(\varphi)$. 
\end{prop}

\begin{proof} The first statement is clear. As for the second, assume $\varphi$ is submodular and $\varphi\geq 0$, and put $\chi\coloneqq\upmin(\varphi)$. Clearly, $\chi$ is nonnegative and $\chi(\emptyset)=0$. It remains to show that
\[
\chi(\ssI_1)+\chi(\ssI_2) \geq \chi(\ssI_1 \cup \ssI_2) + \chi(\ssI_1 \cap \ssI_2) \qquad \forall \, \ssI_1, \ssI_2 \subseteq V.
\]
Now, for each pair of subsets $\ssK_1, \ssK_2\subseteq V$ with $\ssK_1\supseteq \ssI_1$ and $\ssK_2\supseteq \ssI_2$, from the definition of $\chi$,
\[
\varphi(\ssK_1\cup \ssK_2) \geq \chi(\ssI_1 \cup \ssI_2)\quad \textrm{and}\quad \varphi(\ssK_1\cap \ssK_2)\geq   \chi(\ssI_1 \cap \ssI_2).
\]
Applying submodularity of $\varphi$ to $\ssK_1$ and $\ssK_2$, we get 
\[
\varphi(\ssK_1)+\varphi(\ssK_2)\geq\varphi(\ssK_1 \cup \ssK_2) +\varphi(\ssK_1\cap\ssK_2) \geq \chi(\ssI_1 \cup \ssI_2) + \chi(\ssI_1 \cap \ssI_2),
\]
from which we deduce
\[
\chi(I_1)+\chi(I_2)=\min_{\ssK_1\supseteq \ssI_1}\varphi(\ssK_1)+\min_{\ssK_2\supseteq \ssI_2} \varphi(\ssK_2) \geq \chi(\ssI_1\cup \ssI_2) + \chi(\ssI_1 \cap \ssI_2),
\]
as required. 
\end{proof}

\begin{prop}\label{prop:xi-upmin}
Let $\varphi,\xi\colon \ssub{2}!^V\to\R$ be functions with $\xi$ nonincreasing. Then,
\[
\upmin\bigl(\upmin(\varphi)+\xi\bigr) = \upmin\bigl(\varphi+\xi\bigr).
\]
\end{prop}

\begin{proof}
Let $\mu!_1$ be the UpMin transform of $\varphi$, and $\mu!_2$ be the UpMin transform of $\mu!_1 + \xi$, so
 \[
\mu!_1(I)=\min_{K\supseteq I} \varphi(K) \quad \textrm{and}  \quad \mu!_2(I)=\min_{K\supseteq I} \Bigl(\mu!_1(K) + \xi(K)\Bigr)\qquad\textrm{ for each } I\subseteq V.
  \]
We need to show that $\mu!_2$ is the UpMin transform of $\varphi+\xi$, that is, 
\[
\mu!_2(I)=\min_{K\supseteq I}\Big(\varphi(K)+\xi(K)\Big) \qquad \textrm{for each  } I \subseteq V.
\]
From the very definition of the UpMin transform, for $I\subseteq K\subseteq V$ we have
\[
\mu!_2(I)\leq\mu!_1(K)+\xi(K)\leq \varphi(K)+\xi(K),
\]
and so $\mu!_2(I) \leq \upmin\bigl(\varphi+\xi\bigr)(I)$. We need to prove now the reverse inequality.

Given $I\subseteq V$, there is $I \subseteq \ssK_1\subseteq V$ such that $\mu!_2(I)=\mu!_1(\ssK_1)+\xi(\ssK_1)$.  There is similarly $\ssK_1 \subseteq \ssK_2 \subseteq V$ such that $\mu!_1(\ssK_1) =\varphi(\ssK_2)$. Thus, 
\[
 \mu!_2(I)=\varphi(\ssK_2)+\xi(\ssK_1).
\]
Since $\xi$ is nonincreasing, $\xi(\ssK_1)\geq\xi(\ssK_2)$. So
\[
\mu!_2(I) \geq \varphi(\ssK_2)+\xi(\ssK_2),
\]
from which we deduce $\mu!_2(I) \geq \upmin\bigl(\varphi+\xi\bigr)(I)$, finishing the proof.
\end{proof}

\subsection{Geometric meaning of UpMin} \label{sec:upmin-geometric}
Let $\Delta!_g$ be the standard simplex of width $g$, for $g\geq 0$:
\[
\Delta!_g\coloneqq\Bigl\{q\in\R^V_{\geq 0}\,\st \, q(V) =g\Bigr\} \subset \R^V.
\]

\begin{prop}\label{prop:beautiful-geometric-meaning} Let $\varphi\colon 2^V\to\R$ be a nonnegative submodular function of range $g$. Then,
\[
\P_{\upmin(\varphi)}=\P_{\varphi}\cap\Delta!_g = \P_{\varphi}\cap\R_{\geq 0}^V
\]
\end{prop}

\begin{proof} Put $\chi\coloneqq\upmin(\varphi)$. Clearly, $\chi(I)\leq\varphi(I)$ for each $I\subseteq V$, with equality if $I=V$. It follows that  $\P_{\chi}\subseteq\P_{\varphi}$. Also,  $\chi(V\setminus I)\leq\chi(V)$, and thus
\[
q(I)=q(V)-q(V\setminus I)\geq\chi(V)-\chi(V\setminus I)\geq0
\]
for each $q\in\P_{\chi}$, yielding $\P_{\chi}\subseteq\Delta!_g$. Finally, if $q\in\P_{\varphi}\cap\R_{\geq 0}^V$, then we have $q(I)\leq q(K)\leq\varphi(K)$ for each $I\subseteq K\subseteq V$, whence $q\in\P_{\chi}$. 
\end{proof}

\begin{figure}
    \centering
   \scalebox{.25}{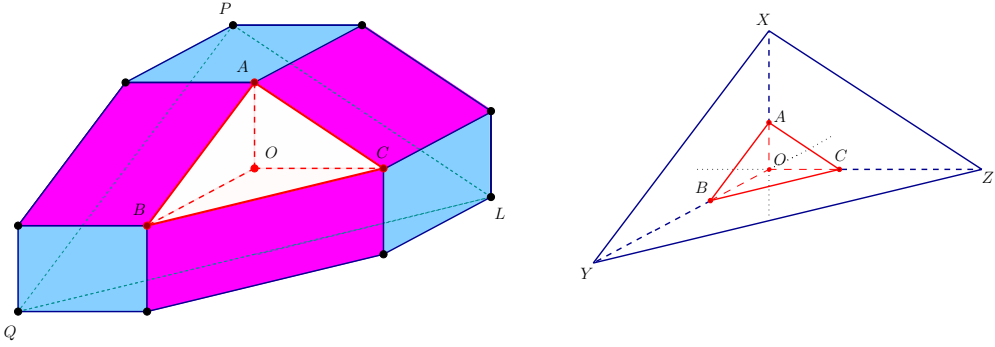}
    \caption{The polytope $\P_\varphi$ on the left is the one associated to the submodular function $\varphi$ defined in Example~\ref{ex:4lines-polytope}. The drawing shows its projection to $\R^3$, using the last three coordinates. It has $7$ visible and $4$ hidden facets. The origin $O$ lies in the interior of $\P_\varphi$. The red tetrahedron with vertices $O, A, B$ and $C$ is the intersection of $\P_\varphi$ with $\Delta!_3 \subset \R_{\geq 0}^4$ (projected down again to $\R^3$): it is the polytope $\P_\chi$ associated to the UpMin transform $\chi = \upmin(\varphi)$. The rescaled figure on the right shows the position of $\P_\chi$ within $\Delta!_3$ (the big tetrahedron in blue with vertices $O, X, Y, Z$).}
    \label{fig:polytope-4-lines}
\end{figure}

\begin{example}\label{ex:4lines-polytope} Consider the set $V = \{\ssu_0, \ssu_1,\ssu_2, \ssu_3\}$ and the submodular function $\varphi \colon \ssub{2}!^V \to \Z$ of range 3 defined by 
\begin{align*}
    &\varphi(\ssu_0)=\varphi(\{\ssu_0,\ssu_i\})=5,\quad \varphi(\{\ssu_0, \ssu_i, \ssu_j\})=4 \quad \text{for all distinct $i, j \in \{1,2,3\}$},\\
    &\varphi(I)=1\quad \text{for all nonempty }I\subseteq\{\ssu_1,\ssu_2,\ssu_3\}.
\end{align*}
The polytope $\P_{\varphi}\subset\R^4$, depicted in Figure~\ref{fig:polytope-4-lines}, is the set of all points $(\ssq_0, \ssq_1, \ssq_2, \ssq_3)$ in $\R^4$ with $\ssq_0 =3 -\ssq_1-\ssq_2-\ssq_3$ that verify, for all pairs $i,j$ of distinct elements in $\{1,2,3\}$, the inequalities
 $   -1 \leq \ssq_i \leq 1,  \,\, -2 \leq \ssq_i +\ssq_j \leq 1 \,\, \text{ and }\,\, -2 \leq  \ssq_1+\ssq_2+\ssq_3 \leq 1.$
The function $\chi=\upmin(\varphi)$ satisfies $\chi(I) = 1$ for each nonempty $I\subseteq \{\ssu_1, \ssu_2, \ssu_3\}$, and $\chi(I)=3$ for each $I\subseteq V$ containing
$\ssu_0$. The polytope $\P_\chi$, also depicted in Figure~\ref{fig:polytope-4-lines}, is the set of all points $(\ssq_0, \ssq_1, \ssq_2, \ssq_3) \in \R_{\geq 0}^4$ with $\ssq_0=3 -\ssq_1-\ssq_2-\ssq_3$ that verify $\ssq_1+\ssq_2+\ssq_3 \leq 1$.
\end{example}

\subsection{The DownSum transform}
Given a function $\varepsilon \colon \ssub{2}!^V \to \R$, define the \emph{DownSum} transform of $\varepsilon$ denoted $\downsum(\varepsilon)\colon \ssub{2}!^V \to \R$ by 
\[
\downsum(\varepsilon)(I) \coloneqq  \sum_{J\subseteq I} \varepsilon(J)   \qquad \text{for each }I \subseteq V.
\]

\begin{prop}\label{prop:min-submodular} Let $\varepsilon \colon \ssub{2}!^V \to \R$. If $\varepsilon\geq 0$, then $\downsum(\varepsilon)$ is supermodular and nondecreasing. 
\end{prop}

\begin{proof} For each three subsets $I_1,I_2,J$ of $V$, we have that $J\subseteq I_1\cup I_2$ if $J\subseteq I_1$ or $J\subseteq I_2$, and $J\subseteq I_1\cap I_2$ if $J\subseteq I_1$ and $J\subseteq I_2$. This simple observation coupled with the fact that $\varepsilon\geq 0$ proves the supermodularity inequality of $\downsum(\varepsilon)$ for $I_1$ and $I_2$. 
\end{proof}

\subsection{The submodular function of a subset}\label{sec:submodular-subset-upmin} 

For each subset $J\subseteq V$, denote by $\xi!_J\colon \ssub{2}!^V\to\Z$ the function satisfying
\[
\xi!_J(I)=\begin{cases}
0&\text{if }I\not\supseteq J,\\
-1&\text{if }I\supseteq J.
\end{cases}
\]
If $J=\{v\}$ we abbreviate $\xi!_v\coloneqq\xi!_J$.

\begin{prop}\label{xi-submodular} For each $J\subseteq V$, the function $\xi!_J$ is submodular and nonincreasing.
\end{prop}

\begin{proof} This follows from Proposition~\ref{prop:min-submodular} and the equality
\[
\xi!_J=-\downsum(\one!_J),
\]
where $\one!_J\colon \ssub{2}!^V \to \R$ is the characteristic function of $J$, taking values $1$ on $J$ and $0$ elsewhere.
\end{proof}

\begin{prop}\label{prop:xi-singleton} Let $\varphi \colon \ssub{2}!^V\to \Z$ be a nonnegative integer-valued submodular function, and let $J\subseteq V$. Then, $\varphi+\xi!_J\geq 0$ if and only if there is $v\in J$ such that $\varphi+\xi!_v\geq 0$.  
\end{prop}

\begin{proof} Since $\xi!_J\geq\xi!_v$ for each $v\in J$, if $\varphi+\xi_v\geq 0$, then  $\varphi+\xi_J\geq 0$. Conversely, assume $\varphi+\xi_J\geq 0$. Suppose for the sake of contradiction that for each $v\in J$, there is $\ssI_v\subseteq V$ such that $\varphi(\ssI_v)+\xi!_v(\ssI_v)< 0$. Since $\varphi \geq 0$, we must have $v\in \ssI_v$. And since $\varphi$ is integer-valued, it follows that $\varphi(I_v)=0$ for each $v\in J$. Let $I\coloneqq \bigcup \ssI_v$. Since $\varphi$ is submodular, we get
\[
0\leq\varphi(I)\leq\sum_{v\in J}\varphi(I_v) = 0,
\]
and hence $\varphi(I)=0$. But $J\subseteq I$, so $\xi!_J(I)=-1$, whence $\varphi(I)+\xi_J(I)=-1$, a contradiction.
\end{proof}

\subsection{Submodular functions of subspaces} \label{sec:submodular-subspace} 

Let $\k$ be a field. For each $v\in V$, let $\VS_v$ be a vector space over $\k$. Put $\VS\coloneqq \bigoplus_{v\in V} \VS_v$. For each subset $I\subseteq V$, let $\VS_I\coloneqq \bigoplus_{v\in I}\VS_v$ and denote by $\proj{I}\colon \VS\to \VS_{I}$ the corresponding projection map. 

Let $\ssW\subseteq \VS$ be a finite-dimensional vector subspace. For each $I\subseteq V$, put $\ssW_I\coloneqq \proj{I}(\ssW)$, and denote by $\nu!^*_{\ssW}(I)$ its dimension, that is, $\nu!_{W}^*(I) \coloneqq \dim_{\k}(\ssW_I)$. 

\begin{prop}\label{prop:submodular-subspace} The function $\nu!_{W}^* \colon \ssub{2}!^V \to \Z$ is submodular, nonnegative and nondecreasing. 
\end{prop}

\begin{proof} Submodularity is direct, see~\cite[Prop.~8.1]{AE-modular-polytopes}. The rest is clear.
\end{proof}

\subsection{Simpleness} A submodular function $\varphi \colon \ssub{2}!^V \to \R$ is called \emph{simple} if for each partition of $V$ into nonempty sets $I$ and $J$, we have $\varphi(I) + \varphi(J)>\varphi(V)$. If $|V|=1$, this is automatic.

\begin{prop}\label{prop:simpleness-min-operation-0}
Let $\varphi\colon 2^V\to\R$ be a submodular function and $\varphi!^*$ its adjoint. Then, $\varphi$ is simple if and only if $\varphi^*<\varphi$. In particular, if $\varphi$ is simple and $\varphi^*\geq 0$, then $\varphi>0$.
\end{prop}

\begin{proof} Immediate.
\end{proof}

\begin{prop}\label{prop:simpleness-min-operation} Let $\varphi\colon 2^V\to\R$ be a nonnegative submodular function. If $\upmin(\varphi)$ is simple, then $\varphi$ is simple. The converse holds if $\varphi$ is positive with positive range.
\end{prop}

\begin{proof}  Put $\chi \coloneqq \upmin(\varphi)$. Since $\varphi$ is submodular and nonnegative, so is $\chi$ by Proposition~\ref{prop:min-operation}. 

Consider a partition of $V$ into nonempty subsets $I$ and $J$. The proposition follows from the following two assertions:
\begin{enumerate}
\item If $\varphi(I)+\varphi(J)=\varphi(V)$, then $\chi(I)+\chi(J)=\chi(V)$. 
\item The converse holds if $\varphi$ is positive with positive range.
\end{enumerate}

The first assertion follows from the inequalities
\[
\varphi(I)+\varphi(J)\geq\chi(I)+\chi(J)\geq\chi(V)=\varphi(V).
\]
As for the second, we have $\chi(I)=\varphi(K)$ and $\chi(J)=\varphi(N)$ for certain subsets $K,N$ of $V$ such that $I\subseteq K$ and $J\subseteq N$. Assume  $\chi(I)+\chi(J)=\chi(V)$. Since $\varphi$ is submodular and nonnegative,
\[
\chi(V)=\varphi(K)+\varphi(N)\geq\varphi(K\cup N)+\varphi(K\cap N)\geq\varphi(K\cup N)=\varphi(V)=\chi(V).
\]
It follows that $\varphi(K)+\varphi(N)=\varphi(V)$ and $\varphi(K\cap N)=0$. If $\varphi$ is positive with positive range, then $K\cap N=\emptyset$. Since $I\subseteq K$ and $J\subseteq N$, we must have $I=K$ and $J=N$, and hence $\varphi(I)+\varphi(J)=\varphi(V)$.
\end{proof}


\section{Graphs}\label{sec:graphs}

Let $G = (V, E)$ be a finite graph with vertex set $V$ and edge set $E$. If $e\in E$ connects $v,w\in V$, we write $e=\{v,w\}$. Given a subset $I\subseteq V$, we denote by $E(I)$ the set of edges in $E$ with both extremities in $I$. The \emph{induced subgraph on $I$}, denoted $G[I]$, is the subgraph of $G$ with vertex set $I$ and edge set $E(I)$. A \emph{connected component} of $G$ is an induced subgraph $G[I]$ which is connected and, in addition, there is no edge in $G$ that connects a vertex in $I$ to one in the complement $\ssI^c$. 
The \emph{genus} of a connected graph $G$ is the quantity $g(G) \coloneqq |E|-|V|+1$.

We denote by $\E$ the set of \emph{arrows} of $G$: each edge $e=\{v,w\}$ gives rise to two arrows $v\to w$ and $w\to v$ in $\E$. If $e$ is a loop, i.e.,~$v=w$, we still have two arrows $v\leftrightarrow w$, corresponding to the two half edges issued from $e$. 

If $a$ is an arrow with orientation $v \to w$, we call  $v$ the \emph{tail} and $w$ the \emph{head} of $a$, and write $\te_a=v$ and $\he_a=w$, as well as $a=vw$. The arrow in $\E$ with the reverse orientation is denoted $\overline a$; it has tail $w$ and head $v$.

Given a set of arrows $A\subseteq\E$ and subsets $I,J \subseteq V$, denote by $A(I,J)$ the subset of arrows in $A$ with tail in $I$ and head in $J$. We abbreviate $A(I)\coloneqq A(I,I)$ and $\ssA_I\coloneqq A(I,V)$; see Figure~\ref{fig:digraph} for an example. If $I=\{v\}$, a singleton, we simply write $\ssA_v$. Clearly, $\ssA_I$ is the disjoint union of the $\ssA_v$ for $v\in I$, that is, $\ssA_I = \bigsqcup_{v\in I}\ssA_v$.

\subsection{Submodular functions arising from digraphs} 

An \emph{orientation} $\orient$ of $G$ is a map $\orient\colon E \to \E$ that associates to each edge $e$ of $G$ one of the two arrows in $\E$ arising from $e$. A \emph{digraph} $H=(V, A)$ is a graph $G=(V,E)$ with an orientation $\orient$ of its edges so that $A = \orient(E)$.

\begin{prop}\label{prop:modular-digraph} 
Let $H=(V, A)$ be a finite digraph. Let $\varphi,\varphi!_1,\varphi!_2 \colon \ssub{2}!^V \to \Z$ be defined by
\begin{align*}
\varphi(I) \coloneqq |A(I)| + |A(\ssI^c, \ssI)|, \qquad \varphi!_1(I) \coloneqq |A(\ssI^c, \ssI)|, \qquad
\varphi!_2(I) \coloneqq - |A(I)|,
\end{align*}
for each $I\subseteq V$. Then, $\varphi$ is modular, and $\varphi!_1$ and $\varphi!_2$ are submodular.
\end{prop}

\begin{proof} All of the assertions can be reduced to the case where $A$ has only one element, and then proved by a case analysis. Also, observe that $\varphi!_1=\varphi+\varphi!_2$; thus, if $\varphi$ is modular, $\varphi!_1$ is submodular if and only if so is $\varphi!_2$.
\end{proof}

\subsection{Level structures on graphs}\label{sec:levelgraphs} 

In this paper, an ordered set is a set with a total order on its elements. We view a subset of $\R$ as ordered with the total order induced by that of $\R$. We add few more terminology to those in Section~\ref{sec:ls-rs}.

Recall that a level structure on the graph $G=(V,E)$ is the data $\pi$ of an ordered partition of the vertex set $V$.  Equivalently, a level structure is the data of a set $\ssV_\pi$ endowed with a total order $\subface_{\pi}$ and a surjection $h=\ssh_\pi\colon V \to \ssV_\pi$. This way, we can view $\ssV_\pi$ as the unordered partition underlying the level structure. The surjection $h=\ssh_\pi\colon V \to \ssV_\pi$ takes then each vertex to the part it belongs to. We refer to either $(G, \pi)$ or $(G,h)$ as a level graph, and call $h$ the associated level function. For each vertex $v\in V$, $h(v)$ is called the level of $v$. For each pair of vertices $u,v\in V$, abusing notation, we write $u\subface_\pi v$ if $h(u) \subface_\pi h(v)$. If the ordered partition $\pi$ or the function $h$ is understood from the context, we abbreviate $(G, \pi)$ and $(G,h)$ to $G$, leaving $\pi$ and $h$ implicit. 

A function $h\colon V\to\R$ induces the ordered partition $\pi=\sspi_h$ of $V$ such that $h=\iota\ssh_{\pi}$ for an order-preserving injection $\iota\colon\ssV_{\pi}\to\R$. Abusing notation, we call $(G,h)$ a level graph, and $h$ its level function.

An edge $e=\{u,v\}$ in $G$ is horizontal if $u$ and $v$ belong to the same part in the partition, that is, if $h(u)=h(v)$; otherwise, we call it vertical. We denote the set of vertical edges by $\ssE_\pi$. Its complement, $\ssE_\pi^c = E \setminus \ssE_\pi$, is the set of horizontal edges. Denote by $\E_\pi$ the subset of $\E$ consisting of arrows on vertical edges, thus called vertical. The complement $\E_\pi^c = \E \setminus \E_\pi$ is the set of arrows on horizontal edges, thus called horizontal. 

We use the level structure to define a bipartition of $\E_\pi$. An arrow $u\rightarrow v$ on a vertical edge in $E_\pi$ is compatible with the level structure if $u \supface_\pi v$. (Visually, the arrow $\rightarrow$ and the order $\supface_\pi$ point to $v$.) Denote by $\A_\pi$ the set of arrows on vertical edges compatible with the level structure. We refer to the elements of $\ssA_\pi$ as upward arrows. (Visually, we think of the level as depth, so the upward arrows on vertical edges point to the vertex with smaller depth.) Its complement in $\E_\pi$ is denoted $\barA_\pi$; the elements of $\barA_\pi$ are called downward arrows. Equivalently, 
$\barA_\pi \coloneqq \{\bar a\st a\in \A_\pi\}$. This way we obtain a bipartition $\E_\pi = \A_\pi \sqcup \barA_\pi$.

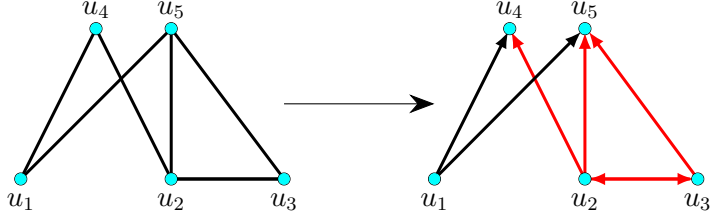
\begin{figure}[ht]
    \centering
\begin{tikzpicture}[scale=.8pt]

\draw[line width=0.4mm] (5,-2) -- (6,0);
\draw[line width=0.4mm] (7,-2) -- (6,0);
\draw[line width=0.4mm] (7,-2) -- (7,0);
\draw[line width=0.4mm] (5,-2) -- (7,0);
\draw[line width=0.4mm] (8.5,-2) -- (7,0);
\draw[line width=0.4mm] (7,-2) -- (8.5,-2);
\draw[line width=0.4mm] (8.5,-2) -- (7,-2);

\draw[line width=0.1mm] (6,0) circle (0.8mm);
\filldraw[aqua] (6,0) circle (0.7mm);
\draw[line width=0.1mm] (7,0) circle (0.8mm);
\filldraw[aqua] (7,0) circle (0.7mm);
\draw[line width=0.1mm] (5,-2) circle (0.8mm);
\filldraw[aqua] (5,-2) circle (0.7mm);
\draw[line width=0.1mm] (7,-2) circle (0.8mm);
\filldraw[aqua] (7,-2) circle (0.7mm);
\draw[line width=0.1mm] (8.5,-2) circle (0.8mm);
\filldraw[aqua] (8.5,-2) circle (0.7mm);

\draw [-{latex}, line width=0.4mm](10.5,-2) -- (11.48,-0.04);
\draw [-{latex},line width=0.4mm, red](12.5,-2) -- (11.52,-0.04);
\draw [-{latex},line width=0.4mm, red](12.5,-2) -- (12.5,-0.05);
\draw [-{latex}, line width=0.4mm](10.5,-2) -- (12.45,-0.04);
\draw [-{latex},line width=0.4mm, red](14,-2) -- (12.55,-0.04);
\draw [-{latex},line width=0.4mm, red](12.5,-2) -- (13.95,-2);
\draw [-{latex},line width=0.4mm, red](13.9,-2) -- (12.55,-2) ;

\draw[line width=0.1mm] (11.5,0) circle (0.8mm);
\filldraw[aqua] (11.5,0) circle (0.7mm);
\draw[line width=0.1mm] (12.5,0) circle (0.8mm);
\filldraw[aqua] (12.5,0) circle (0.7mm);
\draw[line width=0.1mm] (10.5,-2) circle (0.8mm);
\filldraw[aqua] (10.5,-2) circle (0.7mm);
\draw[line width=0.1mm] (12.5,-2) circle (0.8mm);
\filldraw[aqua] (12.5,-2) circle (0.7mm);
\draw[line width=0.1mm] (14,-2) circle (0.8mm);
\filldraw[aqua] (14,-2) circle (0.7mm);

\draw[black] (6,0.25) node{$\ssu_4$};
\draw[black] (7,0.25) node{$\ssu_5$};
\draw[black] (5,-2.3) node{$\ssu_1$};
\draw[black] (7,-2.3) node{$\ssu_2$};
\draw[black] (8.5,-2.3) node{$\ssu_3$};

\draw[black] (11.5,0.25) node{$\ssu_4$};
\draw[black] (12.5,0.25) node{$\ssu_5$};
\draw[black] (10.5,-2.3) node{$\ssu_1$};
\draw[black] (12.5,-2.3) node{$\ssu_2$};
\draw[black] (14,-2.3) node{$\ssu_3$};

\draw [-{Stealth[scale=2]}](8.5,-1) -- (10.5,-1);

\end{tikzpicture} 
\caption{A level graph with two levels on the left, and the corresponding upward and horizontal arrows on the right. The set $\E_I$ for $I=\{\ssu_2,\ssu_3\}$ is depicted in red.}
\label{fig:digraph}
\end{figure}

In figures, we draw horizontal edges horizontally and vertical edges vertically, with the arrows in $\ssA_\pi$ pointing upward. As the directions of the arrows are visually clear, we sometimes drop them.  
Figure~\ref{fig:digraph} illustrates a level graph $G$ on the left, using a level function $h\colon V \to \{1,2\} \subset \R$. The part with smallest level is depicted on top. Thus, $h(\ssu_1) = h(\ssu_2) = h(\ssu_3) > h(\ssu_4) = h(\ssu_5)$, so, for example, $\ssu_1 \supface_\pi \ssu_4$. On the right, we depict the upward and horizontal arrows of $G$

\subsection{Residue spaces and residue conditions}\label{sec:residue-spaces-and-conditions} 

We associate to a level graph $(G,\pi)$ vector spaces $\GlobSp_{\pi}\subseteq\bV$ which will play an important role in our study of canonical series; see Section~\ref{sec:twisted-differentials}. We suppose that $G$ is connected.

Let $h=\ssh_\pi \colon V \to \ssV_\pi$ be the corresponding level function. We fix a field $\k$. Define $\bV\coloneqq\bV_G\coloneqq\k^{\E}$.  For each $\psi \in \bV$ and $a\in \E$, denote by $\sspsi_a$ the $a$-coordinate of $\psi$. For each $v\in V$, we set $\bV_v\coloneqq\k^{\E_v}$. Then, we have the decomposition $\bV =\bigoplus\bV_v$.

\begin{defi} Let $(G, \pi)$ be  a level graph. We let $\GlobSp_{\pi}\subseteq\bV$ be the subspace of the $\psi$ which satisfy Conditions~\ref{cond:R1} through~\ref{cond:R4} in Section~\ref{sec:ls-rs}, and call it the \emph{$(G,\pi)$-residue space}, or \emph{$\pi$-residue space} for short if $G$ is fixed.  
\end{defi}

\begin{thm}\label{thm:dimG} The dimension of $\GlobSp_\pi$ is the genus of $G$, that is, 
\[
\dim_{\k}\GlobSp_\pi=|E|-|V|+1.
\]
\end{thm}
\begin{proof}
In the absence of horizontal edges, this can be obtained from \cite[Prop.~6.3]{BCGGM19}, using Picard--Lefschetz Theory on a degenerating family of Riemann surfaces. The statement in this general form is proved in \cite{AEG-residue-polytope}.
\end{proof}

\begin{defi}
We define $\gamma!_{\pi} \colon \ssub{2}!^V \to \Z$ to be the submodular function associated to the subspace $\GlobSp_{\pi}\subseteq \bV=\bigoplus_{v\in V} \bV_v$, that is, $\gamma!_{\pi}\coloneqq\nu!^*_{\GlobSp_{\pi}}$.
\end{defi}
The range of $\gamma!_{\pi}$ is thus the genus of $G$. For future use, we note the following.

\begin{prop}\label{prop:increasing-gamma} Let $\pi'$ be an ordered partition of $V$ which is a coarsening of $\pi$, that is, each part of $\pi$ is included in a part of $\pi'$. Then, we have $\gamma!_{\pi'}\geq\gamma!_{\pi}$.
\end{prop}
\begin{proof}
This is the second statement in~\cite[Thm.~6.1]{AEG-residue-polytope}. 
\end{proof}

\subsection{Metric graphs}
Given an edge length function $\ell \colon E \to\R_{>0}$ on the edge set of the graph $G=(V,E)$, we denote by $\ell_e$ the real number $\ell(e)$ associated to each $e\in E$. We get a metric graph $\Gamma$ by plugging an interval of length $\ell_e$ between the two extremities of each edge $e$. It is a metric space endowed with the path metric. The pair $(G,\ell)$ is called a \emph{model} for $\Gamma$.

We extend $\ell$ to $\E$ by putting $\ell_a\coloneqq \ell_e$ for each arrow $a\in \E$ on an edge $e\in E$. 

\subsection{Divisor theory on metric graphs}
Let $\Gamma$ be a metric graph with model $(G,\ell)$. Denote by $\Div(\Gamma)$ the group of divisors on $\Gamma$, which is by definition the free Abelian group generated by the points on $\Gamma$. Write $(x)$ for the generator associated to each  point $x\in \Gamma$, so that 
\[
\Div(\Gamma)=\Bigl\{\sum_{x\in S} \ssn_x (x) \,\Big|\, S\subseteq\Gamma\textrm{ finite and }\ssn_x\in \mathbb Z\,\text{ for all } x\in S \Bigr\}.
\]
 
For each divisor $D\in \Div(\Gamma)$ and $x\in \Gamma$, denote by $D(x)$ the coefficient of $(x)$ in $D$, and by $\supp(D)$ the support of $D$, that is, the set of points $x$ with $D(x)\neq 0$. 
The degree of $D$, denoted $\deg(D)$, is defined as 
\[
\deg(D)\coloneqq \sum_{x\in \Gamma} D(x).
\]

A \emph{rational function} on $\Gamma$ is a continuous function $h\colon \Gamma \to \R$ whose restriction to each edge of $\Gamma$ is piecewise affine with integral slopes.  The set of rational functions on $\Gamma$ is denoted $\Rat(\Gamma)$. The \emph{order of vanishing} at $x\in\Gamma$ of a rational function $h\in\Rat(\Gamma)$, denoted $\ord_x(h)$, is the sum of the incoming slopes of $h$ at $x$. The divisor $\div(h)$ associated to $h\in \Rat(\Gamma)$ is then defined by
\[
\div(h)\coloneqq \sum_{x\in \Gamma} \ord_x(h) (x).
\]
Elements of this form in $\Div(\Gamma)$ are called \emph{principal}. Two divisors $\ssD_1$ and $\ssD_2$ are called \emph{linearly equivalent}, and we write $\ssD_1 \sim \ssD_2$, if the difference $\ssD_1 -\ssD_2$ is principal.

\subsection{Admissible divisors on metric graphs}\label{sec:admdivtilpoly}
We briefly review the theory of admissible divisors on metric graphs from our work~\cite{AE-modular-polytopes}.

Let $\Gamma$ be a metric graph with model $(G,\ell)$. A divisor $D$ on $\Gamma$ is \emph{$G$-admissible}, or simply \emph{admissible} if $G$ is clear, if for each edge $e$ of $G$, we have $D(x)=0$ for every point $x$ on the interior of $e$ with at most one possible exception for which we then have $D(x)=1$. 
For each collection $\cl$ of divisors on $\Gamma$, denote by $\stable!_G(\cl)$ the subset of all $G$-admissible divisors $D \in \cl$.

A divisor $D$ on $\Gamma$ supported on the vertex set $V$ is admissible. Moreover, if $D'$ is another admissible divisor on $\Gamma$, then the sum $D+D'$ is also admissible.

\subsubsection{Integer slopes}

For each function $h\colon V \rightarrow \R$, define $\slztwist{\ell}{}h \colon \E \to \Z$ by setting
\[
\slztwist{\ell}{}h(a)\coloneqq \Big\lfloor \frac{h(\te_a)-f(\he_a)}{\ell_a}\Big\rfloor\quad\text{for each }a\in\E.
\]
We call $\slztwist{\ell}{}h(a)$ the \emph{incoming integer slope of $h$ at} $\te_a$ \emph{along $a$}. It is actually the incoming slope of a rational function on $\Gamma$ extending $h$; see Subsection~\ref{sec:adm}.

The following is immediate.

\begin{prop}\label{prop:slope-function} For each function $h\colon V \to \R$, 
the sum $\slztwist{\ell}h(a) + \slztwist{\ell}h(\bar a)$ is either $0$ or $-1$, depending on whether $\ell_a$ divides $h(\te_a)-h(\he_a)$ or not, respectively.
\end{prop}

\subsubsection{Admissible extensions of functions}\label{sec:adm}

For each divisor $D$ on $\Gamma$ and function $h\colon V \rightarrow \R$, \cite[Thm.~5.5]{AE-modular-polytopes} guarantees the existence and uniqueness of a rational function $\hat h \colon \Gamma \to \R$ extending $h$ for which $D+\div(\hat h)$ is admissible. The extension $\hat h$ is the same for all divisors $D$ supported on the vertices of $G$.
We restate the theorem below only in this case, the only one we need in this paper. (We will set $D$ to be $K$, the canonical divisor; see \eqref{eq:canonical}.) 

For each arrow $a=uv \in \E$ of $G$ on each edge $e=\{u,v\}$ and $t\in[0,\ell_e]$, 
denote by $x^a_t$ the point of $\Gamma$ on $e$ at a distance $t$ from $u$ along $e$. Thus, $x^a_t=x^{\bar a}_{\ell_e-t}$.

\begin{thm}[Admissible extension of a function]\label{thm:existence-uniqueness-canonical-extension}
For each real-valued function $h\colon V\to\R$, there is a unique rational function $\hat h\colon \Gamma \rightarrow \R$ extending $h$ such that $\div(\hat h)$ is admissible. Furthermore, $\hat h$ is characterized by the following properties:
\begin{enumerate}
\item[(1)] The incoming slope of $\hat h$ at each $v\in V$ along each arrow $a\in\E_v$ is $\slztwist{\ell}h(a)$. 
\item[(2)] For each arrow $a\in \E$ and $t\in (0, \ell_a)$, the divisor $\div(\hat h)$ takes value $0$ at $x_t^a$, unless $\ell_a-t$ is the remainder of the Euclidean division\footnote{For $n, l \in \R$ with $l>0$, writing $n = \lfloor\frac{n}{l}\rfloor +r$ with $0\leq r<l$, $r$ is the rest of the Euclidean division of $n$ by $l$.} of $h(\te_a) -h(\he_a)$ by $\ell_a$, in which case it takes value $1$.
\end{enumerate}
\end{thm}

\begin{defi}\label{defi:adm-ext} We call the function $\hat h$ described above the \emph{admissible extension} of $h$, and denote by $\ssD_h$ the admissible divisor $\div(\hat h)$.
\end{defi}

The following result gives a description of the principal divisor $\ssD_h=\div(\hat h)$.

\begin{prop}\label{prop:adm}
Let $h\colon V \to\R$ be a function, and $\hat h\colon\Gamma\to \R$ its admissible extension. Then, for each $x\in\Gamma$, we have
\begin{enumerate}
 \item[(1)] If $x\in V$, then $\ssD_h(x) = \sum_{a\in\E_x}\,\slztwist{\ell}h(a)$.

 \item[(2)] If $x=x_r^{\bar a}$ for an arrow $a=uv$, where $r$ is the remainder of the Euclidean division of $h(u)-h(v)$ by $\ell_a$, and $r\neq 0$, then $\ssD_h(x)=1$.
 \item[(3)] Otherwise, $\ssD_h(x)=0$.
\end{enumerate}
\end{prop}

\begin{thm}[Characterization of admissible divisors] \label{thm:admissible-divisors} Let $D$ be a divisor on $\Gamma$ supported on the vertices of $G$. Then, each admissible divisor $D'$ linearly equivalent to $D$ is of the form $D'=D+\div(\hat h)$ for the admissible extension $\hat h$ of a function $h \colon V \to \R$. Furthermore, if $\Gamma$ is connected, $h$ is unique up to addition by constants.
\end{thm}

The upshot is a bijection between admissible divisors in the linear equivalence class of $D$ and level functions $h\colon V \to \R$ up to addition by constants.

\subsection{Slope functions}\label{sec:slope}

The following definition is motivated by Proposition~\ref{prop:slope-function}.

\begin{defi}\label{defi:slopefcn} Let $G=(V, E)$ be a graph. A function $s\colon \E \to \Z$ is called a \emph{slope function} if \[
s(a)+s(\bar a) \in \{-1, 0\}\quad\text{for each arrow }a\in\E.
\]
Given a slope function $s\colon \E \to \Z$, we denote by $\ssA_s$ the set of arrows $a\in \E$ with $s(\bar a)<0$ (equivalently, either $s(a)>0$, or $s(a)=0$ and $s(a)+s(\bar a)=-1$), and call the arrows in $\ssA_s$ \emph{upward} for $s$. The reverse arrows are called \emph{downward} for $s$; their set is denoted $\barA_s$. All of these arrows and their underlying edges are called \emph{vertical} for $s$; the remaining arrows and edges are called \emph{horizontal} for $s$. If $s(a)+s(\bar a)=0$, we say that $a$ and the underlying edge are \emph{integer} for $s$. Denote by $\ssE_s^{\rmint}$ (resp.~$\ssA_s^{\rmint}$) the set of integer vertical edges (resp.~integer upward arrows).
\end{defi}

\begin{remark}\label{rem:slope-level-pair}
By Proposition~\ref{prop:slope-function}, for each edge length function $\ell\colon E \to\R_{>0}$, and function $h\colon V \to \R$, the function $\slztwist{\ell}{}h\colon \E \to \Z$ is a slope function. Furthermore, if $s=\slztwist{\ell}{}h$ and $\pi=\sspi_h$, the level structure induced by $h$, then $\A_s=\A_{\pi}$ and $\barA_s=\barA_{\pi}$. Also, $a=uv\in\A_{\pi}$ is integer for $s$ if and only if $\ell_a$ divides $h(u)-h(v)$.  
\end{remark}

\begin{defi}\label{defi:zetas} For each slope function $s \colon \E \to \Z$, define $\zeta!_s \colon \ssub{2}!^V \to \Z$ by
\[
\zeta!_s(I) \coloneqq |\ssA_s(\ssI^c, \ssI)|+\sum_{a \in \E(I, \ssI^c)}s(a) \qquad \forall\, I \subseteq V.  \qedhere
\]
\end{defi}

\begin{remark} If $s=\slztwist{\ell}h$ for some edge length function $\ell\colon E\to\R_{>0}$ and function $h\colon V\to\R$, then the function $\zeta!_s$ bears some similarity with the function $\mu!_{\ssD_h}^*$ defined in \cite[\S 4.2]{AE-modular-polytopes}, which satisfies
\[
\mu!_{\ssD_h}^*(I) = \frac{|\E(\ssI^c,\ssI)|}{2}+\sum_{a \in \E(I, \ssI^c)}\slztwist{\ell}h(a) \qquad \forall\, I \subseteq V. 
\]
We know $\mu!_{\ssD_h}^*$ is submodular by \cite[Prop.~4.3]{AE-modular-polytopes}. Remarkably, so is $\zeta!_s$ by Proposition~\ref{prop:zetah}.
\end{remark}

\begin{lemma}\label{lem-useful-s} For each slope function $s \colon \E \to \Z$ and $I\subseteq V$,
\[
\sum_{a \in \E(I, \ssI^c)}s(a)=\sum_{a \in \E_I}s(a)+|\ssA_s(I)-\ssA_s^{\rmint}(I)|.
\]
\end{lemma}

\begin{proof} Immediate from the fact that $s(a)+s(\bar a)$ is $0$ if $a$ is integer for $s$, and $-1$ otherwise.
\end{proof}

\begin{prop}\label{prop:zetah} For each slope function $s\colon \E \to \Z$, the functions $\zeta!_s-\varphi!_1$ and $\zeta!_s-\varphi!_2$ are modular, where $\varphi!_1(I)\coloneqq|\ssA^{\rmint}_s(\ssI^c,I)|$ and $\varphi!_2(I)\coloneqq-|\ssA^{\rmint}_s(I)|$ for each $I\subseteq V$. In particular, $\zeta!_s$ is submodular.
\end{prop}

\begin{proof} For each $I\subseteq V$, from Lemma~\ref{lem-useful-s} we get
\[
\zeta!_s(I)-\varphi!_2(I)=
\sum_{a \in \E_I}s(a)+|\ssA_s(I)|+|\ssA_s(\ssI^c, \ssI)|.
\]
Clearly, the function $I\mapsto \sum_{a \in \E_I}s(a)$ is modular. It thus follows from Proposition~\ref{prop:modular-digraph} that $\zeta!_s-\varphi!_2$ is modular. The same proposition yields that $\varphi!_1-\varphi!_2$ is modular, whence $\zeta!_s-\varphi!_1$ is modular. The last statement follows from the fact that $\varphi!_1$ and $\varphi!_2$ are submodular, again by Proposition~\ref{prop:modular-digraph}.
\end{proof}

We state the following monotonicity property, which will be useful later in Section~\ref{sec:bricks-fans}.

\begin{prop}\label{prop:monotonicity-zeta} Let $s,s'\colon \E \to \Z$ be slope functions with $s'\geq s$. Then, $\zeta!_{s'} \geq \zeta!_s$.
\end{prop}

\begin{proof}
Since $s'\geq s$, we only need to prove that for each subset $I\subseteq V$ and arrow $a\in \ssA_s(\ssI^c,I)$ such that $a\not\in A_{s'}(\ssI^c,I)$, we have $s'(\bar a)>s(\bar a)$. But this is clear: since $a\in A_s(\ssI^c,I)$, we have $s(\bar a)<0$, and since $a\not\in A_{s'}(\ssI^c,I)$, we have $s'(\bar a)\geq 0$. Hence, $s'(\bar a)>s(\bar a)$, as claimed.
\end{proof}


\section{Tropicalization}\label{sec:tropicalization}
We discuss in this section tropicalization of differentials and of canonical series. We state in particular the analogue to the theorem proved in our work~\cite{AE-modular-polytopes} that tropicalization of canonical linear series gives rise to tilings of simplices by polytopes.

We fix an algebraically closed field $\K$ which is complete for a nontrivial non-Archimedean valuation $\valuation$. For instance, $\K$ could be obtained from a discrete valuation ring (like the ring of power series in one variable over a field) by taking the completion of the algebraic closure of its field of fractions. Put $\valgroup=\valuation(\K\setminus \{0\})$, the value group of $\valuation$. Denote by $\varR$ the valuation ring of $\valuation$, and by $\varm$ and $\k$ the maximal ideal and the residue field of $\varR$, respectively.

Let $\varX$ be a smooth proper connected curve over $\K$. Denote by $\varX^{\an}$ the Berkovich analytification of $\varX$. Recall that the points on $\varX^{\an}$ are in bijection with the union of the (closed) points on $\varX$ and the valuations on the function field $\K(\varX)$ that extend the valuation of $\K$. For each point $x$ of $\varX^{\an}$, denote by $\valuation_x$ the corresponding valuation, by $\varR_x$ its valuation ring, and by $\k(x)$ its residue field.
Each rational function $\varf$ on $\varX$ gives rise to an evaluation map $\ev_{\varf} \colon \varX^{\an} \to \R\cup\{\pm \infty\}$, that sends $x$ to $\valuation_x(\varf)$.

We fix a \emph{semistable vertex set} $V$ for $\varX^{\an}$. It is a finite set of type 2 points on $\varX^{\an}$ whose complement $\varX^{\an} \setminus V$ is a disjoint union of finitely many open annuli and infinitely many open disks. The semistable vertex set $V$ defines a \emph{skeleton} $\Gamma$ for $\varX^{\an}$, the union in $\varX^{\an}$ of $V$ and the skeleta of the open annuli in $\varX^{\an} \setminus V$. The canonical metric on each of these skeleta identifies it with an open interval whose closure has one or two boundary points in $V$. This way, the semistable vertex set defines a \emph{graph} $G=(V, E)$ and an \emph{edge length function} $\ell\colon E\to\R_{>0}$, and we may view the skeleton $\Gamma$ as the metric graph associated to the pair $(G, \ell)$ canonically embedded in $\varX^{\an}$.  Denote by $\tau\colon \varX^{\an} \to \Gamma$ the retraction map from $\varX^{\an}$ to $\Gamma$. 

For each point $x$ of type 2 on $\varX^\an$, the field $\k(x)$ is of transcendental degree one over $\k$; let $\varC_x$ be the corresponding smooth proper curve over $\k$. Each unit outgoing tangent vector $\nu \in \TT_x\Gamma$ gives a point on $\varC_x$  denoted $\ssp^{\nu}$. If $x\in V$, and $\nu$ is parallel to the arrow $a=xv\in \E$, then we denote this point by $\ssp^{a}$.

Gluing the curves $\varC_v$ for $v\in V$, by identifying the points $\ssp^a$ on $\varC_v$ and $\ssp^{\bar a}$ on $\varC_w$ for all arrows $a=vw\in \E$, gives a nodal curve $X$ whose dual graph is $G$. We say that the pair $(X,\Gamma)$ is obtained as a \emph{tropicalization} of $\varX$.

\subsection{Tropicalization of linear series and tilings of simplices} \label{sec:tiling-divisor} 

We first briefly discuss tropicalization of divisors and spaces of linear series. 

Let $\varD$ be divisor of degree $d$ on $\varX$. Let $\varH$ be a vector subspace of dimension $r+1$ of the Riemann--Roch space of $\varD$:
\[
\varH \subseteq \bigl \{\varf \in\K(\varX) \,\st\, \div(\varf) + \varD \geq 0\bigr\}.
\]
The pair $(\varD, \varH)$ defines a linear series $\g!_{d}^r$ on $\varX$.

The tropicalization of $(\varD, \varH)$ gives rise to a set of rational functions on $\Gamma$, as follows.  Let $\tau(\varD) = \sum_{x\in \varX} \varD(x)\tau(x)$ be the divisor on $\Gamma$ obtained by the tropicalization of $\varD$.   For each nonzero rational function $\varf \in\K(\varX)$, we denote by $\trop(\varf)$ the restriction to $\Gamma$ of the evaluation map $\ev_{\varf} \colon \varX^{\an} \to \R\cup\{\pm \infty\}$. This is a rational function on $\Gamma$ and by Specialization Lemma~\cite[Thm.~4.5]{AB15}, we have 
\[
\tau(\varD) + \div(\trop(\varf)) = \tau(\varD + \div(\varf)).
\]

We enrich the tropicalization by adding a collection of $\k$-vector subspaces $\VS_v\subseteq\k(v)$ associated to the vertices $v$ of $G$ as follows.  First, we choose a homomorphism $\ba\colon\valgroup\to\K\setminus\{0\}$ that is a section of the valuation $\valuation\colon\K\setminus \{0\} \to \valgroup$. Since $\valgroup$ is divisible, a section exists. Put $\ba_{\lambda}:=\ba(\lambda)$ for $\lambda\in\valgroup$. Then, for each vertex $v\in V$, the nonzero elements of the $\k$-vector space $\VS_v$ are the \emph{reductions} of the nonzero functions $\varf \in \varH$, which are the image in the residue field $\k(v)$ of $\ba_{\valuation_v(\varf)}^{-1}\varf \in \varR_v$; see~\cite[\S 4.4]{AB15}. It is proved in~\cite[Lem.~4.3]{AB15} that $\VS_v$ is a $\k$-vector space of dimension $r+1$.  We associate to the linear series $(\varD, \varH)$ the $\k$-vector space 
\[
\VS \coloneqq \bigoplus_{v\in V} \VS_v.
\]

Denote by $\cl$ the linear equivalence class of $\tau(\varD)$ on $\Gamma$. Let $\cl_{\valgroup}$ be the subset of $\cl$ consisting of all $\valgroup$-rational divisors. Recall that a divisor is $\valgroup$-rational if its support consists only of $\valgroup$-rational points, that is, points whose distances to vertices of $V$ are all in $\valgroup$.  

Let $D\in\cl_{\valgroup}$ be a $G$-admissible divisor. In \cite{AE-modular-polytopes}, we associate to $D$ an $(r+1)$-dimensional subspace of $\VS$ as follows. Since $D\in\cl_{\valgroup}$, we can write $D=\div(h)+\tau(\varD)$ for a rational function $h\in \Rat(\Gamma)$ such that $h(V)\subseteq\valgroup$. Denote by $\varM_h$ the space of all functions $\varf\in \varH$ such that $\trop(\varf)(v) \geq h(v)$ for $v\in V$, that is, 
\[
\varM_h \coloneqq \bigl\{\varf\in \varH \,\st\, \trop(\varf)(v)  \geq h(v)\quad \forall v\in V\bigr\}.
\]
Then, we let $\ssW_h$ be the image of $\red!_h\colon\varM_h\to\VS$, defined by 
letting $\red!_h(\varf)_v$ be the reduction of $\ba_{h(v)}^{-1}\varf$ in $\k(v)$ for each $\varf\in\varH$ and $v\in V$; this is \emph{reduction relative to $h$}. The space $\ssW_h$ has dimension $r+1$ by \cite[Lem.~7.9]{AE-modular-polytopes}. It depends only on $D$, rather than on $h$. This is the $(r+1)$-dimensional subspace of $\VS$ associated to $D$. But since in this paper we are viewing the function $h$ as a level function, we prefer to use the notation $\ssW_h$ in the sequel.

We may now associate to $D$ the submodular function $\nu!_h^* \colon \ssub{2}!^V \to \R$ given by $\nu!^{*}_{h}\coloneqq \nu!_{\ssW_{h}}^*$; see Proposition~\ref{prop:submodular-subspace}. Consider as well the corresponding polytope
\[
\Q_h\coloneqq \P_{\nu^{*}_{h}}=\bigl\{q\in\R^V\,\st \, q(I)\leq\nu!_h^*(I)\,\,\forall I\subseteq V, \text{ with equality if }I=V\bigr\}.
\]
It is contained in the standard simplex $\Delta!_{r+1} \subset \R^V$ consisting of all the points $q\in\R_{\geq0}^V$ with $q(V)=r+1$. We proved in~\cite[Thm.~8.4]{AE-modular-polytopes} the following theorem.

\begin{thm}\label{thm:tiling-simplex} The collection of polytopes $\Q_h$ associated to $G$-admissible divisors in $\cl_\valgroup$, where $\cl$ is the linear equivalence class of the divisor $\tau(\varD)$, form a tiling of the standard simplex $\Delta!_{r+1}$.
\end{thm}

\subsection{Tropicalization of differentials}
\label{sec:trop-diff}
For each point $x$ of type 2 on $\Gamma \subset \varX^{\an}$, denote by $\Omega!_x$ the space of meromorphic differentials of $\varC_x$, and by $\g(x)$ its genus. Extending $\g$ by zero over all the other points on $\Gamma$, we get the \emph{genus function} $\g\colon \Gamma \to \Z_{\geq 0}$. We refer to $(\Gamma, \g)$ as an \emph{augmented metric graph}. By an abuse of notation, we denote as well by $\g$ the restriction of the genus function to $V$, and refer to the pair $(G,\g)$ as an \emph{augmented graph}.

Denote by $K$ the canonical divisor of $(\Gamma,\g)$, defined by
\begin{align}\label{eq:canonical}
K \coloneqq \sum_{x \in \Gamma} (2 \g(x) - 2 + \val(x)) \, (x) = \sum_{v \in V} (2 \g(v) - 2 + \val(v)) \, (v),   
\end{align} 
where $\val(x)$ denotes the valence of $x$ in $\Gamma$. Denote by $\clcan$ the linear equivalence class of $K$, and by $\clcan_\valgroup\subseteq\clcan$ the subset consisting of $\valgroup$-rational divisors. Since the edge lengths in $G$ are all in $\valgroup$, the vertices of $G$ are $\valgroup$-rational, and thus $K\in\clcan_\valgroup$.

Denote by $\Omega!_{\varX}$ the space of meromorphic differentials of $\varX$, and by $\varH$ the $g$-dimensional subspace of Abelian differentials. Let $\|{\cdot}\|$ be the K\"ahler norm on $\Omega!_{\varX}$ that at each point $x\in \varX^{\an}$ associates to each $\varalpha\in\Omega!_{\varX}$ a value $\|\varalpha\|_x \in \R_{\geq 0}\cup\{+\infty\}$. We refer to Temkin~\cite[\S 6]{Tem16} for the definition of 
$\|{\cdot}\|$ in the general case, and to~\cite{TT22} in the case of curves. For each $\varf\in \K(\varX)$ and each point $x\in \varX^{\an}$, the norm satisfies
\[
\|\varf \varalpha \|_x = e^{-\valuation_x(\varf)}\|\varalpha\|_x.  
\]
Also, $\|\varalpha\|_x$ is 0, resp. $+\infty$, if $x$ is a closed point of $\varX$ which is a zero, resp.~a pole, of $\varalpha$. 

For each point $x$ of type 2 on $\varX^{\an}$, the meromorphic differentials $\varalpha\in\Omega!_{\varX}$ with $\|\varalpha \|_x\leq 1$ form thus a module over the valuation ring $\varR_x$. This module is endowed with a natural map to $\Omega!_x$, which becomes injective after tensoring with $\k(x)$. The image in $\Omega!_x$ of a differential in the module is called its \emph{reduction}. In any case, for each nonzero $\varalpha\in\Omega!_{\varX}$, we have that $\log \|\varalpha\|_x\in\valgroup$, and hence there is $c\in\K$ such that $\|\varalpha \|_x=\exp(\valuation(c))$, or equivalently $\|c\varalpha \|_x= 1$, and hence $c\varalpha$ has nonzero reduction at $x$.

\begin{defi}[Tropicalization of differentials]
For each nonzero meromorphic differential $\varalpha$ of $\varX$, the \emph{tropicalization} of $\varalpha$ is the function
\[
\trop(\varalpha) \colon \Gamma \longrightarrow \R, \qquad x \longmapsto -\log \|\varalpha\|_x. \qedhere
\]
\end{defi}

\begin{thm}\label{thm:specialization-form}
The tropicalization $\trop(\alpha)$ of  a meromorphic differential $\varalpha\in\Omega!_{\varX}$ is a rational function on $\Gamma$. It verifies the slope formula,
\[ \div(\trop(\varalpha)) + K = \tau(\div(\varalpha)), \]
where $\div(\alpha)$ is the divisor of $\alpha$ on $\varX$.
\end{thm}

\begin{proof}
The first part is a consequence of~\cite[Thm.~8.2.4]{Tem16}. The second is proved in~\cite[Prop.~2.3.3]{TT22} using \cite[Lem.~3.3.2]{BT20} (see as well~\cite{Ami-W, KRZ16, BN16}).
\end{proof}

\begin{defi}\label{defi:reducction-form-x} Let  $x \in \varX^{\an}$ be a point of type 2. Let $\varalpha$ be a meromorphic differential on $\varX$. We denote by  $\wt\varalpha_x$ the reduction of $\varalpha$ at $x$, defined as the  restriction of $\ba_{\lambda}^{-1}\varalpha$ to $\varC_x$ for $\lambda\coloneqq \trop(\varalpha)(x)$.  
 \end{defi}

Note that $\wt\varalpha_x$ is a meromorphic differential on the curve $\varC_x$. It changes by a scalar when the section $\ba$ of $\valuation$ is changed, so it is well-defined up to scaling by an element of $\k^\times$.  For the definition of reduction, we refer to~\cite[\S 2]{TT22}. 
We have the following result \cite[Lem.~3.3.2]{BT20}.

\begin{lemma}[Slope Formula] \label{lem:slope-formula} Let $x$ be a point of type $2$ on $\Gamma$ and $\varalpha$ a meromorphic differential on $\varX$. For each unit outgoing tangent vector $\nu\in \TT_x\Gamma$ at $x$, the order of vanishing of $\wt \varalpha_x$ at the point $\ssp^{\nu}$ is given by $-1 + \slp_{\nu} F$, where $F=\trop(\varalpha)$ and $\slp_{\nu} F$ is the outgoing slope of $F$ at $x$ along $\nu$.
\end{lemma}

For each vertex $v\in V$, let $\VS_v$ be the vector subspace of $\Omega!_v$ spanned by the meromorphic differentials on $\varC_v$ which are reductions $\widetilde \varalpha_v$ at $v$ of Abelian differentials $\alpha\in\varH$.  This is a vector space of dimension $g$, as in~\cite[Lem.~4.3]{AB15}.
Put
\[
\VS \coloneqq \bigoplus_{v\in V} \VS_v \subset \Omega\coloneqq \bigoplus_{v\in V} \Omega!_v.
\]
The spaces $\VS$ and $\Omega$ are endowed with an action of $\k^V\coloneqq \prod_{v\in V}\k$. For each element $\psi\in\k^V$ and $v\in V$, we denote by $\psi_v$ the element in $\k$ given by the $v$-component of $\psi$. Each $\psi\in\k^V$ corresponds to an endomorphism of $\VS$, respectively, $\Omega$, given by componentwise multiplication by $\psi_v$, $v\in V$, which, abusing notation, we will also denote by $\psi$.

\subsection{Reductions of differentials} 
\label{sec:red-diff}

We explain how to associate a $g$-dimensional subspace of $\VS$ to each $\valgroup$-rational $G$-admissible divisor on $\Gamma$ linearly equivalent to  $K$. 

Recall that $\clcan$ is the linear equivalence class of $K$, that $\clcan_{\valgroup}\subseteq\clcan$ is the subset of $\valgroup$-rational divisors, and $\stable_G(\clcan_{\valgroup})\subseteq\clcan_{\valgroup}$ is the subset of $G$-admissible divisors.

Notice that $K$ is $\Lambda$-rational because all the edge lengths in $G$ are in $\valgroup$. Let $\adfun_{G}(K;\valgroup)$ be the set of all functions $h\in \Rat(\Gamma)$ such that $\div(h) +K$ is $G$-admissible and $h(V)\subseteq\valgroup$. Then, $\stable_G(\clcan_{\valgroup})$ is the set of divisors of the form $\div(h) +K$ for $h\in\adfun_{G}(K;\valgroup)$.

By Theorem~\ref{thm:existence-uniqueness-canonical-extension}, there is a bijection from $\adfun_{G}(K;\valgroup)$ to the collection of maps $V \to \valgroup$. It is given by taking $h\in \adfun_{G}(K;\valgroup)$ to $h\rest{V}$, with inverse taking $h\colon V \to \valgroup$ to its admissible extension $\hat h$. Abusing notation, we will not distinguish $h\in \adfun_{G}(K;\valgroup)$ from
$h\rest{V}$.

For each $h\colon V\to\valgroup$, let $\varM_h$ be the space of all Abelian differentials $\varalpha\in\varH$ such that $\trop(\varalpha)(v) \geq h(v)$ for every $v\in V$. Recall that we fixed a section $\ba\colon\valgroup\to\K\setminus\{0\}$ of the valuation $\valuation$. Notice that $\trop(\ba_{h(v)}^{-1}\varalpha)(v)\geq 0$
for each $\varalpha \in \varM_h$.

\begin{defi}[Relative reduction map $\red!_h$] \label{def:rel-red} The \emph{reduction map relative to $h$} is the map 
\[
\red!_h \colon \varM_h \to \VS
\]
defined as follows: For each $\varalpha \in\varM_h$, let $\red!_{h}(\varalpha)$ be the element of $\VS$ whose component $\ssub{(\red!_h(\varalpha))}!_v$ is the reduction of $\ba_{h(v)}^{-1}\varalpha$ to the curve $\varC_v$ for each $v\in V$. 
\end{defi} 

Note that $\ssub{(\red!_h(\varalpha))}!_v=0$ if $\trop(\varalpha)(v)>h(v)$; otherwise, $\ssub{(\red!_h(\varalpha))}!_v = \wt \varalpha_v$.

\begin{defi}[Reduction of $\varH$ relative to $h$]\label{defi:rel-red-space} We define the reduction of the space of Abelian differentials $\varH$ relative to $h$ as the subspace $\ssW_h \subset \VS$ consisting of the reductions of $\varalpha \in \varM_h$ relative to $h$, that is, 
 \[
 \ssW_h \coloneqq \bigl\{\, \red!_h(\varalpha) \,\st\, \varalpha\in \varM_h\,\bigr\}\subseteq\VS\subset\Omega.
 \] 
We call $\ssW_h$ the \emph{limit canonical series} associated to $h$.\qedhere
\end{defi}

Note that for each $\lambda\in\valgroup$, the space $\ssW_{h+\lambda}$ coincides with $\ssW_h$. This implies that $\ssW_h$ depends only on $D=K +\div(h)$. However, since we will view $h$ as a level function, we will mainly use the notation $\ssW_h$ in what follows. 

\begin{remark}\label{rem:trop-diff} If we choose a nontrivial Abelian differential $\varalpha$ of $\varX$, and set $\varD\coloneqq \div(\varalpha)$, then we are in the setting of Section~\ref{sec:tiling-divisor} and we get an identification  of the spaces $\ssW_h$ defined above with the spaces $\ssW_{h'}$ from Section~\ref{sec:tiling-divisor}, with $h'$ an appropriate translation of $h$ using $\trop(\varalpha)$. Nevertheless, it is fundamental that we treat the $\ssW_h$ as spaces of meromorphic differentials for the role residues will play later.
\end{remark} 

A different choice of the section $\ba$ produces a subspace which might be different from $\ssW_h$, but which lies on the same orbit of $\VS$ by the natural action of $\Gm^{\ssV}\coloneqq \prod_{v\in V}\k^\times$. More precisely, let $\pi$ be the ordered partition associated to $h$, and let $\Gm^{\ssV_{\pi}} \coloneqq \prod_{n\in \ssV_\pi}\k^\times \subseteq \Gm^{\ssV}$ be the subgroup consisting of those elements $\psi\colon V\to\k^*$ which factor through the quotient map $V\to\ssV_{\pi}$, that is, $\psi_u=\psi_v$ whenever $u$ and $v$ belong to the same level, that is, whenever $h(u)=h(v)$.

\begin{prop}\label{prop:orbit}
 The orbit $\Gm^{\ssV_{\pi}} \ssW_h$ is well-defined, that is, it does not depend on the choice of the section $\ba$.
\end{prop}

\begin{proof} Let $\ba'$ be another section of the valuation map, so that $\ba_{\lambda}'\ba_{\lambda}^{-1}$ has valuation zero for every $\lambda\in\valgroup$. Let $\red_h' \colon \varM_h \to \VS$ be the reduction map with respect to this new section. For each $\lambda\in\valgroup$, let $\sspsi_{\lambda}$ be the reduction of $\ba_{\lambda}'\ba_{\lambda}^{-1}$ in $\k$. Then, for each $\alpha \in \varM_h$, we have $\ssub{(\red_h'(\alpha))}!_v = \psi_{h(v)}^{-1}\ssub{(\red_h(\alpha))}!_v$ for each $v\in V$, and the statement follows.
\end{proof}

\subsection{The fundamental collection} 
Assume $g>0$. For each $h\colon V\to\valgroup$, we say that the subspace $\ssW_h\subseteq\Omega$ is \emph{simple} if the corresponding submodular function $\nu!_{\ssW_h}^*$ is simple.

\begin{defi}[Fundamental collection of limit canonical series]\label{def:FC}  Notation as above, we denote by $\FC(\varX)$ the collection of simple subspaces $\ssW_h\subseteq\Omega$, and call it the \emph{fundamental collection of limit canonical series} on $X$ associated to $\varX$. 
\end{defi}

By the finite generation property~\cite[Thm.~9.6]{AE-modular-polytopes}, the spaces in the fundamental collection generate all the other limit canonical series $\ssW_h$. 

\begin{thm}\label{thm:W-tiling} The elements $\ssW_h$ of $\FC(\varX)$ are precisely those limit canonical series whose corresponding polytope $\P_{\nu!_{\ssW_h}^*}$ is full-dimensional. They and their faces form a tiling of the standard simplex $\Delta!_g$ in $\R^V$. 
\end{thm}
  
\begin{proof} The first statement follows from \cite[Prop.~2.7]{AE-modular-polytopes} (see \S~3.1 in \emph{loc.~cit.}). The second now follows directly from Theorem~\ref{thm:tiling-simplex}, using Remark~\ref{rem:trop-diff}.
\end{proof}


\section{Residue conditions}\label{sec:twisted-differentials}

We fix a nodal projective connected reduced curve $X$ defined over an algebraically closed field $\k$. Let $G=(V,E)$ be its dual graph and $\E$ its arrow set. Let $\varC_v$ be the normalization of the component of $X$ associated to each $v\in V$. Each edge $e=\{u,v\}$ in $E$ gives a point on $\varC_u$, denoted $p^a$, and a point on $\varC_v$, denoted $p^{\bar a}$, where $a=uv$ and $\bar a=vu$ are the arrows on $e$, whose identification is the node $p^e$ on $X$ corresponding to $e$. For each $v\in V$, denote by $\Omega!_v$ the space of meromorphic differentials on $\varC_v$.  Put $\Omega\coloneqq\bigoplus\Omega!_v$.

\subsection{The residue map}\label{sec:residue-map-Res} For each vertex $v$ of $G$
and each point $p$ on $\varC_v$, denote by 
\[
\res_p \colon \Omega!_v \to \k
\]
the residue map at $p$, associating to a meromorphic differential $\alpha$ on $\varC_v$ its residue at $p$. For each arrow $a \in \E_v$, denote by $\res_a$ the residue map on the curve $\varC_v$ associated to the point $\ssp^{a}$:
\[
\res_a  \coloneqq \res_{p^a}\colon \Omega!_v \to \k.
\]

We put all the residue maps $\res_a$ for $a\in \E$ into a single map denoted by $\Res$ as follows. Consider $\bV=\k^{\E} = \bigoplus_v \bV_v$, $\bV_v=\k^{\E_v}$ for $v\in V$, from Section~\ref{sec:residue-spaces-and-conditions}. The sum of $\res_a$ for $a\in \E_v$ defines a map 
\[
\Res_{v}\coloneqq\sum_{a\in \E_v} \res_a \colon \Omega!_v \to \bV_v.
\]
Finally, the residue maps $\Res_{v}$ for $v\in V$ define altogether
\begin{equation}\label{eq:residue-map-diff}
\Res \coloneqq\sum_{v\in V}\Res_v \colon \Omega \to \bV.
\end{equation}

\subsection{Tropicalization and the residue conditions}

Let $\varX$ be a smooth proper connected curve over an algebraically closed field $\K$ which is complete for a nontrivial non-Archimedean valuation $\valuation$ with value group $\valgroup\subseteq\R$ and residue field $\k$. Assume that $(X,\Gamma)$ is a tropicalization of $\varX$ for a metric graph $\Gamma$ with model $(G,\ell)$, for $\ell\colon E\to\valgroup$ an edge length function. 

For each $h\colon V\to \valgroup$, we consider the corresponding level graph $(G,h)$. Let  $\pi=\sspi_h$ and recall the subspace $\GlobSp_\pi\subseteq\bV$ from Section~\ref{sec:residue-spaces-and-conditions} given by residue conditions~\ref{cond:R1}-\ref{cond:R2}-\ref{cond:R3}-\ref{cond:R4}. Consider as well the subspace $\ssW_h\subseteq\Omega$ arising from reduction relative to $h$ of the space $\varH$ of Abelian differentials of $\varX$; see Section~\ref{sec:red-diff}. We have the following theorem.

\begin{thm}\label{thm:ImWh} Assume $\car\k=0$. Then, for each $h\colon V\to\valgroup$, the residue map $\Res\colon\Omega\to\k^{\E}$ takes $\ssW_h$ into the $\pi$-residue subspace $\GlobSp_\pi\subseteq\k^{\E}$, where $\sspi$ is the level structure on $G$ given by $h$.   
\end{thm}
For the proof, we review some background. The residue map $\Res$ has a counterpart on the smooth curve $\varX$ defined by Temkin and Tyomkin~\cite{TT22}:
\[ 
\bRes \colon \varH \to \K^{\E}. 
\]
This morphism sends each Abelian differential $\alpha \in \varH$ to the map 
\begin{align*}
\bRes_{\bullet}(\alpha) \colon \E &\to \K, \qquad \qquad a \mapsto \bRes_{a}(\alpha),
\end{align*}
defined on each arrow $a$ by the residue of the restriction of $\alpha$ to the oriented annulus in $\varX^{\an}$ associated to $a$. More precisely, we write $\alpha$ in a local coordinate $z$ on the oriented annulus as a formal Laurent series of the form $\sum_{j} c_j z^j \frac{\mathrm dz}z$, and the residue $\Res_a(\alpha)$ is the coefficient $c_0$. This is well-defined by Bojkovi\'c~\cite{Boj19}. 

The residue map $\bRes$ is $\K$-linear, and satisfies the following residue analogue of Slope Formula (Lemma~\ref{lem:slope-formula}): For each vertex $v$ and each arrow $a\in \E_v$,  
\[ 
\res_{p^a}(\wt \alpha_v) = \widetilde{\bRes_a(\alpha)} \qquad \qquad \forall \, \alpha \in \varH. 
\]
Here, $\wt \alpha_v$ is the reduction of $\alpha$ at $v$, and  $\widetilde{\bRes_a(\alpha)}$ is the reduction of $\bRes_a(\alpha) \in \K$, both defined using the section $\ba$ of the valuation map $\valuation$. 

By \cite[Thm.~3.1.1]{TT22}, harmonicity holds at the vertices of $G$, that is $\bRes_{\bullet}(\alpha)$ belongs to the first homology group   $H_1(G, \K)$ of $G$ with coefficients in $\K$. It follows that $\Res$ is a map from $\varH$ to $H_1(G, \K) \subseteq \K^{\E}$.

\begin{proof}[Proof of Theorem~\ref{thm:ImWh}] First, we note that all vectors in $\Res(\ssW_h)$ satisfy Conditions~\ref{cond:R1} through \ref{cond:R3} in Section~\ref{sec:residue-spaces-and-conditions}: \ref{cond:R1} is a consequence of Lemma~\ref{lem:slope-formula},~\ref{cond:R2} is the residue theorem on smooth proper curves, and~\ref{cond:R3} is the classical Rosenlicht condition (see~\cite{TT22}). Furthermore, each $\beta=(\beta_v)_{v\in V}\in\ssW_h$ satisfies $\beta=\red_h(\alpha)$ for some $\alpha\in\varM_h$. If $\beta_v\neq 0$ for every $v\in V$, then 
$\trop(\alpha)(v) = h(v)$ and $\beta_v=\ssub{\bigl(\red_h(\alpha)\bigr)}!_v =\wt\alpha_v$ for every vertex $v\in V$. In this case, the global residue conditions for $\beta$ coincide with the ones established in~\cite{BCGGM18} for the twisted differential $\ssub{(\wt\alpha_v)}!_{v\in V}$, rather obtained in our level of generality from the harmonicity of the residue map $\bRes_{\bullet}(\alpha) \colon \E \to \K$ by \cite[Rmk.~3.1.2]{TT22}. 

For $\ssW_h$ in the fundamental collection $\FC(\varX)$, since $\ssW_h$ is simple, there is $\beta=(\beta_v)_{v\in V}$ with $\beta_v \neq 0$ for all $v\in V$, and thus a general element in $\ssW_h$ has the same property. Therefore, in this case, by linearity of $\Res$, we obtain $\Res(\ssW_h) \subseteq \GlobSp_\pi$. 

By the finite generation property~\cite[Thm.~9.6]{AE-modular-polytopes}, the spaces in the fundamental collection generate all the other spaces $\ssW_h$. That is, denoting by $\bphi^{h'}_h \colon \ssW_{h'} \to \ssW_h$ the $\k$-linear map between the spaces $\ssW_{h'}$ and $\ssW_{h}$ for each $h',h\colon V\to\valgroup$ defined in~\cite[\S 8.3]{AE-modular-polytopes}, each element $\beta \in \ssW_h$ is a linear combination of elements of the form $\bphi^{h'}_h(\beta')$ for $\beta'\in \ssW_{h'}$ with $\ssW_{h'}$ in $\FC(\varX)$.
Now, as we have seen above, denoting by $\pi'=\sspi_{h'}$ the level structure induced by $h'$, we have $\Res(\beta') \in \GlobSp_{\pi'}$ for these $\beta'$. Thus, the theorem follows from Proposition~\ref{prop:compatibility-residue} below.
\end{proof}

\begin{prop}\label{prop:compatibility-residue} Let $h,h'\colon V\to\valgroup$ be functions with induced level structures $\pi$ and $\pi'$, respectively. Let $\beta'\in \ssW_{h'}$ such that $\Res(\beta') \in \GlobSp_{\pi'}$. Then $\Res(\bphi^{h'}_h(\beta')) \in \GlobSp_{\pi}$.  
\end{prop}

\begin{proof} If $\beta'_u=0$ for every minimizer of $h'-h$, then $\bphi^{h'}_h(\beta')=0$, whence the statement is clear; see \cite[\S 8.3]{AE-modular-polytopes}. We may thus suppose there is a minimizer $u$ for $h'-h$ such that $\beta'_u\neq 0$. 

Since $\ssW_{h'} =\ssW_{h'+\lambda}$ for every constant $\lambda\in \valgroup$, we may assume $h'(u)=h(u)$. Then, there is a sequence $\ssv_1,\dots,\ssv_m\in V\setminus \{u\}$ and positive constants $\lambda!_i\in \valgroup$ for $i=1, \dots, m$, such that, putting 
$\ssh_0\coloneqq h$ and $\ssh_i\coloneqq \ssh_{i-1}+\lambda!_i\one_{v_i}$ for $i=1,\dots,m$,  we have
\begin{itemize}
    \item $\ssh_m=h'$, and
    \item for each $i=1, \dots, m$, no vertex $v\in V$ satisfies $h_{i-1}(v_i)<h_{i-1}(v)< h_i(v_i)$.
\end{itemize}
(The second condition means that $v_i$ is moved from its level in $(G,h_{i-1})$ to a level below in $(G,h_i)$ not exceeding any level of $(G,h_{i-1})$ below that of $v_i$.)

 Since $\beta'_u\neq 0$, and since $\ssh_i(u) = h'(u)=h(u)$ for all $i=0, \dots, m$, we have $\phi^{\ssh_j}_{\ssh_{k}}\circ\phi^{\ssh_i}_{\ssh_{j}}\neq 0$ for all $i>j>k$. Then $\phi^{h'}_h=\phi^{\ssh_1}_{\ssh_0}\circ\phi^{\ssh_2}_{\ssh_1}\circ\cdots\circ\phi^{\ssh_m}_{\ssh_{m-1}}$ by  \cite[Prop.~8.10]{AE-modular-polytopes}. Proceeding by induction on $m$, we may thus assume that $\ssh'=\ssh+\lambda\one_w$ for some $w\in V \setminus \{u\}$ and some positive $\lambda\in \valgroup$, and that there is no vertex $v$ with $h(w)<h(v)<h(w)+\lambda$.

Let $\beta=(\beta_v)_{v\in V}\coloneqq \phi^{h'}_{h}(\beta') \in \ssW_h$. Then, $\beta_v = \beta'_v$ for $v\neq w$, and  $\beta_w=0$.

For each $n\in\valgroup$, let $\ssV_{h=n}$ (resp.~$\ssV_{h<n}$) be the set of vertices $v$ with $h(v)=n$ (resp.~$h(v)<n$). Analogously for $\ssV_{h'=n}$ and $\ssV_{h'<n}$. Since there is no $v$ with $h(w)<h(v)<h(w)+\lambda$, for each $n$ in the image of $h$, we have
\[
\bigl(\ssV_{h<n},\ssV_{h=n}\bigr)=
\begin{cases}
\bigl(\ssV_{h'<n},\ssV_{h'=n}\bigr)&\text{if }h(w)>n\text{ or } h(w)<n-\lambda,\\
\bigl(\ssV_{h'<n}\cup\{w\}, \ssV_{h'=n}\setminus \{w\}\bigr)&\text{if } h(w)=n-\lambda,\\
\bigl(\ssV_{h'<n},\ssV_{h'=n}\cup\{w\}\bigr)&\text{if }h(w)=n.
\end{cases}
\]

The vectors in $\Res(\ssW_h)$ satisfy conditions~\ref{cond:R1} through \ref{cond:R3} in Section~\ref{sec:residue-spaces-and-conditions}, so we only need to check \ref{cond:R4}, that is, we need to show that for each connected component $\Xi$ of $G[\ssV_{h<n}]$, we have
\begin{equation}\label{grcell'}
\sum_{a \in \E(\ssV_{h=n},\, \Xi)}\res_a(\beta_{\te_a})=0,
\end{equation}
where $a$ runs through the (upward) arrows connecting a vertex of level $n$ in $\ssV_{h=n}$ to a vertex in $\Xi$. We proceed by a case analysis, according to the three cases for $\bigl(\ssV_{h<n},\ssV_{h=n}\bigr)$ listed above. 

\smallskip

$\bullet$ In the first case, $h(w)>n\text{ or } h(w)<n-\lambda$, we have $\ssV_{h <n}=\ssV_{h' <n}$ and  $\ssV_{h=n}= \ssV_{h'=n}$. In this case, $\Xi$ is a connected component of $G[\ssV_{h'<n}]$ as well, and the arrows $a$ appearing in the sum are the ones connecting $\ssV_{h'=n}$ to $\Xi$. Also, since $h(w)\neq n$, we have $\beta_{\te_a} =\beta'_{\te_a}$ for all these arrows $a$.  Since $\Res_{h'}(\beta')\in\GlobSp_{\pi'}$, we have 
\[
\sum_{a \in \E(\ssV_{h=n},\, \Xi)}\res_a(\beta_{\te_a})= \sum_{a \in \E(\ssV_{h'=n}, \,\Xi)}\res_a(\beta'_{\te_a})=0,
\]
as required.

\smallskip

$\bullet$ In the second case, $h(w)=n-\lambda$ (and hence $h'(w)=n$), we have $\ssV_{h<n} = \ssV_{h'<n}\cup\{w\}$, and $\ssV_{h=n}= \ssV_{h'=n} \setminus \{w\}$. If $w$ is not a vertex of $\Xi$, then 
$\Xi$ is as well a connected component of $G[\ssV_{h'<n}]$, which is moreover not connected to $w$. It follows that the upward arrows connecting $\ssV_{h=n}$ to $\Xi$ are the same as those connecting $\ssV_{h'=n}$ to $\Xi$. In addition, since $h(w)\neq n$, for each such arrow $a$, we have $\beta_{\te_a} =\beta'_{\te_a}$, whence \eqref{grcell'} holds as before.

Suppose now that $w$ is a vertex of $\Xi$. Let $\Xi_1,\dots,\Xi_m$ be the connected components of the graph obtained from $\Xi$ by removing the vertex $w$ (and the edges that connect $w$ to the other vertices of $\Xi$). Note that $w$ belongs to $\ssV_{h'=n}$, and each $\Xi_j$ is a connected component of $G[\ssV_{h'<n}]$.  Let $\ssa_1,\dots,\ssa_k$ be the arrows of $\E_w$ with head another vertex of $\Xi$.  Let $\ssb_1,\dots,\ssb_s$ be the horizontal arrows in $\E_{h'}^c$ with tail $w$ and  heads $\ssw_1,\dots,\ssw_s \in \ssV_{h'=n}$, respectively. Then, the reverse arrows $\barb_1 = \ssw_1w,\dots,\barb_s=\ssw_sw$ become upward for $h$, connecting vertices of $\ssV_{h=n}$ to $w$; see Figure~\ref{fig3}. This gives the following characterization of the arrow set  $\E(\ssV_{h=n},\Xi)$:
\begin{align}\label{eq:upward-xi}
\E(\ssV_{h=n},\Xi) = \E(\ssV_{h'=n},\Xi_1) \sqcup \dots \sqcup \E(\ssV_{h'=n},\Xi_m) \sqcup \{\barb_1, \dots, \barb_s\} \setminus \{\ssa_1, \dots, \ssa_k\}.
\end{align}

\begin{figure}
\centering
\begin{tikzpicture}[scale=.7pt]
\draw [-{Latex}, line width=0.3mm, blue](7,-2) -- (5.5,-1);
\draw [-{Latex}, line width=0.3mm, red](7,-2) -- (5,-2);
\draw [-{Latex}, line width=0.3mm](8.5,-2) -- (8.5,0);
\draw [-{Latex}, line width=0.3mm](7,-2) -- (8.4,0);

\draw[line width=0.1mm] (5.5,-1) circle (0.8mm);
\filldraw[blue] (5.5,-1) circle (0.7mm);

\draw[line width=0.1mm] (8.5,0) circle (0.8mm);
\filldraw[aqua] (8.5,0) circle (0.7mm);

\draw[line width=0.1mm] (5,-2) circle (0.8mm);
\filldraw[aqua] (5,-2) circle (0.7mm);
\draw[line width=0.1mm] (7,-2) circle (0.8mm);
\filldraw[blue] (7,-2) circle (0.7mm);
\draw[line width=0.1mm] (8.5,-2) circle (0.8mm);
\filldraw[aqua] (8.5,-2) circle (0.7mm);

\draw[black] (7,-2.3) node{$w$};

\draw [line width=0.3mm, blue](15,-1) -- (13.5,-1);
\draw [-{Latex}, line width=0.3mm, red] (13,-2) -- (15,-1.05);
\draw [-{Latex}, line width=0.3mm](16.5,-2) -- (16.5,0);
\draw [-{Latex}, line width=0.3mm](15,-1) -- (16.5,0);

\draw[line width=0.1mm] (13.5,-1) circle (0.8mm);
\filldraw[blue] (13.5,-1) circle (0.7mm);

\draw[line width=0.1mm] (16.5,0) circle (0.8mm);
\filldraw[aqua] (16.5,0) circle (0.7mm);

\draw[line width=0.1mm] (13,-2) circle (0.8mm);
\filldraw[aqua] (13,-2) circle (0.7mm);
\draw[line width=0.1mm] (15,-1) circle (0.8mm);
\filldraw[blue] (15,-1) circle (0.7mm);
\draw[line width=0.1mm] (16.5,-2) circle (0.8mm);
\filldraw[aqua] (16.5,-2) circle (0.7mm);

\draw[black] (15.1,-1.3) node{$w$};
\end{tikzpicture} 
\caption{The second case appearing in the proof of Proposition~\ref{prop:compatibility-residue}, with $h'(w)= h(w)+\lambda =n$. The level graph given by $h'$ is depicted on the left and the one by $h$ on the right. We have $s=k=1$, with horizontal $b_1$ in red and upward $a_1$ in blue depicted on the left. Both have tail $w$. The connected component $\Xi$ in $G[\ssV_{h<n}]$ on the right consists of the two blue vertices and the top azure vertex. $\barb_1$ becomes upward for $h$ on the right, and $a_1$ becomes horizontal.}
\label{fig3}
\end{figure}
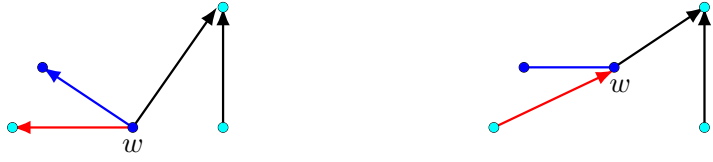

Since $\Res_{h'}(\beta')\in\GlobSp_{\pi'}$, $\beta'$ satisfies the global residue conditions for the level graph $(G,h')$ and the connected components $\Xi_1, \dots, \Xi_m$ of $G[\ssV_{h' <n}]$, that is,
\begin{align*}
    \sum_{a \in \E(\ssV_{h'=n},\, \Xi_j)}\res_a(\beta'_{\te_a})=0, \qquad \qquad \textrm{ for }j =1, \dots, m.
\end{align*}
Summing up these equations, separating the arrows $\ssa_1, \dots, \ssa_k$ with tail $w$, and using that $\beta'_{\te_a} = \beta_{\te_a}$ for every arrow $a$ whose tail is not $w$, gives
\begin{equation}\label{eq:sum-xi}
\sum_{i=1}^k \res_{a_i}(\beta'_{w})+ \sum_{j=1}^m\,\, \sum_{\substack{a \in \E(\ssV_{h'=n}, \, \Xi_j) \\ a \not \in\{a_1,\dots, a_k\}}}\res_a(\beta_{\te_a})=0.
\end{equation}
Given that $\beta'$ has vanishing residues along downward arrows for $h'$ and satisfies the global residue conditions for the remaining components of $G[\ssV_{h'<n}]$, the local residue condition at $w$, which $\beta'$ satisfies as well, yields
\[
\sum_{i=1}^k \res_{a_i}(\beta'_{w}) + \sum_{i=1}^s \res_{b_i}(\beta'_{w})=0.
\]
Moreover, $\beta'$ satisfies Rosenlicht conditions at the horizontal edges $\ssb_1, \dots, \ssb_s$, yielding
\[
\sum_{i=1}^s \res_{b_i}(\beta'_{w}) = -\sum_{i=1}^s \res_{\barb_i}(\beta'_{\te_{\barb_i}})= -\sum_{i=1}^s \res_{\barb_i}(\beta_{\te_{\barb_i}}).
\]
Therefore, we have
\[
\sum_{i=1}^k \res_{a_i}(\beta'_{w}) = \sum_{i=1}^s \res_{\barb_i}(\beta_{\te_{\barb_i}}).
\]
Combining this with Equation~\eqref{eq:sum-xi}, we get 
\[
\sum_{i=1}^s \res_{\barb_i}(\beta_{\te_{\barb_i}})+ \sum_{j=1}^m\,\, \sum_{\substack{a \in \E(\ssV_{h'=n},\, \Xi_j) \\ a \not \in\{a_1,\dots, a_k\}}}\res_a(\beta_{\te_a})=0.
\]
Given the characterization~\eqref{eq:upward-xi} of $\E(\ssV_{h=n}, \Xi)$, this is precisely \eqref{grcell'}.

\smallskip

$\bullet$ Finally, in the third case, $h(w)=n$, and so $h'(w) = n+\lambda$. Then, $\Xi$ is a connected component of $G[\ssV_{h'<n}]$ as well and $w \not\in \ssV_{h'=n}$. Thus, $\beta'$ satisfies the global residue condition for $\Xi$, and $\beta'_v = \beta_v$ for each $v\in \ssV_{h'=n}$, whence 
\[
\sum_{a \in \E(\ssV_{h'=n}, \, \Xi)} \res_a(\beta_{\te_a}) =
\sum_{a \in \E(\ssV_{h'=n},\, \Xi)} \res_a(\beta'_{\te_a})= 0.
\]
Let $\ssa_1,\dots,\ssa_k$ be the upward arrows for $h$ connecting $w$ to $\Xi$. Note that 
\[
\E(\ssV_{h=n},\Xi) = \E(\ssV_{h'=n}, \Xi) \sqcup\{\ssa_1,\dots,\ssa_k\}. 
\]
Since $\beta_w=0$, we have $\res_{\ssa_i}(\beta_w)=0$ for all $i=1, \dots, s$. Therefore,
\[
\sum_{a\in \E(\ssV_{h=n},\, \Xi)} \res_a(\beta_{\te_a}) = \sum_{i=1}^k\res_{\ssa_i}(\beta_w) + \sum_{a \in \E(\ssV_{h'=n},\, \Xi)} \res_a(\beta_{\te_a}) = 0,
\]
which is \eqref{grcell'}. The proof is complete.
\end{proof}

\section{General position conditions}\label{genpos}

Fix a field $\k$. Let $n$ be a nonnegative integer and let $\rmU$ be a vector space of finite dimension over $\k$. A \emph{flag $\filt^\bullet$ of length $n$} in $\rmU$ is a decreasing filtration of vector subspaces of $\rmU$
\[
\filt^\bullet \,:\quad \filt^n \subsetneq \filt^{n-1}\subsetneq \dots \subsetneq\filt^0 \qquad \text{with $\codim(\filt^i,\filt^{i-1})=1$ for each $i$.}
\]  
 We say that $\filt^\bullet$ \emph{complete} if $\dim \rmU =n$. (This implies $\filt^0 =\rmU$ and $\filt^n=0$.)

A collection $\Fl$ of complete flags of $\rmU$ is said to be \emph{in general position} if for each finite collection $\scrC$ of subspaces of $\rmU$, 
there is $\filt^\bullet\in\Fl$ such that:
\smallskip
 
\noindent (\emph{General position property}) For each $0 \leq i\leq \dim\rmU$, we have
\begin{equation}\label{eq:genK}
\codim\Big(Z\cap \filt^{i}, Z\Big)=\min\left(\dim Z,\,i\right)\quad \text{for each $Z\in\scrC$}.
\end{equation}

\subsection{Flags induced by points on curves}\label{sec:flag-curves}  We give an example of a collection of flags in general position that will appear in our study. Let $\varC$ be an irreducible proper curve over $\k$ and let $L$ be an invertible sheaf on $\varC$. Let $\rmU \coloneqq H^0(\varC, L)$ be the space of global sections of $L$ and denote by $n$ its dimension. For each point $p$ on $\varC$ and for each $j=0,\dots,n$, put $\filt^j_p \coloneqq H^0(\varC,L(-jp))$. We have the filtration
\[
\filt^n_p \subseteq \filt^{n-1}_p\subseteq \dots \subseteq\filt^0_p = \rmU.
\]

\begin{lemma}\label{filt1} Notation as above, assume the characteristic of $\k$ is zero. Then, there is a finite set $S\subset \varC$ containing the singular locus of $\varC$ such that 
\begin{enumerate}
\item for each $p\in \varC\setminus S$, the filtration $\filt^\bullet_p$  is a complete flag of $\rmU$, and
\item the collection $\Fl\coloneqq \bigl\{\,\filt^\bullet_p\,\st\,p\in \varC\setminus S\,\bigr\}$ is in general position.
\end{enumerate}
\end{lemma}

\begin{proof} The set $S$ is  the union of the singular locus of the curve $\varC$ and the set of Weierstrass points of the complete linear system of sections of $L$. Since $\car\k$ is $0$, for each $p \not\in S$, the filtration $\filt^\bullet_p$ is a complete flag.

As for the second statement, for each finite collection $\scrC$ of subspaces of $\rmU$, just choose $p\in \varC\setminus S$ such that $p$ is not a Weierstrass point for the linear system $(L,Z)$ for any $Z\in\scrC$, and put $\filt^\bullet\coloneqq \filt^\bul_p$. Then 
\eqref{eq:genK} holds for each $i\leq n$.
\end{proof}

\subsection{Realizability, I}\label{sec:realizabilitythm} 

Let $V$ be a finite set. For each $v\in V$, let $\rmU_v$ be a finite-dimensional vector space over $\k$. Set
\[
\rmU\coloneqq \bigoplus_{v\in V}\rmU_v \quad \textrm{and} \quad \rmU_I\coloneqq \bigoplus_{v\in I}\rmU_v \quad \textrm{for each }I\subseteq V.
\]
Let $\proj{I} \colon \rmU \to \rmU_I$ and $\proj{I\succ J}\colon \rmU_I\to \rmU_J$ denote the natural surjections for $J\subseteq I \subseteq V$.

\smallskip

Recall from Section~\ref{sec:submodular-subspace} that to each subspace $\ssW \subseteq \rmU$ we associate the submodular function $\nu!_{\ssW}^*\colon 2^V\to\Z$, given by $\varphi(I) = \dim_{\k}(\ssW_I)$ for each $I\subseteq V$, with $\ssW_I\coloneqq\proj{I}(\ssW) \subseteq \rmU_I$.

Recall from Section~\ref{sec:submodular-subset-upmin} that to each $J\subseteq V$ we associate the function $\xi!_J\colon 2^V\to\Z$ with
\[
\xi!_J(I)=\begin{cases}
0&\text{if }I\not\supseteq J,\\
-1&\text{if }I\supseteq J.
\end{cases}
\]
By Proposition~\ref{xi-submodular}, the function $\xi!_J$ is nonincreasing and  submodular. 

Given a finite sequence $\calJ$ of subsets  $\ssJ_1,\dots,\ssJ_d$ of $V$, we put 
\[
\xi!_{\calJ}\coloneqq\sum_{i=1}^d\xi!_{\ssJ_i}.
\]
Then, $\xi!_{\calJ}$ is submodular and nonincreasing. Let $\varepsilon!_{\calJ}\colon 2^V\to\Z_{\geq 0}$ be the \emph{counting function} of $\calJ$ defined by letting $\varepsilon!_{\calJ}(I)$ be the number of $\ssJ_i$ in the sequence $\calJ$ which are equal to $I$. We have the tautological identities:
\[
d = |\calJ| = \sum_{I\subseteq V}\varepsilon!_{\calJ}(I)\qquad \textrm{and} \qquad  \xi!_{\calJ}=\sum_{i=1}^d\xi!_{\ssJ_i} = \sum_{I\subseteq V}\varepsilon!_{\calJ}(I)\xi!_I.
\]

Given a subspace $\ssW\subseteq\rmU$ with $\nu!_{\ssW}^*+\xi!_{\calJ}\geq 0$, our next result asserts that the UpMin transform $\chi\coloneqq \upmin(\nu!_{\ssW}^*+\xi!_{\calJ})$, defined by 
\[
\chi(I)\coloneqq\min_{K \supseteq I} \Big(\nu!_{\ssW}^*(K) +\xi!_{\calJ}(K)\Big)\qquad\text{for each }I \subseteq V,
\]
is \emph{realizable} in the sense that it is the submodular function associated to a subspace $\ssW'$ of $\rmU$. Moreover, $\ssW'$ can be found even imposing certain conditions on it, necessary for us later on.

\begin{thm}[First realizability theorem]\label{thm:global} Let $\mathcal W$ be a finite collection of subspaces of $\rmU$. For each $I\subseteq V$, let $\ssm_I\in\Z_{\geq 0}$ and $\Fl_I$ be a collection of complete flags of $\rmU_I$ in general position. Then, there are flags $\filt^\bul_{I,1},\dots,\filt^\bul_{I,\ssm_I}\in\Fl_I$ for each  $I\subseteq V$ such that for each $\ssW\subseteq\rmU$ in $\mathcal W$, each sequence $\calJ$ of subsets $\ssJ_1,\dots,\ssJ_d$ of $V$ such that $\nu!_{\ssW}^*+\xi!_{\calJ}\geq 0$, and each partition $\unde_I=(\sse_{I,1},\dots,\sse_{I,\ssm_I})$ of $\varepsilon!_{\calJ}(I)$ into nonnegative integers for each $I\subseteq V$, we have
\begin{equation}\label{eq:W'nu}
\upmin(\nu!_{\ssW}^*+\xi!_{\calJ}) = \nu!_{\ssW'}^* \qquad \textrm{ for } \qquad   \ssW'\coloneqq \ssW \cap\bigcap_{I\subseteq V} \bigcap_{r=1}^{\ssm_I} \proj{I}^{-1}(\filt^{\sse_{I,r}}_{I,r})\subseteq\rmU.
\end{equation}
In particular, $\codim\bigl(\ssW',\ssW\big)=|\xi!_{\calJ}(V)| = |\calJ|.$
\end{thm}

\begin{remark} In the applications given in this paper, Theorem~\ref{thm:characterization-FC} for instance, and the results that lead to it, we have: $|\ssJ_i|=1$ for all $i$, $\k$ is algebraically closed of characteristic zero, $\rmU_v$ is a space of sections of an invertible sheaf on a smooth proper curve over $\k$, and the flags in $\Fl_{\{v\}}$ are the complete flags $\filt^\bul_p$ for general points $p$ on the curve for each $v\in V$, as in Section~\ref{sec:flag-curves}. Since the proof is the same for the more general statement in the theorem, we treat it in this level of generality.
\end{remark}

The crux of the proof of Theorem~\ref{thm:global} is the following lemma, which concerns the case $|\calJ|=1$, thus $\calJ=\{\ssJ\}$ and $\xi!_{\calJ}=\xi!_J$.

\begin{lemma}\label{lem:eta-gen-F} Let $\ssW\subseteq\rmU$ be a subspace and  $J\subseteq V$ be a subset such that $\ssW_J\neq 0$. Let $\chi = \upmin(\nu!_{\ssW}^*+\xi!_J)$. Let $\rmH \subset \rmU_J$ be a hyperplane. For each $I\subseteq V$, put
  \[
    \ssZ^{I\cup J}_I\coloneqq \proj{I\cup J \succ J}\Bigl(\ssW_{I \cup
      J} \cap\ker(\proj{I\cup J\succ I})\Bigr) = \proj{I\cup J \succ
      J}\Bigl(\proj{I\cup J}(\ssW)\cap\ker(\proj{I\cup J\succ
      I})\Bigr)\subseteq\rmU_J.
  \]
Then, the following statements are equivalent:
\begin{enumerate}
\item $\chi$ is the submodular function associated to $\ssW\cap\proj{J}^{-1}(\rmH) \subseteq \rmU$.
\item For each
$I\subseteq V$ such that $\ssZ^{I\cup J}_I\neq 0$, we have $\ssZ^{I\cup J}_I\not\subseteq\rmH$.
\end{enumerate}
\end{lemma}

\begin{proof} Put $\ssW'\coloneqq \ssW\cap\proj{J}^{-1}(\rmH)$, and let $\varphi\coloneqq\nu!^*_{\ssW}$ and $\varphi'\coloneqq \nu!^*_{\ssW'}$ be the associated submodular functions.  Since $\varphi$ is nondecreasing, we have
\begin{equation}\label{eq:chi-H}
\chi(I)=\min(\varphi(I),\varphi(I\cup J)-1)\quad\text{for each }I\subseteq V.
\end{equation}

If $\ssW_J\subseteq\rmH$, then $\varphi'=\varphi$, whence $\varphi'\neq\chi$. Since $\ssZ^J_\emptyset = \ssW_J\neq 0$, neither (1) nor (2) holds. 

Assume $\ssW_J\not\subseteq\rmH$. Then, for each $I\subseteq V$ we have 
\begin{equation}\label{eq:phi-phi}
\varphi'(I)\leq\varphi'(I\cup J)=\varphi(I\cup J)-1.
\end{equation}
Since also $\varphi'\leq\varphi$, we get from \eqref{eq:chi-H} that $\varphi'\leq\chi$, with equality if and only if, for each $I\subseteq V$,
\begin{equation}\label{eq:etaor}
\varphi'(I)=\varphi(I)\qquad\text{or}\qquad\varphi'(I)=\varphi(I\cup J)-1.
\end{equation}
We now need to verify for each $I\subseteq V$ that \eqref{eq:etaor} holds if and only if 
\begin{equation}\label{eq:Z-H}
\ssZ^{I\cup J}_I=0\quad\text{or}\quad\ssZ^{I\cup J}_I\not\subseteq\rmH.
\end{equation}
Let $I\subseteq V$. Consider the diagram of short exact sequences:
\[
\begin{CD}
0 @>>> \ssZ_I' @>>> \proj{I\cup J}(\ssW') @>\proj{I\cup J\succ I}>> \proj{I}(\ssW') @>>> 0\\
@. @VVV @VVV @VVV @.\\
0 @>>> \ssZ_I @>>> \proj{I\cup J}(\ssW) @>\proj{I\cup J\succ I}>> \proj{I}(\ssW) @>>> 0.
\end{CD}
\]
Clearly, $\ssZ_I'=\ssZ_I\cap\proj{I\cup J\succ J}^{-1}(\rmH)$. Using $\varphi\geq\varphi'$ and \eqref{eq:phi-phi}, and Snake Lemma,
\[
1=\varphi(I\cup J)-\varphi'(I\cup J)\geq \varphi(I)-\varphi'(I)\geq 0,
\]
with equality in the middle if and only if $\ssZ_I'=\ssZ_I$. So, $\proj{I\cup J\succ J}(\ssZ_I)\subseteq\rmH$ if and only if $\varphi'(I)=\varphi(I)-1$. But $\proj{I\cup J\succ J}$ maps $\ssZ_I$ isomorphically onto $\ssZ^{I\cup J}_I$ by definition. Thus, 
\begin{equation}\label{eq:ZIJ-H}
\ssZ^{I\cup J}_I\subseteq\rmH\quad\text{if and only if}\quad\varphi'(I)\neq\varphi(I),\quad\text{if and only if}\quad\varphi'(I)=\varphi(I)-1.
\end{equation}
Furthermore, the exactness of the bottom row in the above diagram yields
\begin{equation}\label{eq:ZIJ=0}
\ssZ^{I\cup J}_I\neq 0\quad\text{if and only if}\quad\varphi(I)\neq\varphi(I\cup J),\quad\text{if and only if}\quad\varphi(I)\leq\varphi(I\cup J)-1.
\end{equation}

If $\ssZ^{I\cup J}_I=0$, then $\varphi'(I)=\varphi(I)-1$ by \eqref{eq:ZIJ-H} and $\varphi(I)=\varphi(I\cup J)$ by \eqref{eq:ZIJ=0}. Thus, \eqref{eq:etaor} holds.

Finally, assume $\ssZ^{I\cup J}_I\neq 0$. Then $\varphi(I)\leq\varphi(I\cup J)-1$ by \eqref{eq:ZIJ=0}. Since $\varphi'\leq\varphi$, it follows that \eqref{eq:etaor} holds if and only if $\varphi'(I)=\varphi(I)$, hence if and only if $\ssZ^{I\cup J}_I\not\subseteq\rmH$ by \eqref{eq:ZIJ-H}. The proof is complete.
\end{proof}

\begin{proof}[Proof of Theorem~\ref{thm:global}]  We do induction on the sum of the $\ssm_I$. If the sum is zero, there is nothing to prove. Otherwise, let $J\subseteq V$ such that $\ssm_J>0$. 
  
For each subspace $\ssW\subseteq\rmU$ and $I\subseteq V$, let
\[
\ssZ^{I\cup J}_I(\ssW)\coloneqq\proj{I\cup J\succ J}\big(\ssW_{I\cup J}\cap\ker(\proj{I\cup J\succ I})\big)\subseteq\rmU_{J}.
\]
Notice that $\ssZ^{J}_{\emptyset}(\ssW)=\ssW_J$. Let $\scrZ_J$ be the collection of the subspaces $\ssZ^{I\cup J}_I(\ssW)$ of $\rmU_J$ for $\ssW\in\mathcal W$ and $I\subseteq V$. By hypothesis, there is a flag $\filt^\bul\in\Fl_J$ such that 
\begin{equation}\label{eq:ZF}
\codim\Bigl(\,Z\cap \filt^i, Z\Bigr)=\min\Bigl(\dim Z, i\Bigr)\quad\text{for each }i\leq\dim_{\k}\rmU_J\text{ and } Z\in\scrZ_J.
\end{equation}

For each $\ssW\subseteq\rmU$ and $i=0,\dots,\dim_{\k}\rmU_J$, define
\[
\ssW^i\coloneqq \ssW\cap\proj{J}^{-1}(\filt^i),
\]
and let $\varphi!_{\ssW}^i$ be its associated submodular function. Clearly, 
\begin{equation}\label{eq:iZ}
\ssZ^{I\cup J}_I(\ssW^i)=\ssZ^{I\cup J}_I(\ssW)\cap\filt^i
\quad\text{for each }I\subseteq V. 
\end{equation}
Then, $\ssW^i_{J}=\ssW_J\cap\filt^{i}\neq 0$ for each $i<\dim_{\k}(\ssW_J)$, allowing us to apply Lemma~\ref{lem:eta-gen-F} to  $\ssW^i$. It follows from \eqref{eq:ZF} and \eqref{eq:iZ} that either we have $\ssZ^{I\cup J}_I(\ssW^i)=0$ or $\ssZ^{I\cup J}_I(\ssW^i)\not\subseteq\filt^{i+1}$ for each $W\in\mathcal W$, and hence $\varphi!_{\ssW}^{i+1}=\upmin(\varphi!_{\ssW}^i+\xi!_J)$ by Lemma~\ref{lem:eta-gen-F}.  Applying Proposition~\ref{prop:xi-upmin} repeatedly, we get
\begin{equation}\label{eq:miI}
\varphi!_{\ssW}^i=\upmin(\nu!_{\ssW}^*+i\xi!_J)\quad\text{for each }W\in\mathcal W\text{ and } i\in[\dim_{\k}(\ssW_J)].
\end{equation}

Let $\wt{\mathcal W}$ be the collection of the $\ssW^i$ for $\ssW\in\mathcal W$ and $i\in[\dim_{\k}(\ssW_J)]$. In particular, $\mathcal W\subseteq\wt{\mathcal W}$. Put $\ssm_I'\coloneqq\ssm_I$ for $I\subseteq V$ distinct from $J$ and $\ssm'_J\coloneqq\ssm_J-1$. By induction hypothesis, there are flags $\filt^\bul_{I,1},\dots,\filt^\bul_{I,\ssm_I'}\in\Fl_I$ for each  $I\subseteq V$ such that for each $\ssW\in\wt{\mathcal W}$, for each sequence $\calJ$ of subsets $\ssJ_1,\dots,\ssJ_d$ of $V$ such that $\nu!_{\ssW}^*+\xi!_{\calJ}\geq 0$, and for each partition $\unde_I=(\sse_{I,1},\dots,\sse_{I,\ssm_I'})$ of $\varepsilon!_{\calJ}(I)$ into nonnegative integers, for each $I\subseteq V$, we have 
\begin{equation}\label{eq:wtW'nu}
 \upmin(\nu!_{\ssW}^*+\xi!_{\calJ}) = \nu!_{\ssW''}^* \qquad \textrm{ for } \qquad   \ssW''\coloneqq \ssW \cap\bigcap_{I\subseteq V} \bigcap_{r=1}^{\ssm'_I} \proj{I}^{-1}(\filt^{\sse_{I,r}}_{I,r})\subseteq\rmU.
\end{equation}

Put $\filt^\bul_{J,\ssm_J}\coloneqq\filt^\bul$. Let $W\in\mathcal W$, and let $\calJ$ be a sequence of subsets $\ssJ_1,\dots,\ssJ_d$ of $V$ such that $\nu!_{\ssW}^*+\xi!_{\calJ}\geq 0$. For each $I\subseteq V$, let $\unde_I=(\sse_{I,1},\dots,\sse_{I,\ssm_I})$ be a partition of $\varepsilon!_{\calJ}(I)$ into nonnegative integers. We need to prove \eqref{eq:W'nu}.

Let $k\coloneqq \sse_{J,\ssm_J}$. Notice that 
$k\leq\varepsilon!_{\calJ}(J)\leq\dim\ssW_J$, because $\nu!_{\ssW}^*+\xi!_{\calJ}\geq 0$. Let $\wt{\calJ}$ be the subsequence obtained from $\calJ$ by removing $k$ of the $\ssJ_i$ equal to $J$. Then, $\xi!_{\calJ}=\xi!_{\wt{\calJ}}+k\xi!_J$.

If $k=0$, then $\wt{\calJ}=\calJ$ and also $(\sse_{J,1},\dots,\sse_{J,\ssm'_J})$ is a partition of $\varepsilon!_{\calJ}(J)$. Then, \eqref{eq:wtW'nu} holds. But $\ssW''=\ssW'$ because $\sse_{J,\ssm_J}=0$, and thus \eqref{eq:W'nu} holds. 

If $k>0$, we have that $\ssW^k$ is in $\wt{\mathcal W}$, and its submodular function $\varphi!_{\ssW}^k$ satisfies, by \eqref{eq:miI} and Proposition~\ref{prop:xi-upmin},
\begin{equation}\label{eq:wtJ-J}
\begin{aligned}
\varphi!_{\ssW}^k+\xi!_{\wt{\calJ}}=&\upmin(\nu!_{\ssW}^*+k\xi!_J)+\xi!_{\wt{\calJ}}=\upmin(\nu!_{\ssW}^*+k\xi!_J+\xi!_{\wt{\calJ}})\\=&\upmin(\nu!_{\ssW}^*+\xi!_{\calJ})\geq 0.
\end{aligned}
\end{equation}
Since $\unde_I=(\sse_{I,1},\dots,\sse_{I,\ssm_I'})$ is a partition of $\varepsilon!_{\wt{\calJ}}(I)$ into nonnegative integers for each $I\subseteq V$, we get that 
\[
 \upmin(\varphi!_{\ssW}^k+\xi!_{\wt\calJ}) = \nu!_{\ssW'''}^* \qquad \textrm{ for } \qquad   \ssW'''\coloneqq \ssW^k \cap\bigcap_{I\subseteq V} \bigcap_{r=1}^{\ssm'_I} \proj{I}^{-1}(\filt^{\sse_{I,r}}_{I,r})\subseteq\rmU.
\]
But $\ssW'''=\ssW'$ because $k=\sse_{J,\ssm_J}$, and $\upmin(\varphi!_{\ssW}^k+\xi!_{\wt\calJ})=\upmin(\nu!_{\ssW}^*+\xi!_{\calJ})$ by \eqref{eq:wtJ-J}. Thus, \eqref{eq:W'nu} holds in this case as well.
\end{proof}

\subsection{Failure of nonnegativity} 
We state what happens if the nonnegativity assumption on $\nu!_{\ssW}^*+\xi!_{\calJ}$ in Theorem~\ref{thm:global} is relaxed. This will be used later in the proof of Proposition~\ref{prop:step1singnot}.

\begin{prop}\label{prop:zeroprojection} Let $\ssW\subseteq\rmU$ be a subspace and $\calJ$ a finite sequence of subsets of $V$ such that 
\[
\nu!_{\ssW}^*(J)+\xi!_{\calJ}(J)\leq 0
\]
for some nonempty subset $J\subseteq V$. 
For each $I\subseteq V$, let $\ssm_I\in\Z_{\geq 0}$ and $\Fl_I$ be a collection of complete flags of $\rmU_I$ in general position. Then there are flags $\filt^\bul_{I,1},\dots,\filt^\bul_{I,m_I}\in\Fl_I$ for all
$I\subseteq V$ such that, for each partition $\unde_I=(\sse_{I,1},\dots,\sse_{I,\ssm_I})$ of $\varepsilon_{\calJ}(I)$ into nonnegative integers for each $I\subseteq V$, putting
\[
\ssW'\coloneqq \ssW \cap\bigcap_{I\subseteq V}\bigcap_{r=1}^{\ssm_I} \proj{I}^{-1}(\filt^{\sse_{I,r}}_{I,r}),
\]
we have $\proj{v}(\ssW')=0$ for some $v\in V$.    
\end{prop}

\begin{proof} Put $\varphi\coloneqq \nu!_{\ssW}^*$. By hypothesis, there is a subsequence $\wt{\calJ}$ of $\calJ$ such that $\varphi+\xi!_{\wt{\calJ}}\geq 0$ and $\varphi(J)+\xi!_{\wt{\calJ}}(J)=0$. By Theorem~\ref{thm:global}, there are $\filt^\bul_{I,1},\dots,\filt^\bul_{I,\ssm_I}\in\Fl_I$ for each $I\subseteq V$ such that, for each partition $\unde_I'=(\sse_{I,1}',\dots,\sse_{I,\ssm_I}')$ of $\varepsilon!_{\wt{\calJ}}(I)$ into nonnegative integers for each $I\subseteq V$, putting
\[
\ssW''\coloneqq \ssW \cap \bigcap_{I \subseteq V}\bigcap_{r=1}^{\ssm_I} \proj{I}^{-1}(\filt^{\sse'_{I,r}}_{I,r}),
\]
we have 
\[
\nu!^*_{\ssW''}=\upmin(\varphi+\xi!_{\wt{\calJ}}).
\]
Then, $\proj{J}(\ssW'')=(0)$. 

Now, let $\unde_I=(\sse_{I,1},\dots,\sse_{I,\ssm_I})$ be a partition of $\varepsilon!_{\calJ}(I)$ into nonnegative integers for each $I\subseteq V$. Since 
$\wt{\calJ}$ is a subsequence of $\calJ$, for each $I\subseteq V$ there is a partition $\unde_I'=(\sse_{I,1}',\dots,\sse_{I,\ssm_I}')$ of $\varepsilon!_{\wt{\calJ}}(I)$ into nonnegative integers with $\sse_{I,r}'\leq\sse_{I,r}$ for each $r\in[\ssm_I]$. Then, $\ssW'\subseteq\ssW''$, and hence  $\proj{J}(\ssW')=(0)$. Taking $v\in J$, we conclude. 
\end{proof}

\subsection{Realizability, II: general hyperplanes} We will combine Theorem~\ref{thm:global} with Theorem~\ref{thm:eta-gen-2-hyper} below in Section~\ref{sec:submodular-Wh} to compute submodular functions of spaces of sections of sheaves obtained by gluings.  Recall $V$ and $\rmU=\bigoplus\rmU_v$ from Section~\ref{sec:realizabilitythm}.

\begin{thm}[Second realizability theorem]\label{thm:eta-gen-2-hyper} Let $\ssW \subseteq \rmU$ be a subspace and $\varphi$ be its submodular function. Let $\calJ$ be a sequence of subsets $\ssJ_1,\dots,\ssJ_d$ of $V$ such that $\varphi+\xi!_{\calJ}\geq 0$. For each $i=1,\dots,d$ and $v\in\ssJ_i$, let $\rmH_{i,v}\subseteq\rmU_v$ be a subspace of codimension at most $1$.  Assume that for each collection $\mathscr S$ of subsets $\ssS_i\subseteq\ssJ_i$ for $i=1,\dots,d$ verifying the inequality
\[
\varphi+\sum_{i=1}^d\sum_{v\in\ssS_i}\xi!_v\geq 0,
\]
the space
\[
\ssW(\mathscr S)\coloneqq\ssW \cap\bigcap_{i=1}^d\bigcap_{v\in \ssS_i}\proj{v}^{-1}(\rmH_{i,v})
\]
has submodular function 
\[
\varphi!^{\mathscr S} \coloneqq \upmin\Bigl(\varphi + \sum_{i=1}^d\sum_{v\in\ssS_i}\xi!_v \Bigr).
\]
Then, $\bigoplus_{v\in\ssJ_i}\rmH_{i,v}$ is a proper subspace of $\rmU_{\ssJ_i}$ for each $i=1,\dots,d$. Moreover, letting $\mathscr H$ be a collection of general hyperplanes $\rmH_i\subset\rmU_{\ssJ_i}$ containing $\bigoplus_{v\in\ssJ_i}\rmH_{i,v}$ for $i=1,\dots,d$, the space
\[
  \ssW(\mathscr H)\coloneqq\ssW \cap\bigcap_{i=1}^d\proj{\ssJ_i}^{-1}(\rmH_i)
\]
has 
submodular function $\upmin(\varphi+\xi!_{\calJ})$.
\end{thm}

We need the following lemma.

\begin{lemma}\label{lem:eta-gen-hyper} Let $\ssW \subseteq \rmU$ be a subspace and $\varphi$ its submodular function. Let $J\subseteq V$ such that $\ssW_J\neq 0$. For each $v\in J$, let $\rmH_v\subseteq\rmU_v$ be a subspace of codimension at most $1$, put $\ssW(v)\coloneqq\ssW\cap\proj{v}^{-1}(\rmH_v)$, and let $\varphi^v$ be its submodular function. Assume that 
\begin{equation}\label{eq:eta-v}
\text{either\quad  
$\proj{v}(\ssW)=0$, \quad or }\quad \varphi^v= \upmin(\varphi+\xi!_v).
\end{equation}
Then, $\bigoplus_{v\in J}\rmH_v$ is a proper subspace of $\rmU_J$. Moreover, letting $\rmH \subset \rmU_J$ be a general hyperplane containing $\bigoplus_{v\in J}\rmH_v$, the space $\ssW'\coloneqq\ssW\cap\proj{J}^{-1}(\rmH)$ has submodular function $\varphi'$ satisfying
\[
\varphi'= \upmin(\varphi+\xi!_J).
\]
\end{lemma}

\begin{proof} First of all, for each $v\in J$ such that $\proj{v}(\ssW)\neq 0$, we have $\ssW(v)\neq \ssW$ by \eqref{eq:eta-v}, whence $\rmH_v\subset\rmU_v$ is a hyperplane. Since $\ssW_J\neq 0$, we must have $\proj{v}(\ssW)\neq 0$ for some $v\in J$, and hence $\bigoplus_{v\in J}\rmH_v$ is a proper subspace of $\rmU_J$.

For each $I\subseteq V$, let
\[
\ssZ^{I\cup J}_I\coloneqq \proj{I\cup J \succ J}\Bigl(\ssW_{I \cup J} \cap\ker(\proj{I\cup J\succ I})\Bigr) \subseteq\rmU_J.
\]
By Lemma~\ref{lem:eta-gen-F}, we only need to prove that $\ssZ^{I\cup J}_I\not\subseteq\rmH$ for any $I\subseteq V$ such that $\ssZ^{I\cup J}_I\neq 0$.

If $\ssZ^{I\cup J}_I\neq 0$, then $\proj{J\succ v}(\ssZ^{I\cup J}_I)\neq 0$ for some $v\in J$. In particular, $\proj{v}(\ssW)\neq 0$. Since
\[
\ssZ^{I\cup v}_v= \proj{J\succ v}(\ssZ^{I\cup J}_I),
\]
and since \eqref{eq:eta-v} holds, it follows from Lemma~\ref{lem:eta-gen-F} that $\ssZ^{I\cup v}_v\not\subseteq \rmH_v$. Then, $\proj{J\succ v}(\ssZ^{I\cup J}_I)\not\subseteq \rmH_v$. For general $\rmH$, we thus have that $\ssZ^{I\cup J}_I\not \subseteq\rmH$ for any $I\subseteq V$ such that $\ssZ^{I\cup J}_I\neq 0$. 
\end{proof}

\begin{proof}[Proof of Theorem~\ref{thm:eta-gen-2-hyper}] 
Assume $\varphi+\xi!_{\calJ}\geq 0$. Then, by Proposition~\ref{prop:xi-singleton}, for each $i=1,\dots,d$ there is $\ssv_i\in\ssJ_i$ such that $\varphi+\xi!_{\ssv_i}\neq 0$. By hypothesis, $\ssW\cap\proj{\ssv_i}^{-1}(\ssH_{i,\ssv_i})$ has submodular function $\upmin(\varphi+\xi!_{\ssv_i})$, whence $\ssH_{i,\ssv_i}\subset\rmU_{\ssv_i}$ is a hyperplane. It follows that $\bigoplus_{v\in\ssJ_i}\rmH_{i,v}$ is a proper subspace of $\rmU_{\ssJ_i}$ for each $i=1,\dots,d$.

Let $\mathscr S$ (resp.~$\mathscr T$) be a collection of subsets $\ssS_i\subseteq\ssJ_i$ (resp.~$\ssT_i\subseteq\ssJ_i$) for $i\in[d]$ such that $\ssS_i\cap\ssT_i=\emptyset$ for every $i$ and
\begin{equation}\label{eq:etasumS}
\varphi+\sum_{i=1}^d\sum_{v\in \ssS_i}\xi!_v+\sum_{i=1}^d\xi!_{\ssT_i}\geq 0.
\end{equation}
Let $\mathscr H$ be a collection of general
hyperplanes $\rmH_i\subset\rmU_{\ssT_i}$ containing $\bigoplus_{v\in\ssT_i}\rmH_{i,v}$ for each $i\in [d]$ such that $\ssT_i\neq\emptyset$. (As above, \eqref{eq:etasumS} guarantees the existence of $\rmH_i$.) We will prove by induction on $\sum|\ssT_i|$ that the space
\[
  \ssW(\mathscr S,\mathscr T;\mathscr H)\coloneqq\ssW(\mathscr S) \cap\bigcap_{i=1}^d\proj{\ssT_i}^{-1}(\rmH_i)
\]
has submodular function $\varphi^{\mathscr S,\mathscr T}$ given by
\begin{equation}\label{eq:etasub}
\varphi!^{\mathscr S,\mathscr T}\coloneqq \upmin\Big(\varphi+\sum_{i=1}^d\sum_{v\in\ssS_i} \xi!_v+\sum_{i=1}^d\xi!_{\ssT_i}\Big).
\end{equation}
The case where $\sum |\ssT_i|=0$, that is, all the $\ssT_i$ are empty, is our hypothesis. The case where $\sum_{i=1}^d |\ssT_i|$ is maximum, that is, $\ssT_i=\ssJ_i$ for every $i$, is our thesis. 

Assume $\sum|\ssT_i|>0$. Without loss of generality, assume $\ssT_d\neq\emptyset$. Let $\mathscr T'$ be the collection of subsets $\ssT_i'\subseteq\ssJ_i$ for $i=1,\dots,d$ with $\ssT_i'=\ssT_i$ for $i<d$ and $\ssT_d'=\emptyset$. Assume \eqref{eq:etasumS} holds. Then. we have
\begin{equation}\label{eq:etasum-1}
  \varphi+\sum_{i=1}^d\sum_{v\in \ssS_i}\xi!_v+\sum_{i=1}^{d-1}\xi!_{\ssT_i}\geq 0.
\end{equation}
By the induction hypothesis, $\ssW(\mathscr S,\mathscr T';\mathscr H')$ has submodular function $\varphi!^{\mathscr S,\mathscr T'}$, where
\begin{equation}\label{eq:etasub-1}
  \varphi!^{\mathscr S,\mathscr T'}\coloneqq \upmin\Big(\varphi+\sum_{i=1}^d\sum_{v\in \ssS_i}\xi!_v+ \sum_{i=1}^{d-1}\xi!_{\ssT_i}\Big),
\end{equation}
and $\mathscr H'$ is the subcollection of $\mathscr H$ where we exclude $\rmH_d$.

Notice that, applying Proposition~\ref{prop:xi-upmin} to \eqref{eq:etasub} and \eqref{eq:etasub-1}, we get 
\begin{equation}\label{eq:STST'}
\varphi!^{\mathscr S,\mathscr T}=\upmin\Big(\varphi!^{\mathscr S, \mathscr T'}+\xi!_{\ssT_d}\Big).
\end{equation}
Notice as well that
\begin{equation}\label{eq:WST=WST'Hd}
\ssW(\mathscr S,\mathscr T;\mathscr H)=\ssW(\mathscr S,\mathscr T'; \mathscr H')\cap\proj{\ssT_d}^{-1}(\rmH_d).
\end{equation}
We would like to apply Lemma~\ref{lem:eta-gen-hyper} to $\ssW(\mathscr S,\mathscr T'; \mathscr H')$. Now,
\begin{equation}\label{eq:WSwT=WST'Hdw}
\ssW(\mathscr S,\mathscr T'; \mathscr H')\cap\proj{w}^{-1}(\rmH_{d,w})=\ssW(\mathscr S^w,\mathscr T';\mathscr H')
\end{equation}
for each $w\in\ssT_d$, where $\mathscr S^w$ is the collection of subsets
$\ssS^w_i\subseteq\ssJ_i$ for $i=1,\dots,d$ with $\ssS^w_i=\ssS_i$ for $i=1, \dots, d-1$, and
$\ssS^w_d=\ssS_d\cup\{w\}$. By induction, for each $w\in\ssT_d$ such that 
\begin{equation}\label{eq:etasum-w}
  \varphi+\sum_{i=1}^d\sum_{v\in \ssS_i}\xi!_v+\sum_{i=1}^{d-1}\xi!_{\ssT_i}+\xi!_w\geq 0,
\end{equation}
the space $\ssW(\mathscr S^w,\mathscr T';\mathscr H')$ has submodular function $\varphi!^{\mathscr S^w,\mathscr T'}$ given by
\begin{equation}\label{eq:etasub-w}
  \varphi!^{\mathscr S^w,\mathscr T'}\coloneqq \upmin\Big(\varphi+\sum_{i=1}^d\sum_{v\in \ssS_i}\xi!_v+ \sum_{i=1}^{d-1}\xi!_{\ssT_i}+\xi!_w\Big).
\end{equation}
Notice that \eqref{eq:etasub-1} and \eqref{eq:etasub-w} yield, by 
Proposition~\ref{prop:xi-upmin},
\begin{equation}\label{eq:usethis}
\varphi!^{\mathscr S^w,\mathscr T'}=\upmin\Big(\varphi!^{\mathscr S, \mathscr T'}+\xi!_w\Big)
\end{equation}
for each $w\in\ssT_d$.

Now, for each $w\in\ssT_d$, it  follows from
\eqref{eq:etasum-1}~and~\eqref{eq:etasub-1} that \eqref{eq:etasum-w}
does not hold if and only if there is $I\subseteq
V$ such that $I$ contains $w$ and $\varphi!^{\mathscr S,\mathscr T'}(I)=0$, hence, if and
only if $\proj{w}(\ssW(\mathscr S,\mathscr T';\mathscr H'))=0$ (because $\varphi!^{\mathscr S,\mathscr T'}$ is the submodular function of $\ssW(\mathscr S,\mathscr T';\mathscr H')$). In other words, \eqref{eq:etasum-w} is equivalent to
\[
\proj{w}(\ssW(\mathscr S,\mathscr T';\mathscr H'))\neq 0.
\]
Using \eqref{eq:WSwT=WST'Hdw} and \eqref{eq:usethis}, we get that  $\ssW(\mathscr S,\mathscr T'; \mathscr H')\cap\proj{w}^{-1}(\rmH_{d,w})$ has submodular function 
\[
\upmin\Big(\varphi!^{\mathscr S, \mathscr T'}+\xi!_w\Big)
\]
for each $w\in\ssT_d$ such that $\proj{w}(\ssW(\mathscr S,\mathscr T';\mathscr H'))\neq 0$.

On the other hand, it follows from \eqref{eq:etasumS} and Proposition~\ref{prop:xi-singleton} that there is $w\in\ssT_d$ such that \eqref{eq:etasum-w} holds, whence $\proj{w}(\ssW(\mathscr S,\mathscr T';\mathscr H'))\neq0$, and in particular, 
\[
\proj{\ssT_d}(\ssW(\mathscr S,\mathscr T';\mathscr H'))\neq 0.
\]

Applying Lemma~\ref{lem:eta-gen-hyper} to $\ssW(\mathscr S,\mathscr T'; \mathscr H')$, and using \eqref{eq:STST'} and \eqref{eq:WST=WST'Hd}, we get that the space $\ssW(\mathscr S,\mathscr T;\mathscr H)$ has submodular function
$\varphi!^{\mathscr S,\mathscr T}$. The proof is complete.
\end{proof}


\section{Expected limits}\label{sec:expected}

We fix a nodal projective connected reduced curve $X$ defined over an algebraically closed field $\k$. Let $g$ be its (arithmetic) genus. Let $G=(V,E)$ be its dual graph and $\E$ its arrow set. Let $\varC_v$ be the normalization of the component of $X$ associated to each $v\in V$.  Each edge $e=\{u,v\}$ in $E$ with $a=uv$ and $\bar a=vu$ the arrows on $e$, gives a point $\ssp^a\in \varC_u$ and a point $\ssp^{\bar a}\in\varC_v$ whose identification is the node $\ssp^e$ on $X$ corresponding to $e$. Let $\g\colon V \to \Z_{\geq 0}$ be the genus function assigning the genus $\g_v$ of $\varC_v$ to each $v\in V$. For each $v\in V$, denote by $\omega!_v=\omega!_{\varC_v}$ the sheaf of differentials on $\varC_v$, by $\Omega!_v$ the space of meromorphic differentials, and and put $\Omega = \bigoplus\Omega!_v$.

We fix an ordered partition $\pi$ of $V$, yielding thus the level graph $(G,\pi)$. We fix as well a slope function $s\colon\E\to\Z$; see Definition~\ref{defi:slopefcn}.

Recall the notation associated to the graph $G$, the level graph $(G,\pi)$ and the slope function $s$ introduced in Section~\ref{sec:graphs}, particularly in Section~\ref{sec:levelgraphs} and Definition~\ref{defi:slopefcn}. 

If $(s,\pi)=(\slztwist{\ell}h,\sspi_h)$ for certain edge length function $\ell\colon E\to\R_{>0}$ and function $h\colon V\to\R$, then an arrow is upward for $\pi$ if and only if it is upward for $s$, that is, $\ssA_s=\ssA_{\pi}$; see Remark~\ref{rem:slope-level-pair}.

\begin{defi}\label{defi:slope-level-pair} We call $(s,\pi)$ a \emph{slope-level pair} on $G$ if $\ssA_s=\ssA_{\pi}$. 
\end{defi}

We assume from now on that $(s,\pi)$ is a slope-level pair. 

\subsection{Sheaves on the curves $\varC_v$}\label{sec:slope-sheaves}

For each vertex $v\in V$, define the sheaf $\ssL_{s,v}$ on $\varC_v$ by
\begin{equation}\label{eq:shf-Lh}
\ssL_{s,v}\coloneqq \omega!_{v}\Big(\sum_{\substack{a\in\E_v}} \big(1+s(a)\big) \ssp^a\Big).
\end{equation}
Notice that $1+s(a) \leq 0$ if and only if $a\in \barA_{s,v}$. Define the divisors $\ssP_{s,v}$  and  $\ssN_{s,v}$ on $\varC_v$ by
\begin{equation}\label{eq:PhNh}
\ssP_{s,v}\coloneqq\sum_{a \in \E_v \setminus \barA_{s,v}} \big(1+s(a)\big)\ssp^{a} \qquad \text{and} \qquad \ssN_{s,v}\coloneqq\sum_{a \in \barA_{s,v}}\big(-1 - s(a)\big)\ssp^{a}.
\end{equation}
The divisors $\ssP_{s,v}$ and $\ssN_{s,v}$ are both effective, and $\ssL_{s,v}=\omega!_v(\ssP_{s,v}-\ssN_{s,v})$. We view $\ssL_{s,v}$ as a subsheaf of the constant sheaf $\Omega!_v$ on $\varC_v$ consisting of the meromorphic differentials with divisors of poles and zeros bounded by $\ssP_{s,v}$ and $\ssN_{s,v}$, respectively.

\subsection{The space $\hatWexp_{s,\pi}$ given by residue conditions}\label{sec:hatWexp} 
 Consider the residue map $\Res \colon \Omega \to \bV=\bigoplus\bV_v$ from Section~\ref{sec:residue-map-Res}, and the $\pi$-residue subspace $\GlobSp_{\pi} \subset \bV$ from Section~\ref{sec:residue-spaces-and-conditions}. Define
\[
\hatWexp_{s,\pi} 
\coloneqq \Res^{-1}(\GlobSp_{\pi}) \bigcap \left(\bigoplus_{v\in V} H^0(\varC_v,\ssL_{s,v})\right) \subseteq \Omega.
\]

As we will see in Section~\ref{sec:limit-canonical-series-FC}, describing $\hatWexp_{s,\pi}$ is an important step in describing limit canonical series on $X$. Theorem~\ref{step1sing} gives the submodular function associated to $\hatWexp_{s,\pi}\subset \Omega$ under certain conditions.

\begin{prop}\label{prop:What-invariance}
The subspace $\hatWexp_{s,\pi}\subset \Omega$ is invariant by the action of $\Gm^{\ssV_{\pi}}$.
\end{prop} 
Recall that $\Gm^{\ssV_{\pi}}\subseteq \Gm^{\ssV}$ is the subgroup of $\psi$ which verify $\psi_u=\psi_v$ for $u, v$ of the same level.
\begin{proof}
The map $\Res\colon \Omega \to \bV$ is $\Gm^{\ssV}$-equivariant. In addition, the subspace $\bigoplus H^0(\varC_v,\ssL_{s,v}) \subset \Omega$ is $\Gm^{\ssV}$-invariant, whereas the subspace $\GlobSp_{\pi} \subset \bV$ is only  $\Gm^{\ssV_{\pi}}$-invariant in general. The result follows.
\end{proof}

\subsection{Pole order and residue conditions}\label{sec:ResP} 

We will need the following basic result on residues. 

\begin{prop}\label{prop:auxiliary} Let $\varC$ be a smooth connected proper curve over an algebraically closed field $\k$ and $\omega$ its canonical sheaf. Let $\ssp_1,\dots,\ssp_m\in\varC$ be distinct points and $\ssn_1,\dots,\ssn_m$ positive integers. Then, the image of the map 
\[
\bigoplus_{i=1}^m \res_{\ssp_i} \colon   H^0\Big(\varC, \omega\Big(\sum_{i=1}^m\ssn_i \ssp_i\Big)\Big) \longrightarrow \k^m, \qquad \qquad \alpha \mapsto (\res_{\ssp_1}(\alpha), \dots, \res_{\ssp_m}(\alpha)).
\]
is the set of all vectors in $\k^m$ whose sum of coordinates is zero. 
\end{prop}

\begin{proof} By the residue theorem, $\sum\res_{p_i}(\alpha)=0$ for all meromorphic differentials $\alpha$ on $C$ with poles at most at the $\ssp_i$. By Riemann--Roch, for each pair of distinct elements $i,j\in [m]$, there is a meromorphic 1-form on $\varC$ for which $\ssp_i$ and $\ssp_j$ are the unique poles, both of order one. This implies that the image of $\bigoplus\res_{\ssp_i}$ contains all the vectors of the form $e_i -e_j$ for $i,j\in [m]$,  and the proposition follows.
\end{proof}

For each $I\subseteq V$, by a slight abuse of notation, we denote  by $\proj{I}$ both projection maps:  
\[
\proj{I}\colon\Omega=\bigoplus_{v\in V} \Omega!_v \longrightarrow\bigoplus_{v\in I} \Omega!_v=\Omega!_I \qquad \textrm{and} \qquad \proj{I}\colon \bV=\bigoplus_{v\in V} \k^{\E_v} \to \bigoplus_{v\in I}\k^{\E_v}=\bV_I.
\] 

Let $\bV^0_\pi \subseteq \bV$ be the vector subspace consisting of those $\psi$ that verify Conditions~\ref{cond:R1} and \ref{cond:R2} in Section~\ref{sec:residue-spaces-and-conditions}. For each $I\subseteq V$, let $\bV^0_{\pi,I} \coloneqq \proj{I}(\bV^0_\pi) \subseteq \bV_{I}$. We have $\bV^0_\pi=\bigoplus\bV^0_{\pi,v}$. It is easy to see that for each $v\in V$, we have 
\[
\dim \bV_{\pi,v}^0 = \max\{0, |\E_v\setminus \barA_\pi|-1\}.
\]

Our next lemma is stated more generally than what we initially need it for, which is the case where $\ssP_v=\ssP_{s,v}$ for all $v$ and $\scrG=\GlobSp_{\pi}$. But this level of generality will be needed later.

\begin{lemma}\label{lem:inverse-residue} For each $v\in V$, let $\ssP_v$ be an effective divisor on $\varC_v$ having the same support as $\ssP_{s,v}$. Let $\scrG \subseteq\bV^0_{\pi}$ be a vector subspace and put
\[
\cW\coloneqq\Res^{-1}(\scrG)\bigcap \Bigl(\bigoplus_{v\in V} H^0(\varC_v,\omega!_v(\ssP_v))\Bigr) \subset\Omega.
\]
Then, for each $I\subseteq V$,
\begin{align*}
\dim \proj{I}(\cW)&=\dim \proj{I}(\scrG) +\sum_{v\in I}\big(\dim H^0(\varC_v,\omega!_v(\ssP_v))-\dim\bV^0_{\pi,v}\big)\\
&= \dim \proj{I}(\scrG)+\g(I)+ \sum_{v\in I}\deg(\ssP_v)-|\E_I\setminus \barA_{\pi}|.
\end{align*}
\end{lemma}

\begin{proof} To simplify, put $\rmH_v\coloneqq H^0(\varC_v,\omega!_v(\ssP_v))$ for each $v\in V$. 
Let $J\coloneqq \ssI^c$, the complement of $I$. Let $\scrG_{I} = \proj{I}(\scrG)$ and denote by $\scrG^J$ the kernel of the surjection 
$\scrG\to \scrG_I$, viewed as a subspace of $\bV^0_{\pi,J}$. We obtain an injective homomorphism of short exact sequences, and consider its cokernel, denoted $\cC$:
\[
\begin{CD}
    0 @>>> \cG^J @>>> \cG @>>> \cG_I @>>> 0\\
    @.@VVV@VVV @VVV @.\\
    0@>>>\bV^0_{\pi,J} @>>> \bV^0_\pi  @>>> \bV^0_{\pi,I} @>>>0\\
    @.@VVV@VVV @VVV @.\\
    0@>>>\cC^J @>>> \cC @>>> \cC_I@>>>0
\end{CD}
\]

Let $\cW_I= \proj{I}(\cW)$ and $\cW^J = \ker(\cW \to \cW_I)$. By Proposition~\ref{prop:auxiliary}, the image under $\Res$ of the sum $\bigoplus\rmH_v$ is $\bV^0_{\pi}$. Thus,
composing $\Res$ with the cokernel, we get a surjective homomorphism of short exact sequences whose kernel is the short exact sequence $0\to\cW^J \to \cW \to \cW_I\to0$:
\[
\begin{CD}
   0 @>>> \cW^J @>>> \cW @>>> \cW_I@>>>0\\
   @. @VVV@VVV @VVV @.\\
   0 @>>> \bigoplus_{v\in J}\rmH_v @>>> \bigoplus_{v\in V}\rmH_v @>>> \bigoplus_{v\in I}\rmH_v @>>> 0\\
   @. @VVV @VVV @VVV @.\\
   0 @>>> \cC^J @>>> \cC @>>> \cC_I @>>> 0.
\end{CD}
\]
The first equation in the lemma follows. To get the second one, we combine Riemann--Roch with the fact that $\dim \bV_{\pi,v}^0 = \max\{0, |\E_v\setminus \barA_\pi|-1\}$ for each $v\in V$.
\end{proof}

\subsection{Submodular function $\sshateta_{s,\pi}$}\label{sec:zeta} 
Let $\gamma!_{\pi}$ be the submodular function associated to the $\pi$-residue space $\GlobSp_\pi$. Let $\zeta!_s$ be the submodular function associated to $s$ from Definition~\ref{defi:zetas}:
\[
\zeta!_s(I) = |\ssA_s(\ssI^c, \ssI)|+\sum_{a \in \E(I, \ssI^c)}s(a) \qquad \forall\, I \subseteq V.
\]

\begin{defi}
 Define the function $\sshateta_{s,\pi}\colon \ssub{2}^V \to \Z$ by 
\begin{align*}
\sshateta_{s,\pi}(I)&\coloneqq\gamma!_{\pi}(I)+ \g(I)+\zeta!_s(I)+|\ssA^{\rmint}_s(I)|\\
&=\gamma!_{\pi}(I)+\g(I)+|\ssA_s(\ssI^c, \ssI)|+ |\ssA_s^{\rmint}(I)|+\sum_{a\in\E(I,\ssI^c)}s(a) \qquad \forall\, I\subseteq V.\qedhere
\end{align*}
\end{defi}

\begin{prop}\label{prop:hateta-submodular} 
The function $\sshateta_{s,\pi}$ is submodular of range $g + |\ssA^\rmint_s|$.
\end{prop}

\begin{proof} The function $\gamma!_{\pi}$ is submodular. Also, $\g$ is modular. Finally, $\sshateta_{s,\pi}-\gamma!_{\pi}-\g$ is modular by Proposition~\ref{prop:zetah}. Then, $\sshateta_{s,\pi}$, being a sum of submodular functions, is submodular. 

The range of $\sshateta_{s,\pi}$ is obtained by evaluating at $I=V$, which gives, using Theorem~\ref{thm:dimG} and the genus formula $g=|E|-|V|+1 + \g(V)$,
\[
\sshateta_{s,\pi}(V) = \gamma!_{\pi}(V)+ \g(V)+ |\ssA_s^\rmint| =|E|-|V|+1 + \g(V) + |\ssA_s^\rmint|=g+|\ssA_s^\rmint|. \qedhere
\]
\end{proof}

We have the following alternate description of $\sshateta_{s,\pi}$, which we will need later.

\begin{prop}\label{prop:connection-hateta-gamma} Put 
\[
\hatWexp_{s,\pi}^+\coloneqq \Res^{-1}(\GlobSp_{\pi})\bigcap\Bigl(\bigoplus_{v\in V}H^0(\varC_v,\omega!_v(\ssP_{s,v}))\Bigr)\subset\Omega.
\]
Then, for each $I \subseteq V$, we have 
\[
\sshateta_{s,\pi}(I)=\dim\proj{I}(\hatWexp_{s,\pi}^+)-\sum_{v\in I}\deg(\ssN_{s,v}).
\]
\end{prop}

\begin{proof} We have
\begin{align*}
\dim\proj{I}(\hatWexp_{s,\pi}^+)&-\sum_{v\in I}\deg(\ssN_{s,v})
=\gamma!_{\pi}(I)+\g(I)+ \sum_{a\in\E_I\setminus \barA_s}s(a)+ \sum_{a\in\E_I\cap\barA_s}(1+s(a))\\
=&\gamma!_{\pi}(I)+\g(I)+\sum_{a\in\E_I}s(a)+
|\E_I\cap\barA_s|\\
=&\gamma!_{\pi}(I)+\g(I)+\sum_{a\in\E(I,\ssI^c)}s(a)-|\ssA_s(I)-\ssA_s^{\rmint}(I)| +|\ssA_s(\ssI^c,I)|+|\ssA_s(I)|\\
=&\gamma!_{\pi}(I)+\g(I)+\sum_{a\in\E(I,\ssI^c)}s(a)+|\ssA_s^{\rmint}(I)| +|\ssA_s(\ssI^c,I)|.
\end{align*}
The first equality above is obtained from Lemma~\ref{lem:inverse-residue} applied to $\scrG=\GlobSp_{\pi}$ and $\ssP_v=\ssP_{s,v}$ for each $v$, and the descriptions of the $\ssP_{s,v}$ and the $\ssN_{s,v}$; the second is trivial; the third follows from Lemma~\ref{lem-useful-s} as well as the trivial equalities $\barA_s(I,\ssI^c)=\ssA_s(\ssI^c,I)$ and 
$|\barA_s(I)|=|\ssA_s(I)|$; the fourth is trivial.  Clearly, the last sum is $\sshateta_{s,\pi}(I)$, 
which concludes the proof. 
\end{proof}

\subsection{Submodular function $\nu!^*_{\hatWexp_{s,\pi}}$}  \label{sec:submodular-hatWexp}
We now express the submodular function 
\[
\nu!^*_{\hatWexp_{s,\pi}} \colon \ssub{2}!^V \to \Z, \qquad  I \mapsto \dim \proj{I}(\hatWexp_{s,\pi}),
\]
assuming that the function $\sshateta_{s,\pi}$ defined in the last section is nonnegative and the branches on the curves $\varC_v$ over the nodes of $X$ are in general position.

\begin{thm}\label{step1sing} 
Assume $\car\k=0$ and that the branches over nodes on $\varC_v$ for each $v\in V$ are in general position. Assume $\sshateta_{s,\pi}(I)\geq 0$ for all $I \subseteq V$. Then, we have
\[
\nu!^*_{\hatWexp_{s,\pi}}= \upmin(\sshateta_{s,\pi}).
\]
In particular, $\dim \hatWexp_{s,\pi}=g+|\ssA_s^\rmint|$.
\end{thm}

\begin{proof} For each $v\in V$ and each point $p$ on $\varC_v$, consider the filtration $\filt_p^\bullet$ of $H^0(\varC_v,\ssL_{s,v})$ given by  $\filt_p^j\coloneqq H^0(\varC_v,\ssL_{s,v}(-jp))$ for $j=0,\dots,\dim H^0(\varC_v,\ssL_{s,v})$. Since $\car\k=0$, Lemma~\ref{filt1} yields that $\filt_p^\bullet$ is a complete flag for each $p$ in general position, and that the collection $\Fl_v$ of these flags of $H^0(\varC_v,\ssL_{s,v})$ is in general position.

As in Section~\ref{sec:submodular-subset-upmin}, for each vertex $v \in V$, denote by $\xi!_v\colon \ssub{2}!^V\to\Z$ the function which takes value $-1$ on each subset $I\subseteq V$ if $I$ contains $v$, and value $0$ otherwise. Put
\begin{equation}\label{eq:xiN}
\xi\coloneqq\sum_{v\in V}\deg(\ssN_{s,v})\xi_v.
\end{equation}
Put $\mathcal W\coloneqq\{\hatWexp_{s,\pi}^+\}$. For each $v\in V$, put $\ssm_v\coloneqq|\bar\A_{s,v}|$, and let $\unde=(\sse_{v,1},\dots,\sse_{v,\ssm_v})$ be the partition of $\deg(\ssN_{s,v})$ with $\sse_{v,j}\coloneqq 1+s(\ssa_j)$ for every $j$, where $\ssa_1,\dots,\ssa_{\ssm_v}$ is any ordering of the arrows in $\bar\A_{s,v}$. Proposition~\ref{prop:connection-hateta-gamma} yields that 
\begin{equation}\label{eq:sshateta-formula}
\sshateta_{s,\pi}=\nu!_{\hatWexp_{s,\pi}^+}^*+\xi,
\end{equation}
whence $\sshateta_{s,\pi}\geq 0$ if and only if $\nu!_{\hatWexp_{s,\pi}^+}^*+\xi\geq 0$. We may thus apply Theorem~\ref{thm:global} to obtain the first part of the theorem. 

The second part follows from the first, by putting $I=V$, which gives $\dim \hatWexp_{s,\pi} = \sshateta_{s,\pi}(V)$, and then applying Proposition~\ref{prop:hateta-submodular}.
\end{proof}

If $\sshateta_{s,\pi}$ is not positive, for future use, we note the following result.

\begin{prop}\label{prop:step1singnot} 
Assume that $\car\k=0$ and the branches over nodes on $\varC_v$ for each $v\in V$ are in general position. If $\sshateta_{s,\pi}(I)\leq 0$ for some nonempty $I\subseteq V$, then $\proj{v}(\hatWexp_{s,\pi})=0$ for some $v\in V$.
\end{prop}

\begin{proof} We follow the proof given to Theorem~\ref{step1sing} but finish by applying Proposition~\ref{prop:zeroprojection}. 
\end{proof}

\subsection{Submodular function $\eta!_{s,\pi}$}\label{sec:etah} Define the function $\eta!_{s,\pi}\colon\ssub{2}!^V\to\Z$ by
\[
\eta!_{s,\pi}(I)\coloneqq\gamma!_{\pi}(I)+\g(I)+ \zeta!_s(I) = \sshateta_{s,\pi}(I) - |\ssA^{\rmint}_s(I)|.
\]

\begin{prop}\label{prop:eta-submodular} The function $\eta!_{s,\pi}$ is submodular of range $g$. It is simple if and only if for each nonempty proper subset $I$ of $V$, we have 
\begin{itemize}
\item either, $\gamma!_{\pi}(I) + \gamma!_{\pi}(\ssI^c) > |E|-|V|+1$,
\item or, there is an integer vertical edge for $s$ connecting $I$ and $\ssI^c$ 
\end{itemize}
\end{prop}

\begin{proof} The first claim is straightforward, as $\eta!_{s,\pi}$ is the sum of two submodular functions, $\gamma!_{\pi}$ and $\zeta!_s$, the latter by Proposition~\ref{prop:zetah}, and the modular function $\g$. As for its range, $\eta!_{s,\pi}(V) = \sshateta_{s,\pi}(V) - |\ssA_s^\rmint| = g$ by Proposition~\ref{prop:hateta-submodular}. This proves the first assertion.

As for the second, we first observe that for $\varphi=\varphi!_1+\varphi!_2$ with $\varphi!_1, \varphi!_2$ submodular, $\varphi$ is split at $I$ (that is, $\varphi(I)+\varphi(\ssI^c)=\varphi(V)$) if and only if both $\varphi!_1$ and $\varphi!_2$ are split at $I$. 

Since $\gamma!_{\pi}(V)=|E|-|V|+1$ by Theorem~\ref{thm:dimG}, for each nonempty proper subset $I\subset V$, we have that $\gamma!_{\pi}$ splits at $I$ if and only if $\gamma!_\pi(I)+\gamma!_\pi(\ssI^c)=|E|-|V|+1$. 

Also, since the function $\zeta!_s-\varphi!_s$ is modular by Proposition~\ref{prop:zetah}, where $\varphi!_s(I)\coloneqq|\ssA^{\rmint}_s(\ssI^c,I)|$ for each $I\subseteq V$, we have that $\zeta!_s$ splits at $I$ if and only if $\varphi!_s$ splits at $I$. Since $\varphi!_s$ is nonnegative with range $0$, this happens if and only if both $\ssA^\rmint_s(\ssI^c,I)$ and $\ssA^\rmint_s(I,\ssI^c)$ are empty, that is, no vertical integer edge for $s$ connects $I$ and $\ssI^c$.
\end{proof}

\subsection{The expected space $\Wexp_{s,\pi}(\varrho)$ associated to a gluing $\varrho$} \label{sec:submodular-Wh}

Viewing each $\ssL_{s,v}$ as a subsheaf of the constant sheaf $\Omega!_v$ on $\varC_v$, and their sum as a subsheaf of the sum on $X$ of the $\Omega!_v$, we consider the subsheaf on $X$
\[
\ssL_{s}^{\varrho}\subseteq\bigoplus_{v\in V}\ssL_{s,v}
\]
defined by imposing
\begin{enumerate}[label=(Exp\arabic*)]
\item\label{Exp1} equality away from the nodes $\ssp^e$ corresponding to integer edges $e$ for $s$; 
\item\label{Exp2} Rosenlicht condition at the node $\ssp^e$ corresponding to a horizontal edge $e$ for $s$, that is, the sum of the residues of the local differentials at the branches over $p^e$ is zero;
\item \label{Exp3} partial gluing of the two sheaves $\ssL_{s,\te_a}$ and $\ssL_{s,\he_a}$ at the points $\ssp^a$ and $\ssp^{\bar a}$ for each $a\in \ssA_s^{\rmint}$, i.~e., the induced maps $\ssL_s^{\varrho}\rest{p^e}\to\ssL_{s,\te_a}\rest{\ssp^{a}}$ and $\ssL_s^{\varrho}\rest{p^e}\to\ssL_{s,\he_a}\rest{\ssp^{\bar a}}$ differ by a given isomorphism
\[
\varrho!_{a}\colon \ssL_{s,\te_a}\rest{\ssp^a} \to \ssL_{s,\he_a}\rest{\ssp^{\bar a}}\qquad  \qquad \textrm{for each }a\in \ssA_s^{\rmint}.
\]
\end{enumerate}
We denote by $\varrho = (\varrho!_a)_{a\in \ssA_s^\rmint}$ the given collection of isomorphisms $\varrho!_a$ for $a\in \ssA_s^{\rmint}$.

\begin{defi}\label{defi:expsp}
Given a collection $\varrho=(\varrho!_a)_{a\in \ssA^\rmint_s}$, the \emph{expected space of limit canonical series} $\Wexp_{s,\pi} = \Wexp_{s,\pi}(\varrho)$ is defined as
  \[
    \Wexp_{s,\pi} = \Wexp_{s,\pi}(\varrho) \coloneqq \hatWexp_{s,\pi} \cap H^0(X,\ssL_s^\varrho)\subset\Omega. \qedhere
  \]
\end{defi} 

As the name suggests, we will see in Section~\ref{sec:limit-canonical-series-FC} that the $\Wexp_{s,\pi}$ are limit canonical series on $X$ under certain conditions.

\begin{thm}\label{step2sing-expected} 
Let $\varrho =(\varrho!_a)_{a\in \ssA_s^\rmint}$ where
$\varrho!_{a}\colon \ssL_{s,\te_a}\rest{\ssp^a} \to \ssL_{s,\he_a}\rest{\ssp^{\bar a}}$ is an isomorphism for each $a\in \ssA_s^{\rmint}$. Assume that
\begin{enumerate}[label=(\arabic*)]
\item \label{expected1} $\car\k=0$ and the branches over nodes on $\varC_v$ for each $v\in V$ are in general position,
\item \label{expected2} $\eta!_{s,\pi} \geq 0$, that is, $\eta!_{s,\pi}$ takes nonnegative values. 
\end{enumerate}
Then, the submodular function $\nu!^*_{\Wexp_{s,\pi}}$ of the subspace $\Wexp_{s,\pi} = \hatWexp_{s,\pi} \cap H^0(X,\ssL_s^\varrho)\subset\Omega$ satisfies
\[
\nu!^*_{\Wexp_{s,\pi}}=\upmin(\eta!_{s,\pi}).
\]
Furthermore, if 
\begin{enumerate}[resume,label=(\arabic*)]
    \item \label{expected3} $\proj{v}(\Wexp_{s,\pi})\neq 0$ for each $v\in V$,
\end{enumerate}
then $\Wexp_{s,\pi}$ generates $\ssL_s^{\varrho}$ at $\ssp^e$ for each integer vertical edge $e$ for $s$.
\end{thm}

\begin{proof} We assume~\ref{expected1}~and~\ref{expected2} hold, and show the first statement. As in the proof of Theorem~\ref{step1sing}, we apply Theorem~\ref{thm:global}, this time to the spaces
\[
\hatWexp_{s,\pi}^+[T]\coloneqq \Res^{-1}(\GlobSp_{\pi})\bigcap\Bigl(\bigoplus_{v\in V}H^0(\varC_v,\omega!_v(\ssP_{s,v}-\sum_{a\in \ssT_v}\ssp^a))\Bigr)\subset\Omega
\]
for $T\subseteq\ssA_s^{\rmint}$, with $\ssT_v = T\cap \E_v$. Since $s(a)>0$ for each $a\in\ssA_s^{\rmint}$, we have that $\ssP_{s,v}-\sum_{a\in \ssT_v}p^a$ is effective with the same support as $\ssP_{s,v}$. Whence, Lemma~\ref{lem:inverse-residue} can be applied to compute the submodular function of $\hatWexp_{s,\pi}^+[T]$, yielding
\[
\nu!_{\hatWexp_{s,\pi}^+[T]}^*=\nu!_{\hatWexp_{s,\pi}^+}^*+\sum_{a\in T}\xi!_{\te_a}.
\]
Thus, using~\ref{expected1} as in the proof of Theorem~\ref{step1sing},  Theorem~\ref{thm:global} yields that, for each pair of subsets $T,T'\subseteq\ssA_s^{\rmint}$ such that
\[
\nu!_{\hatWexp_{s,\pi}^+}^*+\xi+\sum_{a\in T}\xi!_{\te_a}+
\sum_{a\in T'}\xi!_{\he_a}\geq0,
\]
where $\xi=\sum\deg(\ssN_{s,v})\xi_v$ as in \eqref{eq:xiN}, the subspace
\[
\hatWexp_{s,\pi}[T,T']\coloneqq \Res^{-1}(\GlobSp_{\pi})\bigcap\Bigl(\bigoplus_{v\in V}H^0\big(\varC_v,\ssL_{s,v}(-\sum_{a\in \ssT_v}p^a-\sum_{a\in \ssT_v'}p^{\bar a})\big)\Bigr)\subset\Omega
\]
has submodular function
\[
\nu!_{\hatWexp_{s,\pi}[T,T']}^*=
\upmin\Big(\nu!_{\hatWexp_{s,\pi}^+}^*+\xi+\sum_{a\in T}\xi!_{\te_a}+
\sum_{a\in T'}\xi!_{\he_a}\Big).
\]
Now, as stated in \eqref{eq:sshateta-formula}, we have 
$\sshateta_{s,\pi}=\nu!_{\hatWexp_{s,\pi}^+}^*+\xi$. Since $\sshateta_{s,\pi}\geq 0$ from~\ref{expected2}, Theorem~\ref{step1sing} yields
\[
\nu!^*_{\hatWexp_{s,\pi}}= \upmin(\sshateta_{s,\pi})=\upmin(\nu!_{\hatWexp_{s,\pi}^+}^*+\xi).
\]
Using Proposition~\ref{prop:xi-upmin}, we conclude that
\begin{equation}\label{eq:Am-hatWexp}
\nu!_{\hatWexp_{s,\pi}[T,T']}^*=
\upmin\Big(\nu!_{\hatWexp_{s,\pi}}^*+\sum_{a\in T}\xi!_{\te_a}+
\sum_{a\in T'}\xi!_{\he_a}\Big)
\end{equation}
as long as 
\begin{equation}\label{eq:hatWexp-ineq}
\nu!_{\hatWexp_{s,\pi}}^*+\sum_{a\in T}\xi!_{\te_a}+
\sum_{a\in T'}\xi!_{\he_a}\geq0.
\end{equation}

The last conclusion is fundamental for applying Theorem~\ref{thm:eta-gen-2-hyper} to $\hatWexp_{s,\pi}$, which we will do, as follows. For each $a\in\ssA_s^\rmint$, let $\ssJ_a$ be the set of extremities of $a$. Order $\ssA_s^\rmint$ at will, and consider the sequence $\calJ$ of the $\ssJ_a$ under the order. Then, the function $\xi_{\calJ}$ coincides with the function $I \mapsto -|\ssA^\rmint_s(I)|$. So $\eta!_{s,\pi}=\sshateta_{s,\pi}+ \xi_{\calJ}$. Applying 
Proposition~\ref{prop:xi-upmin}, and using~\ref{expected2}, we get
\begin{equation}\label{eq:WJ}
\nu!_{\hatWexp_{s,\pi}}^*+ \xi_{\calJ}=\upmin(\sshateta_{s,\pi})+ \xi_{\calJ}\geq\upmin(\upmin(\sshateta_{s,\pi})+\xi_{\calJ})=\upmin(\eta!_{s,\pi})\geq 0.
\end{equation}

For each $a=uv\in\ssA_s^\rmint$, define the subspaces 
\[
\begin{aligned}
\rmH_{a,u}&\coloneqq H^0(\varC_u,\ssL_{s,u}(-p^a))\subseteq H^0(\varC_u,\ssL_{s,u}),\\
\rmH_{a,v}&\coloneqq H^0(\varC_v,\ssL_{s,v}(-p^{\bar a}))\subseteq H^0(\varC_v,\ssL_{s,v}),
\end{aligned}
\]
which have codimension at most $1$. 
Then, for each collection $\mathscr S$ of subsets $\ssS_a\subseteq J_a$ for $a\in\ssA_s^\rmint$,
\[
\hatWexp_{s,\pi}\cap\bigcap_{a\in\ssA_s^\rmint}\bigcap_{v\in\ssS_a}\proj{v}^{-1}(\rmH_{a,v})=\hatWexp_{s,\pi}[T,T']
\]
and 
\[
\sum_{a\in\ssA_s^\rmint}\sum_{v\in\ssS_a}\xi_v= \sum_{a\in T}\xi!_{\te_a}+\sum_{a\in T'}\xi!_{\he_a}
\]
for $T\coloneqq\{a\in\ssA_s^{\rmint}\,|\,\te_a\in\ssS_a\}$ and $T'\coloneqq\{a\in\ssA_s^{\rmint}\,|\,\he_a\in\ssS_a\}$. Since \eqref{eq:Am-hatWexp} holds for each such $T$ and $T'$, provided \eqref{eq:hatWexp-ineq} holds, and we have \eqref{eq:WJ}, we may apply Theorem~\ref{thm:eta-gen-2-hyper} to $\hatWexp_{s,\pi}$. It yields that there is a space $\rmH_a$ satisfying
\begin{equation}\label{eq:HHH}
\rmH_{a,u}\oplus\rmH_{a,v} \subsetneqq \rmH_a \subsetneqq H^0(\varC_u,\ssL_{s,u})\oplus H^0(\varC_v,\ssL_{s,v})
\end{equation}
for each $a=uv\in\ssA_s^\rmint$. Moreover, letting $\rmH_a$ be general for $a\in\ssA_s^\rmint$, the submodular function of
\[
\hatWexp_{s,\pi}\cap\bigcap_{a\in\ssA_s^\rmint}\proj{J_a}^{-1}(\rmH_a)
\]
is $\upmin(\nu!_{\hatWexp_{s,\pi}}^*+ \xi_{\calJ})$, which is equal to $\upmin(\eta!_{s,\pi})$ by Theorem~\ref{step1sing} and Proposition~\ref{prop:xi-upmin}.

Now, consider the action of the subtorus $\Gm^{\ssV_{\pi}}\subseteq\Gm^V$ on $\bigoplus\ssL_{s,v}$. For each fixed arrow $a=uv\in \ssA_s^\rmint$, since $\ssL_s^\varrho$ is one nondegenerate gluing of $\ssL_{s,u}$ at $\ssp^a$ with $\ssL_{s,v}$ at $\ssp^{\bar a}$, and since $u\supface_{\pi}v$, as $\psi=(\psi_w)_{w\in V}$ runs over all $\k$-points on $\Gm^{\ssV_{\pi}}$, the sheaf $\psi\cdot \ssL_s^{\varrho}$ runs through all such gluings.  Thus, for general $\psi$, the image of the map 
\[
H^0(X,\psi\cdot \ssL_s^\varrho)\longrightarrow H^0(\varC_u,\ssL_{s,u})\oplus H^0(\varC_v,\ssL_{s,v})
\]
is contained in a general space $\rmH_a$ satisfying \eqref{eq:HHH} for every $a=uv\in\ssA_s^\rmint$. In addition, since $\hatWexp_{s,\pi}$ is $\Gm^{\ssV_{\pi}}$-invariant by Proposition~\ref{prop:What-invariance}, we have
\begin{equation}\label{eq:psi-W-inv}
\psi\cdot\Wexp_{s,\pi}=\psi\cdot\hatWexp_{s,\pi}\cap H^0(X,\psi\cdot\ssL_s^\varrho)=
\hatWexp_{s,\pi}\cap H^0(X,\psi\cdot\ssL_s^\varrho).
\end{equation}
It follows that 
\begin{equation}\label{eq:finalW}
\psi\cdot\Wexp_{s,\pi}\subseteq\hatWexp_{s,\pi}\cap \bigcap_{a\in\ssA_s^\rmint}\proj{J_a}^{-1}(\rmH_a).
\end{equation}
Now, on the one hand, by the generality of the $\rmH_a$, the larger space has submodular function $\upmin(\eta!_{s,\pi})$, which has range $g$. And on the other hand, by definition, the dimension of $\Wexp_{s,\pi}$ is at least 
$\dim\hatWexp_{s,\pi}-|\ssA_s^{\rmint}|$, which is equal to $g$ by Theorem~\ref{step1sing}. Thus, equality holds in \eqref{eq:finalW}, and hence the submodular function of $\Wexp_{s,\pi}$ is $\upmin(\eta!_{s,\pi})$, finishing the proof of the first statement of the theorem.

\medskip

Assume now that also \ref{expected3} holds. Then $\nu!_{\Wexp_{s,\pi}}^* + \xi!_{\he_a}\geq 0$ for each $a\in\ssA_s^{\rmint}$. It follows from the first statement of the theorem that $\eta!_{s,\pi}+\xi!_{\he_a}\geq 0$. Since  $\eta!_{s,\pi}=\sshateta_{s,\pi}+ \xi_{\calJ}$, we get
\[
\sshateta_{s,\pi}+\xi_{\calJ}+\xi_{\he_a}\geq 0.
\]
Using $\nu!_{\hatWexp_{s,\pi}}^*=\upmin(\sshateta_{s,\pi})$ and Proposition~\ref{prop:xi-upmin}, we have that~\eqref{eq:hatWexp-ineq} holds and thus~\eqref{eq:Am-hatWexp} holds for $T=\emptyset$ and $T'=\{a\}$, that is,
\[
\nu_{\hatWexp_{s,\pi}[\emptyset,\{a\}]}^*= \upmin(\nu_{\hatWexp_{s,\pi}}^*+\xi_{\he_a}).
\]
Applying Proposition~\ref{prop:xi-upmin} once again, we get
\[
\nu_{\hatWexp_{s,\pi}[\emptyset,\{a\}]}^* +\xi_{\calJ}\geq 0.
\]

We apply now Theorem~\ref{thm:eta-gen-2-hyper} to $\hatWexp_\ssh[\emptyset,\{a\}]$. As before, using \eqref{eq:psi-W-inv}, it yields that for a general  $\psi\in\Gm^{\ssV_{\pi}}$, the intersection 
\[
\hatWexp_{s,\pi}[\emptyset,\{a\}]\cap H^0\bigl(X,\psi\cdot\ssL_s^{\varrho}\bigr)=\big(\psi\cdot \Wexp_{s,\pi}\big)\cap
\ssub{\theta}!_{\he_a}^{-1}\Bigl(H^0\bigl(\varC_{\he_a},\ssL_{s,\he_a}(-\ssp^{\bar a})\bigr)\Bigr)
\]
is contained in a space that 
has submodular function given by 
\[
\upmin\big(\nu!_{\hatWexp_{s,\pi}[\emptyset,\{a\}]}^*+\xi_{\calJ}\big)=\upmin\big(\nu!_{\Wexp_{s,\pi}}^*+\xi_{\he_a}\big),
\]
where we applied Proposition~\ref{prop:xi-upmin} again. Since the range of this function is $g-1$, we must have that 
\[
\proj{\he_a}(\psi\cdot \Wexp_{s,\pi})\not\subseteq H^0\bigl(\ssL_{s,\he_a}(-\ssp^{\bar a})\bigr).
\]
Since $\ssub{\theta}!_{\he_a}(\psi\cdot \Wexp_{s,\pi})=\ssub{\theta}!_{\he_a}(\Wexp_{s,\pi})$, 
it follows that $\Wexp_{s,\pi}$ generates $\ssL_s^{\varrho}$ at the node corresponding to $a$. Since $a\in\ssA_s^{\rmint}$ was arbitrary, the proof is complete.
\end{proof}

\subsection{Failure of nonnegativity} 
We have the following complementary result.

\begin{prop}\label{prop:step2singnot} Assume $\car\k=0$ and that the branches over nodes on $\varC_v$ for each $v\in V$ are in general position. Let $\varrho =(\varrho!_a)_{a\in \ssA_s^\rmint}$ be a collection of isomorphisms $\varrho!_{a}\colon \ssL_{s,\te_a}\rest{\ssp^a} \to \ssL_{s,\he_a}\rest{\ssp^{\bar a}}$ for $a\in \ssA_s^{\rmint}$. If $\eta!_{s,\pi}(I) \leq 0$ for some nonempty subset $I\subseteq V$, then $\proj{v}(\Wexp_{s,\pi})=0$ for some $v\in V$.
\end{prop}

\begin{proof} Consider first the case where $\sshateta_{s,\pi}(I)\leq 0$ for some nonempty $I\subseteq V$. In this case,  Proposition~\ref{prop:step1singnot} yields that $\proj{v}(\hatWexp_{s,\pi})=0$ for some $v\in V$. Since $\Wexp_{s,\pi} \subseteq \hatWexp_{s,\pi}$, it follows that $\proj{v}(\Wexp_{s,\pi})=0$, as required. 

We may thus assume $\sshateta_{s,\pi}\geq 0$. Consider now the case where $\eta!_{s,\pi}\geq 0$. Since $\eta!_{s,\pi}(I)\leq 0$ for some nonempty subset $I\subseteq V$, we get $\eta!_{s,\pi}(I)= 0$.  Theorem~\ref{step2sing-expected} then yields $\proj{I}(\Wexp_{s,\pi})=0$, whence $\proj{v}(\Wexp_{s,\pi})=0$ for each $v\in I$. 

In the remaining case, we have $\sshateta_{s,\pi}\geq 0$ but $\eta!_{s,\pi}(I)<0$ for some subset $I\subseteq V$ (which has to be nonempty). By definition, 
$\eta!_{s,\pi}=\sshateta_{s,\pi}+ \xi!_{\calJ}$, 
where $\ssJ_a$ is the set of extremities of $a$ for each $a\in\ssA_s^{\rmint}$ and $\calJ$ is the sequence of the $\ssJ_a$ for some order for the $a\in\ssA_s^{\rmint}$, as in the proof of Theorem~\ref{step2sing-expected}. Since $\sshateta_{s,\pi}\geq 0$, we may assume, up to replacing $I$ if needed, that there is a proper subset $S\subset\ssA_s^\rmint$ such that denoting by $\mathcal S$ the collection of $\ssJ_a$ for $a\in S$, we have 
\begin{equation}\label{eq:etahS}
\sshateta_{s,\pi}+\xi!_{\mathcal S}\geq 0
\end{equation}
but 
\[
\sshateta_{s,\pi}(I)+\xi!_{\mathcal S}(I)+\xi!_{\ssJ_b}(I)<0
\]
for some $b\in\ssA_s^\rmint\setminus S$. In particular, since $\xi!_{\ssJ_b}$ takes values $-1$ and $0$, we must have
\begin{equation}\label{eq:etahSb-0}
\sshateta_{s,\pi}(I)+\xi!_{\mathcal S}(I)=0
\end{equation}
Let $\ssL_{s}^{S, \varrho}\subseteq\bigoplus\ssL_{s,v}$ be the subsheaf defined similar to $\ssL_s^{\varrho}$ but we now impose only the gluing data $\varrho!_a$ for $a\in S$ (thus, treating the integer vertical edges underlying the arrows in $\ssA_s^\rmint\setminus S$ as if they were noninteger). Then, $\ssL_s^\varrho \subseteq\ssL_s^{S,\varrho}$. As in the proof of Theorem~\ref{step2sing-expected}, applying Theorem~\ref{thm:eta-gen-2-hyper} to $\hatWexp_{s,\pi}$, it follows from \eqref{eq:etahS} that $\hatWexp_{s,\pi}\cap H^0(X,\psi\cdot \ssL_s^{S,\varrho})$ is contained in a space with submodular function equal to 
\[
\upmin\Big(\sshateta_{s,\pi}+ \xi!_{\mathcal S}\Big)
\]
for a general point $\psi \in\Gm^{\ssV_{\pi}}$. But then, $\proj{I}\big(\hatWexp_{s,\pi}\cap H^0(X,\psi\cdot \ssL_s^{S,\varrho})\big)=0$ by \eqref{eq:etahSb-0}. Since
\[
\psi\cdot\Wexp_{s,\pi}=\psi\cdot \big(\hatWexp_{s,\pi} \cap H^0(X,\ssL_s^{\varrho})\big) \subseteq \psi\cdot\big(\hatWexp_{s,\pi} \cap H^0(X,\ssL_s^{S,\varrho})\big) \subseteq \hatWexp_{s,\pi}\cap H^0(X,\psi\cdot\ssL_s^{S,\varrho}),
\]
it follows that $\proj{I}(\psi\cdot\Wexp_{s,\pi})=0$, whence $\proj{I}(\Wexp_{s,\pi})=0$. Again, $\proj{v}(\Wexp_{s,\pi})=0$ for  each $v\in I$, as required. 
\end{proof}


\section{Limit canonical series}\label{sec:limit-canonical-series-FC}

Our aim in this section is to give a characterization of the fundamental collection of limit canonical series $\FC(\varX)$ for a smooth curve $\varX$ whose tropicalization has nodes in general position; see Definition~\ref{def:FC}. This is the subject of Theorem~\ref{thm:characterization-FC}.

We let $\varX$ be a smooth proper connected curve of positive genus over an algebraically closed field $\K$ which is complete for a nontrivial non-Archimedean valuation $\valuation$ with value group $\valgroup\subseteq\R$ and residue field $\k$. Let $(X,\Gamma)$ be a tropicalization of $\varX$. Here, $X$ is a nodal connected reduced proper curve over $\k$ and $\Gamma$ is a metric graph with model $(G,\ell)$, where $G=(V,E)$ is the dual graph of $X$. Recall the notation from the beginning of Section~\ref{sec:expected}.

Let $\varH$ be the space of Abelian differentials on $\varX$. For each $h\colon V\to\valgroup$, we let $\ssW_h\subseteq\Omega$ be the reduction of $\varH$ relative to $h$; see Definition~\ref{defi:rel-red-space}. Recall that  $\ssW_h$ depends on a choice of sections $\ba\colon\valgroup\to\K\setminus\{0\}$ of the valuation $\valuation$ which we now fix.

\subsection{Sheaves on residual curves}\label{sec:residual-sheaves}
Each function $h \colon V \to \valgroup$ gives rise to a slope-level pair $(\pi,s)=(\sspi_h,\slztwist{\ell}h)$ on $G$. Recall the sheaves $\ssL_{s,v}$ and divisors $\ssP_{s,v}$ and $\ssN_{s,v}$ from Section~\ref{sec:slope-sheaves}, and the notation in Definition~\ref{defi:slopefcn}. In particular, $\ssE_s^{\rmint}$ (resp.~$\A_s^{\rmint}$) denotes the set of integer vertical edges (resp.~upward arrows) for $s$.

The following two propositions, proved using Domination Property~\cite[Lem.~5.1]{AE-modular-polytopes}, will allow us to define a gluing of the sheaves $\ssL_{s,v}$ along the nodes $\ssp^a$ and $\ssp^{\bar a}$ for $a\in\A_s^\rmint$, using the data of the space $\ssW_h$.

\begin{prop}\label{prop:WhLh} For each function $h\colon V\to\valgroup$, we have $\ssW_h \subseteq \bigoplus H^0(\varC_v, \ssL_{s,v})$, where $s$ is the slope function corresponding to $h$.
\end{prop}
  
\begin{proof} Let $\varalpha \in \varM_h$ and $F\coloneqq \trop(\varalpha)$. By definition, for each $v\in V$, we have $F(v) =h(v)$ if and only if $(\red_h(\varalpha))_v\neq 0$, in which case  $(\red_h(\varalpha))_v =\wt\varalpha_v$. Recall that we denote by $h$ as well the admissible extension of $h$ to $\Gamma$; see Definition~\ref{defi:adm-ext} and Section~\ref{sec:red-diff}. Since $\varalpha\in\varM_h$, by Domination Property~\cite[Lem.~5.1]{AE-modular-polytopes}, we have $F(x) \geq h(x)$ for every $x\in \Gamma$. Therefore, if $(\red_h(\varalpha))_v\neq 0$, then for each unit outgoing tangent vector $\nu$ to $\Gamma$ at $v$, we have $\slp_{\nu}F\geq \slp_{\nu}h$. Now, Slope Formula, Lemma~\ref{lem:slope-formula}, states that $\ord_{p^{\nu}}(\wt \varalpha_v) = - 1 +\slp_{\nu}F$. In addition, if $\nu$ is parallel to an arrow $a \in \E_v$, we have $p^{\nu}=p^a$ and $s(a) = -\slp_{\nu}h$, the latter by Theorem~\ref{thm:existence-uniqueness-canonical-extension}. Thus,
\[
\div(\wt\varalpha_v) + \ssP_{h,v} - \ssN_{h,v} = \sum_{\nu}\big(-\slp_{\nu}h + \slp_{\nu}F\big)\ssp^{\nu} \geq 0,
\]
which shows that $\wt\varalpha_v \in H^0(\varC_v, \ssL_{h,v})$, as required. 
\end{proof}
 
\begin{prop}\label{prop:sim-vanishing} Let $h \colon V \to \valgroup$ be a function and $s\colon\E\to\Z$ the corresponding slope function. For each $a=uv\in\A_h^\rmint$ and each Abelian differential $\varalpha\in\varM_h$, the following two statements are equivalent:
\begin{enumerate}
    \item $\ssub{(\red!_h(\varalpha))}!_u \neq 0$ and $\ord_{\ssp^a}(\wt\varalpha_u) +s(a)+1=0$.
    \item $\ssub{(\red!_h(\varalpha))}!_v\neq 0$ and $\ord_{\ssp^{\bar a}}(\wt\varalpha_v) +s(\bar a)+1=0$.
\end{enumerate}
In other words, $\ssub{(\red!_h(\varalpha))}!_u$ does not vanish at $\ssp^a$ as a section of $\ssL_{s,u}$ if and only if $\ssub{(\red!_h(\varalpha))}!_v$ does not vanish at $\ssp^{\bar a}$ as a section of $\ssL_{s,v}$
\end{prop}

\begin{proof} Let $F=\trop(\varalpha)$. By symmetry, it will be enough to assume (1) and prove (2).  As in the proof of Proposition~\ref{prop:WhLh}, since $\varalpha\in\varM_h$, we have that $F \geq h$ pointwise. Also, the condition $\ssub{(\red!_h(\varalpha))}!_u \neq 0$ is equivalent to $F(u) = h(u)$. Let $\nu$ be the unit tangent vector at $u$ parallel to $a=uv$. By Slope Formula, Lemma~\ref{lem:slope-formula}, we have $\ord_{\ssp^a}(\wt\varalpha_u) = -1+\slp_{\nu}F$, with $\slp_{\nu}F$ denoting the outgoing slope of $F$ at $u$ along $\nu$. Therefore, the condition $\ord_{\ssp^a}(\wt\varalpha_u) +s(a)+1=0$ implies that $\slp_{\nu}h=-s(a)= \slp_{\nu}F$. Since $h(u) =F(u)$, the two functions $F$ and $h$ have to coincide on a small interval on $e$ containing  $u$. 

Now, since $e\in \ssE_h^\rmint$, the function $h$ is affine on $e$. Also, $\div(F) + K \geq 0$ by Theorem~\ref{thm:specialization-form}, which yields that  $\div(F)\geq 0$ on the interior of $e$. Since $F\geq h$ and since $F$ and $h$ already agree on an interval on $e$ containing $u$, we must have that $F$ and $h$ coincide on $e$. This implies $F(v)=h(v)$, whence $\ssub{(\red!_h(\varalpha))}!_v\neq 0$. Also, letting $\bar \nu$ be the unit tangent vector at $v$ parallel to $\bar a$, applying Lemma~\ref{lem:slope-formula} a second time, we conclude that $\ord_{\ssp^{\bar a}}(\wt\varalpha_v) = -1 + \slp_{\bar \nu}F=-1 + \slp_{\bar \nu}h =-1 - s(\bar a)$.

The last statement is clear from Proposition~\ref{prop:WhLh} and the definition of $\ssL_{s,u}$ and $\ssL_{s,v}$ in Section~\ref{sec:slope-sheaves}. 
\end{proof}

\subsection{Limit canonical series}\label{sec:gluing-sheaves}

We fix $h\colon V\to\valgroup$, and consider the reduction $\ssW_h$ of the space of Abelian differentials of $\varX$. We let $(s,\pi)$ be the slope-level pair associated to $h$ (and $\ell$). 

For each vertical integer edge $e=\{u,v\} \in \ssE_s^\rmint$, with arrows $a=uv$ and $\bar a=vu$, and each $\beta = (\beta_v)_{v\in V}\in\ssW_h$, Proposition~\ref{prop:sim-vanishing} says that $\beta!_u$ vanishes at $\ssp^a$ as a section of $\ssL_{s,u}$ if and only if $\beta!_v$ vanishes at $\ssp^{\bar a}$ as a section of $\ssL_{s,v}$. This implies that $\ssW_h$ satisfies at least one nondegenerate gluing condition at $\ssp^e$. Furthermore, there is a unique such gluing condition, or equivalently, $\ssW_h$ generates $\ssL_{s,u}$ at $\ssp^a$ (and $\ssL_{s,v}$ at $\ssp^{\bar a}$), if and only if there is an element $\beta\in\ssW_h$ such that $\beta!_u(\ssp^a)\neq 0$ (and $\beta!_v(\ssp^{\bar a})\neq 0$).

\smallskip
This being said, consider the subsheaf
\[
\ssL_h\subseteq\bigoplus_{v\in V}\ssL_{s,v}
\]
defined by properties \ref{Exp1} and \ref{Exp2} in Section~\ref{sec:submodular-Wh}, and an analogue of \ref{Exp3}, taking into account the observation made above on $\ssW_h$, namely,
\begin{itemize}
\item equality away from the nodes $\ssp^e$ corresponding to integer edges $e$ for $s$;
\item Rosenlicht condition at the node $\ssp^e$ corresponding to a horizontal edge $e$ for $s$; and 
\item a partial gluing of the two sheaves $\ssL_{s,u}$ and $\ssL_{s,v}$ at the points on $X$ above $p^e$ for each $e=\{u,v\}\in\ssE_s^\rmint$ satisfied by all elements of $\ssW_h$, that is, an isomorphism 
\[
\varrho!_{h,a}=\varrho!_{h,a}(\varX)\colon \ssL_{s,\te_a}\rest{\ssp^a} \to \ssL_{s,\he_a}\rest{\ssp^{\bar a}}\qquad \textrm{for each }a\in \ssA_s^{\rmint}
\]
such that the induced maps $\ssL_h\rest{p^e}\to\ssL_{s,\te_a}\rest{\ssp^{a}}$ and $\ssL_h\rest{p^e}\to\ssL_{s,\he_a}\rest{\ssp^{\bar a}}$ differ by it on $\ssW_h\subseteq H^0(X,\ssL_h)$.
\end{itemize}
Denote $\varrho!_{h}(\varX) \coloneqq (\varrho!_{h,a})_{a\in \ssA_s^\rmint}$. Then, notation as in Section~\ref{sec:submodular-Wh}, we have that $\ssL_h=\ssL_s^{\varrho!_h(\varX)}$. Notice that, as $\ssW_h$ depends on the fixed choice of section $\ba$ of the valuation $\valuation$, so does $\varrho!_h(\varX)$, and in particular, $\ssL_h$.  Theorem~\ref{step2sing} below gives a sufficient condition for having unique such gluings $\varrho!_{h,a}$ for $a\in \ssA_s^\rmint$ determined by $\ssW_h$, in which case, $\ssL_h$ is completely determined by $\ssW_h$.

\smallskip

Assume $\car\k=0$. It follows from Theorem~\ref{thm:ImWh} and Proposition~\ref{prop:WhLh} that $\ssW_h\subseteq\hatWexp_{s,\pi}$ (see Section~\ref{sec:hatWexp}). Since the gluing data $\varrho!_h(\varX)$ is compatible with $\ssW_h$, we have $\ssW_h\subseteq H^0(X,\ssL_h)$.  Thus $\ssW_h \subseteq \Wexp_{s,\pi} = \Wexp_{s,\pi}(\varrho!_h(\varX))\coloneqq \hatWexp_{s,\pi} \cap H^0(X,\ssL_h)$. Equality holds under the conditions of  Theorem~\ref{step2sing-expected}.

\begin{thm}\label{step2sing} Assume that
\begin{enumerate}[label=(\arabic*)]
\item $\car\k=0$ and the branches over nodes on $\varC_v$ for each $v\in V$ are in general position,
\item  $\eta!_{s,\pi} \geq 0$, that is, $\eta!_{s,\pi}$ takes nonnegative values. 
\end{enumerate}
Then, $\ssW_h=\Wexp_{s,\pi}$ and $\nu!^*_{\ssW_h}=\upmin(\eta!_{s,\pi})$. Furthermore, if
\begin{enumerate}[label=(\arabic*), resume]
    \item $\proj{v}(\ssW_h)\neq 0$ for each $v\in V$,
\end{enumerate}
then $\ssW_h$ generates $\ssL_h$ at $\ssp^e$ for each integer vertical edge $e$, and therefore $\ssL_h$ is uniquely determined from $\ssW_h$.
\end{thm}

\begin{proof} It will be enough to show that $\ssW_h = \Wexp_{s,\pi}$, as the remaining statements then follow directly from Theorem~\ref{step2sing-expected}. 

As we observed above, $\ssW_h\subseteq\Wexp_{s,\pi}$. By Proposition~\ref{prop:eta-submodular}, $\eta!_{s,\pi}$ has range $g$. It thus follows from Theorem~\ref{step2sing-expected} that $\Wexp_{s,\pi}$ has dimension $g$. Since $\dim\ssW_h=g$ as well, we get the equality $\ssW_h = \Wexp_{s,\pi}$.
\end{proof}


\subsection{Characterization of the fundamental collection} \label{sec:singsmooth} Recall that a nonzero finite dimensional subspace $\ssW\subset \Omega$ is called simple if the corresponding submodular function $\nu!^*_W$ is simple.  Notice that if $\nu!^*_W$ is simple, then it is positive, that is, $\nu!^*_W(I)>0$ for each nonempty subset $I\subseteq V$. 

\smallskip

By Definition~\ref{def:FC}, the fundamental collection of limit canonical series is the collection of simple $\ssW_h$.  Our next theorem describes this collection. Before stating the result, we need the following lemma.

\begin{lemma}\label{lem:bound-s} There exists a positive constant $B=B(G, \g)$ such that for each edge length function $\ell\colon E \to \R_{>0}$ and each function $h\colon V\to\valgroup$ with $\deg(\ssL_{s,v})\geq 0$ for every $v\in V$, where $s=\slztwist{\ell}h$ is the corresponding slope function, we have that
\[
|\sss(a)|\leq B \qquad \text{for each $a\in \E$}.
\]
\end{lemma}

\begin{proof} The condition in the lemma is equivalent to $K+\div(\hat h)\geq 0$, where $\hat h$ is the admissible extension of $h$ to the metric graph $\Gamma$. This inequality gives bounds on the divisor of poles of $\hat h$. Using  Theorem~\ref{thm:existence-uniqueness-canonical-extension}, the result follows from \cite[Lem.~1.8]{GK08}, which bounds the slopes of $\hat h$ in terms of its divisor of poles.
\end{proof}

For each $h\colon V\to\valgroup$, denote by $(\sss_h,\sspi_h)$ the associated slope-level pair on $G$, and recall that $\eta!_{\sss_h,\sspi_h}=\gamma!_{\sspi_h}+\g+ \zeta!_{\sss_h}$. By Proposition~\ref{prop:eta-submodular}, $\eta!_{\sss_h,\sspi_h}$ is simple if and only if for each bipartition $(I,J)$ of $V$: 
\begin{enumerate}[label=(\arabic*)]
\item either, $\gamma!_{\sspi_h}(I)+\gamma!_{\sspi_h}(J) > \gamma!_{\sspi_h}(V) = |E|-|V|+1$.
\item or, there is an integer vertical edge for $\sss_h$ connecting a vertex of $I$ to one of $J$.
\end{enumerate}
\begin{thm}[Characterization of the fundamental collection] \label{thm:characterization-FC} Assume $\car\k$ is $0$ and the branches over nodes of $X$ that lie on $\varC_v$ are in general position for each $v\in V$. Let $\Gamma$ be a metric graph with model $(G,\ell)$. Then, for each smooth proper curve $\varX$ over a valued field $\K$ which tropicalizes to $(X,\Gamma)$, the fundamental collection $\FC(\varX)$ of simple limit canonical series $\ssW_h$ coincides with the collection of spaces 
\[
\Wexp_{\sss_h,\sspi_h}(\varrho!_h) \coloneqq\Res^{-1}(\GlobSp_{\sspi_h})\cap H^0(X,\ssL_{\sss_h}^{\varrho!_h})\subset\Omega
\]
for $h\colon V \to \valgroup$ such that
\begin{itemize}
    \item $\eta!_{\sss_h,\sspi_h}$ is positive, that is, $\eta!_{\sss_h,\sspi_h}(I)>0$ for each nonempty $I\subseteq V$; and
    \item $\eta!_{\sss_h,\sspi_h}$ is simple,
\end{itemize}
and for a unique partial gluing data $\varrho!_h$. More precisely, for each $h\colon V\to\valgroup$, the subspace $\ssW_h\subset\Omega$ is simple if and only if $\eta!_{\sss_h,\sspi_h}$ is positive and simple, in which case there is a unique partial gluing data $\varrho!_h=\varrho!_h(\varX)$ for the sheaves $\ssL_{\sss_h,v}$ on $\varC_v$ defining the sheaf $\ssL_{\sss_h}^{\varrho!_h}$ on $X$ such that $\ssW_h\subseteq H^0(X,\ssL_{\sss_h}^{\varrho!_h})$ and, moreover, $\ssW_h=\Wexp_{\sss_h,\sspi_h}(\varrho!_h)$.
\end{thm}

\begin{proof} Let $B = B(G, \g)$ be the bound given by Lemma~\ref{lem:bound-s}.
There are only finitely many slope functions $s$ which verify 
\begin{equation}\label{eq:bound-slopes}
|s(a)|\leq B \qquad \qquad \text{ for all } a\in \E.
\end{equation}
Therefore, there are only finitely many divisor $\ssP_{s,v}$ and $\ssN_{\sss,v}$ associated to slope functions which verify the bounds in \eqref{eq:bound-slopes}. For each smooth proper curve $\varX$ over a valued field $\K$ with value group $\valgroup$, and each $h\colon V \to \valgroup$, such that $\ssW_h$ is simple, we have $\deg(\ssL_{\sss_h,v}) \geq 0$, and so $\sss_{h}$ verifies \eqref{eq:bound-slopes}. Therefore, it is possible to assume that the branches over nodes on $\varC_v$ for each $v\in V$ are in general enough position to be able to apply Proposition~\ref{prop:step2singnot} and Theorem~\ref{step2sing} for all such $h$.

Namely, consider $\varX$ over some valued field $\K$ with value group $\valgroup$ which tropicalizes to $(X, \Gamma)$, for some metric graph $\Gamma$. Let  $h\colon V\to\valgroup$ such that $\ssW_h\subseteq\Omega$ is simple. Then its submodular function $\nu!_{\ssW_h}^*$ is both simple and positive. Since $\ssW_h\subseteq\Wexp_{\sss_h,\sspi_h}$, the positivity of $\nu!_{\ssW_h}^*$ implies that of $\nu!_{\Wexp_{\sss_h,\sspi_h}}^*$, and hence the positivity of $\eta!_{\sss_h,\sspi_h}$ by Proposition~\ref{prop:step2singnot}. Theorem~\ref{step2sing} yields the equality $\ssW_h=\Wexp_{\sss_h,\sspi_h}$ and 
\begin{equation}\label{eq:nuthm}
\nu!_{\ssW_h}^*=\upmin(\eta!_{\sss_h,\sspi_h}).
\end{equation}
Furthermore, since $\nu!_{\ssW_h}^*$ is positive, Theorem~\ref{step2sing} yields that $\ssL_h$ is uniquely determined by $\ssW_h$, via the gluing data $\varrho!^h_\varX$. As $\nu!^*_{\ssW_h}$ is simple, so is $\eta!_{\sss_h,\sspi_h}$ by \eqref{eq:nuthm} and  Proposition~\ref{prop:simpleness-min-operation}. We have proved that all the limit canonical series $\ssW_h$ in the fundamental collection are of the form $\Wexp_{\sss_h,\sspi_h}(\varrho!_\varX^h)$ for $h\colon V \to \valgroup$ with $\eta!_{\sss_h,\sspi_h}$ positive and simple. 

Let $\ssh_1, \dots, \ssh_m \colon V \to \valgroup$ be all the functions with $\ssW_{h_1}, \dots, \ssW_{h_m}$ simple, as above. Let $(\sss_j, \sspi_j) = (\sss_{h_j}, \sspi_{h_j})$. Let $\P_{\sss_j, \sspi_j}$ be the polytope associated to $\eta!_{\sss_{j}, \sspi_{j}}$. We know that $\eta!_{\sss_{j}, \sspi_{j}}$ are all simple and positive, and the corresponding polytopes $\P_{\sss_j, \sspi_j}$ are all the full dimensional polytopes in the tiling of $\Delta!_g$ associated to $\varX$.

To finish, it will be enough to show that there is no slope-level pair $(\sss_{h}, \sspi_{h})$ associated to some $h\colon V\to \valgroup$ with $h-\ssh_j$ non-constant for all $j$ such that  $\eta!_{\sss_h,\sspi_h}$ is positive and simple. For the sake of a contradiction, suppose this is not the case, and pick such an $h$.  Changing the curve $X$ to $X'$, with the same dual graph $G=(V,E)$ and same components $\varC_v$, $v\in V$, we may assume that the nodes are in general position on each component so that we can apply Theorem~\ref{step2sing} to $\ssh_1, \dots, \ssh_m$, and $\ssh$.  Let $\varX'$ be a smooth proper curve over $\K$ which tropicalizes to $(X', \Gamma)$, the same metric graph as that of $\varX$, see~\cite[Thm.~3.24]{ABBR}. Applying Theorem~\ref{step2sing}, we find that $\Wexp_{\sss_{1}, \sspi_{1}}(\varrho!_{\ssh_1}'), \dots, \Wexp_{\sss_{m}, \sspi_{m}}(\varrho!_{\ssh_m}'), \Wexp_{\sss_{\ssh}, \sspi_{\ssh}}(\varrho!_{\ssh}')$, for $\varrho!_{\ssh_j}'=\varrho!_{\ssh_j}(\varX')$ for all $j$ and $\varrho!_{\ssh}'=\varrho!_{\ssh}(\varX')$, are in the fundamental collection of $\varX'$. Since the polytopes  $\P_{\sss_j, \sspi_j}$ already tile $\Delta!_g$, $\P_{\eta!_{\sss_h, \sspi_h}}$ should be one of them, and therefore, there exists $j$ with $h-\ssh_j$ constant, leading to a contradiction.
\end{proof}

\begin{remark}
Fixing isomorphisms $\mathcal O_{\varC_v}(\ssp^a)\rest{\ssp^a}\cong\k$ for each $v\in V$ and $a\in\E_v$, and isomorphisms $\omega!_X\rest{\ssp^e}\cong\k$ for each $e\in E$, we define isomorphisms $\ssL_{s,v}\rest{\ssp^a}\cong\k$ for each slope function $s\colon\E\to\Z$, each $v\in V$ and $a\in\E_v$. Using these isomorphisms, the partial gluing $\varrho!_h=\varrho!_h(\varX)$ appearing in the theorem is given by 
\[
\varrho!_{h,a} = \ssrho_e^{\sss_h(a)} \qquad \qquad  \text{for all } a\in \ssA_s^\rmint
\]
where $e$ is the underlying edge of $a$ and $\ssrho_e$ is the leading coefficient of the smoothing of the node $\ssp^e$ in $\varX$, see Appendix~\ref{sec:smoothing-and-gluing}. The theorem states that this data can be also recovered from the fundamental collection itself.
\end{remark}

\section{Fan structure on the cone of edge lengths}\label{sec:fan-structure}

Let $G=(V,E)$ be a finite connected graph. Recall the relevant notation in Section~\ref{sec:graphs}.

Given a general nodal curve $X$ with dual graph $G$, for each edge length function $\ell\colon E\to\R_{>0}$, Theorem~\ref{thm:characterization-FC} describes the set of fundamental collections $\FC(\varX)$ of limit canonical series $\ssW_h$ on $X$ arising from curves $\varX$ whose tropicalization is $(X,\Gamma)$, where $\Gamma$ is the metric graph with model $(G,\ell)$. Different edge length functions may yield the same sets of fundamental collections. Our aim in this section is to define a natural fan structure on the cone of edge lengths $\R_{\geq0}^E$ with relative open cones inside $\R_{>0}^E$ parametrizing edge length functions for which the sets of fundamental collections are the same. 

For given $\ell\in\R_{>0}^E$ and $h\colon V\to\R$, though the expected space $\Wexp_{\sss_h,\sspi_h}(\varrho)$ introduced in Definition~\ref{defi:expsp} depends on the "gluing" $\varrho$, it does not really depend on $\ell$ and $h$ but rather on the slope-level pair $(\sss_h,\sspi_h)$, which defines the residue space $\GlobSp_{\sspi_h}$ and the sheaves $\ssL_{\sss_h,v}$. Also, $\eta!_{\sss_h,\sspi_h}$ depends rather on $\sss_h$ and $\sspi_h$. We recall that it is the property that $\eta!_{\sss_h,\sspi_h}$ be positive and simple that determines whether $\Wexp_{\sss_h,\sspi_h}(\varrho)$ belongs to a fundamental collection (for the right $\varrho$ arising from a tropicalization).

The focus of this section will be on the collection $\PS_{\ell}$ of the pairs $(\sss_h,\sspi_h)$ for $h$ such that $\eta!_{\sss_h,\sspi_h}$ is positive and simple. We will see that the $\ell$ giving the same collection $\PS_{\ell}$ lie in the interior of a polyhedral cone in $\R_{\geq 0}^E$, and the collection of these cones (and their faces in the boundary $\R_{\geq 0}^E\setminus\R_{>0}^E$) is a fan on $\R_{\geq0}^E$.

\subsection{Cone of edge lengths}\label{sec:cone-slope-level}
Let $(s,\pi)$ be a slope-level pair on $G$.

\begin{defi}[Cone of edge lengths associated to a slope-level pair]
We define
\begin{align*}
\intcone_{s,\pi} &\coloneqq \Bigl\{\ell \in \R_{>0}^E \,\st\, \exists\,\, h\colon V \to \R \textrm{ with }\slztwist{\ell}h=s \textrm{ and } \sspi_h=\pi \Bigr\} \text{ and }\\
\cone_{s,\pi} &\coloneqq \textrm{ the closure of $\intcone_{s,\pi}$ in $\R^E_{\geq 0}$}.
\end{align*}
We call $\cone_{s,\pi}$ the \emph{cone of edge lengths} associated to $(s,\pi)$.
\end{defi}
 
 Note that by definition, $\cone_{s,\pi}$ is empty if and only if $\intcone_{s,\pi}$ is empty. We characterize the nonempty cones $\cone_{s,\pi}$ in Theorem~\ref{thm:characterization-slope-level-cone} and show that $\intcone_{s,\pi}$ is the relative interior of $\cone_{s,\pi}$.
 
To this purpose, it will be convenient in the sequel to define certain auxiliary graphs associated to ordered partitions and slope-level pairs. 

\smallskip

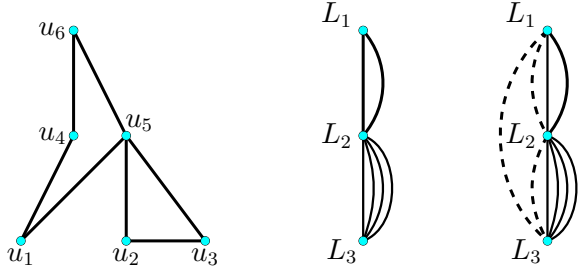
\begin{figure}[ht]
    \centering
\begin{tikzpicture}[scale=.65]
\draw[line width=0.4mm] (5,-2) -- (6,0);
\draw[line width=0.4mm] (7,-2) -- (7,0);
\draw[line width=0.4mm] (5,-2) -- (7,0);
\draw[line width=0.4mm] (8.5,-2) -- (7,0);
\draw[line width=0.4mm] (7,-2) -- (8.5,-2);
\draw[line width=0.4mm] (6,2) -- (6,0);
\draw[line width=0.4mm] (6,2) -- (7,0);

\draw[line width=0.1mm] (6,2) circle (0.8mm);
\filldraw[aqua] (6,2) circle (0.7mm);
\draw[line width=0.1mm] (6,0) circle (0.8mm);
\filldraw[aqua] (6,0) circle (0.7mm);
\draw[line width=0.1mm] (7,0) circle (0.8mm);
\filldraw[aqua] (7,0) circle (0.7mm);
\draw[line width=0.1mm] (5,-2) circle (0.8mm);
\filldraw[aqua] (5,-2) circle (0.7mm);
\draw[line width=0.1mm] (7,-2) circle (0.8mm);
\filldraw[aqua] (7,-2) circle (0.7mm);
\draw[line width=0.1mm] (8.5,-2) circle (0.8mm);
\filldraw[aqua] (8.5,-2) circle (0.7mm);

\draw [line width=0.4mm](11.5,2) -- (11.5,0);
\draw [line width=0.4mm](11.5,2) to [bend left=40] (11.5,0);
\draw [line width=0.3mm](11.5,0) -- (11.5,-2);
\draw [line width=0.3mm](11.5,0) to [bend left=20] (11.5,-2);
\draw [line width=0.3mm](11.5,0) to [bend left=40] (11.5,-2);
\draw [line width=0.3mm](11.5,0) to [bend left=70] (11.5,-2);


\draw[line width=0.1mm] (11.5,2) circle (0.8mm);
\filldraw[aqua] (11.5,2) circle (0.7mm);

\draw[line width=0.1mm] (11.5,0) circle (0.8mm);
\filldraw[aqua] (11.5,0) circle (0.7mm);

\draw[line width=0.1mm] (11.5,-2) circle (0.8mm);
\filldraw[aqua] (11.5,-2) circle (0.7mm);

\draw [line width=0.3mm](15,2) -- (15,0);
\draw [line width=0.4mm](15,2) to [bend left=40] (15,0);
\draw [line width=0.3mm](15,0) -- (15,-2);
\draw [line width=0.3mm](15,0) to [bend left=20] (15,-2);
\draw [line width=0.3mm](15,0) to [bend left=40] (15,-2);
\draw [line width=0.3mm](15,0) to [bend left=70] (15,-2);

\draw[dashed, line width=0.4mm] (15,2) to [bend right=30] (15,0);
\draw[dashed, line width=0.4mm] (15,0) to [bend right=30] (15,-2);
\draw[dashed, line width=0.4mm] (15,2) to [bend right=50] (15,-2);

\draw[line width=0.1mm] (15,2) circle (0.8mm);
\filldraw[aqua] (15,2) circle (0.7mm);
\draw[line width=0.1mm] (15,0) circle (0.8mm);
\filldraw[aqua] (15,0) circle (0.7mm);
\draw[line width=0.1mm] (15,-2) circle (0.8mm);
\filldraw[aqua] (15,-2) circle (0.7mm);

\draw[black] (5.6,2) node{$\ssu_6$};

\draw[black] (5.6,0) node{$\ssu_4$};
\draw[black] (7.25,0.25) node{$\ssu_5$};
\draw[black] (5,-2.3) node{$\ssu_1$};
\draw[black] (7,-2.3) node{$\ssu_2$};
\draw[black] (8.5,-2.3) node{$\ssu_3$};

\draw[black] (11,0) node{$\ssL_2$};
\draw[black] (11.1,-2.2) node{$\ssL_3$};
\draw[black] (11,2.3) node{$\ssL_1$};

\draw[black] (14.5,0) node{$\ssL_2$};
\draw[black] (14.6,-2.2) node{$\ssL_3$};
\draw[black] (14.5,2.3) node{$\ssL_1$};

\end{tikzpicture} 
\caption{A level graph $(G,\pi)$ on the left, $\ssG_\pi=(\ssV_\pi, \ssE_\pi)$ is given in the middle, and $\ssG^\dagger_\pi=(\ssV_\pi,\ssE^\dagger_\pi)$ on the right. Ghost edges are dashed. }
\label{fig:graph-with-ghosts}
\end{figure}

Given an ordered partition $\pi=(\ssV_{\pi},\subface_{\pi})$ of $V$, let $\ssG_{\pi}$ be the graph obtained by removing all the edges of $G$ which are horizontal with respect to $\pi$ and identifying vertices that lie on the same part of $\pi$. In our notation, we naturally identify the vertex and edge sets of $\ssG_{\pi}$ with $\ssV_{\pi}$ and $\ssE_\pi$ (the set of vertical edges for $\pi$), respectively, so $\ssG_\pi=(\ssV_\pi, \ssE_\pi)$. We denote by $\E_\pi$ the set of arrows on $\ssE_{\pi}$; see Figure~\ref{fig:graph-with-ghosts}.

We consider the graph  $\ssG^\dagger_\pi=(\ssV_\pi,\ssE^\dagger_\pi)$
obtained from $\ssG_{\pi}$ by adding an edge between any pair of vertices $u,v$ in $\ssV_{\pi}$. Let $\E_{\pi}^{\dagger}$ be the set of all arrows on edges $\ssE_\pi^\dagger$. We call the edges in $\ssE_\pi^\dagger\setminus\ssE_\pi$ \emph{ghost edges} and refer to arrows supported on them as \emph{ghost arrows}; see Figure~\ref{fig:graph-with-ghosts}. 

Given a slope-level pair $(s,\pi)$, extend $s$ to a function $\E^{\dagger}_{\pi}\to\Z\cup\{-\infty\}$ that we still denote $s$ by abuse of notation, by putting for each ghost arrow $a=uv$,
\[
s(a)\coloneqq \begin{cases}
-\infty&\text{if }u\subface_{\pi}v,\\
0&\text{otherwise.}
\end{cases}
\]
Also, we let $\ssG^{\rmint}_{s,\pi}\coloneqq (\ssV_{\pi},\ssE^{\rmint}_{s})$ be the spanning subgraph of $\ssG_{\pi}$ consisting of the edges in $\ssE_s=\ssE_\pi$ which are integer for $s$. Notice that $\ssG^{\rmint}_{s,\pi}$ contains no ghost edges. Clearly, both $\ssG_\pi$ and $\ssG^{\rmint}_{s,\pi}$ are spanning subgraphs of $\ssG^\dagger_\pi$.

\smallskip
Given an edge length function $\ell\colon E\to\R_{>0}$, its restriction to $\ssE_\pi\subseteq E$ defines an edge length function on the edges of $\ssG_\pi$ that we extend to $\ssE_\pi^\dagger$ by putting $\ell_e=1$ for all ghost edges $e\in\ssE_{\pi}^{\dagger}\setminus \ssE_\pi$. 

\smallskip

\begin{figure}[ht]
  \centering
\begin{tikzpicture}[scale=.27pt]

\draw[rounded corners, -{latex}, line width=0.4mm] (7,-3) to [bend right=40] (8,-1) to [bend left=40] (9,1);

\draw[rounded corners, -{latex}, line width=0.4mm] (9,2) to [bend left=30] (12,2) to [bend left=40] (12.7,-2) to [bend right=40] (14,-5);

\draw[rounded corners, -{latex}, line width=0.4mm] (14,-3) to [bend left=30] (15,-4.5) to [bend left=30] (14,-6) to [bend left=40] (7,-4);

\draw[line width=0.1mm] (7,-3) circle (1mm);
\filldraw[aqua] (7,-3) circle (0.9mm);

\draw[line width=0.1mm] (9,1) circle (1mm);
\filldraw[aqua] (9,1) circle (0.9mm);

\draw[line width=0.1mm] (9,2) circle (1mm);
\filldraw[aqua] (9,2) circle (0.9mm);

\draw[line width=0.1mm] (14,-5) circle (1mm);
\filldraw[aqua] (14,-5) circle (0.9mm);

\draw[line width=0.1mm] (14,-3) circle (1mm);
\filldraw[aqua] (14,-3) circle (0.9mm);

\draw[line width=0.1mm] (7,-4) circle (1mm);
\filldraw[aqua] (7,-4) circle (0.9mm);

\draw[black] (7.2,0) node{$\dpat_1$};

\draw[black] (10,-7) node{$\dpat_3$};
\draw[black] (13.7,-1.8) node{$\dpat_2$};

\end{tikzpicture} 
\caption{A sequence of three oriented paths $\dpat_1,\dpat_2,\dpat_3$ in the level graph $(G,\pi)$ satisfying $\te_{\dpat_{i+1}}\subfaceq_{\pi}\he_{\dpat_{i}}$ for $i=1,2,3$.}
\label{fig:compatible}
\end{figure}
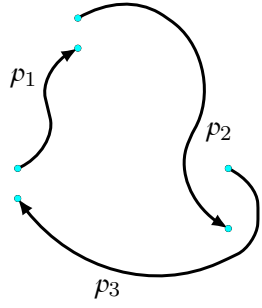

We say that an edge length function $\ell\in\R^E_{>0}$ is \emph{compatible with a slope-level pair $(s,\pi)$ on $G$} if for each circuit $\dq$ in $\ssG_\pi^\dagger$, we have 
\begin{align}\label{eq:slope-level-compatible}
  \sum_{a\in\dq}\ell_as(a)\leq 0,\text{ with equality if and only if the arrows in $\dq$ are all in $\ssG^{\rmint}_{s,\pi}$}.
\end{align}
Equivalently, $\ell$ is compatible with $(s,\pi)$ if for each finite sequence of oriented paths $\dpat_1,\dots,\dpat_k$ in $G$ with $\te_{\dpat_{i+1}}\subfaceq_{\pi}\he_{\dpat_{i}}$ for $i=1,\dots,k$, and $\dpat_{k+1} = \dpat_1$, we have 
\[
\sum_{i=1}^k\sum_{a\in\dpat_i}\ell_as(a)\leq 0,
\]
with equality if and only if each $\dpat_i$ contains only integer arrows for $s$, and $\te_{\dpat_{i+1}}$ and $\he_{\dpat_{i}}$ belong to the same part of the ordered partition $\pi$ for $i=1,\dots,k$; see Figure~\ref{fig:compatible}.

\smallskip

Notice that since $s(a)=0$ for a horizontal arrow $a$, the lengths of horizontal edges for $s$ are not constrained by the inequalities \eqref{eq:slope-level-compatible}. Also, if the arrows in $\dq$ are all in $\ssG^{\rmint}_{s,\pi}$, so is for the circuit $\bar\dq$, defined by the reversed arrows $\bara$ for $a\in q$, and since
\[
\sum_{a\in\dq}\ell_as(a)+\sum_{a\in\bar\dq}\ell_as(a)=\sum_{a\in\dq}\ell_a\big(s(a)+s(\bar a)\big)=0,
\]
equality holds automatically in \eqref{eq:slope-level-compatible} for all circuits in $\ssG^{\rmint}_{s,\pi}$. It is thus the "only if" part in \eqref{eq:slope-level-compatible} which constrains the lengths $\ell_a$.

\smallskip

In the sequel, it will be helpful to define the following auxiliary polyhedral cone
\begin{align*}
\altcone_{s,\pi}&\coloneqq \Bigl\{\ell \in \R_{\geq0}^E \,\st\, \textrm{ inequalities \eqref{eq:slope-level-compatible} hold for all circuits $\dq$ in $\ssG_\pi^\dagger$} \Bigr\}, \quad \text{and set } \\
\intaltcone_{s,\pi} &\coloneqq \Bigl\{\ell \in \R^E_{> 0} \, \st \, \textrm{$\ell$ is compatible with $(s,\pi)$}\Bigr\} \subseteq \altcone_{s,\pi} \cap \R^E_{> 0}.  
\end{align*}
Notice that $0\in \altcone_{s,\pi}$. Moreover, $\intaltcone_{s,\pi}$ might be empty even though $\altcone_{s,\pi}$ is not. Also, if $\intaltcone_{s,\pi}$ is nonempty, it is the relative interior of $\altcone_{s,\pi}$.

\begin{thm}\label{thm:characterization-slope-level-cone} Let $(s,\pi)$ be a slope-level pair on $G$. Then, $\intcone_{s,\pi}$ is the set of edge length functions $\ell\colon E\to\R_{>0}$ which are compatible with $(s,\pi)$, that is, 
$\intcone_{s,\pi}=\intaltcone_{s,\pi}$. In particular, $\intcone_{s,\pi}$ is nonempty if and only if $\cone_{s,\pi}=\altcone_{s,\pi}$, in which case $\intcone_{s,\pi}$ is the interior of $\cone_{s,\pi}$.
\end{thm}

\begin{proof} Let $\ell\in\intcone_{s,\pi}$. Then, there is $h\colon V\to\R$ such that $\slztwist{\ell}h=s$ and $\sspi_h=\pi$. Let $\dpat_1,\dots,\dpat_k$ be paths in $G$ satisfying $\te_{\dpat_{i+1}}\subfaceq_{\pi}\he_{\dpat_{i}}$ for $i=1,\dots,k$, with $\dpat_{k+1} = \dpat_1$. Then,
\begin{align*}
\sum_{i=1}^k\sum_{a\in\dpat_i}\ell_as(a)\leq&\sum_{i=1}^k\sum_{a\in\dpat_i}\big(h(\te_a)-h(\he_a)\big)=\sum_{i=1}^k\big(h(\te_{\dpat_i})-h(\he_{\dpat_i})\big)\\
=& h(\te_{\dpat_1})-h(\he_{\dpat_k})+\sum_{i=1}^{k-1}\big(h(\te_{\dpat_{i+1}})-h(\he_{\dpat_{i}})\big)\leq 0,
\end{align*}
with equality if and only if $\ell_as(a)=h(\te_a)-h(\he_a)$ for each arrow $a$ in each $\dpat_i$, and also $h(\te_{\dpat_{i+1}})=h(\he_{\dpat_{i}})$ for each $i=1,\dots,k$. Equivalently, equality holds if and only if each $\dpat_i$ contains only integer arrows for $s$, and $\te_{\dpat_{i+1}}$ and $\he_{\dpat_{i}}$ belong to the same part of $\pi$ for each $i=1,\dots,k$. Thus, $\ell$ is compatible with $(s,\pi)$, as stated.  

\smallskip

Conversely, let $\ell\colon E\to\R_{>0}$ be a function compatible with $(s,\pi)$.  Let $x\coloneqq\ssx_{\ell,s}\colon\E_{\pi}^\dagger\to\R\cup\{-\infty\}$ be defined by $x(a)\coloneqq \ell_a\sss(a)$ for each $a\in\E_{\pi}^{\dagger}$. Then, $x (a)+x(\bar a)\leq 0$ for all $a\in\E_{\pi}^\dagger$. Furthermore, the compatibility inequalities \eqref{eq:slope-level-compatible} can be rewritten in the form
\[
\sum_{a\in \dq} x(a) \leq 0\quad\text{for each circuit }\dq\text{ in }\ssG_\pi^\dagger.
\]
By Subintegrability Lemma~\ref{lem:key_lemma}, proved in Appendix~\ref{sec:subintegrable}, applied to $\ssG_{\pi}^\dagger$, there is a function $h\colon \ssV_\pi\to\R$ such that 
\begin{equation}\label{eq:xhx}
x(a)\leq h(u)-h(v)
\end{equation}
for each $a=uv \in\E_{\pi}^\dagger$, with equality if and only if either $x(a)+x(\bar a)=0$, or there is a circuit $\dz$ in $\ssG^\dagger_\pi$ containing $a$ such that $\sum_{b\in\dz}x(b)=0$. 

Composing with the projection $V\to\ssV_\pi$, we view $h$ as a function from $V$ to $\R$. We claim that $s=\slztwist{\ell}h$ on $\E$. Indeed, for each edge $e=\{u,v\}$, applying~\eqref{eq:xhx} to the two arrows $a=uv$ and $\bar a=vu \in\E$ over $e$, we get
\begin{equation}\label{eq:shs}
\sss(a)\leq\frac{h(u)-h(v)}{\ell_e}\leq -\sss(\bar a).
\end{equation}
Now, since $s$ is a slope function, $s(a)+s(\bar a)$ is equal to $0$ or $-1$. In the former case, $a$ is an integer arrow for $s$, and equalities holds in \eqref{eq:shs}, whence $s(a)=\slztwist{\ell}h(a)$. In the latter case, we show that $h(u)-h(v)< -\ell_e\sss(\bar a)$, which  implies $\slztwist{\ell}h(a)=s(a)$, as required. Indeed, suppose by contradiction that $h(u)-h(v)=-\ell_a\sss(\bar a)$, that is, equality holds in~\eqref{eq:xhx} for $\bar a$. Since $x(a)+x(\bar a)=-\ell_e<0$, there must exist a circuit $\dz$ in $\ssG^\dagger_\pi$ containing $\bar a$ such that $\sum_{b\in\dz}x(b)=0$. But then, by compatibility of $\ell$ with the pair $(s,\pi)$, the arrow $\bar a$ has to be an integer edge for $s$, a contradiction. We conclude the equality $s=\slztwist{\ell}h$ on $\E$.

It remains to show that $\pi=\sspi_h$. First, since $h$ factors through $\ssV_{\pi}$, it is constant on each part of $\pi$. Conversely, let $u,v\in \ssV_\pi$ be distinct vertices, and suppose without loss of generality that $u\supface_\pi v$. Consider the ghost arrow $uv\in \E^{\dagger}_{\pi}$. By definition of $x$, we have $x(uv)=0$.  Therefore, by \eqref{eq:xhx} applied to the arrow $uv$, we have $h(u)\geq h(v)$. We prove that $h(u) > h(v)$ concluding the proof. 

Suppose by contradiction that $h(u)= h(v)$, and consider now both ghost arrows $uv, vu\in \E_\pi^\dagger$. As $x(uv)=0=h(u)-h(v)$, and $x(uv)+x(vu) = -\infty$, there must exist a circuit $\dz$ in $\ssG_{\pi}^\dagger$ containing the ghost arrow $uv$ such that $\sum_{b\in\dz}x(b)=0$. But then, compatibility of $\ell$ with the pair $(s,\pi)$ implies that $\dz$ contains no ghost arrow, leading to a contradiction.
\end{proof}

We get the following immediate corollary.

\begin{cor} Each nonempty cone $\cone_{s,\pi}$ is a rational polyhedral cone in $\R^E_{\geq 0}$.
\end{cor}

\begin{remark} Though our interest is in the cones $\cone_{s,\pi}$, it is more practical to work with the cones $\altcone_{s,\pi}$, which our Theorem~\ref{thm:characterization-slope-level-cone} allows. In any case, our next main result, Theorem~\ref{thm:face-cone}, yields that $\altcone_{s,\pi}=\cone_{s',\pi'}$ for some slope-level pair $(s',\pi')$. Thus, the collection of the nonempty $\cone_{s,\pi}$ is the collection of the cones $\altcone_{s,\pi}$.
\end{remark}

\subsection{Dimension} The following proposition gives the dimensions of the cones introduced in Section~\ref{sec:cone-slope-level}. 

\begin{prop}\label{prop:dim_count} Let $(s,\pi)$ be a slope-level pair for the graph $G=(V,E)$. Then,
\begin{equation}\label{eq:ineq-dim}
\dim\altcone_{s,\pi} \leq |E|  - g(\ssG_{s,\pi}^{\rmint}).
\end{equation}
If $\intcone_{s,\pi}\neq\emptyset$, then equality holds.
\end{prop}  

\begin{proof} Each inequality \eqref{eq:slope-level-compatible} for $\dq$ in $\ssG_{\pi}^\dagger$ is an equality for all $\ell\in\altcone_{s,\pi}$ if $\dq$ is in $G^{\rmint}_{s,\pi}$. It follows that the codimension of $\altcone_{s,\pi}$ is at least the dimension of the cycle space of the graph $\ssG_{s,\pi}^{\rmint}$, which proves the inequality \eqref{eq:ineq-dim}. On the other hand, if $\intcone_{s,\pi}\neq\emptyset$, then it follows from the equality $\intcone_{s,\pi} =\intaltcone_{s,\pi}$ proved in Theorem~\ref{thm:characterization-slope-level-cone} that equality holds for all $\ell\in\intcone_{s,\pi}$ in \eqref{eq:slope-level-compatible} only for $\dq$ in $\ssG^{\rmint}_{s,\pi}$. Then, in this case, \eqref{eq:ineq-dim} is an equality.
\end{proof}  

\subsection{Essential circuits} Let $(s,\pi)$ be a slope-level pair. Not all circuits $\dz$ in $\ssG_{\pi}^{\dagger}$ are needed to define  $\altcone_{s,\pi}$, but only those verifying properties~\ref{eq:ess1}-\ref{eq:ess2}-\ref{eq:ess3}, given below. Indeed, let $\dz$ be a circuit in $\ssG_{\pi}^{\dagger}$. 

First, if $\dz$ contains a downward ghost arrow, then $\sum_{a\in \dz}\ell_a s(a)=-\infty$ for every $\ell\in\R_{\geq 0}^E$. Thus, to define $\altcone_{s,\pi}$ we need only consider the inequalities \eqref{eq:slope-level-compatible} arising from circuits $\dq$ in $\ssG_{\pi}^{\dagger}$ which verify
\begin{enumerate}[label=(Ess\arabic*)]
\item \label{eq:ess1} \emph{All the ghost arrows in $\dq$, if any, are upward}. 
\end{enumerate}

Also, note that if an arrow $b=uv$ in $\E_\pi = \E_s$ satisfies $s(b)=0$, then, since it is vertical for $s$, we have $s(\bar b)=-1$, and hence $b$ is noninteger for $s$ and $u\supface_\pi v$. Now, if the circuit $\dz$ contains $b$, replacing $b$ with the (upward) ghost arrow $uv\in \E_\pi^{\dagger}$, we obtain a new circuit $\dz'$ in $\ssG_{\pi}^{\dagger}$ for which
\begin{equation}\label{eq:dz=dz}
\sum_{a\in \dz'}\ell_a s(a) = \sum_{a\in \dz}\ell_a s(a)\quad\text{for all }\ell\in\R^E_{\geq 0}.
\end{equation}
Thus, to define $\altcone_{s,\pi}$ we need only consider the inequalities \eqref{eq:slope-level-compatible} arising from those circuits $\dq$ in $\ssG_{\pi}^{\dagger}$ that verify 
\begin{enumerate}[label=(Ess\arabic*), resume]
\item \label{eq:ess2} \emph{The circuit $\dq$ does not contain any arrow $a \in \E_\pi$ with $s(a)=0$}. 
\end{enumerate}

Finally, assume $\dz$ verifies \ref{eq:ess1}. Assume as well that $\dz$ contains two distinct (upward) ghost arrows $\ssu_1\ssv_1$ and $\ssu_2\ssv_2$ which \emph{overlap}, more precisely, which satisfy $\ssu_1  \supface_{\pi} \ssu_2 \supfaceq_{\pi} \ssv_1 \supface_{\pi} \ssv_2$. 

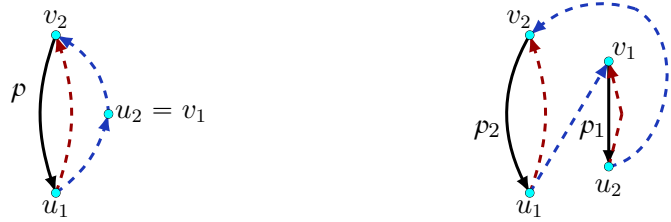
\begin{figure}[ht]
  \centering
\begin{tikzpicture}[scale=.7pt]
\draw[-{latex}, line width=0.4mm] (5,0) to [bend right=20] (5,-3);

\draw[dashed, -{latex}, pblue, line width=0.4mm] (5,-3) to [bend right=20] (6,-1.5);

\draw[dashed, -{latex}, darkred, line width=0.4mm] (5,-3) to [bend right=20] (5,0);


\draw[dashed, -{latex}, pblue, line width=0.4mm] (6,-1.5) to [bend right=30] (5,0);

\draw[line width=0.1mm] (5,-3) circle (0.8mm);
\filldraw[aqua] (5,-3) circle (0.7mm);

\draw[line width=0.1mm] (5,0) circle (0.8mm);
\filldraw[aqua] (5,0) circle (0.7mm);

\draw[line width=0.1mm] (6,-1.5) circle (0.8mm);
\filldraw[aqua] (6,-1.5) circle (0.7mm);


\draw[black] (5,-3.3) node{$\ssu_1$};
\draw[black] (5,0.3) node{$\ssv_2$};

\draw[black] (7,-1.5) node{$\ssu_2 = \ssv_1$};

\draw[black] (4.3,-1.1) node{$\dpat$};


\draw[-{latex}, line width=0.4mm] (14,0) to [bend right=30] (14,-3);

\draw[dashed, -{latex}, darkred, line width=0.4mm] (14,-3) to [bend right=20] (14,0);

\draw[dashed, -{latex}, pblue, line width=0.4mm]  (14,-3) to (15.5,-.5) [bend left=50] ;

\draw[-{latex}, line width=0.4mm] (15.5,-.5) to (15.5,-2.5);

\draw[dashed, darkred, line width=0.4mm]   (15.5,-2.5) to (15.75,-1.5) [bend right=60] ;

\draw[dashed, -{latex}, darkred, line width=0.4mm]   (15.75,-1.5) to (15.5,-.5) [bend right=10] ;

\draw[dashed, pblue, line width=0.4mm] (15.5,-2.5) to [bend right=70] (16,.5);

\draw[dashed, -{latex},  pblue, line width=0.4mm] (16,.5) to [bend right=30] (14,0);

\draw[line width=0.1mm] (14,-3) circle (0.8mm);
\filldraw[aqua] (14,-3) circle (0.7mm);

\draw[line width=0.1mm] (14,0) circle (0.8mm);
\filldraw[aqua] (14,0) circle (0.7mm);

\draw[line width=0.1mm] (15.5,-2.5) circle (0.8mm);
\filldraw[aqua] (15.5,-2.5) circle (0.7mm);

\draw[line width=0.1mm] (15.5,-.5) circle (0.8mm);
\filldraw[aqua] (15.5,-.5) circle (0.7mm);

\draw[black] (14,-3.3) node{$\ssu_1$};
\draw[black] (13.8,0.3) node{$\ssv_2$};

\draw[black] (15.5,-2.9) node{$\ssu_2$};
\draw[black] (15.8,-.3) node{$\ssv_1$};

\draw[black] (15.2,-1.8) node{$\dpat_1$};

\draw[black] (13.2,-1.8) node{$\dpat_2$};

\end{tikzpicture} 
\caption{The replacements used in \ref{eq:ess3}. Blue ghosts in $\dz$ are replaced by red ones.}
\label{fig:replacement1}
\end{figure}

If $\ssu_2=\ssv_1$, we replace $\ssu_1\ssv_1$ and $\ssu_2\ssv_2$ in $\dz$ with the upward ghost arrow $\ssu_1\ssv_2$ to obtain a circuit $\dz'$ satisfying \eqref{eq:dz=dz} and containing less ghost arrows than $\dz$; see Figure~\ref{fig:replacement1} on the left. Clearly, the inequality \eqref{eq:slope-level-compatible} for $\dq=\dz$ is implied by that for $\dq=\dz'$.

If $\ssu_2 \supface_{\pi} \ssv_1$, we remove $\ssu_1\ssv_1$ and $\ssu_2\ssv_2$ from $\dz$, resulting in two oriented paths $\dpat_1$ and $\dpat_2$, with $\dpat_1$ starting at $\ssv_1$ and ending at $\ssu_2$, and $\dpat_2$ starting at $\ssv_2$ and ending at $\ssu_2$. Then, we add the upward ghost arrow $\ssu_2\ssv_1$ to $\dpat_1$ to get a circuit $\dz_1$, and add the upward ghost arrow $\ssu_1\ssv_2$ to $\dpat_2$ to get a circuit $\dz_2$; see Figure~\ref{fig:replacement1} on the right. Then, we get
\[
\sum_{a\in \dz}\ell_a s(a) = \sum_{a\in \dz_1}\ell_a s(a) + \sum_{a\in \dz_2}\ell_a s(a)\quad\text{for all }\ell\in\R^E_{\geq 0}.
\]
Both $\dz_1$ and $\dz_2$ have less ghost arrows than $\dz$. And, the inequality \eqref{eq:slope-level-compatible} for $\dq=\dz$ is implied by those for $\dq=\dz_1$ and $\dq=\dz_2$. 

In any case, it follows by induction on the number of ghost arrows that we need only consider circuits $\dq$ in $\ssG_{\pi}^{\dagger}$ which verify 

\begin{enumerate}[label=(Ess\arabic*), resume]
\item \label{eq:ess3} \emph{The ghost arrows in $\dq$ do not overlap}. 
\end{enumerate}

\begin{defi}[Essential circuit for a slope-level pair] Let $(s,\pi)$ be a slope-level pair. A circuit $\dq$ in $\ssG_{\pi}^{\dagger}$ is called \emph{essential} if it verifies~\ref{eq:ess1}-\ref{eq:ess2}-\ref{eq:ess3}, listed above.
\end{defi}

Notice that circuits $\dq$ in $\ssG^{\rmint}_{s,\pi}$ are all essential. From the above discussion we obtain:

\begin{prop}\label{prop:essential-enough} Let $(s,\pi)$ be a slope-level pair on $G$. Then, the cone $\altcone_{s,\pi}$ is given in $\R^E_{\geq 0}$ by the inequalities \eqref{eq:slope-level-compatible} for all essential circuits $\dq$ in $\ssG_\pi^{\dagger}$.
\end{prop}

\subsection{Active circuits} Let $(s,\pi)$ be a slope-level pair on $G$. For each circuit $\dz$ in $\ssG_\pi^{\dagger}$, define the hyperplane
\[
\ssH_{\dz} \coloneqq\Bigl\{\ell \in \R^E\,\Big|\, \sum_{a\in \dz} \ell_a s(a) =0\Bigr\}.
\]
Since $\sum_{a\in \dz} \ell_a s(a)\leq 0$ for all $\ell\in \altcone_{s,\pi}$, we have that $H_{\dz}$ is a supporting hyperplane for $\altcone_{s,\pi}$, meaning that $\altcone_{s,\pi}$ lies in one of the two closed half-spaces defined by $\ssH_\dz$. Thus, $\altcone_{s,\pi}\cap H_{\dz}$ is a face of $\altcone_{s,\pi}$. 

\begin{defi}[Active circuits for a slope-level pair]
A circuit $\dz$ in $\ssG_\pi^\dagger$ is called \emph{active for a slope-level pair $(s,\pi)$} if 
\begin{itemize}
    \item $\dz$ is not in $\ssG^{\rmint}_{s,\pi}$, and
    \item the intersection $\altcone_{s,\pi}\cap H_{\dz}$  contains an edge length function $\ell\in\R^E_{>0}$. \qedhere
\end{itemize} 
\end{defi}

Note that, in general, the intersection $\altcone_{s,\pi}\cap H_{\dz}$ may not be a proper face of $\altcone_{s,\pi}$. However, if $\dz$ is active for $(s,\pi)$ and $\altcone_{s,\pi}=\cone_{s,\pi}$, then, by Theorem~\ref{thm:characterization-slope-level-cone}, the intersection $\altcone_{s,\pi}\cap H_{\dz}$ will be a proper face of $\altcone_{s,\pi}$ (and it will coincide  with $\cone_{s,\pi}\cap H_{\dz}$).

\subsection{Squashing process} \label{sec:squashing} Let $(s,\pi)$ be a slope-level pair on $G$. Our next result asserts that for an essential circuit $\dz$ in $\ssG_\pi^{\dagger}$ which is active for $(s,\pi)$, the face $\altcone_{s,\pi}\cap H_{\dz}$ is of the form $\altcone_{s',\pi'}$ for a slope-level pair $(s',\pi')$ defined by a \emph{squashing process} that we describe in this section.  First, we need the following lemma.

\smallskip

 Given a circuit $\dz$ in $\ssG_\pi^\dagger$  verifying~\ref{eq:ess1}, we say that \emph{an arrow $a=xy \in \E_{\pi}$ is covered by a ghost of $\dz$} if there is a ghost arrow $uv \in \dz$ such that $u \supfaceq_\pi x$ and $y\supfaceq_\pi v$. (Note that both $a$ and $\bar a$ are covered by a ghost $uv$ of $\dz$ if $u \supfaceq_\pi x,y\supfaceq_\pi v$.)

\begin{lemma}\label{lem:arrows-covered-by-ghosts} Let $\dz$ be a circuit in $\ssG_\pi^\dagger$ which is active for $(s,\pi)$. Let $b\in\E_\pi$ be an arrow covered by a ghost of $\dz$. Then, the following hold:
\begin{itemize}
    \item We have $s(b)\leq 0$ with equality only if $\bar b\not\in\dz$.
    \item If in addition $\dz$ is essential for $(s,\pi)$, and also $\bar b$ is covered by a ghost of $\dz$, then $b,\bar b\not\in\dz$.
\end{itemize} 
\end{lemma}
Note that since $\dz$ is active for $(s,\pi)$, there is $\ell\in\altcone_{s,\pi}$ with $\ell>0$ such that
\begin{equation}\label{eq:l-in-cone}
\sum_{a\in \dz}\ell_as(a)=0.
\end{equation}
In particular, $\dz$ automatically verifies~\ref{eq:ess1}.

\begin{proof} Write $b=xy$ with $x,y\in \ssV_\pi$. Let $uv\in\dz$ be a ghost covering $b$. 
We prove first that $s(b)\leq 0$. Consider the oriented path $\dpat$ from $u$ to $v$ in $\ssG_\pi^\dagger$ composed by the ghost arrow $ux$ if $u\neq x$, the arrow $b$, and the ghost arrow $yv$ if $y\neq v$; see Figure~\ref{fig:path1} on the left. Replace the ghost arrow $uv$ in $\dz$ by $\dpat$ to obtain the circuit $\dq$. Then,
\[
\sum_{a\in \dq} \ell_as(a) = \ell_bs(b) + \sum_{a\in \dz}\ell_a s(a) = \ell_{b}s(b), 
\]
because all ghosts of $\dq$ are upward and \eqref{eq:l-in-cone} holds. Since $\ell\in\altcone_{s,\pi}$, it follows that $\ell_{b} s(b) \leq 0$. Since $\ell>0$, we have $s(b)\leq 0$.

\begin{figure}[ht]
  \centering
\begin{tikzpicture}[scale=.7pt]
\draw[dashed, -{latex}, line width=0.4mm] (5,-3) to [bend left=20] (5,0);

\draw[dashed, -{latex}, pblue, line width=0.4mm] (5,-3) to [bend right=50] (7.5,-.5);

\draw[-{latex}, pblue, line width=0.4mm] (7.5,-.5) to [bend right=20] (6,-1.5);

\draw[dashed, -{latex}, pblue, line width=0.4mm] (6,-1.5) to [bend right=30] (5,0);

\draw[line width=0.1mm] (5,-3) circle (0.8mm);
\filldraw[aqua] (5,-3) circle (0.7mm);

\draw[line width=0.1mm] (5,0) circle (0.8mm);
\filldraw[aqua] (5,0) circle (0.7mm);

\draw[line width=0.1mm] (6,-1.5) circle (0.8mm);
\filldraw[aqua] (6,-1.5) circle (0.7mm);

\draw[line width=0.1mm] (7.5,-.5) circle (0.8mm);
\filldraw[aqua] (7.5,-.5) circle (0.7mm);

\draw[black] (5,-3.3) node{$u$};
\draw[black] (5,0.3) node{$v$};

\draw[black] (6,-1.8) node{$y$};
\draw[black] (7.8,-.5) node{$x$};

\draw[black] (6.6,-1.1) node{$b$};


\draw[dashed, -{latex}, line width=0.4mm] (14,-3) to [bend left=30] (14,0);

\draw[-{latex}, pblue, line width=0.4mm] (15.5,-2.5) to [bend left=30]  (14,-3);

\draw[dashed, -{latex}, pblue, line width=0.4mm]  (14,-3) to [bend left=30] (15.5,-2.5);

\draw[-{latex}, line width=0.4mm] (15.5,-2.5) to [bend right=20] (15.5,-.5);

\draw[-{latex}, darkred, line width=0.4mm]  (14,0) to [bend left=30] (15.5,-.5);

\draw[dashed, -{latex}, darkred, line width=0.4mm] (15.5,-.5) to [bend left=30] (14,0);

\draw[line width=0.1mm] (14,-3) circle (0.8mm);
\filldraw[aqua] (14,-3) circle (0.7mm);

\draw[line width=0.1mm] (14,0) circle (0.8mm);
\filldraw[aqua] (14,0) circle (0.7mm);

\draw[line width=0.1mm] (15.5,-2.5) circle (0.8mm);
\filldraw[aqua] (15.5,-2.5) circle (0.7mm);

\draw[line width=0.1mm] (15.5,-.5) circle (0.8mm);
\filldraw[aqua] (15.5,-.5) circle (0.7mm);

\draw[black] (14,-3.3) node{$u$};
\draw[black] (14,0.3) node{$v$};

\draw[black] (15.8,-2.5) node{$x$};
\draw[black] (15.8,-.5) node{$y$};

\draw[black] (15.4,-1.5) node{$b$};

\draw[black] (14.8,.3) node{$\dpat_1$};

\draw[black] (14.8,-3.3) node{$\dpat_2$};

\end{tikzpicture} 
\label{fig:path1}
\caption{The replacements used in the proof of the first statement of Lemma~\ref{lem:arrows-covered-by-ghosts}. Ghost arrows are dashed. The path $\dpat$ which replaces the ghost arrow $uv$ in $\dz$ is depicted in blue on the left. The circuits $\dq_1$ and $\dq_2$ are depicted in red and blue on the right, respectively.}
\end{figure}
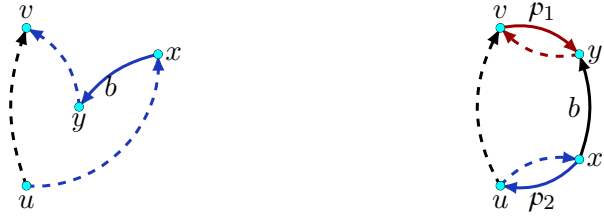

Assume now that $s(b)=0$. Since $b$ is in $\E_\pi$, it is vertical and so we have $s(\bar b)=-1$. This implies that $b$ is upward, i.e., $x\supface_\pi y$. Assume for the sake of a  contradiction that $\bar b\in\dz$. Let $\dpat_1$ and $\dpat_2$ be the paths on $\dz$ from $v$ to $y$ and from $x$ to $u$, respectively. Let $\dq_1$ be the circuit in $\ssG^{\dagger}_{\pi}$ obtained by adding to $\dpat_1$ the upward  ghost arrow $yv$ if $y\neq v$. Let $\dq_2$ be the circuit obtained by adding to $\dpat_2$ the upward ghost arrow $ux$ if $u\neq x$; see Figure~\ref{fig:path1} on the right.  Since $\ell\in\altcone_{s,\pi}$, we get the inequalities
\[
 \sum_{a \in \dpat_1}\ell_as(a) =\sum_{a \in \dq_1}\ell_as(a)  \leq 0 \qquad \textrm{and} \qquad \sum_{a \in \dpat_2}\ell_as(a)  = \sum_{a \in \dq_2}\ell_as(a) \leq 0.
\]
But then
\[
0=\sum_{a\in \dz}\ell_as(a) = \sum_{a \in \dpat_1}\ell_as(a) + \ell_{b}s(\bar b)+\sum_{a \in \dpat_2}\ell_as(a) \leq - \ell_{b} <0,
\]
leading to a contradiction. This proves the first statement.

As for the second statement, if $b$ and $\bar b$ are covered by ghosts of $\dz$, it follows from the first part of the lemma that $s(b),s(\bar b)\leq 0$. But, $s(b)+s(\bar b)\geq -1$. Thus, either $s(b)=0$ or $s(\bar b)=0$. 

If $s(b)=0$, then $\bar b\not\in\dz$ by the first statement.  If in addition $\dz$ is essential for $(s,\pi)$, then $b\not\in\dz$ by \ref{eq:ess2}, as required. 

In the other case, when $s(\bar b)=0$, we conclude in the same way.
\end{proof}

Given a circuit $\dz$ of $\ssG_\pi^{\dagger}$ which is both essential and active for $(s,\pi)$, we define a slope-level pair $(s',\pi')$ on $G$ such that $\altcone_{s',\pi'}$ is the face $\altcone_{s,\pi}\cap H_{\dz}$ of $\altcone_{s,\pi}$; see Proposition~\ref{thm:face1}.

\smallskip

The function $s' \colon \E \to \Z$ is defined by 
\[
s'(a) \coloneqq \begin{cases} s(a)+1 & \textrm{if $\bar a$ is a noninteger arrow for $s$ in $\dz$},\\
s(a)+1 & \textrm{if $a$ and $\bar a$ are covered by the same ghost of $\dz$ and $s(\bar a) = 0$},\\ 
s(a) & \textrm{otherwise}.
\end{cases}
\]
In particular, $s'(a)=s(a)$ if $a\not\in\E_\pi$. Notice as well that Lemma~\ref{lem:arrows-covered-by-ghosts} yields that the two conditions implying $s'(a)=s(a)+1$ in the above definition do not simultaneously hold. 

\smallskip

We define the ordered partition $\pi'$ of $V$ as a the following coarsening of $\pi$: For each ghost arrow $uv \in \dz$, we replace all parts of $\pi$ containing vertices $w$ with $u \supfaceq_\pi w\supfaceq_\pi v$ by their union. Since the ghost arrows in $\dz$ do not overlap by our assumption \ref{eq:ess3}, we get a well-defined ordered partition of $V$ that we denote by $\pi'$.

\smallskip

We say that the pair $(s',\pi')$ is obtained by \emph{squashing $(s,\pi)$ relative to the active essential circuit $\dz$}.

\subsection{Faces, I} Notation as in the previous section.
\begin{prop}\label{thm:face1} Let $(s,\pi)$ be a slope-level pair on $G$. Let $\dz$ be a  circuit of $\ssG_{\pi}^{\dagger}$ which is essential and active for $(s,\pi)$. Let $(s', \pi')$ be the pair obtained by squashing $(s,\pi)$ relative to $\dz$. Then, the following holds.
\begin{enumerate}[label=(\arabic*)]
    \item  \label{thm:face1-1} The pair $(s',\pi')$ is a slope-level pair on $G$. 

    \item \label{thm:face1-3} The cone $\altcone_{s',\pi'}$ is the intersection of $\altcone_{s,\pi}$ with the hyperplane $H_{\dz}$ in $\R^E$, that is,
    \[
    \altcone_{s',\pi'} = \altcone_{s,\pi} \cap \Bigl\{\ell \in \R^E \st \sum_{a\in \dz} \ell_a s(a) =0\Bigr\}.
    \]
\end{enumerate}
In particular, $\altcone_{s',\pi'}$ is a face of $\altcone_{s,\pi}$.
\end{prop}
The proof of this result is given in Sections~\ref{sec:proof-face1-1} and \ref{sec:proof-face1-2}, below. 

\subsubsection{$(s',\pi')$ is a slope-level pair.}\label{sec:proof-face1-1} We prove first that $s'$ is a slope function. Indeed, let $a\in\E$ such that $s'(a)\neq s(a)$. Then, $a\in\E_\pi$ and $s'(a)=s(a)+1$. We need only  prove that $a$ is not integer for $s$, and $s'(\bar a)=s(\bar a)$. 

By definition, $s'(a)=s(a)+1$ only in two cases: 

- First, we may have that $a$ is a noninteger arrow for $s$ and $\bar a\in\dz$. But then, $a\not\in\dz$ because $\dz$ is a circuit, and $a$ and $\bar a$ are not covered by ghosts of $\dz$ by  Lemma~\ref{lem:arrows-covered-by-ghosts}. Thus, 
$s'(\bar a)=s(\bar a)$ by definition. 

- Second, we may have that $a$ and $\bar a$ are covered by a ghost of $\dz$ and $s(\bar a)=0$. Then, $a,\bar a\not\in\dz$ by Lemma~\ref{lem:arrows-covered-by-ghosts}. Also, $s(a)=-1$ because $a\in\E_\pi$, and hence $a$ is not integer for $s$, and $s'(\bar a)=s(\bar a)$ by definition.

This proves that $s'$ is a slope function. 

It remains to prove that $(s',\pi')$ is a slope-level pair. For this, we must prove that  an arrow $a\in\E$ is upward for $s'$ if and only if it is upward for $\pi'$. In other words, for each arrow $a=xy$ of $G$, we must prove that $s'(\bar a)<0$ if and only if $x\supface_{\pi'}y$.

\smallskip

First, if $s'(\bar a)<0$, then also $s(\bar a)<0$, and hence $x\supface_{\pi}y$. Then, $x\supfaceq_{\pi'}y$, with equality if and only if $a$ and $\bar a$ are covered by the same ghost of $\dz$. Thus, if equality were the case, then Lemma~~\ref{lem:arrows-covered-by-ghosts} would yield $s(a),s(\bar a)\leq 0$. Since $s(\bar a)<0$, we would get $s(a)=0$ and $s(\bar a)=-1$. But then, we would obtain $s'(\bar a)=s(\bar a)+1=0$ by definition, leading to a contradiction. We conclude that $x\supface_{\pi'}y$, as required.

\smallskip

Conversely, we suppose $x\supface_{\pi'}y$, and show that $s'(\bar a)<0$. From the construction of $\pi'$, we get $x\supface_{\pi}y$, and hence $s(\bar a)<0$. Also, $a$ and $\bar a$ are not covered by the same ghost of $\dz$. Now, since $s(\bar a)<0$, we must have $s'(\bar a)\leq 0$. If we had equality, $s'(\bar a)=0$, then we would get $s(\bar a)=-1$, and, by definition of $(s',\pi')$, the arrow $a$ would be noninteger for $s$ and in $\dz$. But then, since $s$ is a slope function, we would get $s(a)=0$, and $\dz$ would not verify \ref{eq:ess2}. We conclude that $s'(\bar a)<0$, as required.

\subsubsection{$\altcone_{s',\pi'}=\altcone_{s,\pi}\cap H_{\dz}$.} \label{sec:proof-face1-2} We  prove now the second statement of the proposition. 

For each circuit $\dq$ in $\ssG_{\pi}^{\dagger}$, let $\dq'$ denote the contraction of $\dq$ in $\ssG_{\pi'}^{\dagger}$, defined as follows. A ghost arrow $a=uv$ in $\dq$ remains in $\dq'$ if the two vertices $u,v$ remain in different parts of $\pi'$; otherwise, it is removed. Moreover, each arrow $a\in \E_\pi$ which is in $\dq$ remains an arrow in $\dq'$ provided $a$ and $\bar a$ are not covered by the same ghost of $\dz$; otherwise, $a$ is removed. In the latter case, we have $s(a)\leq 0$ by Lemma~\ref{lem:arrows-covered-by-ghosts}. In any case, for each $\ell\in\R^E_{>0}$, we get 
\begin{equation}\label{eq:dq<dq'}
\sum_{a\in \dq} \ell_a s(a)\leq \sum_{a \in \dq'} \ell_a s'(a).
\end{equation}

If $\dq=\dz$, then all ghost arrows in $\dz$ are removed. By Lemma~\ref{lem:arrows-covered-by-ghosts}, every nonghost arrow $a$ in $\dz$ remains an arrow in $\dz'$. In addition, since $\bar a\not\in\dz$, because $\dz$ is a circuit, $s'(a)=s(a)$ by definition. In this case, we have equality in \eqref{eq:dq<dq'} for each $\ell\in\R^E_{>0}$:
\begin{equation}\label{eq:dz=dz'}
\sum_{a\in \dz}\ell_as(a) = \sum_{a\in \dz'}\ell_as'(a).
\end{equation}
Furthermore, whether a nonghost $a\in\dz$ is integer for $s$ or not, $a$ will be integer for $s'$. Thus, $\dz'$ will be in $\ssG_{s',\pi'}^{\rmint}$.

\smallskip

We prove first that $\altcone_{s',\pi'} \subseteq \altcone_{s,\pi}\cap H_{\dz}$. Let $\ell\in\altcone_{s',\pi'}$ and $\dq$ be a circuit in $\ssG_{\pi}^{\dagger}$. Since $\ell\in\altcone_{s',\pi'}$ it follows from \eqref{eq:dq<dq'} that $\sum_{a\in \dq} \ell_a s(a)\leq 0$. Thus, 
$\ell\in\altcone_{s,\pi}$. In addition, since $\dz'$ is in $\ssG_{s',\pi'}^{\rmint}$, it follows from \eqref{eq:dz=dz'} that 
$\sum_{a\in \dz}\ell_as(a)=0$, that is, $\ell\in H_{\dz}$. This proves the inclusion $\altcone_{s',\pi'}\subseteq\altcone_{s,\pi}\cap H_{\dz}$.

\smallskip

We prove now the reverse inclusion. Let $\ell \in \altcone_{s,\pi}\cap H_{\dz}$.  We need only prove that 
\begin{equation}\label{eq:dq'<0}
\sum_{a\in\dq'}\ell_as'(a)\leq0
\end{equation}
for every circuit $\dq'$ in $\ssG^{\dagger}_{\pi'}$ that verifies \ref{eq:ess1} in $\ssG^{\dagger}_{\pi'}$.

Since $\ell\in H_{\dz}$, it follows from \eqref{eq:dz=dz'} that
\begin{equation}\label{eq:dz'=0}
\sum_{a\in \dz'}\ell_as'(a)=0.
\end{equation}
Viewing each $\dq'$ and $\dz'$ as cycles, we may write the sum $\dq'+\dz'$ as a sum of pairwise compatible circuits $\dq'_1,\dots,\dq'_s$ in $\ssG^{\dagger}_{\pi'}$. No new ghost arrows are introduced, so all the $\dq'_i$ verify \ref{eq:ess1}. If we prove \eqref{eq:dq'<0} for each such $\dq'_i$ instead of $\dq'$, we will have it for $\dq'$ as well, because of \eqref{eq:dz'=0}. We may thus assume that $\bar a\not\in\dz'$ for each $a\in\dq'$. 

\smallskip

The circuit $\dq'$ is a contraction of a circuit $\dq$ in $\ssG^{\dagger}_{\pi}$: We keep all the nonghost arrows of $\dq'$, and add ghost arrows which are either contracted to the ghost arrows of $\dq$ or removed. The nonghost arrows $a$ in $\dq$ are not contracted and satisfy $\bar a\not\in\dz$; thus $s'(a)=s(a)$ by definition. Let $N$ be the number of downward ghost arrows of $\dq$. These are removed when contracting to $\dq'$ as $\dq'$ satisfies \ref{eq:ess1}.

If $N=0$, then equality holds in \eqref{eq:dq<dq'}. Since $\ell\in\altcone_{s,\pi}$, we get \eqref{eq:dq'<0}, as required. Assume $N>0$. We will proceed by descending induction on $N$. Let $uv$ be a downward ghost arrow in $\dq$. Since it does not appear in $\dq'$, and $\dz$ satisfies \ref{eq:ess1}, there is a ghost arrow $xy$ in $\dz$ such that $x\subfaceq_{\pi}v\subface_{\pi}u\subfaceq_{\pi}y$. Consider the closed path $\dpat$ in $\ssG^{\dagger}_{\pi}$ made up of all the arrows of $\dq$ except $uv$, all the arrows of $\dz$ except $xy$, and the upward ghost arrows $xv$ if $x\neq v$ and $uy$ if $u\neq y$. Its contraction $\dpat'$ to $\ssG^{\dagger}_{\pi'}$, viewed as a cycle, is $\dq'+\dz'$. Since $\dz$ has only  upward ghost arrows, $\dpat$ has less downward ghost arrows than $\dq$. Writing $\dpat$ as a sum of pairwise compatible circuits $\dq_1,\dots,\dq_s$ in $\ssG^{\dagger}_{\pi}$, as before, we need only prove \eqref{eq:dq'<0} for the contractions $\dq'_i$ of the $\dq_i$ instead of $\dq'$. Since each $\dq_i$ has less downward arrows that $\dq$, we apply induction to finish the proof. 
\begin{proof}[Proof of Proposition~\ref{thm:face1}] We have established the two statements in the proposition.
\end{proof}

\subsection{Faces, II} Here is our main theorem on the faces of the cone of edge lengths.
\begin{thm}\label{thm:face-cone} Let $(s,\pi)$ be a slope-level pair on $G$ and $F$ a face of $\altcone_{s,\pi}$ intersecting $\R^E_{>0}$. Then, there is a slope-level pair $(s',\pi')$ on $G$ such that  $F=\altcone_{s',\pi'}=\cone_{s',\pi'}$, and the following properties are verified:
\begin{enumerate}[label=(F\arabic*)]
    \item\label{item:coarse} $\pi'$ is a coarsening of $\pi$,
    \item\label{item:s<s'} $s\leq s'$, and
    \item\label{item:intE} $\ssE_{s'}^\rmint\supseteq\ssE_{s}^\rmint$.
 \end{enumerate}
In particular, if $\intcone_{s,\pi}$ is nonempty, then each face of $\cone_{s,\pi}$ intersecting $\R^E_{>0}$ is $\cone_{s',\pi'}$ for some slope-level pair $(s',\pi')$ on $G$ satisfying \ref{item:coarse}-\ref{item:s<s'}-\ref{item:intE}.
\end{thm}

\begin{proof} If $F$ intersects $\intcone_{s,\pi}$, then  by Theorem~\ref{thm:characterization-slope-level-cone}, $\altcone_{s,\pi}=\cone_{s,\pi}$ and $\intcone_{s,\pi}$ is the interior of $\altcone_{s,\pi}$. Since 
$F$ is a face of $\altcone_{s,\pi}$, we must have $F=\altcone_{s,\pi}$, and hence we may put $(s',\pi')\coloneqq (s,\pi)$. 

Assume now that $F$ does not intersect $\intcone_{s,\pi}$. Then, there is an essential circuit $\dz$ in $\ssG_\pi^\dagger$ which is not in $\ssG^\rmint_{s,\pi}$ such that $F\subseteq H_{\dz}$. Since $F$ contains $\ell\in\R^E_{>0}$, the circuit $\dz$ is active for $(s,\pi)$. Then, Proposition~\ref{thm:face1} yields that $F$ is a face of $\altcone_{s',\pi'}$ for the pair $(s',\pi')$ associated to $\dz$ by the squashing process.

By definition of the pair $(s',\pi')$, properties \ref{item:coarse} and \ref{item:s<s'} hold. As for \ref{item:intE}, an arrow $a\in\E_{\pi}$ which is integer for $s$ is present in $\ssG_{\pi'}$. Indeed, if not, then $a$ and $\bar a$ would be covered by the same ghost of $\dz$, and hence $s(a),s(\bar a)\leq 0$ by Lemma~\ref{lem:arrows-covered-by-ghosts}. Since $a$ is integer for $s$, we would get $s(a)=s(\bar a)=0$, that is, $a$ would be  horizontal for $s$, contradicting $a\in\E_{\pi}$. Also, since 
\[
0=s(a)+s(\bar a)\leq s'(a)+s'(\bar a)\leq 0,
\]
we must have $a \in \E_{s'}^\rmint$, showing \ref{item:intE}.

\smallskip

Now, if $F$ intersects $\intcone_{s',\pi'}$, as before we get 
$F=\altcone_{s',\pi'}=\cone_{s',\pi'}$. Otherwise, we proceed as above, getting a slope-level pair $(s'',\pi'')$ such that $F$ is a face of $\altcone_{s'',\pi''}$, and such that $\pi''$ is a coarsening of $\pi'$, $s'\leq s''$, and $\ssE_{s''}^\rmint\supseteq\ssE_{s'}^\rmint$. Clearly, $\pi''$ is a coarsening of $\pi$, and we have $s\leq s''$ and $\ssE_{s''}^\rmint\supseteq\ssE_{s}^\rmint$. And, we restart this cycle of reasoning. This process will eventually end, and we obtain the theorem.
\end{proof}

We obtain the following refined characterization of the facets.

\begin{thm}[Characterization of facets]\label{thm:facet}
Let $(s,\pi)$ be a slope-level pair on $G$ with nonempty $\cone_{s,\pi}$. Let $F$ be a facet of $\cone_{s,\pi}$ meeting $\R_{>0}^E$. Then, there exists an essential circuit $\dz$ in $\ssG_{\pi}^\dagger$ which is active for $(s,\pi)$ such that, denoting by $(s',\pi')$ the slope-level pair associated to $\dz$ by squashing, we have 
\begin{itemize}
\item $F = \cone_{s',\pi'}$, and 
\item $g(\ssG^{\rmint}_{s',\pi'}) = g(\ssG^{\rmint}_{s,\pi})+1$.
\end{itemize}
\end{thm}

Theorem~\ref{thm:face-cone} implies that there is a sequence of circuits whose consecutive squashings lead to a slope-level pair that defines the facet. The content of the above theorem is that this can be done by a single squashing. 

\begin{proof} Since $\cone_{s,\pi}=\altcone_{s,\pi}$ by Theorem~\ref{thm:characterization-slope-level-cone}, it follows from Proposition~\ref{prop:essential-enough} that $\cone_{s,\pi}$ is given by the inequalities \eqref{eq:slope-level-compatible} for essential circuits. Thus $F = \cone_{s,\sigma} \cap H_{\dz}$ for an essential circuit $\dz$ in $\ssG_{\pi}^\dagger$ not in $\ssG_{s,\pi}^{\rmint}$. Since there is $\ell\in F$ with $\ell>0$, it follows that $\dz$ is active for $(s,\pi)$. Proposition~\ref{thm:face1} then implies that $F=\altcone_{s',\pi'}$ for the slope-level pair $(s',\pi')$ obtained by squashing $(s,\pi)$ relative to $\dz$.

We claim that $g(\ssG^{\rmint}_{s',\pi'}) > g(\ssG^{\rmint}_{s,\pi})$. Indeed, $\dz$ is not in $\ssG^{\rmint}_{s,\pi}$, but the contracted circuit $\dz'$ in 
$\ssG^\dagger_{\pi'}$ is in $\ssG^{\rmint}_{s', \pi'}$, as we observed in the proof of Proposition~\ref{thm:face1}. Since $F = \cone_{s,\sigma} \cap H_{\dz}$ and $F$ is a facet of $\cone_{s,\sigma}$, it follows that $\dz'$ is not in the cycle space of $\ssG^{\rmint}_{s', \pi'}$ generated by the circuits contracted from $\ssG^{\rmint}_{s, \pi}$. As no arrow of any circuit in $\ssG^{\rmint}_{s, \pi}$ is actually contracted, that is, $\ssE_{s}^\rmint\subseteq\ssE_{s'}^\rmint$, we get our claim. 

We show now that the cone $\cone_{s', \pi'}$ in not empty. Assume by contradiction that $\cone_{s', \pi'}=\emptyset$. Then, the proof of Theorem~\ref{thm:face-cone} yields that $F\subseteq\altcone_{s'',\pi''}$ for some slope-level pair $(s'',\pi'')$ derived by squashing $(s',\pi')$ relative to a circuit $\dz'$ in $\ssG_{\pi'}^\dagger$ which is essential and active for $(s',\pi')$. As proved above, $g(\ssG^{\rmint}_{s'',\pi''}) > g(\ssG^{\rmint}_{s',\pi'})$. Now, Proposition~\ref{prop:dim_count} implies that
\[
\dim F\leq |E|-g(\ssG_{s'',\pi''}^{\rmint})\quad\text{and}\quad
\dim\cone_{s,\pi}=|E|-g(\ssG_{s,\pi}^{\rmint}).
\]
But, $g(\ssG_{s'',\pi''}^{\rmint})>g(\ssG_{s',\pi'}^{\rmint})>g(\ssG_{s,\pi}^{\rmint})$, contradicting the fact that $F$ has codimension one in $\cone_{s,\pi}$. 

We conclude that $\cone_{s', \pi'}\neq\emptyset$, whence $F=\cone_{s',\pi'}$ by Theorem~\ref{thm:characterization-slope-level-cone}. Applying Proposition~\ref{prop:dim_count}, we get
\[
1=\codim(F,\cone_{s,\pi})=\dim\cone_{s,\pi}-\dim\cone_{s', \pi'}
=g(\ssG_{s',\pi'}^{\rmint})-g(\ssG_{s,\pi}^{\rmint})\geq 1,
\]
from which we deduce that $g(\ssG^{\rmint}_{s',\pi'}) = g(\ssG^{\rmint}_{s,\pi})+1$, as required.
\end{proof}

\subsection{Fans associated to bricks}\label{sec:bricks-fans} 

Let $g\in\Z_{\geq 0}$. Consider the standard simplex $\Delta!_g$ in $\R^V$ of width $g$. Each $\beta \in \Delta!_{g}$ defines the polytope $\B_\beta$ consisting of the set of all points $q \in \Delta!_g$ that verify the inequalities 
\[
\lfloor \beta(S)\rfloor \leq q(S) \leq \lceil \beta(S)\rceil\qquad \qquad \forall \,\, S \subseteq V.
\]
We refer to $\B_\beta$ as the \emph{brick} defined by $\beta$; see \cite[\S 10.2]{AE-modular-polytopes} and Appendix~\ref{sec:torusgrass}.

By \cite[Thm.~10.4]{AE-modular-polytopes}, the collection of bricks $\B_\beta$ for $\beta \in \Delta!_{g}$ is a tiling of $\Delta!_{g}$. Furthermore, the base polytope $\P$ of any polymatroid with integral vertices in $\Delta!_{g}$ is tiled by the set of bricks included in $\P$. Denote by $\Br_g$ the collection 
\[
\Br_g \coloneqq \Bigl\{\,\B_\beta \st \beta \in \Delta!_g \quad \text{and} \quad \B_\beta  \text{ is full-dimensional}\,\Bigr\}.
\]
 
For each slope-level pair $(s,\pi)$ on $G$, recall the submodular function $\gamma!_{\pi}$ associated to the $\pi$-residue space $\GlobSp_\pi\subset \bV$; see Section~\ref{sec:residue-spaces-and-conditions}. Recall as well the submodular function $\zeta!_s$ associated to $s$ by Definition~\ref{defi:zetas}. Whereas $\gamma!_{\pi}$ has range $g(G)$, the function $\zeta!_s$ has range 0. Assume $g\geq g(G)$, and fix a \emph{genus function} $\g\colon V \to \Z_{\geq 0}$ of range $g-g(G)$.

Consider the submodular function 
\[
\eta!_{s,\pi}\coloneqq\eta!_{s,\pi}(\g)\coloneqq\gamma!_\pi+\g+ \zeta!_s,
\]
and let $\P_{s,\pi}$ denote the base polytope associated to the submodular function $\upmin(\eta!_{s,\pi})$.

\begin{defi}[Permissible slope-level pairs]\label{defi:permissible}
A slope-level pair $(s,\pi)$ is \emph{permissible} if 
\begin{itemize}
\item the set $\intcone_{s,\pi} \subset \R_{>0}^E$ is nonempty, and
\item the  submodular function $\eta!_{s,\pi}$ is positive and simple.
\end{itemize}
Denote by $\PSL=\PSL(G)$ the collection of all permissible slope-level pairs on $G$. 
\end{defi}

Since there are finitely many ordered partitions $\pi$, and the slopes taken by $s$ with $\eta!_{s,\pi}$ positive and simple are all bounded by a constant by Lemma~\ref{lem:bound-s}, the set $\PSL$ is finite. Since each $\P_{s,\pi}$ is full-dimensional, we can write 
\[
\PSL=\bigcup_{\B\in\Br_g}\PSL(\B) \quad \text{where} \quad \PSL(\B)\coloneqq\bigl\{(s,\pi)\in\PSL\,\st\, \P_{s,\pi}\supseteq\B\bigr\}.
\]

\begin{thm}\label{thm:fan-sigmaB} Let $\B \in \Br_g$ be a full-dimensional brick in $\Delta!_g$. The collection of cones $\cone_{s,\pi}$ for $(s,\pi) \in \PSL(\B)$ and all of their faces which lie on the boundary $\R_{\geq 0 }^E \setminus \R_{>0}^E$ form a rational fan with support $\R_{\geq 0}^E$.
\end{thm}
 
\begin{defi}\label{def:sigmaB} For each full-dimensional brick $\B \in \Br_g$, we denote by $\ssSigma_{\B}$ the fan produced by the above theorem.
\end{defi} 
We need the following lemma.
\begin{lemma}\label{lem:PSLB} Let $\B \in \Br_g$ be a full-dimensional brick in $\Delta!_g$. For each $\ell\in\R^E_{>0}$ there is a unique $(s,\pi)\in\PSL(\B)$ such that $\ell\in\intcone_{s,\pi}$. In particular, the cones $\intcone_{s,\pi}$ for  $(s,\pi) \in \PSL(\B)$ are pairwise disjoint.
\end{lemma}

\begin{proof} We derive this as a consequence of our tiling theorem,  Theorem~\ref{thm:W-tiling}. Thus, let $X$ be a nodal connected reduced curve proper over an algebraically closed field $\k$ whose associated dual graph is $G=(V,E)$. Assume the branches over its nodes are in general position on the normalizations of its components. 

Let $\ell\in\R^E_{>0}$ and denote by $\Gamma$ the metric graph associated to $(G,\ell)$. There exists a complete algebraically closed non-Archimedean field $\K$ with residue field $\k$ and value group $\valgroup=\R$, and a smooth proper curve $\varX$ over $\K$ whose tropicalization is $(X,\Gamma)$, see~\cite[Thm.~3.24]{ABBR}.  For each function $h\colon V \to \R$, let $\ssW_h$ be the associated reduction of the space of Abelian differentials on $\varX$. By Theorem~\ref{thm:W-tiling}, for each $\B\in\Br_g$, there is a function $h\colon V \to \R$ such that the base polytope of $\ssW_h$ contains $\B$. Then, $\ssW_h$ is simple, and thus $\eta!_{s,\pi}\geq 0$ by Proposition~\ref{prop:step2singnot}, where $(s,\pi)\coloneqq (\slztwist{\ell}h,\sspi_h)$. Then, Theorem~\ref{step2sing} yields that the submodular function associated to $\ssW_h$ is $\upmin(\eta!_{s,\pi})$, and hence $\P_{s,\pi} \supseteq \B$. Since $\ssW_h$ is simple, so is $\eta!_{s,\pi}$ by Proposition~\ref{prop:simpleness-min-operation}. Therefore, $(s,\pi)\in\PSL(\B)$. Since $(s,\pi)= (\slztwist{\ell}h,\sspi_h)$, we have  $\ell\in\intcone_{s,\pi}$ by definition, as required. This proves the existence.

To prove the uniqueness, note that if $\ell\in\intcone_{s',\pi'}$ for another slope-level pair $(s',\pi')\in\PSL(\B)$, then there is a function $h'\colon V\to\R$ with $(s',\pi')= (\slztwist{\ell}h',\sspi_{h'})$ and $\P_{s',\pi'} \supseteq \B$. But then the base polytope of $\ssW_{h'}$ contains $\B$ as well. By Theorem~\ref{thm:W-tiling}, the base polytopes of $\ssW_h$ and $\ssW_{h'}$ will then coincide, and hence $h'-h$ is constant by \cite[Prop.~8.9]{AE-modular-polytopes}. It follows that $(s,\pi)=(s',\pi')$.
\end{proof}
 
\begin{proof}[Proof of Theorem~\ref{thm:fan-sigmaB}] Lemma~\ref{lem:PSLB} yields that 
\[
\bigcup_{(s,\pi) \in \PSL(\B)} \intcone_{s,\pi}= \R_{>0}^E,
\]
and that the union is disjoint.

Let $(s,\pi)\in\PSL(\B)$. Then $\cone_{s,\pi}\neq\emptyset$. By Theorem~\ref{thm:face-cone}, a face $F$ of $\cone_{s,\pi}$ meeting $\R_{>0}^E$ is equal to $\cone_{s',\pi'}$ for a slope-level pair $(s', \pi')$ satisfying Conditions~\ref{item:coarse}-\ref{item:s<s'}-\ref{item:intE} therein. Now, since $\pi'$ is a coarsening of $\pi$, it follows from Proposition~\ref{prop:increasing-gamma} that $\gamma!_{\pi'}\geq\gamma!_{\pi}$. Also, since $s'\geq s$, by Proposition~\ref{prop:monotonicity-zeta}, we have $\zeta!_{s'}\geq\zeta!_{s}$. We infer that $\eta!_{s',\pi'}\geq\eta!_{s,\pi}$. It follows that  the polytope $\P_{s',\pi'}$ contains the polytope $\P_{s,\pi}$, and hence contains $\B$. We conclude that $(s',\pi') \in \PSL(\B)$. This shows that, for each $(s,\pi)\in\PSL(\B)$, each face of the cone $\cone_{s,\pi}$  that meets $\R^E_{>0}$ is equal to $\cone_{s',\pi'}$ for a certain $(s',\pi') \in \PSL(\B)$.

Furthermore, let $(s,\pi), (s',\pi')\in\PSL(\B)$ be such that the intersection
$Z\coloneqq\cone_{s,\pi}\cap\cone_{s',\pi'}\cap\R^E_{>0}$ is nonempty. For each $\ell\in Z$, let $F$ and $F'$ be the faces of $\cone_{s,\pi}$ and $\cone_{s',\pi'}$ that contain $\ell$ in their relative interiors, respectively. As we showed above, $F=\cone_{s'',\pi''}$ and $F'=\cone_{s''',\pi'''}$ for certain $(s'',\pi''),(s''',\pi''')\in \PSL(\B)$. But then Lemma~\ref{lem:PSLB} yields $(s'',\pi'')=(s''',\pi''')$, whence $F=F'$. It follows that $Z$ is contained in the union of the common faces of $\cone_{s,\pi}$ and $\cone_{s',\pi'}$, whence contained in a common face $F$, since  $Z$ is convex. As seen above, $F=\cone_{s'',\pi''}$ for some $(s'',\pi'')\in \PSL(\B)$. Since $Z=F\cap\R^E_{>0}$, it follows that $\cone_{s,\pi}\cap\cone_{s',\pi'}=\cone_{s'',\pi''}$.

The upshot is that the collection $\mathscr C$ of cones $\cone_{s,\pi}$ for $(s,\pi) \in \PSL(\B)$ verifies the following three properties:
\begin{itemize}
    \item If $\sigma$ is in $\mathscr C$ and $F$ is a face of $\sigma$ that meets $\R_{>0}^E$, then $F$ belongs to $\mathscr C$.
    \item If $\sigma$ and $\tau$ are in $\mathscr C$ and the intersection $\sigma\cap \tau \cap \R_{>0}^E$ is nonempty, then $\sigma\cap \tau$ is a common face of $\sigma$ and $\tau$, and belongs to $\mathscr C$. 
    \item The collection $\mathscr C$ covers $\R_{\geq 0}^E$.
\end{itemize}
It is easy to conclude now. 
\end{proof}

\subsection{Permissible collections of slope-level pairs}
Given a subset $S\subseteq \PSL$, we define 
\[
\cone_S \coloneqq \bigcap_{(s, \pi)\in S} \cone_{s,\pi}. 
\]
For $S=\{(s,\pi)\}$, a singleton, we obtain $\cone_S = \cone_{s,\pi}$. The following proposition is immediate from Theorem~\ref{thm:characterization-slope-level-cone}.

\begin{prop} For each subset $S\subseteq \PSL$, the set $\cone_S$ is a rational polyhedral cone in $\R_{\geq 0}^E$ provided that it is nonempty. 
\end{prop}

\subsection{The meet of the fans associated to all full-dimensional bricks}\label{sec:fan-structure-fc} For each edge length function $\ell\colon E \to\R_{>0}$, let $\PS_\ell \subseteq \PSL$ be defined by
\[
\PS_\ell \coloneqq \Bigl\{(s,\pi) \in \PSL \, \st \, \ell \in \intcone_{s,\pi} \Bigr\} =  \Bigl\{(\slztwist{\ell}h, \pi_h)  \,\st\, \eta!_{\slztwist{\ell}h,\sspi_h}\text{ is positive and simple}\Bigr\}.
\] 
Denote by 
\[
\slsfunc \colon \R_{>0}^E \longrightarrow \ssub{2}!^{\PSL}, \qquad \qquad 
\ell \mapsto \PS_\ell
\]
the corresponding map. Notice that $\cone_{\PS_\ell}$ is nonempty, as it contains $\ell$. Consider the collection of all the rational cones $\cone_{\PS_\ell}$ for $\ell\in\R^E_{>0}$, and denote by $\Sigma=\Sigma(G,\g)$ the collection of these cones and all of their faces on the boundary $\R_{\geq 0}^E \setminus\R_{>0}^E$.

For fans $\ssSigma_1,\dots,\ssSigma_m$ with the same support in $\R^V$, we define their \emph{meet}, denoted $\bigwedge\ssSigma_i$, to be the fan obtained by taking all the nonempty intersections $\cone_1 \cap\cdots\cap\cone_m$ for cones $\cone_i \in \ssSigma_i$ for $i=1,\dots,m$. Note that the meet has the same support as the $\ssSigma_i$.

\begin{thm}\label{thm:Sigma} Notation as above, $\Sigma$ is a rational fan with support $\R_{\geq 0}^E$. It is given by
\begin{equation}\label{eq:sigma-sigmaB}
\Sigma =  \Bigl\{\bigcap_{\B \in \Br_g} \sssigma_\B \, \st\, \sssigma_\B \text{ in } \ssSigma_\B \text{ for } \B\in \Br_g\Bigr\}.
\end{equation}
In other words, $\Sigma=\bigwedge_{\B \in \Br_g} \ssSigma_{\B}$ is the meet of the fans $\ssSigma_{\B}$ associated to the full-dimensional bricks $\B$ in $\Delta!_g$.
\end{thm}  

\begin{proof} In order to prove the theorem, it is enough to show the equality \eqref{eq:sigma-sigmaB}, as its right-hand side is a fan. It will be thus enough to show that the collection of cones in $\Sigma$ meeting $\R^E_{>0}$ is the collection of cones in the meet $\bigwedge\ssSigma_{\B}$ meeting $\R^E_{>0}$.
  
A cone in $\Sigma$ meeting $\R^E_{>0}$ is $\cone_{\PS_\ell} = \bigcap_{(s,\pi) \in \PS_{\ell}} \cone_{s,\pi}$ for some $\ell \in \R_{>0}^E$. 
Each cone $\cone_{s,\pi}$ above belongs to $\ssSigma_{\B}$ for each brick $\B \in \Br_g$ contained in the base polytope $\P_{s,\pi}$. Furthermore, Lemma~\ref{lem:PSLB} yields that for each brick $\B\in\Br_g$, there is $(s,\pi) \in \PS_{\ell}$ such that $\P_{s,\pi}\supseteq\B$. Thus the cone $\cone_{\PS_\ell}$ belongs to $\bigwedge\ssSigma_{\B}$.

Conversely, consider a cone $\cone \in \bigwedge \ssSigma_{\B}$ meeting $\R_{>0}^E$. Then, $\cone = \bigcap_{\B \in \Br_g} \cone_{\sss_{\B}, \sspi_{\B}}$ for certain $(\sss_{\B}, \sspi_{\B})\in \PSL(\B)$. Let $\ell \in \intcone$ be in the relative interior of $\cone$. Then the collection $\PS_\ell$ coincides with the set of all $(\sss_{\B}, \sspi_{\B})$. Indeed, for each $\B\in\Br_g$, we have $\intcone\subseteq \intcone_{\sss_{\B},\pi_{\B}}$, and hence $(\sss_{\B}, \sspi_{\B})\in\PS_\ell$. On the other hand, given $(s,\pi)\in\PS_\ell$, there is $\B\in\Br_g$ such that $\P_{s,\pi}\supseteq\B$. But then $(s,\pi)\in\PSL(\B)$, and hence  $(\sss_{\B},\sspi_{\B})=(s,\pi)$ by Lemma~\ref{lem:PSLB}. It follows that $\cone = \bigcap_{\B \in \Br_g} \cone_{\sss_{\B}, \sspi_{\B}} = \cone_{\PS_\ell}$, 
finishing the proof of the theorem.
\end{proof}

\begin{figure}
  \centering
 \scalebox{.3}{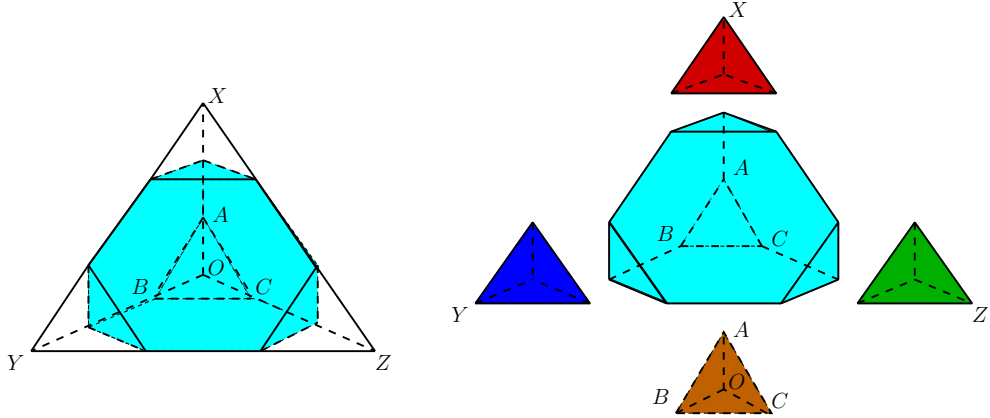}
  \caption{The polytope $\P_{\sss_0, \sspi_0}$ is the one in azure, in the simplex $\Delta!_3=OXYZ$ of width three in $\R^4$, with vertices $O=(3,0,0,0)$, $X=(0,3,0,0)$, $Y=(0,0,3,0)$ and $Z=(0,0,0,3)$. The figure shows the projection to $\R^3$ given by the last three coordinates. The polytopes in brown, red, blue, and green are the four bricks which contain each one of the four vertices of $\Delta!_3$.}
  \label{fig:k4}
\end{figure}

\begin{example}\label{ex:K4} Let $\ssK_4=(V,E)$ be the complete graph on $4$ vertices with zero genus function, $\g=0$. Thus, as in Example \ref{ex:4lines-polytope}, we have $V=\{\ssu_0,\ssu_1,\ssu_2,\ssu_3\}$  and $E=\{\sse_{0,1},\sse_{0,2},\sse_{0,3},\sse_{1,2},\sse_{1,3},\sse_{2,3}\}$, where $\sse_{i,j}$ is the unique edge connecting $\ssu_i$ and $\ssu_j$ for each $i,j$. We will simplify by putting $\ell_{i,j}\coloneqq\ell_{\sse_{i,j}}$ for each $i,j$. We let $\ssa_{i,j}$ be the unique arrow from $u_i$ to $u_j$ for all $i,j$; thus $\overbar{\ssa}_{i,j}=\ssa_{j,i}$, and $\ssa_{i,j}$ and $\ssa_{j,i}$ are the arrows on $\sse_{i,j}$ for all $i,j$ such that $0\leq i<j\leq 3$.

Let $\k$ be a field. The first slope-level pair $(\sss_0,\sspi_0)$ we consider is the trivial one: $\sss_0=0$ and $\sspi_0\coloneqq\{V\}$. Its associated cone is the full $\R^E_{\geq 0}$. The 
slope function $\sszeta_{\sss_0}$ is zero, whence $\eta!_{\sss_0,\sspi_0}=\gamma!_{\sspi_0}$. The residue space $\GlobSp_{\sspi_0}\subseteq\k^{\E}$ is the one given by local and Rosenlicht residue conditions and corresponds to the first homology group of $K_4$ with $\k$-coefficients; there are no vanishing or global residue conditions. We have $\gamma!_{\sspi_0}(\{u_i\})=2$ and $\gamma!_{\sspi_0}(\{u_i,u_j\})=3$ for all distinct $i,j$. Clearly, $\upmin(\eta!_{\sss_0,\sspi_0})=\eta!_{\sss_0,\sspi_0}$, and thus $\P_{\sss_0,\sspi_0}\subseteq\Delta!_3$ is the polytope that contains all bricks except those 4 which contain each a vertex of $\Delta!_3$. It is depicted in azure in Figure~\ref{fig:k4}.

For each $i=0,1,2,3$, let $\B_i$ be the brick of $\Delta!_3$ which contains the point with $i$-th coordinate equal to $3$ (and vanishing remaining coordinates). These are depicted in brown, red, blue, and green, respectively, in Figure~\ref{fig:k4}.  From what we discussed above, $\PSL(\B)=\{(\sss_0,\pi_0)\}$ for each $\B$ distinct from $\B_0,\B_1, \B_2,\B_3$. Also, for each such $\B$, the fan $\ssSigma_\B$ is that defined by $\R^E_{\geq 0}$ and its faces.  

The second slope-level pair $(\sss,\sspi)$ we consider satisfies the following: 
\begin{align*}
\ssV_{\pi}=&\bigl\{\{\ssu_0\},\{\ssu_1,\ssu_2,\ssu_3\}\bigr\},\quad\text{with}\quad \{\ssu_0\}\supface_{\pi}\{\ssu_1,\ssu_2,\ssu_3\}; \quad \text{and}\\
&\sss(\ssa_{i,j})=0 \text{ if $i,j\in\{1,2,3\}$ or $i=0$ and $j=2,3$},\quad\text{and}\\ &\sss(\ssa_{1,0})=\sss(\ssa_{2,0})=\sss(\ssa_{3,0})=-\sss(\ssa_{0,1})=-1.
\end{align*}
Then, $\ssA_{\sss}=\{\sse_{0,1},\sse_{0,2},\sse_{0,3}\}$ and $\ssA_{\sss}^{\rmint}=\{\sse_{0,1}\}$. 
We denote the associated cone by $\sssigma^{0}_{1}\subseteq\R^E_{\geq 0}$; it is given by two inequalities $\ell_{0,1}\leq\ell_{0,2}$ and $\ell_{0,1}\leq\ell_{0,3}$. The associated submodular functions $\gamma!_{\pi}$ and $\sszeta_{\sss}$ satisfy
\begin{align*}
\gamma!_{\pi}(\{\ssu_1,\ssu_2,\ssu_3\})=1,\qquad  \gamma!_{\pi}(\{\ssu_0\})=2, \qquad \gamma!_{\pi}(\{\ssu_0,\ssu_i\})=3 \text{ for $i=1,2,3$},
\end{align*}
from which the other values of $\gamma!_{\pi}$ can be derived, and
\begin{align*}
\zeta!_{\sss}(\{\ssu_0\})=1,\qquad \zeta!_{\sss}(I)=0 \text{ for each $I\subseteq V$ with $I\neq\{\ssu_0\}$}.
\end{align*}
Then, $\eta!_{\sss,\pi}=\gamma!_{\pi}+\zeta!_{\sss}$ and $\upmin(\eta!_{\sss,\pi})$ can be derived; in fact, one checks that  $\eta!_{\sss,\pi}=\upmin(\eta!_{\sss,\pi})$ and 
\begin{align*}
\eta!_{\sss,\pi}(\{\ssu_0\})=3,\qquad \eta!_{\sss,\pi}(\{\ssu_{i}\})=\eta!_{\sss,\pi}(\{\ssu_1,\ssu_2,\ssu_3\})=1 \text{ for each $i=1,2,3$},
\end{align*}
from which the other values of $\eta!_{\sss,\pi}$ can be derived. Finally, it is clear that the associated polytope is defined in $\Delta!_3$ by a single condition: $q(\ssu_1)+q(\ssu_2)+q(\ssu_3)\leq 1$, and hence it is the brick $\B_0$, depicted in brown, with vertices $OABC$ in Figure~\ref{fig:k4}.

\begin{figure}
  \centering
 \scalebox{.3}{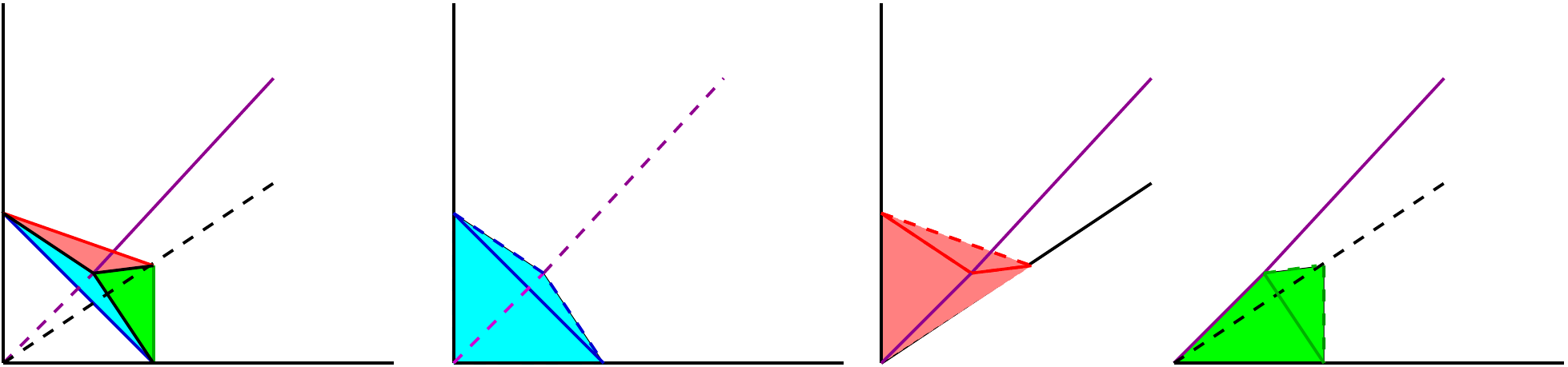}
  \caption{The fan $\ssSigma_{\B_0}$ is the product of the fan $\R_{\geq 0}^{E\setminus \ssE_{0}}$, and the fan given on the left, obtained as the barycentric subdivision of $\R_{\geq 0}^{\ssE_{0}}$ in three 3-dimensional cones $\sssigma^0_1$, $\sssigma^0_2$, and $\sssigma^0_3$, depicted in red, blue, and green.}
  \label{fig:sigma-B0}
\end{figure}

By symmetry, choosing $\ssu_2$ and $\ssu_3$ instead of $\ssu_1$, we get two other slope-level pairs whose associated polytope is $\B_0$, with associated slope-level cones denoted $\sssigma^0_2$ and $\sssigma^0_3$ given by the inequalities $\ell_{0,2}\leq\ell_{0,1},\ell_{0,3}$ and $\ell_{0,3}\leq\ell_{0,1},\ell_{0,2}$, respectively. Notice that $\sssigma^0_1,\sssigma^0_2,\sssigma^0_3$ and their faces form a fan with support $\R^E_{\geq 0}$, therefore, we get a description of the fan $\ssSigma_{\B_0}$ as the product of the fan $\R_{\geq 0}^{E \setminus \ssE_0}$ and the fan on $\R_{\geq 0}^{\ssE_0}$ obtained by a barycentric subdivision, where $E_0\coloneqq\{\sse_{0,1},\sse_{0,2},\sse_{0,3}\}$; see Figure~\ref{fig:sigma-B0}. 
The smallest cone in $\ssSigma_{\B_0}$ that meets $\R^E_{>0}$, which is also the common face of the cones $\sssigma^0_i$, is the cone given by the equalities $\ell_{0,1}=\ell_{0,2}=\ell_{0,3}$. The corresponding slope-level pair can be obtained by successively squashing the slope-level pair $(\sss,\pi)$ above relative to the active circuits $\dz_1=\ssa_{2,0}\ssa_{0,1}$ and $\dz_2=\ssa_{3,0}\ssa_{0,1}$ in any order: it yields $(\sss^0,\pi)$, where $\sss^0(\ssa_{0,i})=1$ for $i=1,2,3$ and all arrows are integer for $\sss^0$. The slope function associated to $\sss^0$ is larger than $\zeta!_{\sss}$; we have
\[
\begin{aligned}
&\zeta!_{\sss^0}(\{\ssu_0\})=3, \qquad \zeta!_{\sss^0}(\{\ssu_0,\ssu_i\})=2\text{ for $i=1,2,3$},\\
&\zeta!_{\sss^0}(\{\ssu_0,\ssu_i,\ssu_j\})=1\text{ for $1\leq i<j\leq 3$},
\end{aligned}
\]
and $\zeta!_{\sss^0}(I)=0$ for other subset $I \subseteq V$. Then, $\eta!_{\sss^0,\pi}$ is the function $\varphi$ in Example~\ref{ex:4lines-polytope}, which, as pointed out there, is quite different from its UpMin transform.

By symmetry, the remaining fans $\ssSigma_{\B_1}$, $\ssSigma_{\B_2}$, and $\ssSigma_{\B_3}$ have a similar description. This gives a description of the fan $\ssSigma(K_4)$. It is not difficult to see that all the cones in $\ssSigma(K_4)$ that meet $\R^E_{>0}$ have a face in common, which is the ray consisting of constant edge length functions. 
\end{example}

\section{The variety of limit canonical series}\label{sec-var-lim-can-series} 

We fix a nodal connected proper curve $X$ over an algebraically closed field $\k$. Let $g$ be its (arithmetic) genus and $G=(V,E)$ be its dual graph. Assume $g>0$. Recall the relevant notation from the beginning of Section~\ref{sec:expected}. In this section, we construct the variety of (fundamental collections of) limit canonical series on $X$

Let $\Sigma$ be the canonical fan of $(G, \g)$. The set  $\PSL$ of permissible slope-level pairs is finite is finite. Therefore, there is a finite-dimensional subspace $\scrU_v \subset \Omega!_v$ for each vertex $v\in V$ such that $H^0(\varC_v, \ssL_{s,v}) \subseteq \scrU_v$ for each $(s,\pi)\in \PSL$ and each $v\in V$. Put $\scrU\coloneqq\bigoplus\scrU_v$.

\subsection{Grassmannians and bricks} 
We will consider the Grassmannian of $g$-dimensional subspaces of $\scrU$, denoted by $\varGr\coloneqq\grass\big(g,\scrU\big)$.

As in Section~\ref{sec:bricks-fans}, let $\Br_g$ denote the collection of all full-dimensional bricks $\B_\beta$ for $\beta\in \Delta!_g$. 
For each full-dimensional brick $\B \in \Br_g$, we define the \emph{brick locus} 
\[
\varGr^{\B}\coloneqq\grass^{\B}(g,\scrU)\coloneqq\Bigl\{\ssW \subseteq \scrU\,\st \,\P_{\nu!^*_\ssW}\supseteq \B\Bigr\},
\]
where $\P_{\nu!^*_\ssW}\subseteq \Delta!_g$ is the base polytope associated to the submodular function $\nu!^*_{\ssW}$; see Section~\ref{sec:submodular} and Appendix~\ref{sec:torusgrass}. By Theorem~\ref{GIT-projective}, the brick locus $\varGr^{\B}$ is in an open subvariety of $\varGr$.

Furthermore, the component-wise action of the torus $\Gm^V$ on $\scrU =\bigoplus_{v\in V}\scrU_v$ induces a natural action of $\Gm^V$ on $\varGr$, with the diagonal $\Gm \hookrightarrow \Gm^V$  acting trivially. Let $\varT \coloneqq\rquot{\Gm^V}{\Gm}$ be the quotient torus. The brick locus $\varGr^{\B}$ is invariant under the action of $\varT$. 

\begin{prop}\label{prop:proj} The categorical quotient $\rquot{\varGr^{\B}}{\varT}$ exists and is projective.
\end{prop}

\begin{proof} This is Theorem~\ref{GIT-projective}.
\end{proof}

\subsection{Strata associated to cones in $\Sigma$ and the variety of limit canonical series}\label{sec:strata-cones}
We assume the characteristic of $\k$ is zero. We assume the branches $\ssp^a$ for $a\in \E_v$ are in general position on the curve $\varC_v$, for each $v\in V$, and fix isomorphisms $\mathcal O_{\varC_{v}}(\ssp^a)\rest{\ssp^a}\cong\k$. In addition, 
for each edge $e\in E$, we fix an isomorphism $\omega\rest{\ssp^e}\cong\k$, where $\omega$ is the canonical sheaf on $X$.

For each permissible slope-level pair $(s,\pi)$, recall the submodular function $\eta!_{s,\pi} = \gamma!_{\pi}+\mathfrak g+\zeta!_{s}$, where $\gamma!_{\pi}$ is the submodular function of the $\pi$-residue space $\GlobSp_{\pi}\subseteq\bV$, $\g$ is the genus function of $X$, and $\zeta!_{s}$ is the slope submodular function; see Section~\ref{sec:etah}. The associated polytope to $\upmin(\eta_{s,\pi})$ in $\Delta!_g$ is denoted $\P_{s,\pi}$.

To each full-dimensional brick $\B\in\Br_g$, we associated the fan $\ssSigma_{\B}$, defined by the cones $\sssigma_{s,\pi}$ for permissible slope-level pairs $(s,\pi)$ in $\PSL(\B)$, verifying $\intcone_{s,\pi}\neq\emptyset$ and $\P_{s,\pi}\supseteq\B$; see Definitions~\ref{defi:permissible}~and~\ref{def:sigmaB}.

By Lemma~\ref{lem:PSLB}, for each cone $\sigma\in\ssSigma_{\B}$ meeting $\R^E_{>0}$, there is a unique $(s,\pi)\in\PSL(\B)$ such that $\sigma=\cone_{s,\pi}$. We can thus put 
\[
\GlobSp_{\sigma}\coloneqq \GlobSp_{\pi}, \quad \ssE^\rmint_{\sssigma} \coloneqq \ssE^\rmint_s \quad \text{and} \quad \ssA^\rmint_{\sssigma} \coloneqq \ssA^\rmint_s,
\]
the $\pi$-residue space, the set of edges and the set of upward arrows which are integer for $s$, respectively. There is a unique upward arrow in $\ssA^\rmint_{\sssigma}$ over each $e\in\ssE^\rmint_{\sssigma}$.

For each $v\in V$, put
\[
\ssL_{\sigma,v}\coloneqq\ssL_{s,v}= \omega!_{v}\Big(\sum_{a\in\E_v} \big(1+ s(a) \big) \ssp^a\Big),
\]
and note that, by the choice of the $\scrU_v$, we have $\bigoplus H^0(\varC_v,\ssL_{\sigma,v}) \subseteq \scrU$. 

Using the isomorphisms $\omega\rest{\ssp^e}\cong\k$ and $\cO_{\varC_{\te_a}}(\ssp^a)\rest{\ssp^a}\cong\k$ to obtain  isomorphisms $\ssL_{\sigma,\te_a}\rest{\ssp^a}\cong \k$ for all arrows $a\in \E$, we may define, for each function $y\colon\ssA^{\rmint}_\sigma\to\k^\times$,
\[
\ssL_{\sigma}^{y}\coloneqq\ssL_{s}^y =\bigcap_{a\in\ssA^\rmint_{\sigma}}\mathrm{Ker}\Bigl(\bigoplus_{v\in V} \ssL_{\sigma,v}\longrightarrow\ssL_{\sigma,\te_a}\rest{p^a}\oplus\ssL_{\sigma,\he_a}\rest{p^{\bar a}}\cong\k\oplus\k \xrightarrow{(-\ssy_a^{-1},1)}\k\Bigr).
\]
Also, for each point $\ssrho$ on $\Gm^E$, we define $\ssy(s, \rho)\colon \ssA^\rmint_\sigma \to \k^\times$ by $y(s,\rho)_a\coloneqq\ssrho_{e}^{s(a)}$ for each $a\in\ssA^{\rmint}_\sigma$ with $e$ the edge underlying $a$. Define the map
\[
\ssboldPhi^{\B}_{\sigma}\colon\Gm^E\longrightarrow   \frac{\varGr^{\B}}{\varT} \quad \text{by} \quad \ssboldPhi^{\B}_{\sigma}(\rho)\coloneqq\Gm^V\ssW_\sigma(y(s,\rho))\quad \text{where}
\]
\begin{equation}\label{eq:orbitW}
\ssW_\sigma(y(s,\rho))\coloneqq\Res^{-1}(\GlobSp_{\sigma}) \bigcap H^0(X,\ssL_{\sigma}^{y(s,\rho)})\subseteq \bigoplus_{v\in V} H^0(\varC_v, \ssL_{\sigma,v})\subseteq \scrU.
\end{equation}

Since $\sigma$ is rational, there is $\ell\in\Z^E_{>0}$ such that $\ell\in\intcone$. As explained in Appendix~\ref{sec:smoothing-and-gluing}, for each $\rho\in\Gm^E$, there is a degeneration to $X$ with singularity degrees given by $\ell$ and leading coefficients given by $\rho$. Consider the associated tropicalization of the geometric generic fiber of the degeneration. It follows from the discussion in Appendix~\ref{sec:smoothing-and-gluing} and Theorem~\ref{thm:characterization-FC} that the limit canonical series whose associated base polytope contains $\B$ is $\ssW_\sigma(y(s,\rho))$. Since the associated submodular function is $\upmin(\eta!_{s,\pi})$, thus independent of $\rho$, it follows that $\ssboldPhi^{\B}_{\sigma}$, which does not depend on the choice of $\ell\in\intcone$, is a well-defined morphism of schemes for each $\B\in\Br_g$ and $\sigma\in\ssSigma_{\B}$.

Let $\Sigma=\ssSigma(G,\g)$ be the meet of the $\ssSigma_{\B}$ as $\B$ runs in $\Br_g$. For each cone $\sigma$ in $\Sigma$ meeting $\R^E_{>0}$, and each full-dimensional brick $\B\in\Br_g$, denote by $\sssigma_{\B}$ the smallest cone in $\ssSigma_{\B}$ that contains $\sigma$. Define
\[
\ssboldPhi_{\sigma}\colon\Gm^E\, \xrightarrow{\ssub{\bigl( \ssPhi^{\B}_{\sssigma_{\B}}\bigr)}!_{\mathrm B\in\Br_g}}\, \prod_{\B\in\Br_g}\frac{\varGr^{\B}}{\varT}, \quad \text{and put} \quad 
\varV_\sigma\coloneqq\mathrm{Im}(\ssboldPhi_{\sigma})\subset
\prod_{\B\in\Br_g}\frac{\varGr^{\B}}{\varT}.
\]

For each $\ell\in\R^E_{>0}$, there is a unique $\sigma\in\Sigma$ such that $\ell\in\intcone$. Conversely, for each $\sigma\in\Sigma$, since $\sigma$ is rational, there is $\ell\in\Z^E_{>0}$ such that $\ell\in\intcone$. As before, combining Theorem~\ref{thm:characterization-FC} with the discussion in Appendix~\ref{sec:smoothing-and-gluing}, we see that a point on $\varV_\sigma$ parameterizes the collection of $\Gm^V$-orbits of the spaces in the fundamental collection $\FC(\varX)$ of limit canonical series associated to the tropicalization of a smooth curve $\varX$ to $(X, \Gamma)$, where $\Gamma$ is the metric graph associated to the pair $(G, \ell)$, or equivalently, to a degeneration to $X$ having singularity degrees given by $\ell$, for all $\ell\in\Z^E_{>0}$ with $\ell\in\intcone$. Conversely, each such collection is parameterized by a point on $\varV_\sigma$.

Thus, the union $\varV$ of the $\varV_\sigma$ is the parameter space for the collections of $\Gm^V$-orbits of the spaces in the fundamental collections of limit canonical series associated to all tropicalizations yielding $X$, or equivalently, all smoothings of $X$. We call $\varV$ the variety of limit canonical series.

\subsection{Projectivity}\label{sec:projectivity} Notation as before.  We show that $\varV$ is closed in the product of the $\varGr^{\B}/\varT$ for $\B\in\Br_g$, hence projective by Proposition~\ref{prop:proj}.

\begin{thm}\label{thm:variety-lcs} Let $X$ be a nodal curve over an algebraically closed field $\k$ of characteristic zero with branches over nodes in general position on the components of its normalization, and let $G$ be its dual graph. Then:
\begin{enumerate}
    \item The variety of limit canonical series $\varV$ parametrizes all fundamental collections associated to any smooth proper curve  $\varX$ over a valued field $\K$ that tropicalizes to $(X, \Gamma)$ for any metric graph $\Gamma$ over $G$. 
    \smallskip
   \item For each cone $\tau$ of $\Sigma$ meeting $\R_{>0}^E$, we have that $\compvarV_{\tau}=\bigcup_{\sigma\supseteq\tau}\varV_{\sigma}$.  In particular, $\varV$ is a projective variety.
\end{enumerate}
\end{thm}

\begin{proof} The partial gluing $\varrho!_h=\varrho!_h(\varX)$ appearing in Theorem~\ref{thm:characterization-FC} is given by 
\[
\varrho!_{h,a} = \ssrho_e^{\sss_h(a)} \qquad \qquad  \text{for all } a\in \ssA_s^\rmint
\]
where $e$ is the underlying edge of $a$ and $\ssrho_e$ is the leading coefficient of the smoothing of the node $\ssp^e$ in $\varX$, see Appendix~\ref{sec:smoothing-and-gluing}.

\smallskip

It follows from this description and Theorem~\ref{thm:characterization-FC} that the data 
\[
(\varrho!_{h,a})_{h\in\R^V,a\in\ssA_{\sss_h}^{\rmint}}
\]
appearing in the description of the fundamental collection
for $\varX$ varying among all smooth curves over valued fields $\K$ tropicalizing to $(X,\Gamma)$ is the same as the data
\[
(\ssrho_a^{\sss_h(a)})_{h\in\R^V,a\in\ssA_{\sss_h}^{\rmint}}
\]
given by $y(\sss_h,\rho)$, 
for $\rho$ varying in $\Gm^E(\k) = \ssub{\bigl(\k^\times\bigr)}!^E$. This proves the first statement.

To prove the second, it is enough to prove the equality $\compvarV_{\tau}=\bigcup_{\sigma\supseteq\tau}\varV_{\sigma}$. Let $\sigma$ and $\tau$ be cones in $\Sigma$ meeting $\R^E_{>0}$ with $\sigma\supseteq\tau$. We show first that $\compvarV_{\tau}\supseteq\varV_{\sigma}$.  We use the notation in Section~\ref{sec:strata-cones}. In particular, for each $\B\in\Br_g$, we denote by $\sssigma_{\B}$ and $\sstau_{\,\B}$ the smallest cones in $\ssSigma_{\B}$ that contain $\sigma$ and $\tau$, respectively.

Let $\rho\in\Gm^E(\k)$ and $\ell\in\ring\sigma$ with $\ell_e\in\N$ for each $e\in E$. Let $\rho'\in\Gm^E(\k[[t]])$ be such that $\ssrho'_e=\ssrho_e t^{\ell_e}$ for each $e\in E$. We only need to prove that
\begin{equation}\label{eq:limW}
\lim_{t\to 0}\ssboldPhi_{\tau}(\rho')=\ssboldPhi_{\sigma}(\rho).
\end{equation}

Take $\B\in\Br_g$ and let $(s,\pi)\in\PSL(\B)$ be the slope-level pair such that $\cone_{\B}=\cone_{s,\pi}$. Up to changing $\ell$ to an integer multiple, we may assume there is $h\colon V\to\Z$ be such that $s=\slztwist{\ell}{h}$. Since $\tau$ is a face of $\cone$, we have that $\sstau_{\B}$ is a face of $\cone_{\B}$, and hence it follows from Theorem~\ref{thm:face-cone} that there is a slope-level pair $(s',\pi')\in \PSL(\B)$ such that $\sstau_{\,\B}=\cone_{s',\pi'}$. Furthermore, $\pi'$ is a coarsening of $\pi$, and we have that $s'\geq s$ and $\ssE_{\sstau_{\B}}^\rmint\supseteq\ssE_{\sssigma_{\B}}^\rmint$. Since $s=\slztwist{\ell}{h}$, we have
\begin{equation}\label{eq:lss1}
\ell_a s'(a)=\ell_a s(a)+\ell_a(s'(a)-s(a))\geq h(\te_a)-h(\he_a)
\end{equation}
for each $a\in\ssA_{\sstau_{\B}}^\rmint$, with equality if and only if $a\in\ssA_{\sssigma_{\B}}^\rmint$.

Finally, let $x\in\Gm^V(\k\{\{t\}\})$ given by $\ssx_v\coloneqq t^{-h(v)}$ for each $v\in V$. Then,
\[
x\Big(\Res^{-1}(\GlobSp_{\sstau_{\,\B}}) \bigcap H^0(X,\ssL_{\sstau_{\,\B}}^{y(s',\ssrho')})\Big)=
\Res^{-1}(x\GlobSp_{\sstau_{\,\B}}) \bigcap H^0(X,\ssL_{\sstau_{\B}}^{y})
\]
where $y\colon \ssA_{\sstau_{\B}}^\rmint\to\k$ is the function given by
\[
\ssy_a\coloneqq\ssrho_{e}^{s'(a)}\,t^{\ell_a s'(a)-h(\te_{a})+h(\he_{a})}
\]
for each $a\in\ssA_{\sstau_{\B}}^\rmint$ with $e$ the edge underlying $a$.
On the one hand, since \eqref{eq:lss1} holds for each $a\in\ssA_{\sstau_{\,\B}}^\rmint$ with equality if and only if $a\in\ssA_{\sssigma_{\,\B}}^\rmint$, we have
\[
\lim_{t\to 0}H^0(X,\ssL_{\sstau_{\,\B}}^{y})=
H^0(X,\ssL_{\sssigma_{\,\B}}^{y(s,\rho)}).
\]
On the other hand, it follows from \cite[Thm.~6.1]{AEG-residue-polytope} that
\[
\lim_{t\to 0}\Res^{-1}(x\GlobSp_{\sstau_{\,\B}})=\Res^{-1}(\GlobSp_{\sssigma_{\,\B}}).
\]
Then, \eqref{eq:limW} holds.

We show that the closure of $\varV_{\tau}$ is contained in the union of the $\varV_{\sigma}$, $\sigma\supseteq\tau$.

Let $\ell'\in\ring\cone$ with $\ell'_e\in\N$ for each $e\in E$. Let $\rho'\in\Gm^E(\k[[t]])$. Let $n\colon E\to\Z_{\geq 0}$, $e\mapsto \ssn_e$, such that for each $e\in E$, we have $\ssrho_e'=t^{\ssn_e}\mu_e$ for $\mu!_e\in\k[[t]]$ with $\mu!_e(0)\neq 0$. Let $\rho\in\Gm^E(\k)$ be given by $\ssrho_e=\mu!_e(0)$ for each $e\in E$. Let $\sigma$ be the cone of $\Sigma$ containing $\ell\coloneqq N\ell'+n$ for all large enough $N\in\N$. Clearly $\sigma\supseteq\tau$. We claim \eqref{eq:limW}.

Indeed, for each $\B\in\Br_g$, as before, let $(s,\pi)$ and $(s',\pi')$ be slope-level pairs in $\PSL(\B)$ such that  $\sssigma_{\B}=\sssigma_{s,\pi}$ and $\sstau_{\,\B}=\sssigma_{s',\pi'}$, and such that $\pi'$ is a coarsening of $\pi$, and we have that $s\leq s'$ and $\ssE_{\sstau_{\B}}^\rmint\supseteq\ssE_{\sssigma_{\B}}^\rmint$. Up to changing $\ell'$ to an integer multiple, and replacing $t$ by a high enough power, we may assume there are $h,h'\colon V\to\Z$ such that $s=\slztwist{\ell}{h}$ and $s'=\slztwist{\ell'}(h')$. Then,
\[
\ssn_e s'(a)=\ell_a s(a)+\ell_a (s'(a)-s(a))-N\ell'_a s'(a)\geq 
h(\te_a)-h(\he_a)-N(h'(\te_a)-h'(\he_a))
\]
for each $a\in\ssA_{\sstau_{\B}}^\rmint$ with underlying edge $e$, with equality if and only if $a\in\ssA_{\sssigma_{\B}}^\rmint$.

Finally, let $x,y\in\Gm^V(\k\{\{t\}\})$ given by $\ssx_v\coloneqq t^{-h(v)}$ and $\ssy_v\coloneqq t^{-h(v)+Nh'(v)}$ for each $v\in V$. Then,
\[
y\Big(\Res^{-1}(\GlobSp_{\sstau_{\,\B}}) \bigcap H^0(X,\ssL_{\sstau_{\,\B}}^{\ssrho'})\Big)=
\Res^{-1}(x\GlobSp_{\sstau_{\,\B}}) \bigcap H^0(X,\ssL_{\sstau_{\B}}^{z}),
\]
where $z\colon \ssA_{\sstau_{\B}}^\rmint\to\k$ is the function satisfying 
\[
\ssz_a\coloneqq\mu!_{e}^{s'(a)}t^{\ssn_e s'(a)-h(\te_a)+h(\he_a)+N(h'(\te_a)-h'(\he_a))}
\]
for each $a\in\ssA_{\sstau_{\B}}^\rmint$ with $e$ the edge underlying $a$. As before,
\[
\lim_{t\to 0}\Res^{-1}(x\GlobSp_{\sstau_{\B}})=\Res^{-1}(\GlobSp_{\sssigma_{\,\B}})\quad\text{and}\quad
\lim_{t\to 0}H^0(X,\ssL_{\sstau_{\,\B}}^{z})=
H^0(X,\ssL_{\sssigma_{\B}}^{y(s,\rho)}).
\]
Then, \eqref{eq:limW} holds, and the theorem follows.
\end{proof}

\begin{example} Let $X$ be the union of four general lines $\ssL_0,\ssL_1,\ssL_2,\ssL_3$ on the plane. Each component is smooth and rational, and contains three branches over nodes, which are thus in general position. Its dual graph, $G=(V,E)$, is the complete graph on $4$ vertices with genus function zero, $\g=0$. In Example~\ref{ex:K4}, we described the canonical fan $\ssSigma$ of $\ssK_4$, and we observed it has a ray $\tau$ which is a face of any cone in the meet which intersects $\R^E_{>0}$. It thus follows from Theorem~\ref{thm:variety-lcs} that the variety of limit canonical series $\varV$ is irreducible, the closure of $\varV_{\tau}$. 

The ray $\tau$ is the one where all lengths are equal. Smoothings with edge length function $\ell$ constant equal to $1$ are obtained by pencils $\ssL_0\ssL_1\ssL_2\ssL_3-tF$ where $F$ is a quartic meeting each $\ssL_i$ away from the branches over the nodes for each $i$. Notation as in Example~\ref{ex:K4}, we let $p_{i,j}\coloneqq\ssp^{\ssa_{i,j}}$ for each distinct $i,j\in\{0,1,2,3\}$.

Fix a line $L_j$. Then $\ssL_i/\ssL_r$ is a local equation for $\ssL_i$ at $p_{i,j}$ for each distinct $i,j,r$. We have
\[
\frac{\ssL_i}{\ssL_r}\frac{\ssL_j}{\ssL_s}=t\frac{F}{\ssL_r^2\ssL_s^2},
\]
for each $i,j,r,s$ such that $\{i,j,r,s\}=\{0,1,2,3\}$. Let $\ssx_1,\ssx_2,\ssx_3$ be homogeneous coordinates for $\mathbf P^2_{\k}$. We may assume $\ssL_i=\ssx_i$ for $i=1,2,3$ and $\ssL_0=\ssx_1+\ssx_2+\ssx_3$. Write $F=\sum\ssa_{i,j,r}\ssx_1^i\ssx_2^j\ssx_3^r$. With these choices, the leading coefficients for the smoothing  $\ssL_1\ssL_2\ssL_3\ssL_4-tF$ are 
\begin{align*}
    \ssrho_{0,1}=&\ssa_{0,4,0}+\ssa_{0,2,2}+\ssa_{0,0,4}-(\ssa_{0,3,1}+\ssa_{0,1,3}), \qquad   \ssrho_{1,2}=\ssa_{0,0,4}\\
    \ssrho_{0,2}=&\ssa_{4,0,0}+\ssa_{2,0,2}+\ssa_{0,0,4}-(\ssa_{3,0,1}+\ssa_{1,0,3}) ,\qquad 
    \ssrho_{1,3}=\ssa_{0,4,0}, \\
    \ssrho_{0,3}=&\ssa_{4,0,0}+\ssa_{2,2,0}+\ssa_{0,4,0}-(\ssa_{3,1,0}+\ssa_{1,3,0}),  \qquad 
    \ssrho_{2,3}=\ssa_{4,0,0}. 
\end{align*}
So we see that we get the whole $\Gm^E$ of leading coefficients as we vary $F$. 
\end{example}

\bibliographystyle{alpha}
\bibliography{bibliography}

\appendix

\section{Torus actions on Grassmannians}\label{sec:torusgrass}

Let $\k$ be an algebraically closed field and $V$ be a finite set. For each $v\in V$, let $\rmU_v$ be a finite-dimensional vector space over $\k$ and denote by $\ssd_v$ its dimension. Let $d\colon V\to\Z$ the function $v\mapsto\ssd_v$ and put $\rmU\coloneqq\bigoplus\rmU_v$. 
Let $n$ be any nonnegative integer and $\varGr\coloneqq\mathrm{Grass}(n,\rmU)$, the Grassmannian of vector subspaces of $\rmU$ of dimension $n$. 

For each subspace $\ssW\subseteq\rmU$ and subset $S\subseteq V$, let $\ssW_S\subseteq\bigoplus_{v\in S}\rmU_v$ denote the image of $\ssW$ under the corresponding projection $\proj{S}\colon \rmU \to \rmU_S \coloneqq \bigoplus_{v\in S}\rmU_v$. To each subspace $\ssW\subseteq\rmU$, we associate its adjoint modular pair $(\nu!_{\ssW},\nu!_{\ssW}^*)$, satisfying the equations
\[
\nu!_{\ssW}^*(S)=\dim\ssW_S=\dim\ssW-\nu!_{\ssW}(V\setminus S)\quad\text{ for all }S\subseteq V,
\]
and the corresponding polytope 
\[
\Q_{\ssW}=\Big\{q\in\R^V\,\Big|\,\nu!_{\ssW}(S)\leq q(S)\leq
\nu!_{\ssW}^*(S)\quad \forall S\subseteq V\Big\};
\]
see \cite[\S 8]{AE-modular-polytopes}. 
The largest possible polytope is achieved for $\ssW$ with projection maps $\ssW\to\rmU_S$ of maximal rank, equal to $\min(\dim\ssW,\ssd(S))$, with $d(S)=\sum_{v\in V}\ssd_v$. So, if $\dim\ssW=n$, then $\Q_{\ssW}\subseteq\Q_{n,d}$ where
\[
\Q_{n,d}\coloneqq\Big\{q\in\R^V\,\Big|\,q(V)=n\text{\,\, and \,\,}q(S)\leq
\min(n,\ssd(S))\,\,\,\forall \,\, S\subseteq V\Big\}.
\]

\subsection{Pl\"ucker coordinates.} Consider the Pl\"ucker embedding of $\varGr$, that is, the association of the one-dimensional subspace $\extp^n\ssW$ of $\extp^n\rmU$ to each $n$-dimensional subspace $\ssW\subseteq\rmU$. To describe $\extp^n\rmU$, let $\Q^{\Z}_{n,d}\subseteq\Q_{n,d}$ be the subset of integer-valued functions. For each $r\in\Q^\Z_{n,d}$, sending $v$ to $r_v$, let $\extp^r\rmU$ denote the tensor product over $v\in V$ of the $\extp^{r_v}\rmU_v$. Then,
\[
\extp^n\rmU=\bigoplus_{r\in\Q^{\Z}_{n,d}}\extp^r\rmU.
\]
Thus, for each $\k$-point $\ssW$ of $\varGr$, there are $\ssx_r=\ssx_r^{\ssW}\in\extp^r\rmU$ for $r\in\Q^{\Z}_{n,d}$ whose sum is a generator of the subspace 
$\extp^n\ssW$ of $\extp^n\rmU$. The $\ssx_r$ are the ``Pl\"ucker coordinates'' of $W$.

\begin{prop}\label{integerP} If $\ssx_r^{\ssW}\neq 0$, then $r\in\Q_{\ssW}$. Conversely, if $r$ is a vertex of $\Q_{\ssW}$, then $\ssx_r^{\ssW}\neq 0$.
\end{prop}

\begin{proof} Clearly, $\ssx_r^{\ssW}\neq 0$ only if $\nu!_{\ssW}^*(S)\geq r(S)$ for each $S\subseteq V$, and hence $r\in\Q_{\ssW}$. 

Conversely, by \cite[Prop.~2.7]{AE-modular-polytopes}, a vertex $r$ of $\Q_{\ssW}$ is the unique point of the base polytope of a $\pi$-splitting of $\nu!_{\ssW}$ for a maximal ordered partition $\pi$ of $V$. Let $v_1,\dots,v_b$ be the ascending sequence of the elements of $V$ compatible with the order in which the singleton $\{v_i\}$ appear in $\pi$. Then, there is a filtration $\ssW=\ssW_b\supseteq \ssW_{b-1}\supseteq\cdots\supseteq \ssW_1\supseteq \ssW_0=0$ such that the kernel of the projection of $\ssW_i$ to $\rmU_{v_i}$ is $\ssW_{i-1}$ and the image has dimension $\ssr_i$ for each $i=1,\dots,b$. It is now easy to show that $\ssx_r^{\ssW}\neq 0$. 
\end{proof}

It will be necessary to compose the Pl\"ucker embedding with the $m$-th Veronese embedding, for each $m\in\mathbb N$. We call the composition the \emph{$m$-th Pl\"ucker embedding} of $\varGr$. The resulting ``Veronese--Pl\"ucker coordinate'' of $\ssW$ are  
\begin{equation}\label{coordW}
\prod_{r\in\Q^{\Z}_{n,d}}\ssx_r^{s_r}
\end{equation}
for all nonnegative functions $s\colon\Q^{\Z}_{n,d}\to\Z_{\geq 0}$, $r\mapsto \sss_r$ with $s(\Q^{\Z}_{n,d})=m$, that is, the sum of $\sss_r$ for $r\in\Q^{\Z}_{n,d}$ equals $m$.

\subsection{The torus action.}

The torus $\Gm^V$ acts naturally on the Grassmannian $\varGr$, with the action of the diagonal $\Gm \hookrightarrow \Gm^V$  being trivial. Let $\varT \coloneqq\rquot{\Gm^V}{\Gm}$ be the quotient torus. 

For each $\beta\in\R^V$ with $\beta(V)\in\Z$, let 
$\B_\beta$ be defined as the set of all functions $q\in\R^V$ verifying
\[
\lfloor \beta(S)\rfloor \leq q(S) \leq \lceil \beta(S)\rceil\qquad \qquad \forall \,\, S \subseteq V.
\]
We call $\B_\beta$ the \emph{brick} associated to $\beta$. 

The maximum dimension of a brick $\B_\beta$ is $|V|-1$, achieved if and only if $\beta(S)\not\in\Z$ for any nonempty proper $S\subset V$. If so, we call $\B_\beta$ \emph{full-dimensional}.

By (a direct generalization of) \cite[Thm.~10.4]{AE-modular-polytopes}, the collection of $\B_\beta$ for $\beta\in\Q_{n,d}$ form a tiling of $\Q_{n,d}$. For each full-dimensional brick $\B\subseteq\Q_{n,d}$, let
\[
\varGr^{\B}=\Bigl\{\, \ssW\subseteq\rmU\,\st \,\Q_{\ssW}\supseteq\B\,\Bigr\} \subseteq\varGr.
\]
(If $\Q_{\ssW}\supseteq\B$, then $\dim\ssW=n$.)

\begin{thm}\label{GIT-projective} For each full-dimensional brick $\B\subseteq\Q_{n,d}$,  $\varGr^{\B}$ is a $\varT$-invariant open subset of $\varGr$ whose categorical quotient $\rquot{\varGr^{\B}}{\varT}$ exists and is projective.
\end{thm}

\begin{proof} We will show that $\rquot{\varGr^{\B}}{\Gm^V}$ is a GIT-quotient which is an orbit space. Let $q\in\mathring\B$, the interior of $\B$, be such that $\ssq_v\in\mathbb Q$ for each $v\in V$. Let $m\in\N$ be such that $\ssp_v\coloneqq\ssq_vm$ is an integer for each $v\in V$.

Let $\ssx_r\in\extp^r\rmU$ for $r\in\Q^{\Z}_{n,d}$ be the ``Pl\"ucker coordinates'' of $\ssW$. The corresponding ``coordinates'' for the $m$-th Pl\"ucker embedding are the products \eqref{coordW} for all nonnegative integer-valued functions $s\colon\Q^{\Z}_{n,d}\to\mathbb Z$ with $s(\Q^{\Z}_{n,d})=m$.

Consider the following linearization to the symmetric product $S^m(\extp^n\rmU)$ of the action of $\Gm^V$ on $\varGr$: The action of $t\in\Gm^V(\k)$ changes the ``coordinates'' \eqref{coordW} of $\ssW$ to the ``coordinates''
\[
\prod_{v\in V}\sst_v^{\sum \sss_r\ssr_v-\ssp_v}\prod_{r\in\Q^{\mathbb Z}_{n,d}}\ssx_r^{\sss_r}.
\]
Notice that the inverse of $\prod\sst_v^{\ssp_v}$ is multiplying all the ``coordinates,'' hence the action of $\Gm^V$ on $S^m(\bigwedge^n\rmU)$ induces the natural action of $\Gm^V$ on $\varGr$. 

It is now enough to prove that the semistable locus of $\varGr$ for this linear action is equal to the stable locus, and also equal to $\varGr^{\B}$. We will use Hilbert--Mumford Criterion. We must argue that, given $\ssW\in\varGr(\k)$, we have  $\ssW\in\varGr^{\B}(\k)$ if and only if for each nonconstant $\mu\colon V\to\Z$, there is a map $s\colon\Q^{\Z}_{n,d}\to\Z_{\geq 0}$ such that 
\begin{enumerate}[label=(\arabic*)] 
    \item \label{st:1} $s(\Q^{\Z}_{n,d})=m$;
    \item \label{st:2} $\ssx_r\neq 0$ for each $r\in\Q^{\Z}_{n,d}$ such that $\sss_r>0$;
    \item \label{st:3} $\sum\sss_r\ssr_v\mu!_v\leq\sum\ssp_v\mu! _v$ (with strict inequality if $\ssW\in\varGr^{\B}(\k)$).
\end{enumerate}

Assume first that $\ssW\in\varGr^{\B}(\k)$. Then, we have $\B\subseteq\Q_{\ssW}$. Therefore, $q$ is a convex linear combination of the vertices of $\Q_{\ssW}$, that is, $q=\sum\ssy_r r$ for a certain  $y\colon\Q^{\Z}_{n,d}\to\mathbb Q_{\geq 0}$ with $y(\Q^{\Z}_{n,d})=1$ such that $\ssy_r\neq 0$ only if $r$ is a vertex of $\Q_{\ssW}$.

Fix a nonconstant $\mu\colon V\to\Z$. Suppose by sake of contradiction that for each $r\in\Q^{\Z}_{n,d}$ such that $\ssx_r\neq 0$, we have $m\sum\ssr_v\mu!_v\geq\sum\ssp_v\mu!_v$. (This is a negation of \ref{st:3} for $s\colon\Q^{\Z}_{n,d}\to\Z_{\geq 0}$ with $\sss_r=m$.) Then, $\sum\ssr_v\mu!_v\geq\sum\ssq_v\mu!_v$ for each vertex $r$ of $\Q_{\ssW}$ by Proposition~\ref{integerP}, and hence 
\begin{equation}\label{eq:st}
\sum_{r,v}\ssy_r\ssr_v\mu!_v\geq\sum_v\ssq_v\mu!_v,
\end{equation}
with equality if and only if $\sum\ssr_v\mu!_v=\sum\ssq_v\mu!_v$ for every vertex $r$ of $\Q_{\ssW}$. But the equality holds in \eqref{eq:st}, since $q=\sum\ssy_r r$. Therefore, $\Q_{\ssW}$ is contained in the hyperplane $(\mu,\cdot)=\sum\ssq_v\mu!_v$. Since $\B\subseteq\Q_{\ssW}$, this is only possible if $\mu$ is constant, a contradiction.

Conversely, if for each nonconstant $\mu\colon V\to\Z$, there is $s\colon\Q^{\Z}_{n,d}\to\Z_{\geq 0}$ satisfying \ref{st:1}, \ref{st:2} and \ref{st:3}, then consider the case where $\mu=\one!_S$ for a nonempty proper subset $S\subset V$. Then, there is $s\colon\Q^{\Z}_{n,d}\to\Z_{\geq 0}$ satisfying \ref{st:1} and \ref{st:2} such that
\[
\sum_r\sss_r r(S)\leq m q(S).
\]
Now, if $\sss_r>0$, then $\ssx_r\neq 0$ by (2), thus $r\in\Q_{\ssW}$ by Proposition~\ref{integerP}, whence $\nu!_{\ssW}(S)\leq r(S)$. It follows that $\nu!_{\ssW}(S)\leq q(S)$. As this holds for each nonempty proper subset $S\subset V$, it follows that $q\in\Q_{\ssW}$, and hence $\ssW\in\varGr^{\B}$.
\end{proof}

\section{Smoothing and gluing}\label{sec:smoothing-and-gluing}

Let $X$ be a nodal reduced connected proper curve over an algebraically closed field $\k$ of characteristic zero whose singularities are nodes.  Let $G$ be the dual graph of $X$, with vertex set $V$ and edge set $E$. Let $\E$ be the arrow set of $G$. For each $v\in V$, let $\ssX_v$ be the corresponding irreducible component of $X$ and $\varC_v$ its normalization. For each $a=uv\in\E$, let $\ssp^a$ be the branch on $\varC_{u}$ above the node $\ssp^e$ of $X$ corresponding to the edge $e\in E$ underlying $a$.

We describe all the fractional ideal limits $\ssT_h$ of the structure sheaf along families of smooth curves degenerating to $X$, by computing the sheaves of fractional ideals $\ssT_{h,v}$ each $\ssT_h$ generates on the $\varC_v$ and the gluing of the $\ssT_{h,v}$ that defines $\ssT_h$.

\subsection{Smoothings}\label{sec:smoothings} 
Let $\scrC/\varD$ be the versal deformation of $X$. The versal deformation comes with an identification of the special fiber of $\scrC/\varD$ with $X$, thus we may view $X\subseteq\scrC$. As explained in \cite[pp.~79--81]{DM69} and reviewed in \cite[pp.~288--9]{EM}, we have that $\varD$ is the power series ring over $\k$ in variables $\sst_e$, for $e\in E$, and $\sss_1,\dots,\sss_l$, for a certain integer $l$, chosen in such a way that for each $e=\{u,v\}\in E$, letting $a=uv$ and $\bar a=vu$ be the arrows supported on $e$, we have an isomorphism of $\varD$-algebras,
\[
\psi_e\:\hatcO_{\scrC,\ssp^e}\xrightarrow{\,\cong\,} \varD[[\ssz_a,\ssz_{\bar a}]]/(\ssz_a \ssz_{\bar a}-\sst_e).
\]
We will assume that $\ssz_a=0$ (resp.~$\ssz_{\bar a}=0$) is the formal local equation of the component $\ssX_v$ (resp.~$\ssX_u$). Letting $\hatz_a$ and $\hatz_{\bar a}$ denote the elements of $\hatcO_{\scrC,\ssp^e}$ corresponding to $\ssz_a$ and $\ssz_{\bar a}$, we have that $\hatz_a$ restricts to a local parameter $\barz_a$ of $\hatcO_{\varC_u,\ssp^a}$, whereas, $\hatz_{\bar a}$ restricts to a local parameter $\barz_{\bar a}$ of 
$\hatcO_{\varC_v,\ssp^{\bar a}}$. 

Let $\ell\colon E\to\N$ be an edge length function. For convenience, for each arrow $a\in\E$, we put $\ell_a\coloneqq\ell_e$, where $e$ is the underlying edge. Let $\rho\in\Gm^E(\k)$. Fix a local homomorphism $\varD\to\k[[t]]$ taking $t_e$ for each $e\in E$ to a power series $\xi_e\in\k[[t]]$ with leading term $\ssrho_e t^{\ell_e}$. Consider the induced Cartesian diagrams:
\[
\begin{CD}
X @>>> \scrX @>>> \scrC\\
@VVV @V{\frakp} VV @VVV\\
\text{Spec}(\k) @>>> \text{Spec}(\k[[t]]) @>>> \text{Spec}(\varD).
\end{CD}
\]
Denoting by $\wtz_a$ and $\wtz_{\bar a}$ the pullbacks under the map $\scrX\to\scrC$ of $\hatz_a$ and $\hatz_{\bar a}$, respectively, we have $\wtz_a\wtz_{\bar a}=\xi_e$ in $\hatcO_{\scrX,\ssp^e}$ for each $e\in E$.

The map $\frakp$ is a flat projective map to $\text{Spec}(\k[[t]])$ whose generic fiber $\scrX_\eta$ is smooth and whose special fiber $\scrX_0$ is isomorphic to $\ssX$. We call such a map a  \emph{smoothing} of $X$. Conversely, by the versal property of $\scrC/\varD$, every smoothing of $X$ arises the way we described. We call $\ell_e$ the \emph{singularity degree} of $\frakp$ at $\ssp^e$ for each $e\in E$, and say the resulting function $\ell\colon E\to\Z_{>0}$ is the associated edge length function. And we call the $\ssrho_e$ the \emph{leading coefficients} of $\frakp$. Clearly, though $\ell$ does not depend on the choices we made, the leading coefficients do. We describe how the limits $\ssT_h$ of the structure sheaf along $\frakp$ depend on the $\ssrho_e$ for $\ell$ fixed.

The components $\ssX_v$ are Cartier divisors of $\scrX$ away from the nodes $\ssp^e$. If $\frakp$ is \emph{regular}, that is, $\scrX$ is regular, equivalently, if all the $\ell_e$ are equal to $1$, then the components $\ssX_v$ are Cartier divisors of $\scrX$.

\subsection{Gluings}\label{sec:twisters} Keep the setup above. Let $\frakq\colon\wtscrX\to\scrX$ be the semistable model of $\frakp$, and $\widetilde\frakp\coloneqq\frakp\frakq$. The  special fiber $\wtscrX_0$ of $\widetilde\frakp$ is obtained from $\scrX_0$ by splitting apart the branches $\ssp^a$ and $\ssp^{\bar a}$ over each node $\ssp^e$, and connecting them by a chain of $\ell_e-1$ rational curves. The model $\wtscrX$ is regular, so all the components of $\wtscrX_0$ are Cartier divisors of $\wtscrX$. 

Let $\Gamma$ be the metric graph associated to the pair $(G,\ell)$. Let $V(\Gamma)$ be the collection of all points on $\Gamma$ at integer distances from vertices of $G$. Clearly, $V\subseteq V(\Gamma)$. For each $v\in V(\Gamma)$, let $\ssX_v$ be the corresponding component of $\wtscrX_0$ (abusing notation, if $v\in V$).

For each function $h\colon V\to\Z$, let $\hat h\colon\Gamma\to\R$ be its unique admissible extension; see Section~\ref{thm:existence-uniqueness-canonical-extension}. Put $\ssD_h\coloneqq\div(\hat h)$. Then $\ssD_h$ is $G$-admissible and is supported on $V(\Gamma)$ by Proposition~\ref{prop:adm}. Put
\[
\widetilde\cT_h\coloneqq\cO_{\wtscrX}\Big(-\sum_{v\in V(\Gamma)}\hat h(v)\ssX_v\Big).
\]
For each $v\in V(\Gamma)$, the degree of the restriction $\widetilde\cT_h\rest{\ssX_v}$ is $\ssD_h(v)$. Also, by Proposition~\ref{prop:adm}, for each edge $e\in E$, the restriction $\widetilde\cT_h\rest{\ssX_v}$ has degree $0$ for each $v\in V(\Gamma)\setminus V$ on $e$, with the possible exception of a single $v$, for which the degree is $1$. By \cite[Thm.~3.1]{EP16}, the sheaf 
\[
\cT_h\coloneqq\frakq_*\widetilde\cT_h
\]
is $\frakp$-flat, and its restriction to each fiber of $\frakp$ has degree $0$. Moreover, its restriction to the generic fiber $\scrX_{\eta}$ is trivial.

Clearly, we may and will identify $\cT_0$ with the structure sheaf $\cO_{\scrX}$, and $\cT_n$ with the sheaf of ideals of $n\scrX_0$ in $\scrX$ for each function with constant values $n$. For each $h\colon V\to\Z$, we may and will pick $\ssh_{\min}\coloneqq\min\{h(v)\,|\,v\in V(G)\}$, and view $\cT_h$ as a subsheaf of $\cT_{\ssh_{\min}}$, thus as a fractional ideal.

For each $h\colon V\to\Z$ and arrow $a=uv\in\E$, put
\begin{equation}\label{eq:s-C}
s(a)\coloneqq\sss_h(a)\coloneqq\slztwist{\ell}h(a) \coloneqq \left\lfloor\frac{h(u)-h(v)}{\ell_{a}}\right\rfloor.
\end{equation}
Notice that $s(a)+s(\bar a)=0$ if $\ell_a$ divides $h(u)-h(v)$ and $s(a)+s(\bar a)=-1$ otherwise. Let $e$ be the edge underlying $a$. We claim that the formal stalk of $\cT_h$ at $\ssp^e$ is the fractional ideal
\begin{equation}\label{fracid}
\big(\sst^{h(u)}\wtz_a^{-s(a)}\,,\, \sst^{h(v)}\wtz_{\bar a}^{s(a)+1}\big).
\end{equation}
Indeed, the formal stalk of $\cT_h$ at $\ssp^e$ is a fractional ideal contained in the possibly larger fractional ideal $(\sst^{\ssh_e})$, where $\ssh_e\coloneqq \min(h(u),h(v))$. Assume $h(v)\leq h(u)$. As shown in \cite[\S 3, p.~14]{CEG}, the formal stalk of $\cT_h$ at $\ssp^e$ is the fractional ideal 
\begin{equation}\label{fracid2}
\big(\wtz_{\bar a}^{s(a)+1}\sst^{h(v)},\wtz_{\bar a}^{s(a)}\sst^{h(v)+\ssr_\ell(a)}\big),
\end{equation}
where $\ssr_\ell(a)\coloneqq h(u)-h(v)-s(a)\ell_e$. It is now easy to check, using $\wtz_a\wtz_{\bar a}=\xi_e$, that \eqref{fracid2} is equal to \eqref{fracid}. An analogous argument works when $h(v)\leq h(u)$. 

Notice that 
\[
\big(\sst^{h(u)}\wtz_a^{-s(a)}\big) = \big(\sst^{h(v)}(\sst^{\ell_e})^{s(a)} \wtz_a^{-s(a)}\big) = \big(\sst^{h(v)}\wtz_{\bar a}^{s(a)}\big)=\big(\sst^{h(v)}\wtz_{\bar a}^{-s(\bar a)}\big)
\]
if $\ell_e$ divides $h(u)-h(v)$, or equivalently, $s(a)+s(\bar a)=0$; in particular, $\cT_h$ is principal at $\ssp^e$ in this case.

For each $u\in V$, the sheaf of fractional ideals $\sst^{-h(u)}\cT_h$ generates a sheaf of fractional ideals $\ssT_{h,u}$ of $\varC_u$. Given the above description of $\cT_h$, for each $a=uv\in\E$, the element $\wtz_a^{-s(a)}$ of the formal stalk of $\sst^{-h(u)}\cT_h$ at $\ssp^e$ is sent to $\barz_a^{-s(a)}$, whereas $\sst^{h(v)-h(u)}\wtz_{\bar a}^{s(a)+1}$ is sent to 0. It follows that, as sheaves of fractional ideals,
\begin{equation}\label{eq:Thu}
\ssT_{h,u}=\cO_{\varC_u}\Big(\sum_{a\in\E_u}\sss_h(a)\ssp^a\Big).
\end{equation}

As in \cite[p.~292]{EM}, we let $\ssT_h$ be the image of the composition
\[
\cT_h\xrightarrow{(\sst^{-h(v)}\,|\,v\in V)} \bigoplus_{v\in V}\sst^{-h(v)}\cT_h \xrightarrow{\qquad\qquad} \bigoplus_{v\in V}\ssT_{h,v}.
\]
Then, $\ssT_h$ is a sheaf of fractional ideals of $X$ isomorphic to $\cT_h\rest{X}$. It is thus a fractional ideal limit, that is, a limit which is a fractional ideal, of the structure sheaf along $\frakp$. As in \emph{loc.~cit.}, it can be argued that all such limits arise this way. Clearly, $\ssT_h$ generates $\ssT_{h,v}$ for each $v\in V$. As the  $\ssT_{h,v}$ are described above in \eqref{eq:Thu}, to describe the subsheaf $\ssT_h$ we need only describe it formally locally at the $\ssp^e$ where $\ssT_h$ is principal, that is, for $e=\{u,v\}\in E$ such that $\ell_e$ divides $h(u)-h(v)$.

Let $\k(v)$ be the field of rational functions of $\varC_v$ for each $v\in V$. For each arrow $a=uv\in\E$, let $\hat\k_{a}$ denote the field of fractions of $\hatcO_{\varC_u,\ssp^a}$; it contains the field of fractions of $\cO_{\varC_u,\ssp^a}$, which is $\k(u)$. 

Thus, for each $e\in E$, we may and will view any local description of a sheaf of fractional ideals of $X$ at $\ssp^e$ in $\hat\k_a\times\hat\k_{\bar a}$, where $a$ and $\bar a$ are the arrows supported on $e$. In particular, for each $h\colon V\to\Z$, the formal stalk of $\ssT_h$ at $\ssp^e$ in $\hat\k_a\times\hat\k_{\bar a}$ is the fractional ideal $(\bar\ssz_a^{-s(a)},\ssrho_e^{-s(a)}\bar\ssz_{\bar a}^{s(a)})$ if $s(a)+s(\bar a)=0$ and $((\bar\ssz_a^{-s(a)},0),(0,\bar\ssz_{\bar a}^{-s(\bar a)}))$ otherwise.

For each arrow $a=uv\in\E$, we use the local parameter $\barz_a$ to establish an isomorphism $\cO_{\varC_u}(m\ssp^a)\rest{\ssp^a}\cong \k$ for each $m\in\Z$, by viewing $\cO_{\varC_u}(m\ssp^a)$ as the sheaf of fractional ideals of $\varC_u$ whose formal stalk at $\ssp^a$ is the fractional ideal generated by $\barz_a^{-m}$, and taking $\barz_a^{-m}$ to 1. Then, for each $e=\{u,v\}\in E$, using these isomorphisms for the arrows $a=uv$ and $\bar a=vu$ supported on $e$, if $s(a)+s(\bar a)=0$ then the stalk of $\ssT_h$ at $\ssp^e$ is the kernel of the surjective composition
\begin{equation}\label{eq:gluing}
\cO_{\varC_u}(s(a)\ssp^a)_{\ssp^a}\oplus \cO_{\varC_v}(s(\bar a)\ssp^{\bar a})_{\ssp^{\bar a}} \xrightarrow{\qquad\quad} \k\oplus\k \xrightarrow{(-\ssrho_e^{-s(a)},1)} \k.
\end{equation}
Notice the symmetry above: if we exchanged $a$ and $\bar a$ the kernel would not change. 

As the $\ssrho_e$ vary, the collection of sheaves $\ssT_h$ for $h\colon V\to\Z$ varies, ultimately depending thus on the free choice
of $\rho\in\Gm^E(\k)$. Also, by Nakayama Lemma, an isomorphism $\cO_{\varC_u}(\ssp^a)|_{\ssp^a}\cong \k$ 
corresponds to choosing a formal local parameter of $\varC_u$ at $\ssp^a$ for each arrow $a=uv\in\E$. We have chosen above the formal local parameters $\bar\ssz_a$ for $\varC_u$ at $\ssp^a$ and  $\bar\ssz_{\bar a}$ for $\varC_v$ at $\ssp^{\bar a}$. Different choices $\barz_a'$ and $\barz_{\bar a}'$ can be expressed as power series $\barz_a'=\ssrho_a\barz_a+\cdots$ and $\barz_{\bar a}'=\ssrho_{\bar a}\barz_{\bar a}+\cdots$ for
$\ssrho_a,\ssrho_{\bar a}\in \k^*$, 
and we obtain the same subsheaf $\ssT_h$ of 
$\bigoplus\ssT_{h,v}$ by replacing $\ssrho_e$ by
$\ssrho_a\ssrho_e\ssrho_{\bar a}$ for each $e\in E$. 

When $\scrX$ is regular, the sheaves $\ssT_h$, rather their isomorphism classes, have been called ``twisters." In this case, the $\ssT_{-\one!_v}$ for the characteristic functions $\one!_v$ of the vertices $v$ of $V$ form what Main\`o called an enriched structure on $X$ in \cite{Maino}, where she described their moduli. To summarize, we state:

\begin{thm}\label{thm:appendix-C} Let $\ell\colon E\to\Z_{>0}$ be an integer-valued edge length function. For each integer-valued function $h\colon V\to\Z$, let $\sss_h\colon\E\to\Z$ be defined by \eqref{eq:s-C} and let $\ssT_{h,u}$ be the sheaf of fractional ideals of $X$ given by \eqref{eq:Thu} for each $u\in V$. Then there is a collection of morphisms
\[
\ssvarphi_{\ell}=(\ssvarphi_{\ell,h})_{h\in\Z^V}\colon\Gm^E\longrightarrow\prod_{h\in\Z^V}\mathrm{Quot}_{\bigoplus\ssT_{h,v}},
\]
where $\mathrm{Quot}$ stands for Grothendieck ``Quot'' scheme, such that, for each element $\rho=(\ssrho_e)_{e\in E}\in\Gm^E(\k)$, the sheaves of fractional ideals of $X$ given by the $\ssvarphi_{\ell,h}(\rho)$ are the fractional ideal limits of the trivial sheaf along a smoothing of $X$ with edge length function $\ell$ and leading coefficients $\rho_e$. More precisely, for each $\rho\in\Gm(\k)$, the subsheaf $\ssT_h\subseteq\bigoplus\ssT_{h,v}$ associated to $\ssvarphi_{\ell,h}(\rho)$ is defined by imposing that the stalk of $\ssT_h$ at $p^e$ is the kernel of \eqref{eq:gluing} for each $e=\{u,v\}\in E$ such that $\ell_e$ divides $h(u)-h(v)$. 
\end{thm}

The definition of $\ssvarphi_{\ell}$ depends on the choices made but its image does not. 


Finally, for each $v\in V$, let $\omega!_v$ be the canonical sheaf of $\varC_v$ for each $v\in V$. Let $\Omega!_v$ be the space of meromorphic differentials of $\varC_v$. We view $\omega!_v$ as a sheaf of meromorphic differentials of $\varC_v$, that is, as a subsheaf of the constant sheaf associated to $\Omega!_v$. Let $\omega$ be the canonical sheaf of $X$. Put $\Omega\coloneqq\bigoplus\Omega!_v$. We view $\omega$ as the sheaf of meromorphic differentials --- a subsheaf of the sum of the constant sheaves associated to the $\Omega!_v$ --- satisfying Rosenlicht residue conditions at the nodes, namely, the sum of the residues at the two branches over each node is zero. Putting $\ssL_h\coloneqq\omega\otimes\ssT_h$ for each $h\colon V\to\Z$, we may and will view 
$\ssL_h$ as a sheaf of meromorphic differentials as well, a subsheaf of $\bigoplus\ssL_{h,v}$ by \eqref{eq:Thu}, where
\[
\ssL_{h,v}\coloneqq\omega!_v\Big(\sum_{a\in\E_u} (1+\slztwist{\ell}h(a))\ssp^a\Big)\quad\text{for each }v\in V,
\]
the sheaf of meromorphic differentials satisfying zero and pole conditions determined by the divisor $\sum (1+\slztwist{\ell}h(a))\ssp^a$. The collection of the $\ssL_h$ can be described by the gluing conditions \eqref{eq:gluing} encoded in $\ssvarphi_{\ell}$.

In Section~\ref{sec:tropicalization-introduction},  we explained how to pass from the above setup of a smoothing $\frakp$ of $X$ to the associated tropicalization of the geometric generic fiber $\varX$ of $\frakp$. The above sheaves $\ssL_h$ are the one introduced in Section~\ref{sec:gluing-sheaves}. The gluings there were carried out by the data $\varrho!_h{(\varX)}$ of certain isomorphisms. These are induced by the leading coefficients $\rho$ after trivializations at branches of nodes are chosen.


\section{Subintegrability lemma}\label{sec:subintegrable} Let $G=(V,E)$ be a finite graph, and $\E$ its set of arrows. Recall the notation in Section~\ref{sec:graphs}. A \emph{path in $G$} is a collection $\dpat$ of arrows of $G$ that can be ordered, $\ssa_1,\dots,\ssa_k$, in such a way that $\he_{\ssa_{i}}=\te_{\ssa_{i+1}}$ for $i=1,\dots,k-1$. If in addition the tails of the $\ssa_i$ are pairwise distinct, and $\he_{\ssa_k}=\te_{\ssa_1}$, then we call $\dpat$ a (elementary) \emph{circuit in $G$}, and call $a_1,\dots,a_k$ an \emph{admissible order} for $\dpat$.

A \emph{cycle} $\dz \in H_1(G, \Z)$ is a function $\dz\colon\E\to\Z$ satisfying
$\dz(a)=-\dz(\bar a)$ for each $a\in\E$ and 
\[
\sum_{a\in\E_v}\dz(a)=0\quad\text{for each }v\in V.
\]
Its \emph{positive support}, which we denote by $\supp_+(\dz)$, is the set of arrows $a\in\E$ such that $\dz(a)>0$. Given a spanning subgraph $G'$ of $G$, we view the cycles $\dz\in H_1(G',\Z)$ as forming a subgroup of $H_1(G,\Z)$ by extending them by zero to all $\E$.

A circuit $\dq$ in $G$ can and will be identified with its corresponding cycle $\dz\in H_1(G,\Z)$, the unique cycle $\dz$ satisfying $\supp_+(\dz)=\dq$ and $\dz(a)=1$ for each $a\in\dq$.

For each cycle $\dz\in H_1(G,\Z)$, let $\bar\dz\colon\E\to\Z$ denote the function satisfying $\bar\dz(a)=\dz(\bar a)$ for each $a\in\E$; it is also a cycle. Moreover, if $\dz$ is a circuit, then also $\bar\dz$ is a circuit; if $a_1,\dots,a_k$ is an admissible order for $\dz$, then $\bar a_k,\dots,\bar a_1$ is an admissible order for $\bar\dz$.

Two circuits $\dz_1$ and $\dz_2$ are called \emph{consistent in their orientations} if there is no $a\in \E$ with $a\in\dz_1$ and $\bar a\in\dz_2$, that is, $\dz_1\cap\bar\dz_2=\emptyset$. Each cycle $\dz\in H_1(G, \Z)$ can be written as a sum of circuits pairwise consistent in their orientations. Moreover, each circuit contained in $\supp_+(\dz)$ is a summand of one of these sums.

In this section we prove the following lemma.

\begin{lemma}[Subintegrability in graphs]\label{lem:key_lemma} Let $G = (V,E)$ be a finite connected graph. Let $x \colon \E \to \R \cup\{-\infty\}$ be a function such that
  \[
  x(a) + x(\bar a) \leq 0\text{ for each }a\in\E \quad \text{and} \quad \sum_{a \in\dz} x(a)\leq 0 \text{ for each circuit $\dz$ in $G$}.
  \]
Then, there is a function $h\colon V\to\R$ such that 
\[
x(a) \leq h(u) -h(v)
\]
for each $a=uv \in \E$, with equality if and only if either $x(a)+x(\bar a)=0$ or there is a circuit $\dz$ containing $a$ such that $\sum_{b \in\dz}  x(b)= 0$.
\end{lemma}

Note that the sum $\sum_{a \in\dq} x(a)$ is nonpositive for every circuit $\dq$ in $G$ if and only if the sum $\sum x(a)\dz(a)$ over the positive support $\supp_+(\dz)$ of $\dz$ is nonpositive for every cycle $\dz\in H_1(G,\Z)$. If this is the case, $\sum x(a)\dz(a)=0$ for a cycle $\dz$ if and only if $\sum_{a \in\dq} x(a)=0$ for each circuit $\dq$ contained in $\supp_+(\dz)$.

\begin{proof}[Proof of Lemma~\ref{lem:key_lemma}] 
  Let $\dz_1$ and $\dz_2$ be two circuits in $G$ with  
  \[
  \sum_{a \in\dz_1}  x(a) = \sum_{a \in\dz_2}  x(a) = 0. 
  \]
 We claim that for each arrow $a\in\E$ 
 such that $a \in\dz_1$ and $\bar a\in \dz_2$, we have $x(a)+x(\bar a)=0$. Indeed, let $\dz \coloneqq\dz_1+\dz_2$.
 We have 
 \[
 \begin{aligned}
   0 \geq \sum_{a\in\dz} x(a)\dz(a) =& \sum_{a\in\dz_1} x(a) + \sum_{a\in\dz_2} x(a) 
  - \sum_{a\in\dz_1\cap\bar\dz_2} (x(a)+x(\bar a))\\
  =&- \sum_{a\in\dz_1\cap\bar\dz_2} (x(a)+x(\bar a)) \geq 0.
 \end{aligned}
  \]
Then the inequalities are all equalities, that is, $x(a) + x(\bar a) =0$ for each $a\in \dz_1\cap\bar\dz_2$.

Let $F(x)$ be the set of all arrows $a\in\E$ with $x(a)+x(\bar a)<0$. Let $S(x)$ be the set of all circuits $\dz$ in $G$ with 
  \[
  \sum_{a \in \dz} x(a) <0.
  \]  
We proceed by induction on the size $s(x)\coloneqq |F(x)|+|S(x)|$. If $s(x) =0$, we get $x(a) = -x(\bar a)$ for each $a\in \E$, and   $\sum_{a \in \dz} x(a) =0$ for each circuit $\dz$ in the graph. It is a standard fact that in this case we can find a function $h \colon V \to \R$ such that $h(u)-h(v) = x(a)$ for each $a=uv\in\E$. 

Let $s$ be a positive integer. 
Suppose the statement in the lemma holds for each $x$ with $|F(x)|+|S(x)|$ smaller than $s$. 
Let $x\colon \E \to \R\cup\{-\infty\}$ be a function verifying the conditions of the lemma with  $s(x) = s$. We claim that $F(x)$ is nonempty. Indeed, otherwise, we would have $x(a)+x(\bar a)=0$ for all $a\in \E$. But then, for each circuit $\dz$ in $G$, we would have
\[
0 \geq \sum_{a \in \dz} x(a) = -\sum_{a \in \dz}x(\bar a) = -\sum_{a\in \bar\dz} x(a) \geq 0.
\]
This would imply that $S(x) =\emptyset$, contradicting the positivity of $s$. 

Let $a\in F(x)$. Let $S_a$, resp.~$S_{\bar a}$, be the set of circuits in the graph $G$ that contain $a$, resp.~$\bar a$.  By the claim in the beginning of the proof, at least one of the two sets, $S_a$ or $S_{\bar a}$, say $S_{a}$ without loss of generality, is entirely included in $S(x)$. Indeed, otherwise, we would find two circuits $\dz_1$ and $\dz_2$ excluded from $S(x)$ with $a\in \dz_1$ and $\bar a\in \dz_2$. This would imply $x(a)+x(\bar a)=0$, contradiction.

Let $-\epsilon$ be the maximum of $x(a)+x(\bar a)$ and the quantities $\sum_{b \in \dz} x(b)$ for circuits $\dz \in S_a$. We have $\epsilon>0$. Define $x'\colon \E \to \R \cup\{-\infty\}$ by $x'(a) \coloneqq x(a)+\epsilon$, and  $x'(b) = x(b)$ for every $b\in \E$ different from $a$.
By the choice of $\epsilon$, we have $x'(b)+x'(\bar b) \leq 0$ for each $b\in\E$, and for each circuit $\dz$ in the graph, 
\[
  \sum_{b \in \dz} x'(b) \leq 0.
\]
Furthermore, since $x'\geq x$, we have $F(x')\subseteq F(x)$ and $S(x')\subseteq S(x)$. 
On the other hand, either $x'(a) +x'(\bar a) =0$, or there exists $\dz\in S_a$ with   $\sum_{b\in\dz} x'(b) = 0.$ This means $|F(x')|+|S(x')|<s$. 

Applying the hypothesis of our induction to $x'$, we get a function $h\colon V \to \R$ such that for each edge $b=uv \in \E$, we have $x'(b) \leq h(u)-h(v)$. Since $x\leq x'$, it follows that
\begin{equation}\label{eq:xhx-0}
    x(a)\leq h(u)-h(v)\quad\text{for each }a=uv\in\E.
\end{equation}

To prove the remainder of the lemma, assume \eqref{eq:xhx-0} holds. If  $x(a)+x(\bar a)=0$ for $a=uv\in\E$, then
\[
0=x(a)+x(\bar a)\leq\big(h(u)-h(v))+(h(v)-h(u)\big)=0,
\]
and hence $x(a)=h(u)-h(v)=-x(\bar a)$. Also, if $\sum_{a\in\dz} x(a)=0$ for a circuit $\dz$ in $G$, then
\[
0=\sum_{a\in\dz} x(a)\leq\sum_{a\in\dz}\big(h(\te_a)-h(\he_a)\big)=0,
\]
and hence $x(a)=h(u)-h(v)$ for each $a=uv\in\dz$.

On the other hand, given $a=uv\in\E$ such that $a\in F(x)$ and $\ssS_a\subseteq S(x)$, there is $\epsilon>0$ such that, defining $x'\colon \E \to \R \cup\{-\infty\}$ by $x'(a) \coloneqq x(a)+\epsilon$ and $x'(b)\coloneqq x(b)$ for $b\neq a$, we have $x'(b)\leq h(y)-h(w)$ for each $b=yw\in\E$ but $a\not\in F(x')$ or $\ssS_a\not\subseteq S(x')$. Then, as shown above, $x'(a)= h(u)-h(v)$, and hence $x(a)<h(u)-h(v)$, as required.
\end{proof} 

\end{document}